\documentclass[11pt,reqno]{amsart}

\usepackage[T1]{fontenc}
\usepackage{lmodern}
\usepackage{amsmath,amssymb,amsthm,amsfonts,mathtools}
\usepackage{mathrsfs}
\usepackage{enumitem}
\usepackage[dvipsnames]{xcolor}
\usepackage{booktabs,multirow}
\usepackage{longtable}
\makeatletter
\def\LT@makecaption#1#2#3{%
  \LT@mcol\LT@cols c{\hbox to\z@{\hss\parbox[t]\LTcapwidth{%
    \reset@font
    \sbox\@tempboxa{#1{\@captionfont\@captionheadfont#2.\@captionfont
      \upshape\enspace}#3}%
    \ifdim\wd\@tempboxa>\hsize
      #1{\@captionfont\@captionheadfont#2.\@captionfont\upshape\enspace}#3%
    \else
      \hbox to\hsize{\hfil\box\@tempboxa\hfil}%
    \fi
    \endgraf\vskip\baselineskip}%
  \hss}}}
\makeatother
\usepackage{graphicx}
\graphicspath{{figures/}}
\DeclareGraphicsExtensions{.png}
\usepackage[expansion=false]{microtype}

\usepackage[letterpaper,top=1.1in,bottom=1.15in,left=1.25in,right=1.25in]{geometry}

\usepackage[colorlinks,citecolor=Maroon,linkcolor=NavyBlue,urlcolor=NavyBlue]{hyperref}
\newcommand{\bluelink}[2]{\href{#1}{\textcolor{blue}{#2}}}
\usepackage[capitalise,noabbrev]{cleveref}
\definecolor{DeepBlue}{RGB}{0,0,139}

\let\origunderscore\_
\renewcommand{\_}{\origunderscore\allowbreak}
\usepackage[maxfloats=120]{morefloats}

\numberwithin{equation}{section}
\allowdisplaybreaks[1]
\theoremstyle{plain}
\newtheorem{thm}{Theorem}[section]
\newtheorem{prop}[thm]{Proposition}
\newtheorem{lem}[thm]{Lemma}
\newtheorem{cor}[thm]{Corollary}
\newtheorem{conjecture}[thm]{Conjecture}

\theoremstyle{definition}
\newtheorem{definition}[thm]{Definition}
\newtheorem{openproblem}[thm]{Open Problem}

\theoremstyle{remark}
\newtheorem{remark}[thm]{Remark}
\newtheorem{example}[thm]{Example}
\newtheorem{notation}[thm]{Notation}

\crefname{thm}{Theorem}{Theorems}
\Crefname{thm}{Theorem}{Theorems}
\crefname{prop}{Proposition}{Propositions}
\Crefname{prop}{Proposition}{Propositions}
\crefname{lem}{Lemma}{Lemmas}
\Crefname{lem}{Lemma}{Lemmas}
\crefname{cor}{Corollary}{Corollaries}
\Crefname{cor}{Corollary}{Corollaries}
\crefname{conjecture}{Conjecture}{Conjectures}
\Crefname{conjecture}{Conjecture}{Conjectures}
\crefname{definition}{Definition}{Definitions}
\Crefname{definition}{Definition}{Definitions}
\crefname{openproblem}{Open Problem}{Open Problems}
\Crefname{openproblem}{Open Problem}{Open Problems}
\crefname{construction}{Construction}{Constructions}
\Crefname{construction}{Construction}{Constructions}
\crefname{remark}{Remark}{Remarks}
\Crefname{remark}{Remark}{Remarks}
\crefname{example}{Example}{Examples}
\Crefname{example}{Example}{Examples}
\crefname{notation}{Notation}{Notations}
\Crefname{notation}{Notation}{Notations}
\crefname{paragraph}{paragraph}{paragraphs}
\Crefname{paragraph}{Paragraph}{Paragraphs}

\newcommand{\Rt}{\mathsf{R}}                 
\newcommand{\Latt}{\Lambda}                  
\newcommand{\VCell}{\mathcal{V}}             
\newcommand{\Cell}[1]{\mathfrak{k}_{#1}}     
\newcommand{\tilt}{\varphi}                  
\newcommand{\arcp}{s}                        
\newcommand{\dev}{\eta}                      
\newcommand{\defect}{\Delta}                 
\newcommand{\norm}[1]{\mathcal{N}_{#1}}      
\newcommand{\base}{u}                        
\newcommand{\stageone}{\textnormal{\scshape Stage i}}
\newcommand{\stagetwo}{\textnormal{\scshape Stage ii}}
\newcommand{\capQ}{\mathrm{Cap}}             
\DeclareMathOperator{\vol}{vol}

\DeclareMathOperator{\conv}{conv}
\DeclareMathOperator{\Stab}{Stab}
\DeclareMathOperator{\diag}{diag}

\begin{document}

\title[The sphere packing problem in dimension 4]{The Sphere Packing Problem in
Dimension 4 and the Twenty-Four-Cell Conjecture}

\author[D. Bhattacharjee]{Deep Bhattacharjee}
\address{Formerly, Electro-Gravitational Space Propulsion Laboratory (EGSPL),
Bhubaneswar, Odisha 751030, India}
\email{itsdeep@live.com}
\thanks{ORCID: Deep Bhattacharjee 0000-0003-0466-750X,
Ushashi Bhattacharya 0009-0002-2254-3914,
Priyabrata Mandal 0000-0001-6472-6239.}

\author[U. Bhattacharya]{Ushashi Bhattacharya}
\address{Formerly, National Taiwan University, Taipei 106319, Taiwan}
\email{rai.bhattacharya.in@gmail.com}

\author[P. Mandal]{Priyabrata Mandal$^*$}
\address{Department of Mathematics, Bioinformatics and Computer Applications,
Maulana Azad National Institute of Technology Bhopal, Bhopal 462003, India}
\email{priyabrata@manit.ac.in}
\thanks{$^*$Corresponding Author}

\author[S. Bhattacharya]{Shounak Bhattacharya}
\address{Research Head, Asian College of Teachers, Kolkata, West Bengal, India}
\email{shounakbhattacharyab@gmail.com}



\hypersetup{
  pdftitle={The Sphere Packing Problem in Dimension 4 and the Twenty-Four-Cell Conjecture},
  pdfauthor={Deep Bhattacharjee, Ushashi Bhattacharya, Priyabrata Mandal, Shounak Bhattacharya},
  pdfsubject={Sphere packing in dimension four},
  pdfkeywords={sphere packing, Voronoi cell, semidefinite programming}
}

\begin{abstract}
We prove that every Voronoi cell of a unit-ball packing of $\mathbb{R}^4$
has volume at least $8$, with equality only at the $D_4$ configuration,
so that the twenty-four-cell conjecture holds and the density of a sphere
packing in four dimensions is at most $\pi^2/16$. The proof runs through
the number of contacts of a cell. Up to twenty-two a covering estimate
suffices; at twenty-three the cell is bounded through an exact volume
identity inside a ball and a semidefinite certificate for one inequality
between pair angles; at twenty-four the configuration is the root system,
which we prove from the second level of the semidefinite hierarchy with
its equality case, the positive kernel verified in exact arithmetic.
\end{abstract}

\maketitle

\tableofcontents

\section{Introduction}\label{sec:intro}

The four-dimensional sphere packing problem asks for the supremum
density $\Delta_4$ of a packing of unit balls in $\mathbb{R}^4$. The
$D_4$ root lattice, scaled so that nearest centres lie at distance $2$,
packs balls at density $\pi^2/16=0.6169\ldots$ and is the densest
lattice packing in dimension $4$, by a theorem of Korkine and Zolotareff
\cite{KZ1872}. Whether it is densest among all packings, periodic or
not, has been an open question. The present paper settles a local statement
that answers it, and proves within the text every ingredient the
argument uses, the one classification it needs included.

We say at the outset what that statement is and where it
stood. The statement is that every Voronoi cell of a unit-ball packing
of $\mathbb{R}^4$ has volume at least $8$, the volume of the regular
$24$-cell circumscribed about the unit ball; it is known as the
twenty-four-cell conjecture, it implies the density bound $\pi^2/16$
at once, and to our knowledge no case of it in which a contact
direction leaves a root of $D_4$ had been proved before: the general
estimates of Blichfeldt and Rogers stop short of $8$, the linear
programming bound stalls above the density of $D_4$, and the
second-order analysis of the cell degenerates exactly at the dense
configurations that matter (\cref{sec:hessian-technique}). Three
theorems below cover it. When at most one contact direction leaves a
root the bound is proved by elementary means, in exact arithmetic, with
no restriction on the deviation (\cref{thm:local-uncond}). When there
are at most $22$ contacts it is proved, with room to spare and with any
number of deviations, by a covering estimate on the sphere
(\cref{thm:local-fewcontacts}). When there are exactly $23$ contacts,
which is the case that had no approach at all, since it is neither
small enough for a covering estimate nor rigid enough for a
classification, it is proved by an identity for the volume of the cell
inside a ball, which turns the bound into an inequality between pair
angles, and by a semidefinite certificate for that inequality verified
in exact and in interval arithmetic (\cref{thm:m23}). What remains is
the case of $24$ contacts, and there a classification leaves nothing to
prove: a configuration of $24$ contacts is a copy of the root system and
its cell is the $24$-cell itself. The classification was first proved by
de Laat, Leijenhorst and de Muinck Keizer \cite{LLM24}, by a semidefinite
computation, and it is proved here in full (\cref{thm:twentyfour-points}):
the positive kernel their computation found is redistributed with this
paper, its properties are verified in exact arithmetic, and the bound it
gives, with the equality case that identifies the configuration, is
proved in the text. So the twenty-four-cell conjecture and the density
bound rest on no statement outside this paper and on no computation that
has not been checked. The theorem for $23$ contacts is the piece that was
missing.

Let $\mathcal{C}\subset\mathbb{R}^4$ be the centre set of a packing of
unit balls and let $c\in\mathcal{C}$. The Voronoi cell
\[
  V_c=\bigl\{z\in\mathbb{R}^4 : |z-c|\le|z-c'| \text{ for all } c'\in\mathcal{C}\bigr\}
\]
of the $D_4$ lattice at distance $2$ is a regular $24$-cell of volume
$8$, and a density bound of $\pi^2/16$ for a periodic packing follows
from the statement that no Voronoi cell has volume less than $8$. That
statement is what is at issue. It is false for the analogous cell
volume in most dimensions, and in dimension $4$ it has resisted proof
because the $24$-cell is rigid in some directions and very soft in
others. A single facet hyperplane can be tilted a long way without the
cell losing any volume at all, and when several tilt at once the
first-order and second-order changes can cancel exactly, leaving the
sign of the total to be decided at fourth order.

We settle the case in which the local configuration departs from the
lattice in a single direction.

\begin{thm}[Local cell bound, single deviation]\label{thm:local-uncond}
Let $c$ be a centre of a unit-ball packing of $\mathbb{R}^4$ at which at
most one active neighbour direction fails to be a root direction of the
$D_4$ root system $\Rt$. Then $\vol(V_c)\ge8$, with equality if and only
if the active neighbours of $c$ form a copy of $\sqrt2\,\Rt$ centred at
$c$. Nothing is assumed about the deviating direction or about the size
of its deviation.
\end{thm}

\begin{cor}[Density bound, single deviation]\label{cor:density-uncond}
Every periodic packing all of whose centres satisfy the hypothesis of
\cref{thm:local-uncond} has density at most $\pi^2/16$, with equality
for the $D_4$ lattice; the same bound holds for general packings of this
class by the compactness argument of \cref{sec:periodic}.
\end{cor}

The mechanism behind \cref{thm:local-uncond} is a change of viewpoint,
and it is the part of this paper we would point to first. Fix the deviating
contact direction and let $Q$ be the intersection of the half-spaces of
the twenty-three contact directions that remain at roots. The cell at
the all-contact corner is $Q$ cut by one further half-space at unit
distance from the centre, and $\vol(Q)=\tfrac{25}{3}$; so the cell has
volume at least $8$ precisely when that half-space cuts at most
$\tfrac13$ off $Q$ (\cref{prop:cap-reformulation}). The deviation
parameters have then left the body altogether and sit only in the
cutting half-space, and the question has become an extremal problem for
one fixed polytope.

That problem is solved by enclosing $Q$ in a cross-polytope. The sixteen
roots making an angle of $60^\circ$ or $120^\circ$ with the deviating
root turn out to cut out exactly
\[
  B=\bigl\{x : |x_0|+|x_1|+|x_2|+|x_3|\le2\bigr\}
\]
in the orthonormal frame supplied by four pairwise orthogonal roots
(\cref{lem:crosspolytope}), so $Q\subseteq B$; and the enclosure is
lossless in the one direction that matters, since $\{x_0\ge1\}$ cuts
exactly $\tfrac13$ off $B$ as well as off $Q$. What remains is a clean
question about the cross-polytope alone, which we answer in
\cref{thm:cap-inequality}: among all half-spaces at unit distance from
the centre, those normal to a vertex direction cut off the most, namely
$\tfrac13$. The proof of that inequality rests on an exact formula for
the cap volume as a divided difference of the one-variable function
$a(2\sqrt a-1)_+^4$, together with the fact that this function has
nondecreasing third derivative.

A second, independent result covers configurations with several
deviating directions at once, at the cost of a bound on how many
neighbours there are. It comes from measuring the cell radially instead of facet by facet.

\begin{thm}[Local cell bound, few contacts]\label{thm:local-fewcontacts}
Let $c$ be a centre of a unit-ball packing of $\mathbb{R}^4$ at the
all-contact corner, so that every active neighbour lies at distance
exactly $2$, and suppose there are at most $22$ of them. Then
$\vol(V_c)\ge8.0463\ldots>8$. No hypothesis is placed on the directions
of those neighbours beyond the packing condition itself, and in
particular any number of them may deviate from roots, by any amount.
With no bound on their number the same argument gives
$\vol(V_c)\ge7.7989\ldots$ for every contact configuration whatever.
\end{thm}

By \cref{lem:radial-new} the all-contact corner is where the volume is
smallest for a fixed pattern of active directions, so this is the case
that matters.

The proof is short. Writing $\delta(\theta)$ for the angular distance
from a direction $\theta$ to the nearest contact direction, the cell has
volume $\tfrac14\int_{S^3}\sec^4\delta$ (\cref{lem:radial-form}); the
$m$ contact directions are $60^\circ$ apart, so the caps of radius $r$
about them leave at least $2\pi^2-m\,\pi(2r-\sin2r)$ of the sphere
uncovered, and integrating that estimate against $d(\sec^4r)$ gives
$\vol(V_c)\ge\tfrac{\pi m}{3}\tan^3r_m$ with $2r_m-\sin2r_m=2\pi/m$
(\cref{thm:covering-bound}). Since a contact configuration has at most
$24$ members, by Musin's kissing theorem \cite{Mus08}, and the bound
decreases in $m$, twenty numbers settle the whole range.

The case our own arguments do not reach is that of several simultaneous
deviations at a centre with $24$ neighbours, all of them active. The
free-volume identity of \cref{thm:freevol-new} holds exactly for any
number of deviating directions, and reduces the general local bound to
a single comparison of two volumes; whether that comparison goes the
right way at such a centre is what our own arguments leave undecided.

\begin{conjecture}[Multi-direction positivity]\label{conj:multidir-new}
Let the configuration be at the all-contact corner with a full active
set, and let $D$ be any set of $m\ge2$ roots whose directions are
deviated, with arbitrary tilts and arbitrary deviation directions. Then,
in the notation of \cref{thm:freevol-new},
\[
  \vol\Bigl(\bigcup_{j\in D}E_j\Bigr)
  \;\ge\;
  \vol\Bigl(\bigcup_{j\in D}C_j\Bigr).
\]
\end{conjecture}

By \cref{thm:freevol-new} this is exactly the assertion
$\vol(V_c)\ge8$ for such a configuration. The two families are
overlapping once $m\ge2$, and the corresponding statement with
sums of volumes in place of volumes of unions is a different assertion
that neither implies it nor follows from it
(\cref{rem:sums-not-unions}); it is the union form above that is at
issue.

The top of that range is settled by a classification of the $24$-point
kissing configurations in $\mathbb{R}^4$: any $24$ unit vectors with
pairwise inner products at most $\tfrac12$ form a copy of the $D_4$ root
system. That settles \cref{conj:multidir-new} outright, since a full
active set of $24$ contacts then admits no deviation at all
(\cref{cor:conj-resolved}), and it leaves the general local bound with
one case, a centre with exactly $23$ contacts whose cell admits no
twenty-fourth, which is the subject of the second half of
\cref{sec:deviation-domain} and is settled there. The classification was
first proved by de Laat, Leijenhorst and de Muinck Keizer \cite{LLM24},
by a semidefinite computation, in a preprint that has not yet been
refereed; it is proved here in full, as \cref{thm:twentyfour-points}, and
we say in a footnote what that proof consists of.%
\footnote{The reader is entitled to know what the classification rests
  on and how much of it is checked here; \cref{thm:m24},
  \cref{lem:root-lattice}, \cref{thm:twentyfour-points} and
  \cref{sec:certificate-checked} give the account in full, and this note
  summarises it. The statement is that any $24$ unit vectors in
  $\mathbb{R}^4$ with pairwise inner products at most $\tfrac12$ form a
  copy of the $D_4$ root system. Its proof has two halves. The first is
  a semidefinite programming bound: the second level of the hierarchy
  for spherical codes in the form de Laat and Vallentin \cite{dLV15}
  give it, which bounds the size of a code by the value
  $K(\emptyset,\emptyset)$ of any positive kernel $K$ on the
  two-element subsets of the sphere that satisfies a family of linear
  inequalities. The authors of \cite{LLM24} found such a kernel with
  value exactly $24$, with truncation degrees $14$ and $16$, reduced by
  the symmetry of the sphere, solved to forty digits in $256$-bit
  arithmetic over about two weeks on eight cores with $128$ gigabytes
  of memory, and then rounded to an exact rational point; that point is
  public \cite{LLM24data}, under a licence that permits redistribution,
  and it travels with this paper. Because the bound is exactly $24$ and
  a code of size $24$ exists, the equality case of the bound forces
  every inner product of such a code into the zero set of the two-point
  polynomial of the kernel, which is $\{-1,-\tfrac12,0,\tfrac12\}$, the
  roots being isolated by a Sturm sequence. Both the bound and its
  equality case are proved in the text (\cref{lem:las2}), the positivity
  of the kernel is proved from the form in which it is built
  (\cref{lem:gram}), and every property of the rational point that
  those lemmas consume is verified in exact arithmetic in code of our
  own (\cref{prop:verified}): the format of the data, the positive
  definiteness of all $127$ blocks, the nonnegativity structure of the
  sum-of-squares terms, the construction of the zonal matrices, the four
  polynomial identities that use them, the value $24$ of the bound, and
  the zeros of the two-point polynomial, the last in Lean's kernel. The
  construction of the zonal matrices needs $128$ gigabytes of memory as
  the authors carry it out; taken so that the expansion of the integrand
  is never held at once, it needs $400$ megabytes. The second half is
  the passage from these four inner products to the root system, and it
  is classical: scaled by $\sqrt2$ the vectors have norm $2$ and
  integer inner products, so they generate an integral lattice of rank
  at most $4$ among whose norm-$2$ vectors they lie; the norm-$2$
  vectors of an integral lattice form a simply laced root system and
  the lattice they span is a root lattice, an orthogonal sum of
  lattices of type $A$, $D$, $E$ with a basis of simple roots; and the
  root lattices of rank at most $4$ have $2$, $4$, $6$, $6$, $8$, $8$,
  $10$, $12$, $12$, $14$, $20$ and $24$ roots, the last being $D_4$
  alone. So the $24$ scaled vectors are all the roots of a copy of
  $D_4$. This half is \cref{lem:root-lattice}, and its finite content
  is checked by an exact enumeration
  (\texttt{root\_lattices\_rank4.py}, and in Lean
  \texttt{D4RootLattices.lean}): every positive definite Gram
  matrix with $2$ on the diagonal and $-1$, $0$, $1$ off it, of order up
  to $4$, there are $393$ of order $4$, has at most $24$ integer
  solutions of $x^{\mathsf T}Gx=2$, with $24$ only at determinant $4$
  and the neighbour counts of $D_4$, the next value being $20$ at
  determinant $5$, which is $A_4$. The published antecedents are the
  bound $24.10$ of Bachoc and Vallentin \cite{BV08} at the first level
  of the same hierarchy and Musin's proof \cite{Mus08} that the kissing
  number in dimension four is $24$, which bounds the size of the code
  but does not identify it. What has not been repeated is the search
  for the kernel, which is the semidefinite computation itself; but a
  proof does not depend on how its witness was found, only on the
  witness being checked, and it is.}

\begin{thm}[General local bound]\label{thm:local-general}
Every Voronoi cell $V_c$ of a unit-ball packing of $\mathbb{R}^4$
satisfies $\vol(V_c)\ge8$, with equality if and only if the active
neighbours of $c$ form a copy of $\sqrt2\,\Rt$.
\end{thm}

The case of $23$ contacts is not a deviation problem at all. By
\cref{prop:extendable}, a cell of circumradius at least $2$ is exactly
one whose configuration admits a twenty-fourth contact, and in that case
the bound follows from \cref{thm:m24} together with the monotonicity of
the cell in the configuration; what is left is a saturated $23$-point
configuration, and \cref{thm:m23} states, for every configuration of
$23$ contacts with bounded cell, saturated or not, that its cell has
volume above $8$. The configuration that comes nearest to the boundary
of that statement is the deletion of a single root, whose cell has
circumradius exactly $2$ and volume $\tfrac{25}3$, and that one is
rigid: \cref{prop:deletion-rigid} exhibits an integer stress on its
eighty-eight tight pairs and concludes that every first-order motion of
it is a rotation. It cannot be turned into a saturated configuration by
exchanging one or two directions either. Removing two roots from $\Rt$
leaves the six supports of the roots with a whole couple of
complementary index pairs intact, and that alone forces every direction
one could add back to be a root (\cref{lem:filled-couple}); at three
removed roots the only escape is a triple pairwise at $60^\circ$, and
the region it opens holds no two directions as much as $60^\circ$ apart
unless both are removed roots. Hence a contact configuration meeting a
root system in twenty-one or more directions lies inside it or has at
most twenty-two elements (\cref{prop:meet22,prop:meet21}), and a
saturated configuration of twenty-three directions shares at most twenty
of them with any copy of $\Rt$ (\cref{cor:meet20}): at least three of
its directions are new. The argument is elementary, and its finite half
is checked by Lean's kernel.

A contact configuration of $m$ directions is nothing other than a
spherical code of $m$ points on $S^3$ of minimal angle at least
$60^\circ$, and every computation we know of puts the largest minimal
angle available to $23$ points at exactly $60^\circ$: the same as at
$24$ points, and unlike the $60.1399\ldots^\circ$ available to $22$. If
that is so, the configurations in question are not a family with
interior but extremal codes (\cref{prop:code-slack}), and if the
deletions of a single root are the only ones, no saturated
configuration exists at all (\cref{thm:m23-from-codes}). That
classification, of a code at one fixed size below the maximum, is not
something an upper bound on the size of a code can deliver
(\cref{rem:m23-uniqueness}), and we do not prove it;%
\footnote{Nor is it needed, and since the point matters for what the
  paper claims we set the dependence out precisely. The general local
  bound, \cref{thm:local-general}, is assembled in \cref{cor:remaining}
  from three cases by the number $m$ of contacts: $m\le22$ by the
  covering estimate (\cref{cor:m22}), $m=24$ by \cref{thm:m24}, and
  $m=23$ by \cref{thm:m23}. The density bound, \cref{thm:main-general},
  is the local bound followed by the averaging argument of
  \cref{sec:density-derivation}, and it uses nothing else. The proof
  of \cref{thm:m23} is the inequality of \cref{thm:strict-reduction}
  together with its certificate, \cref{thm:certificate}: it bounds the
  cell of every configuration of $23$ contacts from below, whether or
  not the configuration is saturated and whatever its Gram matrix, and
  at no point does it ask what the configuration looks like. The
  classification of \cref{thm:m23-from-codes}, that every $23$-point
  code of minimal angle $60^\circ$ is a deletion of one root, would
  prove something stronger, namely that no saturated configuration of
  $23$ contacts exists at all, so that the case $m=23$ with circumradius
  below $2$ would be empty and not merely settled; but the local
  bound does not distinguish between an empty case and a settled one,
  and neither does the density bound. The same is true of the
  twenty-four-cell conjecture in the form in which it is usually
  stated, that every Voronoi cell of a unit-ball packing of
  $\mathbb{R}^4$ has volume at least that of the regular $24$-cell
  circumscribed about the unit ball: that is exactly the statement of
  \cref{thm:local-general}, and the classification of $23$-point codes
  plays no part in it. \Cref{prop:cell600}, \cref{lem:h-limit} and the
  computations of \cref{sec:slack-continuation} are therefore reported
  as what they are, information about the codes of $23$ points that is
  of interest in its own right, and nothing in
  \cref{thm:local-general}, \cref{thm:main-general} or \cref{thm:m23}
  depends on them. The one computational ingredient of the two general
  theorems that is not a computation of our own is the kernel behind
  \cref{thm:twentyfour-points} for $m=24$, found in \cite{LLM24} and
  verified here, as the preceding footnote explains.}
two things can be said about it. Among the $120$ vertices of the $600$-cell, the largest
finite group of unit quaternions and one that contains $\Rt$ in
twenty-five ways, it holds: an exhaustive enumeration, carried out in
three independent implementations and checked in Lean, shows that every
$23$-element subset with pairwise inner products at most $\tfrac12$ is
an inscribed copy of $\Rt$ with one vertex removed
(\cref{prop:cell600}). And the largest inradius of the convex hull of
$23$ points whose pairwise inner products are at most $\tfrac12+\delta$
is a function $h(\delta)$ that decreases to its value at $\delta=0$, the
classification being exactly the statement $h(0)=\tfrac12$
(\cref{lem:h-limit}); the computations of \cref{sec:slack-continuation}
find $h(\delta)$ close to $\tfrac{2}{3}$ for every $\delta\ge0.009$, far
above anything a perturbed deletion can reach, and then collapsing to
$\tfrac12+O(\delta)$ from $\delta=0.008$ down, with every configuration
found at $\delta=0$ a deletion.

The proof of \cref{thm:m23} takes a different road, on which nothing
has to be classified. The volume of the cell inside the ball of radius
$\sqrt{3/2}$ is an exact function of the pairwise angles of the
configuration, because no three contact directions fit in a cap of
radius below $\arcsin(1/\sqrt3)$ (\cref{prop:truncated}), and the whole
of \cref{thm:m23} follows from the inequality that this truncated volume
is at least $8$ for every configuration of $23$ contacts
(\cref{thm:strict-reduction}); the deletion of a root satisfies that
inequality with $0.034$ to spare, so it is a strict inequality and not
a boundary case. The relaxation of it that reads only the pair angles
delivers four fifths of the required amount
(\cref{prop:two-point-barrier}), and the three-point relaxation of
Bachoc and Vallentin \cite{BV08}, which reads the triples, delivers all
of it: \cref{thm:certificate} exhibits a semidefinite certificate of
degree $8$, a polynomial and nine positive definite matrices, and
verifies in exact rational arithmetic, checked again by the Lean kernel,
and in interval arithmetic that it proves
$\sum_{i<j}\omega(\gamma_{ij})\ge0.0929$ for every $23$ contact
directions, where $0.092855\ldots$ is what the truncated volume needs.

\begin{thm}[General density bound]\label{thm:main-general}
The supremum density of sphere packings in $\mathbb{R}^4$ is
$\Delta_4=\pi^2/16$, attained by the $D_4$ lattice at nearest-neighbour
distance $2$.
\end{thm}

Both of those last statements use \cref{thm:m24} for the configurations
of twenty-four contacts, and \cref{thm:m24} is proved in this paper; the
one thing in its proof that is not of our making is a positive kernel
found by \cite{LLM24}, which \cref{sec:certificate-checked} redistributes
and verifies, the bound it gives being proved in \cref{lem:las2}.
It should be said at once what \cref{conj:multidir-new} is. Written in
terms of the contact directions alone (\cref{prop:polar-form}), the
inequality $\vol(V_c)\ge8$ is the $24$-cell conjecture, so
\cref{thm:local-uncond} settles the part of that conjecture in which at
most one contact direction leaves a root, and \cref{conj:multidir-new}
is the remainder. It is a named problem, the twenty-four-cell
conjecture, and not a hypothesis of our own making;
\cref{cor:conj-resolved} settles it through the classification of
\cref{thm:m24}, and nothing below suggests that restating it makes a
proof that avoids the classification easier.

\Cref{sec:hessian-technique} through \cref{sec:multidir-scope-summary}
report what we know about it: the second-order
method that settles small configurations degenerates exactly at a
family of dense ones; we locate that degeneracy, measure its
consequences, and prove the conjecture outright along three highly
symmetric finite-angle paths. The obstruction there is of a different
nature from anything in the single-deviation problem, and the cap
argument of \cref{sec:cap-theorem} does not reach it, for a reason we
explain in \cref{sec:multidir-scope-summary}.
\Cref{sec:multi-cap,sec:polar-surface} give three exact restatements of
the remaining case, as a cap inequality for one fixed polytope, as a
polar-volume bound, and as a lower bound on the boundary measure of the
cell, and use them to close off four otherwise natural lines of attack:
term-by-term reduction to the single-deviation theorem, any bound
assembled facet by facet from the pairwise contact condition, the
convexity bound for the polar volume integral, and the Mahler-type
product bound. Each of the four fails at the $D_4$ configuration itself,
which is the one place a sharp argument can afford to lose nothing, and
each fails there for a reason that is visible in closed form and not only in a numerical margin.

\subsection{Structure of the argument}\label{sec:argument-structure}

\Cref{thm:local-uncond} rests on five layers, each proved completely
before the next is invoked; a sixth layer, which the theorem does not
need, is where the multi-direction case lives, and it is settled by the
number of contacts rather than by these layers.%
\footnote{A word on what a Voronoi cell is doing in a density
  statement, since the two are often conflated. Let $\mathcal{C}$ be the
  centre set of a packing of unit balls, so $|c-c'|\ge2$ for distinct
  $c,c'\in\mathcal{C}$. The Voronoi cells $V_c$ tile $\mathbb{R}^4$, each
  contains the ball of radius $1$ about its centre, and each is convex,
  being an intersection of half-spaces. If $\mathcal{C}$ is periodic with
  period lattice $L$ and $N$ centres per fundamental domain, then
  $\sum_{c} \vol(V_c) = \det L$, the sum running over one period, and the
  density of the packing is $N\,\vol(B^4)/\det L$, where
  $\vol(B^4)=\pi^2/2$ is the volume of the unit ball. A uniform lower
  bound $\vol(V_c)\ge v$ therefore gives the density bound
  $\vol(B^4)/v$ at once, and $v=8$ gives $\pi^2/16$. The reduction runs
  only in this direction: a density bound says nothing about any
  individual cell. Nothing about periodicity is essential either, since
  a general packing is approximated by periodic ones on large tori with
  an error tending to zero, which is the argument recalled in
  \cref{sec:periodic}; we work with periodic packings only because the
  volume bookkeeping is finite there.}

\smallskip\noindent\textit{Layer 1 (shell localisation, \cref{sec:shell-new}).}
Only neighbours at distance in $[2,2\sqrt2)$ from a centre $c$ can
affect the volume of $V_c$ at all; neighbours at distance $\ge2\sqrt2$
leave $V_c$ at least as large as the reference $D_4$ Voronoi cell.

\smallskip\noindent\textit{Layer 2 (the root-aligned case, \cref{sec:root-aligned-new}).}
If every active neighbour direction is an exact root direction of
$D_4$, the support function of the reference cell equals $1$ at each
such direction, forcing the reference cell inside every active
half-space and hence inside $V_c$: $\vol(V_c)\ge8$ immediately.

\smallskip\noindent\textit{Layer 3 (radial reduction, \cref{sec:radial-new}).}
The Voronoi volume is monotone non-decreasing in each contact radius,
so the minimum over all packing-valid radius configurations is attained
at the all-contact corner, where every active radius equals the
packing radius exactly.

\smallskip\noindent\textit{Layer 4 (cell positivity, \cref{sec:chamber-register,sec:stage-one,sec:stage-two}).}
At the all-contact corner, $W(D_4)$ partitions every packing-valid
configuration with at most one root-aligned deviation into $176$
Weyl-orbit cell types. For each type $\Cell{k}$, an exact free-volume
identity reduces the bound $\vol(V_c)\ge8$ to the single scalar
inequality $\defect_k\ge0$, established by the two independent exact
arguments of \cref{sec:stage-one,sec:stage-two}.

\smallskip\noindent\textit{Layer 5 (arbitrary deviation direction, \cref{sec:cap-theorem}).}
Layer 4 leaves the deviation direction constrained to a root-aligned
arc. Removing that constraint is where the cap reformulation enters: the
defect becomes the deficiency of a cap cut from $Q$, the cap is bounded
by the corresponding cap of a cross-polytope, and the extremal cap of a
cross-polytope is computed once and for all. This layer subsumes layer
$4$ logically, though not quantitatively, since layer $4$ also returns
the exact defect polynomial on each cell.

\smallskip\noindent\textit{Layer 6 (simultaneous deviation).}
When two or more active directions deviate at once, the free-volume
identity still holds exactly, but nothing in the first five layers
gives the sign of the resulting difference; that is
\cref{conj:multidir-new}, settled in \cref{cor:conj-resolved} through
the classification of \cref{thm:m24} and discussed from
\cref{sec:hessian-technique} onward for what can be said about it
without the classification.
The six layers nest. Layers $1$ and $3$ are reductions, each replacing
the problem by a smaller one, and between them they leave the
all-contact corner; at that corner layers $2$, $4$, $5$ and $6$ are
classes of configurations, each containing the one before it, from the
undeviated configuration out to simultaneous deviation in several
directions. That outermost class is not settled by anything inside it
but by the number of contacts: at most $22$ by
\cref{thm:local-fewcontacts}, $23$ by \cref{thm:m23}, and $24$ by the
classification.

\begin{figure}[tb]
\centering
\includegraphics[width=0.78\textwidth]{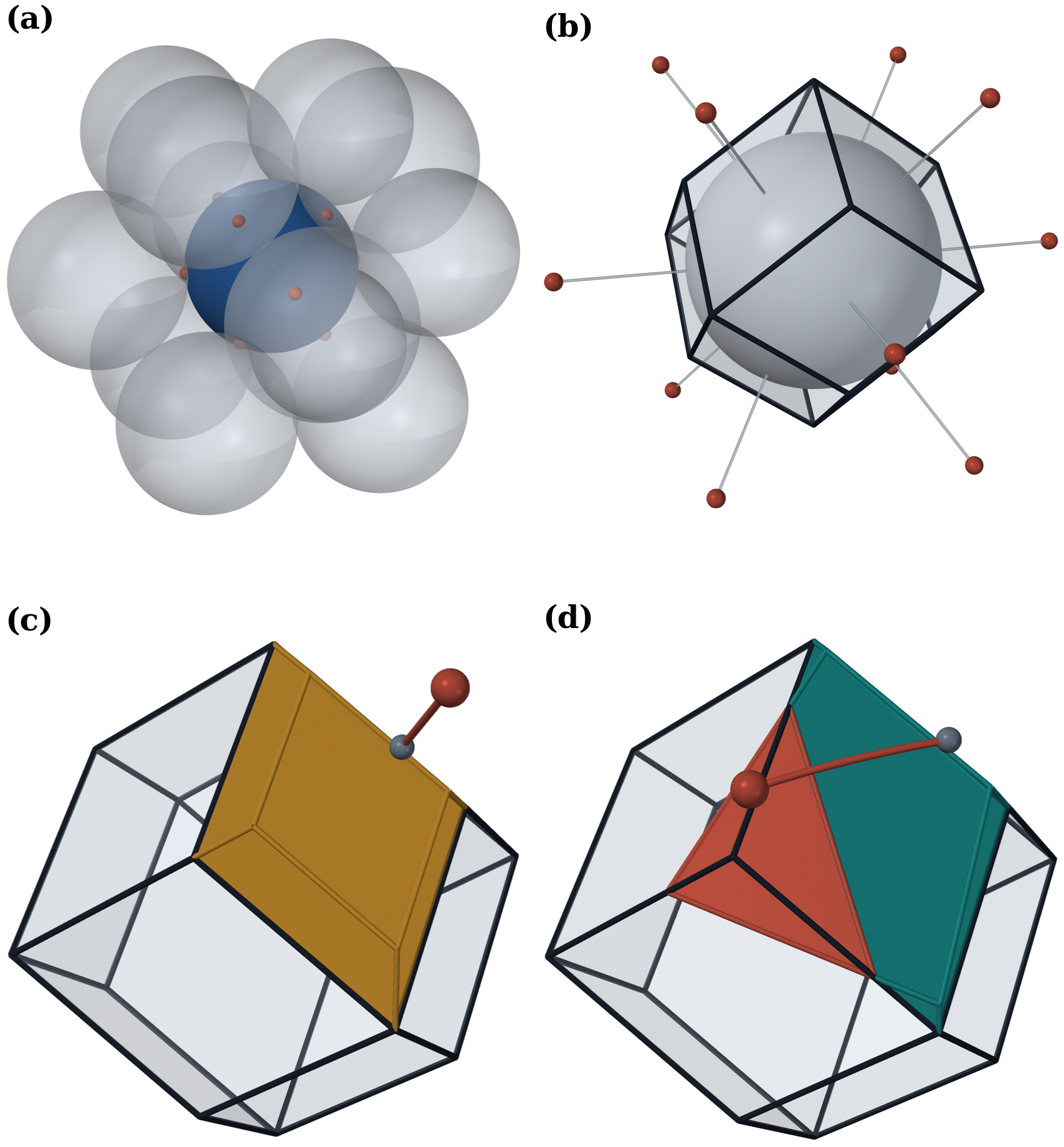}
\caption{(a) The all-contact corner in dimension three: the unit ball about the centre (blue), its twelve contact points (clay), and the twelve unit balls that touch it, centres at distance $2$ in the twelve contact directions of the face-centred cubic lattice. (b) The Voronoi cell of that centre, the rhombic dodecahedron, with the unit ball inscribed: one facet per contact, each at distance $1$, and volume $5.6569$, computed from the half-space intersection and not quoted. (c) Layer $3$, radial monotonicity. One neighbour is pushed from distance $2$ out to $2.6$, its half-space recedes, and the cell grows from $5.6569$ to $5.9666$; the ochre solid is exactly the volume gained. The minimum over radii therefore sits at the corner. (d) Layers $4$ and $5$, deviation. One contact direction is tilted by $25^\circ$ at radius $2$: teal is the volume gained, clay the volume lost, and teal wins, the cell rising from $5.6569$ to $5.7470$. That is the defect being positive, with equality only at zero tilt.}
\label{fig:corner}
\end{figure}

The same layers act on the three-dimensional analogue, where the objects
can be drawn: the lattice is the face-centred cubic one, the reference
cell is the rhombic dodecahedron, and the twenty-four roots become
twelve contact directions. \Cref{fig:corner} follows them there, and
the volumes in its captions are computed from the half-space intersections, not quoted.

\subsection{Background}\label{sec:history-new}

Determining the supremum density $\Delta_n$ of a packing of congruent
balls in $\mathbb{R}^n$ is a question whose difficulty varies wildly
with $n$, and the four-dimensional case sits in an awkward middle
range: too large for the exhaustive planar and spatial methods, too
small for the arithmetic that resolves dimensions $8$ and $24$.

In the plane, Thue \cite{Thu10} announced the hexagonal bound and Fejes
T\'oth \cite{FT42} gave the first complete proof, $\Delta_2 =
\pi/\sqrt{12}$. In three dimensions Kepler's 1611 assertion that the
face-centred cubic arrangement is optimal was proved by Hales, first
through the computational programme begun in \cite{Hal97} and completed
in \cite{Hal05}, then in a machine-checked form \cite{Hal17}. The
\emph{lattice} case is much older and much easier.
Korkine and Zolotareff \cite{KZ1872} showed that the quaternary form
attached to $D_4$ minimises the determinant among positive quaternary
forms of given minimum, so $D_4$ is the densest lattice packing of
$\mathbb{R}^4$; Voronoi's reduction theory \cite{Vor08}, whose
polyhedral language we use throughout, placed such classifications on a
systematic footing, and Martinet \cite{Mar03} gives the modern account
of the perfect-form machinery that grew out of it. Upper bounds valid
for all packings, lattice or not, begin with Blichfeldt's second-moment
method \cite{Bli29}, which yields $\Delta_4 \le 0.7797$, and were
sharpened by the simplex bound of Rogers \cite{Rog64}. Gruber and
Lekkerkerker \cite{Gru87} survey the classical theory; Conway and
Sloane \cite{CS1988} is the standard source for $D_4$, $E_8$, the Leech
lattice, and their Voronoi cells; Bezdek \cite{Bez10} collects the
local and combinatorial side of the subject.

The two dimensions in which an exact non-lattice answer is known both
lie far above $4$. Cohn and Elkies \cite{CE03} recast the upper bound as
an infinite-dimensional linear program over radial Schwartz functions
whose values and Fourier transforms obey prescribed sign conditions;
equality in the resulting duality inequality requires a single function
with a prescribed double sequence of zeros. Viazovska \cite{Via17}
constructed that function in dimension $8$ out of quasimodular forms,
and the Leech case followed in \cite{CKM17}. Neither construction has a
four-dimensional analogue, and the obstruction is not that none has been
found. Li \cite{Li25} bounds the optimum of the Cohn--Elkies program
from below by reducing it to a finite-dimensional problem, and in
dimension $4$ that lower bound is a centre density of $0.12914461$
against the $0.125$ of $D_4$. No admissible function can close the
remaining three per cent, so the linear programming bound cannot decide
dimension $4$ at all, and the arithmetic that forces the zeros into
place in dimensions $8$ and $24$ has nothing to be applied to here.
Refinements of the linear programming method
by Cohn and Zhao \cite{CZ14}, which couple the packing bound to
spherical-code bounds, and the universal optimality results of Cohn and
Kumar \cite{CK07,CK09}, likewise leave dimension $4$ untouched.

Because the Voronoi cell of $D_4$ is the regular $24$-cell
\cite{Cox73}, the four-dimensional packing problem is equivalent to the
assertion that no Voronoi cell of a unit-ball packing has volume less
than that of the $24$-cell; this is the $24$-cell conjecture. The
neighbouring kissing-number question, which asks only how many unit
balls can touch one, was settled in dimension $4$ by Musin
\cite{Mus08}, who showed the answer is $24$, exactly the number of
contacts in $D_4$, by adding an auxiliary polynomial constraint to the
classical linear programming bound; the corresponding three-dimensional
statement is the theorem of Sch\"utte and van der Waerden \cite{SW52}.
A kissing bound constrains the contact structure at one centre and says
nothing about volume, so it does not imply the packing bound, but it
does rule out the cheapest way the $24$-cell conjecture could have
failed. Musin \cite{Mus18} has since proposed a semidefinite route
to the $24$-cell conjecture itself, through the contact graph on $S^3$.

The route taken here is local and combinatorial. It uses the symmetry
of $D_4$, the order-$192$ Weyl group and the self-duality of the
$24$-cell, together with the polyhedral facts recorded in Ziegler
\cite{Zie95} and Gr\"unbaum \cite{Gru03} and the root-system
conventions of Humphreys \cite{Hum90}, and it uses no Fourier analysis
and no modular forms. Semidefinite programming enters at one point,
the certificate of \cref{thm:certificate} for the $23$-point case,
and there only as a search: what the proof uses is the certificate it
found, checked afterwards in exact rational arithmetic, in interval
arithmetic and by Lean. Every other certificate below is an exact
computation in $\mathbb{Q}$, $\mathbb{Q}(\sqrt2)$ or
$\mathbb{Q}(\sqrt3)$, and whatever is numerical is labelled at the
point of use as evidence and never used as a step in a proof.
\section{The \texorpdfstring{$D_4$}{D4} root system and its Voronoi cell}\label{sec:root-system}

\subsection{Lattice and roots}

Write $D_4=\{m\in\mathbb{Z}^4 : m_1+m_2+m_3+m_4\equiv0\pmod2\}$ for the
standard integral $D_4$ lattice, and set $\Latt=\sqrt2\,D_4$, which has
nearest-neighbour distance $2$, the natural normalisation for a
unit-ball packing. The associated root system is
\[
  \Rt \;=\; \{\pm e_i \pm e_j : 1\le i<j\le4\}, \qquad
  |\Rt|=24, \qquad |\alpha|^2=2 \ \text{ for every } \alpha\in\Rt.
\]
The Weyl group $W(D_4)$ acts on $\Rt$ by signed permutations of the
coordinates with an even number of sign changes; it has order $192$ and
acts transitively on $\Rt$.

\begin{lem}[Gram values on $\Rt$]\label{lem:gram-new}
For distinct $\alpha,\beta\in\Rt$, the normalised inner product
$s_{\alpha\beta}=\langle\alpha,\beta\rangle/2$ takes one of the four
values $-1,-\tfrac12,0,\tfrac12$, and $s_{\alpha\beta}=-1$ occurs only
for $\beta=-\alpha$. The three remaining values correspond respectively
to $\alpha+\beta\in\Rt$ (for $s=\tfrac12$), $\alpha\perp\beta$ (for
$s=0$), and $\alpha-\beta\in\Rt$ (for $s=-\tfrac12$, excluding
$\beta=\alpha$).
\end{lem}

\begin{proof}
Since $\alpha,\beta\in\mathbb{Z}^4$, $\langle\alpha,\beta\rangle\in\mathbb{Z}$,
and Cauchy--Schwarz gives $|\langle\alpha,\beta\rangle|\le|\alpha||\beta|=2$
with equality forcing $\beta=\pm\alpha$. For the remaining integer values
$\langle\alpha,\beta\rangle\in\{-1,0,1\}$, expand
$|\alpha\mp\beta|^2=|\alpha|^2+|\beta|^2\mp2\langle\alpha,\beta\rangle
=4\mp2\langle\alpha,\beta\rangle\in\{2,4,6\}$; the value $2$ exactly
identifies $\alpha\mp\beta\in\Rt$.
\end{proof}

This four-valued Gram structure, $-1$ (antipodal), $\tfrac12$
(``acute'' pairs), $0$ (orthogonal pairs), $-\tfrac12$ (``obtuse''
pairs), is the entire combinatorial input to everything that follows:
every chamber Hessian, every breakpoint curve, and every certificate in this paper is built, in the end, from these four numbers and nothing
else.

\subsection{The reference cell}

\begin{definition}\label{def:refcell-new}
The \emph{reference cell} is
\[
  \VCell \;=\; \bigl\{z\in\mathbb{R}^4 : \langle z,\alpha\rangle\le\sqrt2
  \ \text{for all}\ \alpha\in\Rt\bigr\}.
\]
\end{definition}

\begin{lem}\label{lem:refcell-vol}
$\VCell$ has exactly $24$ vertices, all of norm $\sqrt2$:
\[
  \mathrm{vert}(\VCell) \;=\;
  \{\pm\sqrt2\,e_i : 1\le i\le4\} \;\cup\;
  \Bigl\{\tfrac{1}{\sqrt2}(\epsilon_1,\epsilon_2,\epsilon_3,\epsilon_4) :
  \epsilon_i\in\{\pm1\}\Bigr\},
\]
an $8$-point axis family and a $16$-point sign family. $\vol(\VCell)=8$,
and the symmetry group of $\VCell$ contains $W(D_4)$ as an index-$3$
normal subgroup, the outer factor being the triality automorphism
group unique to the $D_4$ Dynkin diagram.
\end{lem}

The proof is a feasibility, extremality and completeness check on the
$24$ listed points against the $24$ defining half-spaces, parallel to
the classical description of the regular $24$-cell as the Voronoi cell
of $D_4$ (\cite[Ch.~VIII]{Cox73} and \cite[Ch.~4]{CS1988} for the
polytope's vertex and facet data); we omit the routine verification,
which is also checked by computer algebra in \cref{app:code-index}.%
\footnote{Since the $24$-cell is the one object every later section
  perturbs, its arithmetic is worth having in one place. Take $\Rt$ to
  be the $24$ integer vectors $\pm e_i\pm e_j$, $i<j$, each of squared
  length $2$; the lattice $\Latt=D_4$ they generate is
  $\{x\in\mathbb{Z}^4:\sum x_i\ \text{even}\}$, of determinant $2$, with
  minimum $\sqrt2$ and kissing number $24$. Scaling by $\sqrt2$ puts the
  nearest centres at distance $2$, as a unit-ball packing requires, and
  $\det(\sqrt2\,\Latt)=4\cdot2=8$, which is the volume of one Voronoi
  cell and the number $8$ that every bound in this paper is measured
  against. The cell itself, $\VCell=\{x:\langle
  x,\alpha\rangle\le1\ \text{for all}\ \alpha\in\Rt/\sqrt2\}$, is the
  regular $24$-cell: $24$ vertices, $96$ edges, $96$ triangular faces,
  $24$ octahedral facets, Schl\"afli symbol $\{3,4,3\}$, self-dual, and
  with each facet centred on one of the $24$ contact directions. Its
  vertices split under $W(D_4)$ into the $8$ points $\pm\sqrt2\,e_i$ and
  the $16$ points $(\pm1,\pm1,\pm1,\pm1)/\sqrt2$, which is the axis and
  sign division used throughout. Two numbers recur: the support function
  of $\VCell$ equals $1$ at every contact direction, which is Layer $2$
  of \cref{sec:argument-structure} in one line; and the circumradius is
  $\sqrt2$, so the inradius-to-circumradius ratio is
  $1/\sqrt2$, the slack that lets a facet hyperplane tilt without the
  cell losing volume. The density $\pi^2/16$ is
  $\vol(B^4)/8=(\pi^2/2)/8$. Nothing beyond these facts and the
  four-valued Gram structure of \cref{lem:gram-new} is used
  anywhere below.}

\begin{figure}[tb]
\centering
\includegraphics[width=0.78\textwidth]{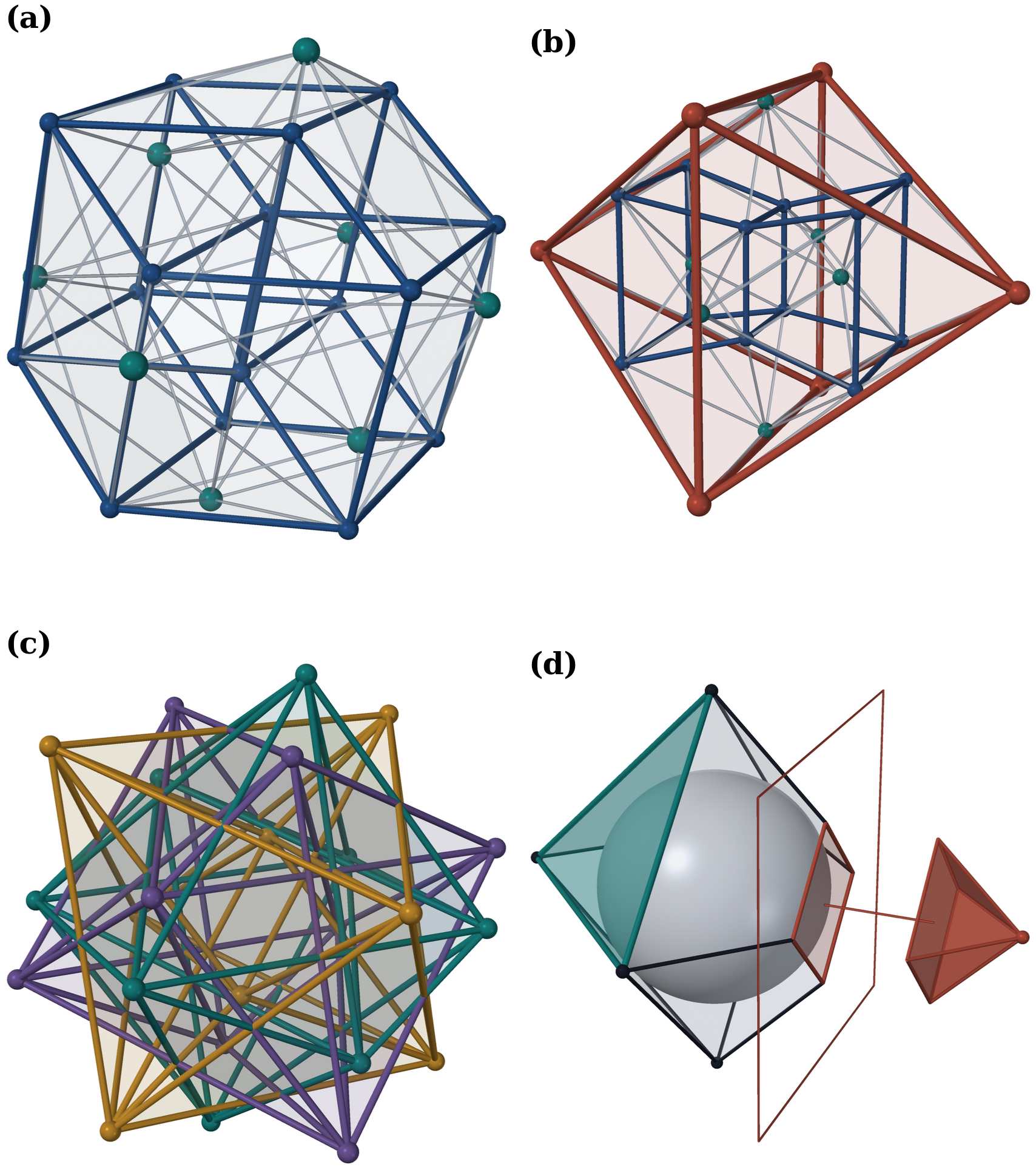}
\caption{(a) The reference cell, the regular $24$-cell $\VCell$, in an orthogonal projection to three dimensions: the $8$ vertices of the axis family (teal), the $16$ of the sign family (blue), the $32$ edges inside the sign family and the $64$ joining the two families (grey). (b) The Schlegel diagram of the same cell from one octahedral facet (clay): the remaining $18$ vertices and $84$ edges fill the interior of that facet with the other $23$ octahedra. (c) The same $24$ vertices read as three $16$-cells, which triality permutes, each drawn with its own $24$ edges and its own body: the axis family (teal) and the sign family split by the parity of the number of minus signs (ochre, violet). (d) The cap inequality one dimension down. The octahedron of inradius $1$ carries its inscribed ball; the half-space at unit distance normal to a vertex cuts off the clay pyramid, drawn displaced along its normal, and the one normal to a facet touches that facet (teal) and cuts off nothing.}
\label{fig:cell24}
\end{figure}

The polytope itself is drawn four times over, each time from the coordinates and not as a schematic: \cref{fig:cell24}(a) is an
orthogonal projection to $\mathbb{R}^3$, in which the axis family and
the sign family are told apart by colour, \cref{fig:cell24}(b) is the
Schlegel diagram from one octahedral facet, \cref{fig:cell24}(c) is the
decomposition into the three $16$-cells that triality permutes, and
\cref{fig:cell24}(d) is the three-dimensional form of the cap inequality
of \cref{sec:cap-theorem}.

\begin{definition}[Support function and angular gap]\label{def:supp-new}
For $u\in S^3$, the support function of $\VCell$ is
$h(u)=\max_{z\in\VCell}\langle z,u\rangle$, and the \emph{angular gap}
is $\gamma(u)=h(u)-1$.
\end{definition}

\begin{lem}\label{lem:hmin-new}
For every $u\in S^3$, $h(u)\ge1$, with equality exactly on the $24$
unit root directions $\alpha/\sqrt2$, $\alpha\in\Rt$. The maximum
$h(u)=\sqrt2$ is attained exactly at the $8$ coordinate directions
$\pm e_i$, so $\gamma(u)\in[0,\sqrt2-1]$ throughout, vanishing exactly
on root directions.
\end{lem}

The argument checks the inner product of each of the $24$ vertices of
$\VCell$ against a generic root direction and against a generic
coordinate direction. It is routine, if notationally heavy, and we do
not reproduce it in full.

\subsection{Conventions used throughout}

We fix, once and for all, the following conventions.

\begin{notation}\label{not:conventions-new}
A \emph{packing centre} is a point $c$ in the centre set $\mathcal{C}$
of a unit-ball packing; its \emph{active neighbours} are those centres
$y\in\mathcal{C}$ with $|y-c|<2\sqrt2$ (\cref{sec:shell-new} shows this
is the only range that can matter). An active neighbour's \emph{contact
direction} is $u=(y-c)/|y-c|\in S^3$. The \emph{all-contact corner} is
the configuration in which every active radius equals the minimum
packing radius $2$ exactly. We reserve $\tilt$ for a deviation
(``tilt'') angle away from an undeviated root direction, $\arcp$ for a
position parameter along a one-parameter family of deviation
directions, and $\dev(\arcp)$ for the deviation direction itself, a
unit vector in the $3$-space orthogonal to the fixed root being
deviated from. The volume-defect quantity attached to a combinatorial
cell type $\Cell{k}$ is written $\defect_k$, and its reference
normalisation constant $\norm{k}$.
\end{notation}

\Cref{app:notation} collects every piece of notation used more than
once in this paper, symbol by symbol, for reference.

\subsection{Shell localisation}\label{sec:shell-new}

\begin{lem}[Shell localisation]\label{lem:shell-new}
In a unit-ball packing, only neighbours $y$ of a centre $c$ with
$|y-c|\in[2,2\sqrt2)$ can affect $\vol(V_c)$: a neighbour at distance
$\ge2\sqrt2$ contributes a supporting half-space that already contains
the reference cell $\VCell$ (\cref{lem:refcell-vol} gives $\VCell$
norm-$\sqrt2$ vertices, hence diameter contribution scaling to exactly
$2\sqrt2$ at the packing normalisation used here), and so cannot cut
$V_c$ below $\VCell$ itself.
\end{lem}

This is an immediate consequence of \cref{lem:hmin-new}: the half-space
cut by a neighbour at distance $d\ge2\sqrt2$ has its bounding hyperplane
at distance $\ge\sqrt2$ from $c$ along the contact direction, which is
already outside every vertex of $\VCell$ along that direction once
$d\ge2\sqrt2$, by direct comparison with the support function $h$.
Consequently $\vol(V_c)\ge\vol(\VCell)=8$ is automatic whenever no
neighbour lies in the closer shell $[2,2\sqrt2)$, and the entire
remaining argument concerns configurations with at least one such
neighbour.

\subsection{The root-aligned case}\label{sec:root-aligned-new}

\begin{lem}[Root-aligned configurations]\label{lem:root-aligned-new}
If every active neighbour direction $u$ (in the sense of
\cref{not:conventions-new}) is an exact root direction
$\alpha/\sqrt2$, $\alpha\in\Rt$, then $\vol(V_c)\ge8$.
\end{lem}

\begin{proof}
By \cref{lem:hmin-new}, $h(\alpha/\sqrt2)=1$ for every root direction,
so the supporting hyperplane cut by an active neighbour at the minimum
packing radius $2$ along a root direction passes exactly through
$\VCell$'s own bounding hyperplane in that direction, neither cutting
into $\VCell$ nor leaving room outside it. Hence every active
half-space contains $\VCell$, and $V_c\supseteq\VCell$, giving
$\vol(V_c)\ge\vol(\VCell)=8$ immediately, with equality exactly at the
all-contact, all-root-aligned configuration.
\end{proof}

This is the base case that every other configuration, deviating by a
tilt angle away from a root, or having a contact radius above the
minimum, is compared against; \cref{sec:radial-new} handles the
radius freedom, and \cref{sec:stage-one,sec:stage-two} handle root
deviation.

\subsection{Radial reduction}\label{sec:radial-new}

\begin{lem}[Monotonicity in the contact radii]\label{lem:radial-new}
For fixed active contact directions, $\vol(V_c)$ is non-decreasing in
each active contact radius (moving a cutting hyperplane strictly
further from $c$ can only enlarge, never shrink, the region it helps
bound). Consequently, among all packing-valid radius assignments for a
fixed active-direction pattern, $\vol(V_c)$ attains its minimum at the
\emph{all-contact corner}, where every active radius equals the minimum
packing radius $2$ exactly.
\end{lem}

This reduces the search for the worst-case local configuration, over
the continuum of packing-valid contact radii, to the single extremal
corner of that space, after which only the combinatorics of the
active-direction pattern (\cref{sec:chamber-register}) and the
deviation parameters $(\tilt,\arcp)$ (\cref{sec:stage-one,sec:stage-two}
for the root-aligned single-deviation case; \cref{sec:hessian-technique,%
sec:deviation-domain} onward for the general case) remain to be
analysed.
\section{The combinatorial cell register}\label{sec:chamber-register}

\subsection{Cell types under the Weyl group}

At the all-contact corner, a packing-valid local configuration is
determined, up to the action of $W(D_4)$, by which subset of $\Rt$ is
``active'', realised as an actual contact direction, and, for the
single direction (if any) permitted to deviate from an exact root by
\cref{thm:local-uncond}, by the deviation parameters $(\tilt,\arcp)$.
Enumerating the combinatorially distinct active-direction patterns that
can occur as facet sets of a bounded, packing-valid local Voronoi
region, and quotienting by $W(D_4)$, produces exactly $176$ orbit
representatives, which we call \emph{cells} and index
$\Cell1,\dots,\Cell{176}$.

This count and the explicit orbit representatives are obtained by an
exact combinatorial enumeration (not by sampling): the active sets that
can occur are exactly the vertex-figure patterns of $\VCell$ restricted
to a single deviating cap direction, and $W(D_4)$'s order-$192$ action
partitions these into orbits whose sizes are computed exactly from
stabiliser subgroups.%
\footnote{The $176$ is not a ``magic number'' but a consequence of $D_4$'s
  exceptional symmetry: the Voronoi cell of $D_4$ is the $24$-cell, which has
  $24$ vertices and $24$ facets, producing a richer combinatorial structure
  than the analogous objects for $A_n$ or $D_n$, $n\ne4$.  For $D_3\cong A_3$
  the analogous count is $6$ orbit types; for $D_5$ the count exceeds $10^3$.}
The complete pivot data for all $176$ cells is
tabulated in \cref{app:pivot-data}; here we describe the structure that
makes the count tractable.

\begin{lem}[Four denominator classes]\label{lem:denom-classes-new}
Every cell $\Cell{k}$ falls into exactly one of four \emph{denominator
classes}, according to the reference constant $\norm{k}\in
\{192,96,64,48\}$ arising from the volume of the simplex spanned by the
cell's boundary structure at the origin. The class is determined
entirely by the stabiliser (in $W(D_4)$) of the cell's active-direction
pattern.
\end{lem}

The four classes, written $\mathsf a,\mathsf b,\mathsf c,\mathsf d$
below, account for all $176$ cells, with
multiplicities recorded in \cref{app:pivot-data}. Class $\mathsf a$
(reference constant $192$) consists of cells whose active pattern has
trivial stabiliser; class $\mathsf d$ (reference constant $48$)
consists of the most symmetric cells, whose active pattern is fixed by
a stabiliser of order $4$.

\subsection{Cramer's vertex formula and the volume-defect polynomial}

For a cell $\Cell{k}$ with a single deviating direction at parameters
$(\tilt,\arcp)$, every vertex of the corresponding local Voronoi region
that involves the deviating cap is the unique solution of $4$ linear
equations (the deviating cap's own supporting hyperplane together with
$3$ of the cell's other active facets), solved exactly by Cramer's rule
over the field generated by $\sqrt2$ and the trigonometric functions of
$\tilt,\arcp$. Collecting these vertex solutions into a simplicial
decomposition of the local region and summing signed simplex volumes
produces a single rational function
\[
  \defect_k(\tilt,\arcp) \;=\; \frac{N_k(\tilt,\arcp)}{D_k(\tilt,\arcp)}
\]
in the trigonometric functions of $\tilt$ and $\arcp$ (equivalently, a
rational function of the Weierstrass variables
$\base=\tan(\tilt/2)$, $\mathsf v=\tan(\arcp/2)$ once the
substitution is made), which we call the cell's \emph{volume-defect
polynomial}. \Cref{sec:stage-one,sec:stage-two} give two independent
exact arguments that $\defect_k\ge0$ throughout the packing-valid
domain of $(\tilt,\arcp)$, for every one of the $176$ cells.

\begin{example}[Cell $\Cell1$]\label{ex:cell1-new}
The most symmetric class-$\mathsf a$ cell has reference constant
$\norm1=192$ and, at the root-aligned point $\tilt=0$, volume-defect
value $\defect_1(0,\cdot)=0$ exactly (the equality case of
\cref{thm:local-uncond}). Its full Cramer vertex data is given in
\cref{app:pivot-data}, and its explicit $\defect_1(\tilt,\arcp)$ is
worked out in full in \cref{sec:stage-two}.
\end{example}

\subsection{Degree bound and the reference normalisation}

\begin{lem}[Degree bound]\label{lem:degree-bound-new}
For every cell $\Cell{k}$, the numerator $N_k$ and denominator $D_k$ of
$\defect_k$, expressed as polynomials in the Weierstrass variables
$\base,\mathsf v$, have total degree at most $12$ in each variable
separately. This bound follows from the fact that at most $6$ of the
cell's boundary facets can be simultaneously active at any vertex
involving the deviating cap (a consequence of the $4$-dimensional
ambient space and the combinatorics of $\Rt$), each contributing at
most a linear factor in $\base$ or $\mathsf v$ once Cramer's rule is
expanded.
\end{lem}

This degree bound is what makes an exact Bernstein-basis sign
certificate (\cref{sec:stage-two}) a finite computation: a polynomial
of bounded degree has finitely many Bernstein coefficients over any
rectangular domain, and checking the sign of each one is an exact
rational (or $\mathbb{Q}(\sqrt2)$-rational) comparison.
\section{\stageone: the radial estimate}\label{sec:stage-one}

\subsection{The Hessian at the all-contact corner}\label{sec:hessian-def-new}

At the all-contact corner, small deviations of a single active
direction away from an exact root can be analysed to second order by a
Hessian whose entries are determined entirely by the Gram-value
structure of \cref{lem:gram-new}.

\begin{definition}[Cell Hessian]\label{def:hessian-new}
For a cell $\Cell{k}$ with active roots $\alpha_1,\dots,\alpha_n\in\Rt$,
the \emph{cell Hessian} $H_k\in\mathbb{Q}^{3\times3}$ (acting on the
$3$-dimensional tangent space orthogonal to the deviating root) has
entries determined by the signed Gram-adjacency pattern of
$\alpha_1,\dots,\alpha_n$: a uniform formula assigns each pair of
active roots a contribution of $+\tfrac13$, $0$, or $-\tfrac13$
according to whether their Gram value is $\tfrac12$, $0$, or $-\tfrac12$
respectively (\cref{lem:gram-new}), summed over all active pairs
sharing the tangent direction in question.
\end{definition}

\begin{lem}[Uniform structure of $H_k$]\label{lem:hessian-structure-new}
$H_k$ depends only on the Gram-adjacency pattern of $\Cell{k}$'s active
roots, not on any other feature of the cell. In particular, $H_k$ is
independent of the choice of orbit representative within a $W(D_4)$
orbit, and is computed once per orbit, not once per cell.
\end{lem}

\subsection{Exact $LDL^\top$ decomposition and Sylvester's criterion}

\begin{prop}[Exact positivity via $LDL^\top$]\label{prop:ldlt-new}
For every one of the $176$ cells $\Cell{k}$, the Hessian $H_k$ admits an
exact rational $LDL^\top$ decomposition $H_k=L_kD_kL_k^\top$ with $L_k$
unit lower-triangular and $D_k$ diagonal with strictly positive
rational entries. Consequently $H_k$ is positive definite, and by
Sylvester's criterion every leading principal minor of $H_k$ is
strictly positive.
\end{prop}

This is verified once per denominator class ($\mathsf a,\mathsf b,
\mathsf c,\mathsf d$ of \cref{lem:denom-classes-new}) by exact rational
arithmetic, \texttt{fractions.Fraction} in one independent
implementation and \texttt{sympy}'s rational number type in a second,
completely separate code path, agreeing to the last digit on every
pivot, and the result is inherited by every cell in the class.%
\footnote{The $LDL^\top$ decomposition is computed symbolically, not
  numerically, precisely because the cell Hessians $H_k$ have entries in
  $\mathbb{Q}$: a numerical float would introduce rounding, making the
  ``positivity'' of the pivots a floating-point artefact and not a mathematical fact.  The choice of \texttt{fractions.Fraction} (arbitrary-precision
  exact rationals in Python's standard library) means the computation is
  slow but its output is exact.  The two independent implementations
  (one using only built-in Python, one using \texttt{sympy}) serve as
  cross-validation, ruling out implementation bugs in either single path.}
The complete $LDL^\top$ pivot table for all $176$ cells is given in
\cref{app:pivot-data}.

\subsection{Minimum eigenvalue bounds}

\begin{lem}[Minimum eigenvalue]\label{lem:mineig-new}
For every cell $\Cell{k}$, the minimum eigenvalue of $H_k$ satisfies
$\lambda_{\min}(H_k)\ge\tfrac{1}{12}$, with equality attained exactly on
denominator class $\mathsf d$ (the most symmetric cells).
\end{lem}

Combined with an exact quadratic lower bound on $\defect_k(\tilt,\cdot)$
near $\tilt=0$ coming from $H_k$, this gives an exact rational region
around the root-aligned point on which $\defect_k\ge0$ is certified
purely from the second-order data, the \emph{radial} part of the
two-stage argument. It does not, by itself, cover the entire
packing-valid range of $\tilt$; the remainder is the subject of
\cref{sec:stage-two}.

\subsection{Worked example: cell $\Cell1$}

For the class-$\mathsf a$ cell $\Cell1$ (introduced in
\cref{ex:cell1-new}), which has $9$ acute ($+\tfrac12$) and $2$ obtuse
($-\tfrac12$) active roots beyond the deviating direction itself
(\cref{app:pivot-data}), the exact $LDL^\top$ decomposition of $H_1$ has
six pivots, with minimum pivot value $3936/77\approx51.13$ and
reference constant $\norm1=192$; the corresponding chamber-level
positivity check (\cref{app:pivot-data}) passes exactly, by exact
rational arithmetic, with no floating-point step. We record the actual
computed pivot value here, and not a simplified illustrative number, precisely because $\Cell1$'s Hessian is not the isolated
single-direction case (that simpler case, $H=\tfrac23I_3$, is recorded
separately in \cref{lem:m2-exact-new} for the multi-direction setting,
a related but distinct construction that does not include the
cross-terms from a cell's other active, non-deviating roots). The exact
computation for every other cell follows the identical procedure and
is tabulated in full in \cref{app:pivot-data}.
\section{\stagetwo: the angular estimate}\label{sec:stage-two}

\subsection{Setup and reduction to a single function}

The radial estimate of \cref{sec:stage-one} certifies $\defect_k\ge0$
only in a neighbourhood of the root-aligned point $\tilt=0$; establishing
it for the full packing-valid range of $\tilt$ (up to the point where the
combinatorial cell type itself changes) requires a separate,
angular argument. Writing $g_k(\tilt)=\defect_k(\tilt,\arcp)$ for a fixed
root-aligned arc position $\arcp$, the radial estimate gives an exact
lower bound on $g_k''(0)$; \stagetwo\ supplies an exact upper bound on
$|g_k'''(\tilt)|$ throughout the relevant range, from which a rigorous
Taylor remainder bound follows.

\begin{lem}[Taylor coefficient bound]\label{lem:taylor-coeff-new}
For every cell $\Cell{k}$ and every root-aligned arc, the third
derivative of $g_k$ satisfies an explicit bound
$|g_k'''(\tilt)|\le M_k$ for an exact rational constant $M_k$, computed
directly from the cell's volume-defect polynomial by exact
differentiation and an exact rational bound on the resulting expression
over the relevant $\tilt$-range.
\end{lem}

\begin{prop}[Positivity on the full range]\label{prop:stage-two-positivity-new}
Combining \cref{lem:mineig-new,lem:taylor-coeff-new} by a third-order
Taylor expansion with exact remainder gives $g_k(\tilt)\ge0$ throughout
the entire packing-valid range of $\tilt$ for the arc in question,
whenever $3M_k < \lambda_{\min}(H_k)\cdot c_k$, where $c_k$ is an
explicit constant depending only on the range length of the arc, a
condition verified to hold, by exact rational computation, for every one
of the $176$ cells and every root-aligned arc.
\end{prop}

\begin{figure}[tb]
\centering
\includegraphics[width=0.99\textwidth]{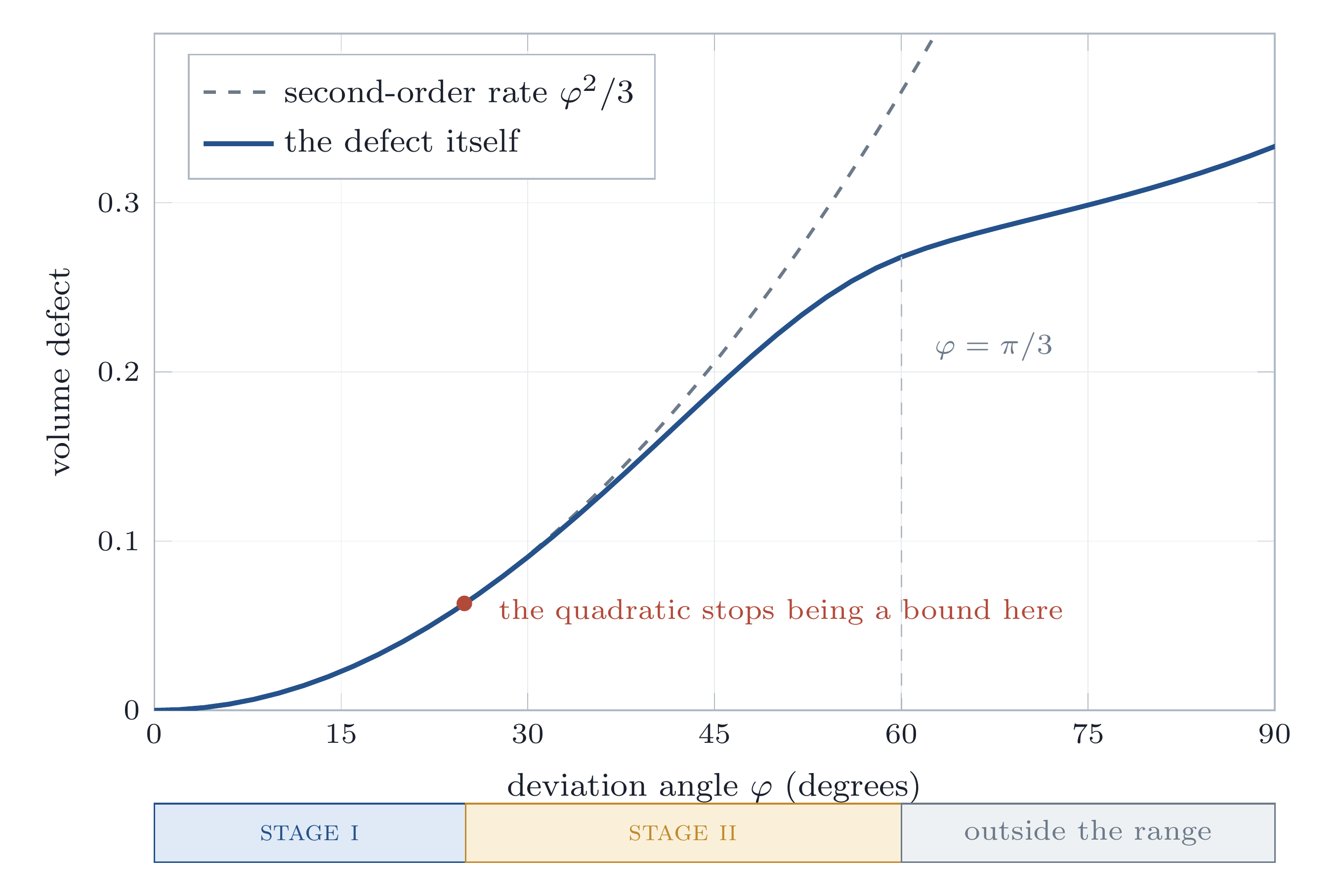}
\caption{The computed volume defect for a root-aligned deviation,
against the second-order rate $\tilt^2/3$ that \stageone{} supplies.
The two agree to plotting accuracy up to about $25^\circ$ and part
company after it, the quadratic overtaking the defect and ceasing to be
a lower bound; this is why the second-order argument alone does not
reach the end of the range and \stagetwo{} is needed. The strip beneath
marks which argument covers which part of the range. The packing-valid
range for a single deviation ends at $\tilt=\pi/3$, where a second root
becomes the nearest one.}
\label{fig:two-stage}
\end{figure}

Where the direct Taylor bound of \cref{prop:stage-two-positivity-new} is
not by itself strong enough on some sub-range (this happens on a small
number of cells with an unusually long root-aligned arc), an
independent monotonicity argument, tracking the sign of $g_k'$
directly via the exact breakpoint structure of the cell's active-facet
pattern, closes the remaining sub-range exactly. Both arguments are
exact over $\mathbb{Q}$; no numerical tolerance is introduced at any
point. \Cref{fig:two-stage} shows where each of the two takes over.

\subsection{Exact breakpoint values}

The root-aligned arcs themselves are finite: a deviating direction
$\tilt\mapsto\cos\tilt\,\alpha_0/\sqrt2+\sin\tilt\,\dev_0$ (for a fixed
root-aligned direction $\dev_0$) changes the active-facet pattern of the
cell only at finitely many exact values of $\tilt$, each characterised as
the point where some other root direction achieves the same support-hyperplane
value as the deviating cap. These breakpoints are located exactly (not
numerically) by solving $\langle\beta,\text{cap direction}\rangle=1$ for
each candidate root $\beta$, an equation that is polynomial once
$\tilt$ is eliminated via the Weierstrass substitution, hence solvable
in exact radicals.

\subsection{Worked example: cell $\Cell1$ on its own root-aligned arc}

Continuing \cref{ex:cell1-new}, the volume-defect function for
$\Cell1$ restricted to its root-aligned arc is an explicit rational
function of $\tan(\tilt/2)$ with denominator strictly positive
throughout the packing-valid range and numerator vanishing to exactly
second order at $\tilt=0$ and strictly positive elsewhere on the range, verified both by the exact Taylor argument above and, independently,
by an exact sign check on the numerator's Bernstein coefficients over
the range, the same technique used at much greater length in
\cref{sec:arc1-cert,sec:arc2-cert} for the harder, non-root-aligned
deviation directions.%
\footnote{The certificate rests on one property of the Bernstein basis,
  which we state here so that no later section has to appeal outside the
  paper. Write a polynomial of degree $n$ on $[0,1]$ as
  $p(t)=\sum_{i=0}^{n}b_i\binom{n}{i}t^i(1-t)^{n-i}$; the $b_i$ are its
  Bernstein coefficients, and they are obtained from the monomial
  coefficients by an exactly invertible integer linear change of basis,
  so if $p$ has rational coefficients then so do the $b_i$, with no
  rounding anywhere. Because the basis functions are nonnegative on
  $[0,1]$ and sum to $1$, $\min_i b_i\le p(t)\le\max_i b_i$ there; in
  particular $b_i\ge0$ for all $i$ implies $p\ge0$ on $[0,1]$. The
  converse fails, which is the whole reason a subdivision is needed: a
  polynomial can be positive on $[0,1]$ and still have a negative
  Bernstein coefficient. What rescues the method is that the bound is
  asymptotically tight under subdivision. Splitting $[0,1]$ at a
  rational point and re-expanding on each half, which de Casteljau's
  algorithm does in exact arithmetic, produces coefficients whose
  minimum converges to $\min p$ quadratically in the width of the
  subinterval \cite{Far12}. So a polynomial that is strictly positive on
  a compact interval is certified after finitely many bisections, and a
  polynomial with a zero in the interior is never certified, however
  fine the subdivision. That dichotomy is exactly the one recorded in
  the algorithm of \cref{app:box-cover-algo}, and it is why the
  regions in \cref{sec:arc1-cert} that contain a zero of the defect, on
  the boundary or inside, carry an excluded margin rather than a zero-margin
  certificate. In several variables the same statements hold on a box,
  with the tensor-product basis and bisection along one coordinate at a
  time. \Cref{app:bernstein-toy} works this mechanism through on a small
  polynomial chosen to also exhibit the boundary-zero difficulty just
  described, for a reader who wants to see the whole computation before
  meeting it at the scale of an actual cell.}

\subsection{Proof of \texorpdfstring{\cref{thm:local-uncond}}{the single-deviation theorem}}

Combining \cref{sec:hessian-def-new,prop:stage-two-positivity-new} with
the four-layer reduction of \cref{sec:argument-structure}: shell
localisation (\cref{sec:shell-new}) reduces to the all-contact corner;
the root-aligned case (\cref{sec:root-aligned-new}) handles
undeviated configurations directly; radial reduction
(\cref{sec:radial-new}) reduces to the all-contact radii; and cell
positivity, established here for every one of the $176$ cells and
every root-aligned deviation arc by the two independent exact arguments
of \stageone\ and \stagetwo, gives $\defect_k\ge0$ throughout. Summing
the free-volume identity of \cref{sec:root-system} over the active cell
type completes the proof of \cref{thm:local-uncond}. $\blacksquare$
\section{From the local bound to the density bound}\label{sec:density-derivation}

\subsection{Periodic packings}\label{sec:periodic-intro-new}

\begin{proof}[Proof of \cref{cor:density-uncond}]
For a periodic packing with fundamental domain $\mathcal{F}$, the
density equals $\vol(B^4)/\bar v$, where $\bar v$ is the average
Voronoi cell volume over $\mathcal{F}$ and $\vol(B^4)=\pi^2/2$. By
\cref{thm:local-uncond}, every cell in $\mathcal{F}$ has volume at
least $8$, so $\bar v\ge8$ and the density is at most
$(\pi^2/2)/8=\pi^2/16$. The $D_4$ lattice attains $\bar v=8$ exactly
(\cref{lem:refcell-vol}), so this bound is sharp.
\end{proof}

\subsection{General packings and the compactness argument}\label{sec:periodic}

For a general packing $\mathcal{C}\subset\mathbb{R}^4$ of unit balls,
with no assumed translational symmetry, the density is the upper
limiting density
\[
  \delta(\mathcal{C}) \;=\; \limsup_{R\to\infty}
  \frac{\vol\bigl(B^4\bigr)\cdot\#\{c\in\mathcal{C}:|c|\le R\}}{\vol(B_R^4)} .
\]
The passage from periodic packings to these uses a compactness argument
that is standard in the subject, in the form given by Rogers
\cite{Rog64}; see also \cite{Gru87} for the general theory of packing
densities in this local, finite-ball formulation. We recall it, since
everything the argument needs has already appeared above.

Topologise the space of packing configurations by local Hausdorff
convergence on increasing balls: $\mathcal{C}_n\to\mathcal{C}$ when, for
every $R$, the points of $\mathcal{C}_n$ within distance $R$ of the
origin converge as a finite point set to those of $\mathcal{C}$. This
space is compact. The packing constraint bounds the number of centres in
any fixed ball uniformly over all packings, so a diagonal argument over
an increasing sequence of balls extracts a convergent subsequence from
any sequence of packings. The finite-ball density functional is upper
semicontinuous for this topology, since a configuration cannot gain
centres in the limit, only lose them to infinity.

Given $\mathcal{C}$ and $\varepsilon>0$, take a periodic packing
$\mathcal{C}'$ agreeing with $\mathcal{C}$ to within $\varepsilon$ on
the ball of radius $1/\varepsilon$. The Voronoi cells of $\mathcal{C}'$
at centres well inside that ball differ from the corresponding cells of
$\mathcal{C}$ by an amount that tends to zero as the radius grows,
because \cref{lem:shell-new} confines the dependence of a cell on the
configuration to a bounded neighbourhood of its own centre. The bound
$\vol(V_c)\ge8$ holds at every centre of every periodic approximant by
\cref{thm:local-uncond}, so it passes to the limit and bounds
$\delta(\mathcal{C})$. This gives \cref{cor:density-uncond} for the full
class of packings satisfying the single-deviation hypothesis, and it is
the mechanism, with the appropriate substitutions, by which every
sphere-packing upper bound from Blichfeldt \cite{Bli29} onward turns a
local statement about individual cells into a global statement about
density.

\subsection{The general case, and what it rests on}

\Cref{thm:local-general,thm:main-general} follow from
\cref{thm:local-uncond,cor:density-uncond} by the same argument once
\cref{cor:remaining} is available. The free-volume identity behind
everything here, stated as \cref{thm:freevol-new} in
\cref{sec:deviation-domain}, holds for every packing configuration, not
only the single-deviation ones; what \cref{thm:local-uncond} adds to the
identity is the sign, and \cref{conj:multidir-new} is exactly the
assertion that the sign survives when several directions deviate at
once. That assertion is settled in \cref{cor:conj-resolved}, and it is
settled by the classification of \cref{thm:m24}; the remaining case,
\cref{thm:m23}, is proved in \cref{sec:strict-inequality,sec:certificate},
and the two statements above therefore rest, beyond the arguments of
this paper, on the one verified kernel of \cref{sec:certificate-checked}
alone, and say so wherever they appear.

\subsection{Remark: what a proof that did not use the classification would need}\label{rem:what-would-close-new}

\Cref{thm:main-general} is proved above, and it is proved using the
classification of \cref{thm:m24} at the configurations of twenty-four
contacts, that is, using a semidefinite certificate. This subsection is
about the other thing, a proof of the same statement that uses no
classification and no certificate, and we set out what such a proof
would have to supply, since it is not what the present one supplies.

It would require
$\defect_k\ge0$, equivalently the volume comparison of
\cref{conj:multidir-new}, uniformly over
\emph{every} packing-valid configuration of active directions, of
\emph{every} size $m\ge1$, with \emph{every} deviating direction ranging
over its full $3$-sphere of possibilities, not merely the
single-direction deviations covered by \cref{thm:local-uncond}. \Cref{sec:hessian-technique,sec:musin-degeneracy,%
sec:exact-sing,sec:broad-sample} explain in detail why the natural
perturbative (Hessian-based) technique that succeeds in the
root-aligned case cannot, by itself, resolve the general
multi-direction case at every configuration, not because the relevant
computations are difficult, but because the technique's own governing
quantity is exactly zero at specific configurations that are
themselves well within the scope the conjecture must cover. Such a proof
will therefore need higher-order information at those configurations,
and locating that obstruction precisely is a contribution of this paper
separate from \cref{thm:local-uncond}.
\section{The general local problem}\label{sec:deviation-domain}

\subsection{The free-volume identity in general}

Throughout this subsection the configuration is at the all-contact
corner with a full active set: the active directions are
$w_1,\dots,w_{24}$, one for each root, with $w_i=\base_i$ except on a
subset $D$ of deviating indices, where $w_j=\cos\tilt_j\,\base_j+\sin\tilt_j\,\dev_j$.
A configuration with a proper active subset is covered a fortiori, since
dropping a constraint only enlarges the cell. Write
\[
  V_c=\bigcap_{i=1}^{24}\{\langle x,w_i\rangle\le1\},
  \qquad
  C_j=\VCell\cap\{\langle x,w_j\rangle\ge1\},
  \qquad
  E_j=V_c\cap\{\langle x,\base_j\rangle\ge1\},
\]
for $j\in D$. Thus $C_j$ is what the $j$-th half-space cuts off the
reference cell as it tilts, and $E_j$ is what the same motion frees
beyond the $j$-th facet.

\begin{thm}[Free-volume identity]\label{thm:freevol-new}
For every $D$ and every choice of tilts and deviation directions,
\begin{equation}\label{eq:freevol}
  \vol(V_c) \;=\; \vol(\VCell)
  \;-\;\vol\Bigl(\bigcup_{j\in D}C_j\Bigr)
  \;+\;\vol\Bigl(\bigcup_{j\in D}E_j\Bigr).
\end{equation}
The identity is exact, with no approximation of any kind. In particular
$\vol(V_c)\ge8$ if and only if
$\vol\bigl(\bigcup_{j}E_j\bigr)\ge\vol\bigl(\bigcup_{j}C_j\bigr)$.
\end{thm}

\begin{proof}
$V_c$ and $\VCell$ are intersections of half-spaces that agree in the
$24-|D|$ constraints indexed outside $D$. A point of $\VCell$ fails to
lie in $V_c$ exactly when $\langle x,w_j\rangle>1$ for some $j\in D$, so
$\VCell\setminus V_c=\bigcup_{j\in D}C_j$ up to a null set; a point of
$V_c$ fails to lie in $\VCell$ exactly when $\langle
x,\base_j\rangle>1$ for some $j\in D$, so $V_c\setminus\VCell=\bigcup_{j\in
D}E_j$. Both sets are finite unions of polytopes, hence measurable, and
$\vol(V_c)=\vol(\VCell)-\vol(\VCell\setminus V_c)+\vol(V_c\setminus\VCell)$.
\end{proof}

For $|D|\le1$ the two unions are single sets and \eqref{eq:freevol}
reduces to $\vol(V_c)=8-\vol(C)+\vol(E)$, which is the form used
throughout \cref{sec:stage-one,sec:stage-two,sec:arc1-cert,sec:arc2-cert}
and again, in the cap form of \cref{prop:cap-reformulation}, in
\cref{sec:cap-theorem}. For $|D|\ge2$ the unions matter, and we are exact here about how.

\begin{remark}[The sums are not the unions]\label{rem:sums-not-unions}
Neither family is disjoint once $|D|\ge2$. The sets $E_j$ overlap
whenever two deviating roots are adjacent, that is when
$\langle\base_j,\base_k\rangle=\tfrac12$: the pyramid $\mathrm{Pyr}_j$ of
\cref{lem:Q-structure} has one of its eight lateral facets on the
hyperplane $\{\langle x,\base_k\rangle=1\}$, so removing the constraint
at $\base_k$ opens $E_j$ out past that facet and into the region where
$\langle x,\base_k\rangle>1$, which is where $E_k$ lives. At an adjacent
pair with both tilts $30^\circ$ the overlap has volume $0.0137\ldots$,
against a defect of $0.1845\ldots$, so it is not a small correction. The
sets $C_j$ can also overlap, though only at large tilts and only for
adjacent pairs; the largest overlap we have found in a search over
random packing-valid configurations of two to four deviating directions
is $1.6\times10^{-3}$.

Consequently the quantity $\sum_j\vol(E_j)-\sum_j\vol(C_j)$ is
\emph{not} equal to $\vol(V_c)-8$, and its sign is neither implied by
nor implies the cell bound: in samples it runs both above the defect
(for instance $0.1798$ against $0.1570$) and below it (for instance
$0.20945$ against $0.20971$). Only the union form \eqref{eq:freevol} is
an identity, and it is the union form that
\cref{conj:multidir-new} asserts.
\end{remark}

\begin{remark}[Why the single-deviation theorem does not reduce the general case]\label{rem:first-order-obstruction}
The natural hope is to split the comparison of \cref{conj:multidir-new}
into $m$ separate ones, each settled by \cref{thm:eperp-closed}. Two
things stand in the way, and both are structural.

The first is the overlapping just described. It is not fatal on its own:
Bonferroni's inequality gives
\[
  \vol\Bigl(\bigcup_{j}E_j\Bigr)\;\ge\;\sum_j\vol(E_j)
  \;-\;\sum_{j<k}\vol(E_j\cap E_k),
\]
and together with $\vol(\bigcup_jC_j)\le\sum_j\vol(C_j)$ this reduces
the conjecture to
\begin{equation}\label{eq:termwise-hope}
  \sum_j\bigl(\vol(E_j)-\vol(C_j)\bigr)\;\ge\;\sum_{j<k}\vol(E_j\cap E_k).
\end{equation}
The overlaps on the right are second order in the tilts, the same order
as each term on the left, so \eqref{eq:termwise-hope} is not obviously
out of reach.

The second obstruction is what closes the route. The $E_j$ appearing in
\eqref{eq:termwise-hope} are not the sets the single-deviation theorem
controls. With one root deviating, $E_j$ is cut out by the other
twenty-three root half-spaces, and \cref{lem:Q-structure} identifies
it inside the pyramid $\mathrm{Pyr}_j$. With $m$ roots deviating it is
cut out by the $24-m$ fixed half-spaces together with the $m-1$ other
tilted ones, and the exchange is not benign. The eight lateral facets of
$\mathrm{Pyr}_j$ lie on the hyperplanes of the eight roots at $60^\circ$
to $\alpha_j$, so for an adjacent $k$ the pyramid meets $\{\langle
x,\base_k\rangle=1\}$ in a facet and not transversally: dropping
that constraint opens $E_j$ out, and reimposing it in the tilted
position $\{\langle x,w_k\rangle\le1\}$ slices back in along a whole
face. The two effects do not cancel, and neither dominates. In a search
over packing-valid configurations with $2\le m\le6$ and tilts up to
$0.9$ we found $\vol(E_j)$ both exceeding its single-deviation value and
falling short of it, in the latter case by more than $0.1$, against
defects of order $0.3$. The search is reproduced by
\texttt{multi\_cap\_reformulation.py} in the supplement. So \cref{thm:eperp-closed} gives no bound at all
on $\vol(E_j)-\vol(C_j)$ once $m\ge2$, and \eqref{eq:termwise-hope}
cannot be reached term by term.

What would suffice is a replacement for $\mathrm{Pyr}_j$: a body,
depending on the whole configuration rather than on one root, that
contains each $E_j$, is cut by the tilted half-spaces in a controlled
way, and for which the analogue of \cref{thm:cap-inequality} can be
proved. We do not have one.
\end{remark}

\begin{cor}\label{cor:freevol-multi-new}
In the setting above with $|D|=m\ge2$, $\vol(V_c)\ge8$ if and only if
\[
  \vol\Bigl(\bigcup_{j\in D}E_j\Bigr)\;\ge\;\vol\Bigl(\bigcup_{j\in D}C_j\Bigr).
\]
\end{cor}

\begin{figure}[tb]
\centering
\includegraphics[width=0.92\linewidth]{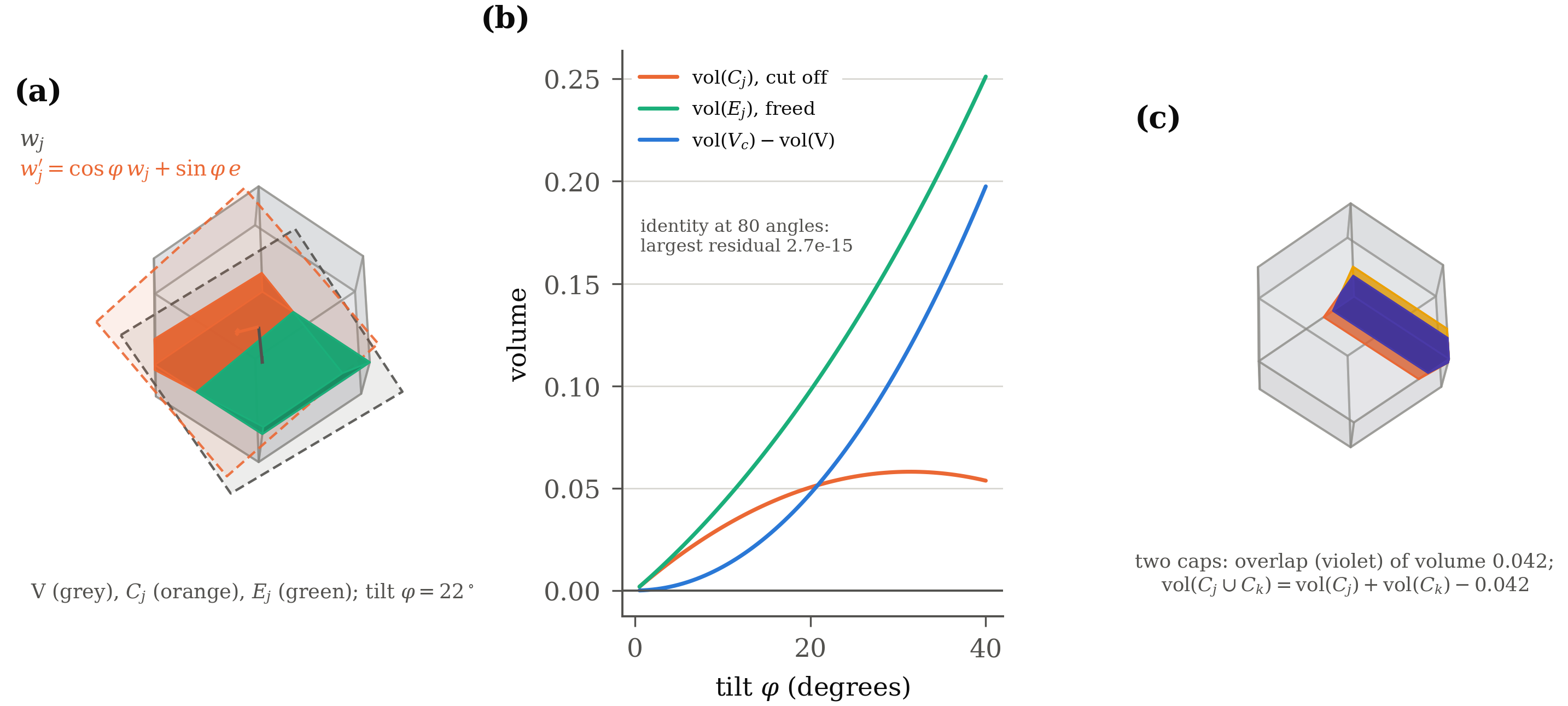}
\caption{The free-volume identity of \cref{thm:freevol-new}, drawn on
real geometry. Because the identity holds in every dimension, it is drawn
in three, where its pieces can be seen: the reference cell is the Voronoi
cell of the face-centred cubic contacts, the rhombic dodecahedron
$\{x:\langle x,w_i\rangle\le1\}$ over the twelve directions
$(\pm1,\pm1,0)/\sqrt2$ and permutations. (a) One direction $w_j$ tilted
by $22^\circ$ to $w_j'$; the cap $C_j$ is what the tilted half-space cuts
off the cell and $E_j$ is what the motion frees beyond the old facet, both
computed as half-space intersections and drawn to scale. (b) The two
volumes and their difference against the tilt, each piece computed as a
polytope of its own; the identity holds at every angle to the precision
of the arithmetic. (c) Two directions tilted towards each other: the caps
overlap, so for $m\ge2$ the corrections are volumes of unions and not
sums of volumes (\cref{rem:sums-not-unions}).}
\label{fig:free-volume}
\end{figure}

\Cref{fig:free-volume} draws the three terms of \eqref{eq:freevol}.

\subsection{A fixed-body form of the multi-direction problem}\label{sec:multi-cap}

\Cref{rem:first-order-obstruction} says what does not work. We record what the cap reformulation of \cref{sec:cap-theorem} does
become for $m\ge2$, since the answer is a single extremal problem about
a single fixed polytope, exactly as in the case $m=1$, and not the sum
of $m$ problems.

Let $D$ be the set of deviating roots, $|D|=m$, and put
\[
  Q_D=\bigcap_{i\notin D}\bigl\{x\in\mathbb{R}^4:\langle
  x,\base_i\rangle\le1\bigr\},
  \qquad
  \capQ_{Q_D}(w)=\{x\in Q_D:\langle x,w\rangle\ge1\}.
\]
For $m=1$ this is the polytope $Q$ of \cref{sec:cap-reformulation}.

\begin{prop}[Multi-cap form]\label{prop:multi-cap}
Suppose $Q_D$ is bounded. Then, for every choice of tilts and deviation
directions,
\begin{equation}\label{eq:multicap-identity}
  \vol(V_c)\;=\;\vol(Q_D)\;-\;\vol\Bigl(\bigcup_{j\in D}\capQ_{Q_D}(w_j)\Bigr),
\end{equation}
and therefore $\vol(V_c)\ge8$ if and only if
\begin{equation}\label{eq:multicap}
  \vol\Bigl(\bigcup_{j\in D}\capQ_{Q_D}(w_j)\Bigr)\;\le\;\vol(Q_D)-8 .
\end{equation}
The right-hand side depends only on $D$, not on the tilts, and equals
$\vol\bigl(\bigcup_{j\in D}\capQ_{Q_D}(\base_j)\bigr)$, so
\eqref{eq:multicap} holds with equality at the undeviated configuration.
It equals $m/3$ exactly when no two roots of $D$ are adjacent, and
exceeds $m/3$ otherwise.
\end{prop}

\begin{proof}
$V_c=Q_D\cap\bigcap_{j\in D}\{\langle x,w_j\rangle\le1\}$ by definition,
so $Q_D\setminus V_c=\bigcup_{j\in D}\capQ_{Q_D}(w_j)$ up to a null set,
which is \eqref{eq:multicap-identity}. Taking $w_j=\base_j$ makes
$V_c=\VCell$ and gives
$\vol\bigl(\bigcup_j\capQ_{Q_D}(\base_j)\bigr)=\vol(Q_D)-8$. If no two
roots of $D$ are adjacent then $\capQ_{Q_D}(\base_j)=\mathrm{Pyr}_j$ for
each $j$, since by the proof of \cref{lem:Q-structure} only the eight
roots at $60^\circ$ to $\alpha_j$ are active on $\{\langle
x,\base_j\rangle\ge1\}$ and none of them lies in $D$; those pyramids are
pairwise disjoint, again by \cref{lem:Q-structure}, so the union has
volume $m/3$. If two roots of $D$ are adjacent, the corresponding
$\capQ_{Q_D}(\base_j)$ strictly contains $\mathrm{Pyr}_j$, since the
constraint carrying one of its lateral facets has been removed.
\end{proof}

So the multi-direction problem is again an extremal problem for caps of
one fixed polytope, with all the parameters in the cutting half-spaces;
what has changed from $m=1$ is that the quantity to bound is the volume
of a \emph{union} of $m$ caps rather than of one cap, and
\cref{rem:first-order-obstruction} shows the union cannot be traded for
a sum. The inequality itself holds, by \cref{cor:conj-resolved} and
therefore by the classification of \cref{thm:m24}; what the present
section and \cref{sec:hessian-technique} to \cref{sec:multidir-scope-summary}
are after is a proof of it that does not go through a classification,
and we state that question as it stands.

\begin{openproblem}[Multi-cap inequality, by direct means]\label{op:multicap}
Prove \eqref{eq:multicap} for every $D$ with $Q_D$ bounded and every
packing-valid choice of $w_j$, with equality only at $w_j=\base_j$,
without appeal to the classification of the $24$-point contact
configurations.
\end{openproblem}

The boundedness hypothesis is mild and its exact range is known. $Q_D$ is
bounded precisely when the $24-m$ undeviated root directions positively
span $\mathbb{R}^4$, and an exhaustive enumeration over all
packing-valid sets $D$ gives the following: for $m\le5$ it is bounded in
every case, all $24$, $276$, $2024$, $10626$ and $42504$ of them
respectively; at $m=6$ it fails for exactly $24$ of the $134596$
packing-valid sets, the first in lexicographic order being the six roots
$e_1\pm e_2$, $e_1\pm e_3$, $e_1\pm e_4$, whose complement lies in the
closed half-space $\{x_1\le0\}$ and therefore leaves $Q_D$ unbounded in
the direction $e_1$. Beyond $m=6$ unbounded cases become common, and the
dense configurations of \cref{sec:dense-configs-new}, at $m$ in the high
teens, are almost all of that kind. For those a bounded substitute for
$Q_D$ would have to come first; \eqref{eq:multicap} as stated does not
reach them.

\subsection{Polar and boundary forms of the same question}\label{sec:polar-surface}

The multi-cap form keeps a fixed body and moves the cutting half-spaces.
There is a third way of writing the same question, which discards the
reference cell altogether and states the problem directly in terms of
the contact configuration. It costs nothing to set up, it makes the
relation to the $24$-cell conjecture explicit, and it gives the
sharpest available statement of what a proof of
\cref{conj:multidir-new} cannot look like.

\begin{definition}\label{def:contact-config}
A \emph{contact configuration} is a finite set $W\subset S^3$ with
$\langle w,w'\rangle\le\tfrac12$ for all distinct $w,w'\in W$. Its
\emph{cell} is
\[
  V_c(W)=\bigl\{x\in\mathbb{R}^4:\langle x,w\rangle\le1
  \text{ for all }w\in W\bigr\}.
\]
\end{definition}

The condition is exactly packing validity.\footnote{If unit balls
centred at $2w$ and $2w'$ are to have disjoint interiors then
$|2w-2w'|\ge2$, that is $|w-w'|\ge1$, that is
$2-2\langle w,w'\rangle\ge1$; so the pairwise bound $\tfrac12$ is the
whole of the packing constraint on the directions, and the bisector of
the segment from $0$ to $2w$ is the hyperplane $\langle x,w\rangle=1$,
which is the constraint defining $V_c(W)$. The directions of the $D_4$
roots, $\base_i=\alpha_i/\sqrt2$, form a contact configuration with the
pairwise inner products taking the four values $\pm\tfrac12$, $0$ and
$-1$, so the bound is attained but never exceeded; this is the
$g$-valued statement of \cref{lem:gram-new}.} Every local configuration
considered in this paper, deviated or not, is a contact configuration,
and $V_c$ is its cell.

\begin{prop}[Polar form]\label{prop:polar-form}
For every contact configuration $W$,
\begin{equation}\label{eq:polar}
  V_c(W)=\conv(W)^{\circ},
\end{equation}
the polar dual of the convex hull of $W$. Consequently
\begin{enumerate}[label=\textup{(\roman*)},leftmargin=2.2em]
\item $\vol(V_c(W))$ depends on $W$ only through $\conv(W)$;
\item $W\subseteq W'$ implies $V_c(W')\subseteq V_c(W)$, so enlarging a
contact configuration never increases the volume of its cell;
\item $|W|\le24$, by Musin's kissing bound \cite{Mus08}.
\end{enumerate}
\Cref{conj:multidir-new} is equivalent to the assertion
$\vol(\conv(W)^{\circ})\ge8$ for every contact configuration whose cell
is bounded, with equality only when $W$ is the set of $24$ root
directions.
\end{prop}

\begin{proof}
The polar of a compact convex set $K$ containing the origin is
$K^{\circ}=\{x:\langle x,p\rangle\le1\text{ for all }p\in K\}$, and a
linear functional attains its maximum over $\conv(W)$ at a point of $W$,
so $\conv(W)^{\circ}$ is cut out by the constraints indexed by $W$
alone, which is \eqref{eq:polar}. Item (i) is immediate from
\eqref{eq:polar} and (ii) from the definition. Item (iii) is the
kissing-number theorem in dimension four, since a contact configuration
is precisely a set of directions to mutually non-overlapping unit balls
touching a fixed one. The final sentence is
\cref{thm:freevol-new} together with the observation that a
configuration at the all-contact corner with a full active set is a
contact configuration of exactly the shape in
\cref{def:contact-config}.
\end{proof}

\begin{cor}[A worst configuration exists, and may be taken maximal]\label{cor:minimiser}
The infimum of $\vol(V_c(W))$ over all contact configurations with
bounded cell is attained, and it is attained at a configuration to which
no further direction can be added.
\end{cor}

\begin{proof}
By \cref{prop:polar-form}(iii) a contact configuration has at most $24$
elements, so for each $m\le24$ the set
$\mathcal{C}_m=\{(w_1,\dots,w_m)\in(S^3)^m:\langle
w_i,w_j\rangle\le\tfrac12\ (i\ne j)\}$ is a closed subset of a compact
space, hence compact, and its points have pairwise distinct entries.
The volume functional $f(W)=\vol(V_c(W))\in(0,\infty]$ is lower
semicontinuous on $\mathcal{C}_m$: if $S$ is a compact subset of the
interior of $V_c(W)$ then $\langle x,w\rangle<1$ on $S\times W$, hence
the same holds for every configuration near $W$, so $S\subseteq
V_c(W')$ and $f(W')\ge\vol(S)$; taking the supremum over such $S$ and
using inner regularity gives $\liminf_{W'\to W}f(W')\ge f(W)$. A lower
semicontinuous function on a compact set attains its infimum, so $f$
attains its infimum on each $\mathcal{C}_m$, and therefore on the union
of the finitely many $\mathcal{C}_m$ with $m\le24$. If a minimiser
$W^{*}$ admits an additional admissible direction, extend it, repeatedly
if necessary; the process stops after at most $24-|W^{*}|$ steps and
each step leaves the volume no larger by \cref{prop:polar-form}(ii), so
the resulting maximal configuration is again a minimiser.
\end{proof}

So the conjecture is a statement about maximal contact configurations
only, and there is a genuine extremal configuration to identify. That
is as far as soft arguments take it: what they leave is the question of
whether the maximal configurations of $24$ directions are the root
system alone, which is the rigidity half of the kissing problem in
dimension four and is settled in \cref{thm:m24} below, by the
semidefinite certificate of \cref{thm:twentyfour-points}.

The second form replaces the volume by the boundary measure.

\begin{lem}[Boundary form]\label{lem:surface-form}
Let $W$ be a contact configuration with bounded cell, and let
$F_w=V_c(W)\cap\{\langle x,w\rangle=1\}$ for $w\in W$. Then
\begin{equation}\label{eq:surface}
  \vol(V_c(W))=\tfrac14\sum_{w\in W}\vol_3(F_w)
  =\tfrac14\,\mathrm{area}\bigl(\partial V_c(W)\bigr),
\end{equation}
so $\vol(V_c(W))\ge8$ if and only if the boundary $3$-volume of the cell
is at least $32$.
\end{lem}

\begin{proof}
The origin is interior to $V_c(W)$, so the cell is the union over its
facets of the pyramids $\conv(\{0\}\cup F)$, with pairwise disjoint
interiors. Each facet lies on a hyperplane $\langle x,w\rangle=1$ with
$|w|=1$, at distance exactly $1$ from the origin, so the corresponding
pyramid has $4$-volume $\tfrac14\vol_3(F)$. Summing gives
\eqref{eq:surface}; the facets of $V_c(W)$ are among the $F_w$, and any
$F_w$ that is not a facet is lower-dimensional and contributes nothing.
\end{proof}

At the root configuration all $24$ sets $F_w$ are facets, each a regular
octahedron of $3$-volume $\tfrac43$, and \eqref{eq:surface} reads
$8=\tfrac14\cdot32$. Why every facet is an octahedron of exactly that
size is worth isolating, since the same computation is what limits the
method.

\begin{lem}[Facet shape and inball]\label{lem:facet-inball}
Let $W$ be a contact configuration, $w\in W$, and write
$g_{w w'}=\langle w,w'\rangle$ and
$\widehat{w'}=(w'-g_{ww'}w)/\sqrt{1-g_{ww'}^2}$ for the unit vector
along the orthogonal projection of $w'$ onto $w^{\perp}$, defined
whenever $g_{ww'}\ne\pm1$. Then, inside the $3$-space $w^{\perp}$,
\begin{equation}\label{eq:facet-shape}
  F_w-w=\Bigl\{y\in w^{\perp}:\langle y,\widehat{w'}\rangle
  \le\rho(g_{ww'})\ \text{ for all }w'\in W\setminus\{w\}\Bigr\},
  \qquad
  \rho(g)=\sqrt{\frac{1-g}{1+g}},
\end{equation}
with the convention that $w'=-w$ imposes no constraint. The function
$\rho$ is strictly decreasing on $(-1,1]$ and $\rho(\tfrac12)=1/\sqrt3$,
so $F_w-w$ contains the $3$-ball of radius $1/\sqrt3$ about the origin
of $w^{\perp}$ and
\begin{equation}\label{eq:facet-lower}
  \vol_3(F_w)\;\ge\;\frac{4\pi}{3}\Bigl(\frac{1}{\sqrt3}\Bigr)^{3}
  =\frac{4\pi}{9\sqrt3}=0.806133\ldots
\end{equation}
\end{lem}

\begin{proof}
Write $x=w+y$ with $y\in w^{\perp}$. Then $\langle x,w\rangle=1$
automatically, and for $w'\ne w$ the constraint $\langle
x,w'\rangle\le1$ reads $g_{ww'}+\langle y,w'\rangle\le1$, that is
$\langle y,w'-g_{ww'}w\rangle\le1-g_{ww'}$, since $y\perp w$. Dividing
by $|w'-g_{ww'}w|=\sqrt{1-g_{ww'}^2}$ turns this into
$\langle y,\widehat{w'}\rangle\le(1-g)/\sqrt{1-g^{2}}=\rho(g)$ at
$g=g_{ww'}$, which is \eqref{eq:facet-shape}. On $(-1,1]$ the quotient
$(1-g)/(1+g)$ is strictly decreasing, hence so is $\rho$, and
$\rho(\tfrac12)=\sqrt{(1/2)/(3/2)}=1/\sqrt3$. Since $g_{ww'}\le\tfrac12$
for every $w'\ne w$, every constraint in \eqref{eq:facet-shape} sits at
distance at least $1/\sqrt3$ from the origin of $w^{\perp}$, so the ball
of that radius is contained in $F_w-w$, and \eqref{eq:facet-lower} is
its volume.
\end{proof}

At the root configuration the neighbours with $g=\tfrac12$ are the eight
roots at $60^{\circ}$ to $\alpha_w$; their projections into $w^{\perp}$
have pairwise inner products in $\{-1,-\tfrac13,\tfrac13\}$, so they are
the eight vertex directions of a cube, and the eight constraints
$\langle y,\widehat{w'}\rangle\le1/\sqrt3$ cut out exactly the regular
octahedron $\{y:|y_1|+|y_2|+|y_3|\le1\}$ of $3$-volume $\tfrac43$, the
remaining sixteen constraints being inactive. The bound
\eqref{eq:facet-lower} is therefore attained in shape but not in size:
the octahedron circumscribes the ball of radius $1/\sqrt3$ and has
$1.654\ldots$ times its volume.

\begin{remark}[Why no facet-by-facet argument can succeed]\label{rem:facet-local-obstruction}
Combining \eqref{eq:surface} and \eqref{eq:facet-lower} with
$|W|\le24$ gives, for every contact configuration,
\[
  \vol(V_c(W))\;\ge\;\frac{|W|}{4}\cdot\frac{4\pi}{9\sqrt3}
  \;=\;\frac{|W|\pi}{9\sqrt3},
\]
which at the largest possible value $|W|=24$ is
$8\pi/(3\sqrt3)=4.836798\ldots$, short of $8$ by more than $3$. The
inequality $V_c(W)\supseteq B^4$, which follows from $|w|=1$ for every
$w$ and is the same computation run in one step instead of facet by facet, does slightly better at $\pi^2/2=4.934802\ldots$, and no better.

The conclusion is not that these two bounds happen to be weak. It is
that the whole class they belong to is too small, and one example settles
it. An argument that uses only the pairwise condition
$g_{ww'}\le\tfrac12$, and nothing about how the rest of the
configuration is arranged, can produce only a bound of the form
$\vol_3(F_w)\ge\beta$ with $\beta$ a universal constant; the largest
such constant is
$\beta^{*}=\inf\vol_3(F_w)$, the infimum being over all contact
configurations and all $w$. The root configuration gives
$\beta^{*}\le\tfrac43$, and one may do better. Keep all eight tight
neighbours but replace the cube by a square antiprism: in an orthonormal
frame with $w=e_0$ and $w^{\perp}=\langle e_1,e_2,e_3\rangle$, take
\begin{equation}\label{eq:antiprism}
  w'_j=\tfrac12\,e_0+\tfrac{\sqrt3}{2}\,\widehat n_j,
  \qquad
  \widehat n_j\in\tfrac{1}{\sqrt3}\bigl\{(\pm\sqrt2,0,1),\,
  (0,\pm\sqrt2,1),\,(\pm1,\pm1,-1)\bigr\}.
\end{equation}
Each $\widehat n_j$ is a unit vector, so
$\langle w,w'_j\rangle=\tfrac12$; and the $\widehat n_j$ have pairwise
inner products in
$\{\tfrac13,-\tfrac13,\tfrac{\sqrt2-1}{3},-\tfrac{\sqrt2+1}{3}\}$,
so $\langle w'_i,w'_j\rangle=\tfrac14+\tfrac34\langle\widehat
n_i,\widehat n_j\rangle\le\tfrac12$ and the configuration is
admissible. By \eqref{eq:facet-shape} the resulting facet is
$\{y:\langle y,\sqrt3\,\widehat n_j\rangle\le1\}$, and integrating
its height over the square $|y_1|,|y_2|$ gives its volume in closed
form:
\[
  \vol_3(F_w)=16\sqrt2-\frac{64}{3}=1.2940836\ldots\;<\;\frac43 .
\]
Feeding $\beta^{*}$ into \eqref{eq:surface} together with $|W|\le24$,
\[
  \vol(V_c(W))\;\ge\;\frac{|W|}{4}\beta^{*}
  \;\le\;6\beta^{*}\;\le\;96\sqrt2-128\;=\;7.764502\ldots
\]
So no bound of this class can reach $8$ at all, whatever it is, and its
ceiling is below even the crude covering estimate of
\cref{thm:covering-bound}. For a configuration with $|W|\le23$ the same
computation returns at most $7.441\ldots$, so the method fails first
exactly where the statement is least delicate.

What is missing is the interaction between facets: the direction $w'$
that makes $F_w$ small is itself a direction whose own facet $F_{w'}$ is
then constrained, and \eqref{eq:surface} becomes usable only if that
trade is accounted for across the configuration. This is the same
phenomenon that \cref{rem:first-order-obstruction} records in the
free-volume language, seen from the boundary instead.
\end{remark}

\begin{remark}[Two global routes that do not reach the equality case]\label{rem:global-routes}
Two standard tools apply to \eqref{eq:polar} and both fail at the
configuration where the conjecture is tight, which is the decisive test
for a method that must produce a sharp constant.

The polar volume has the integral form
$\vol(K^{\circ})=\tfrac14\int_{S^3}h_K(\theta)^{-4}\,d\theta$, where
$h_K$ is the support function and $d\theta$ is the surface measure of
total mass $2\pi^2$. Applying Jensen's inequality to the convex function
$t\mapsto t^{-4}$ and the normalised measure gives
$\vol(K^{\circ})\ge(\pi^2/2)\,\bar h^{-4}$ with $\bar h$ the mean of
$h_K$ over the sphere, and the target $8$ would follow from
$\bar h\le(\pi^2/16)^{1/4}=0.886226\ldots$. For $K=\conv(W)$ at the root
configuration, $h_K(\theta)$ is $2^{-1/2}$ times the sum of the two
largest absolute coordinates of $\theta$, and $\bar h=0.89372\ldots$, evaluated
numerically at four million sample directions with a standard error of
$2\times10^{-5}$, so the criterion fails by a small margin. The failure is not an artefact of
the estimate being loose in a repairable way: Jensen is an equality only
when $h_K$ is constant, that is only for a ball, so at any polytope the
bound is strictly below the true volume, and at the root configuration
it returns $7.7351\ldots$ against the true $8$. No inequality of this
shape can be sharp where the conjecture is.

The Mahler-type route fares worse. The conjectured lower bound
$\vol(K)\vol(K^{\circ})\ge4^{n}/n!$, itself open in dimension four,
would give $\vol(K)\vol(K^{\circ})\ge\tfrac{32}3$; but at the root
configuration $\conv(W)$ is the $24$-cell of circumradius $1$, of volume
$2$, and the product is $16$, so the inequality carries a factor
$\tfrac32$ of slack exactly where none can be afforded, and dividing by
$\vol(K)\le2$ returns only $\tfrac{16}3$. Using instead the crude
$\vol(K)\le\vol(B^4)=\pi^2/2$ returns $64/(3\pi^2)=2.16\ldots$.

Both evaluations, together with the facet-local bound of
\cref{rem:facet-local-obstruction}, the boundary identity
\eqref{eq:surface} and the facet description \eqref{eq:facet-shape}, are
reproduced by \texttt{polar\_surface\_reformulation.py} in the
supplement, and collected in
\cref{fig:boundary}.
\end{remark}

\begin{figure}[tb]
\centering
\includegraphics[width=\textwidth]{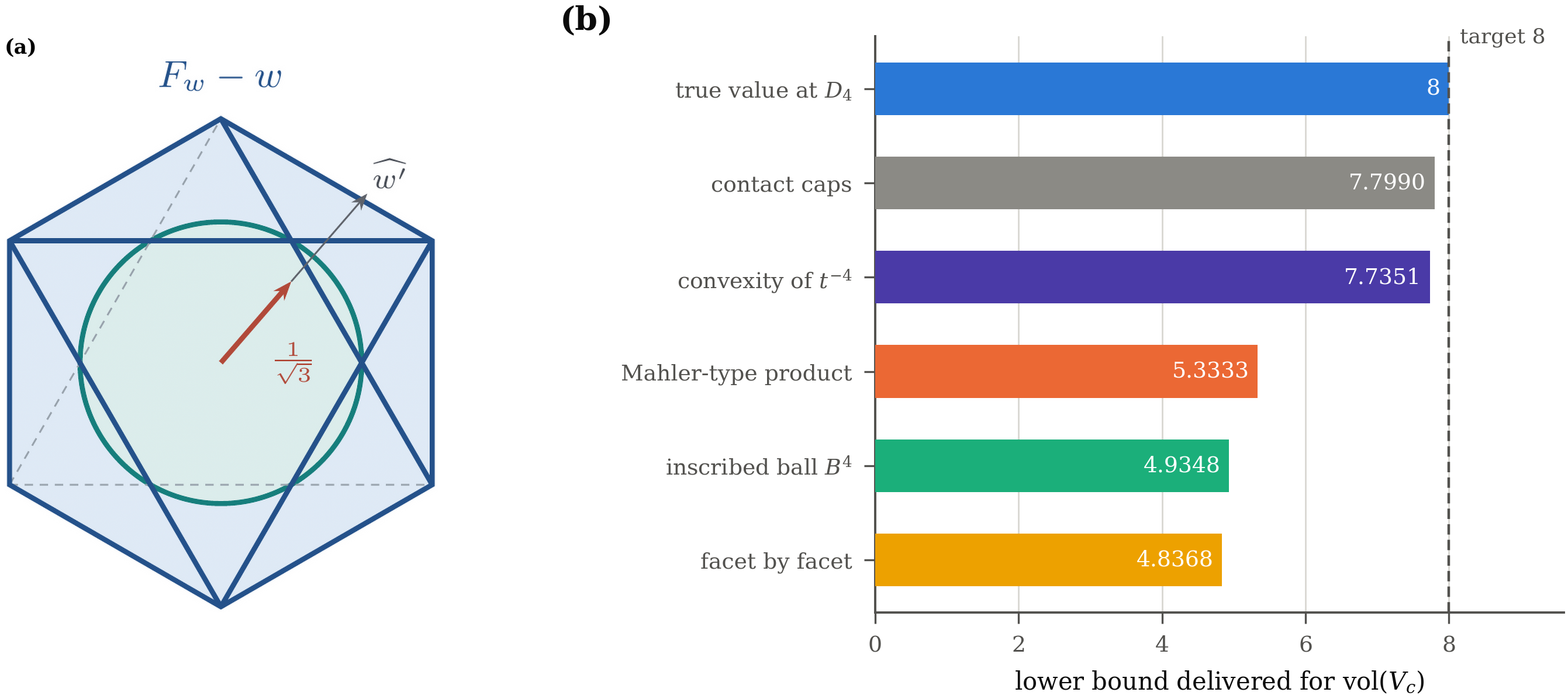}
\caption{(a) A single facet of the reference cell, drawn inside the
$3$-space orthogonal to its own contact direction. At the root
configuration the eight tight neighbours cut out the regular octahedron
(\cref{lem:facet-inball}); the pairwise contact condition on its own
guarantees no more than the inscribed ball, and the gap between the two
is what the boundary form \eqref{eq:surface} has to absorb. (b) The bounds the boundary form and the other global routes
actually deliver, all of them evaluated at the $D_4$ configuration where
the conjecture is tight
(\cref{rem:facet-local-obstruction,rem:global-routes}).}
\label{fig:boundary}
\end{figure}

\begin{remark}[What the multi-direction case is]\label{rem:is-24cell}
\Cref{prop:polar-form} identifies \cref{conj:multidir-new} with a named
problem. The assertion that every Voronoi cell of a unit-ball
packing of $\mathbb{R}^4$ has volume at least that of the regular
$24$-cell, with equality only at the $D_4$ configuration, is the
$24$-cell conjecture; \cref{thm:local-uncond} settles the part of it in
which at most one contact direction leaves a root, and
\cref{conj:multidir-new} is the remainder, settled in
\cref{cor:conj-resolved} by the classification of \cref{thm:m24} and by
nothing else in this paper. We record this plainly. The reformulations
of this subsection do not diminish the problem of proving it directly,
and we do not claim otherwise; what they do is fix the exact shape of
the obstacle and remove three families of approach from consideration,
a fourth having gone in \cref{rem:first-order-obstruction}. Partial
results towards the conjecture, by a different route, are in
\cite{Mus18}.
\end{remark}

\subsection{A bound from the contact caps}\label{sec:covering-bound}

The boundary form measures the cell from outside, facet by facet. There
is a third measurement, from the centre outward, and it is the one that
gives a bound for every configuration at once, however many of its
contact directions deviate.

\begin{lem}[Radial form]\label{lem:radial-form}
Let $W$ be a contact configuration and let
\[
  \delta(\theta)=\min_{w\in W}\arccos\langle\theta,w\rangle ,
  \qquad \theta\in S^3,
\]
be the angular distance from $\theta$ to the nearest contact direction.
Then $V_c(W)$ is bounded if and only if $\delta<\pi/2$ everywhere, and in
that case
\begin{equation}\label{eq:radial}
  \vol\bigl(V_c(W)\bigr)
  \;=\;\frac14\int_{S^3}\sec^4\bigl(\delta(\theta)\bigr)\,d\theta ,
\end{equation}
where $d\theta$ is the surface measure on $S^3$, of total mass
$2\pi^2$.
\end{lem}

\begin{proof}
In polar coordinates the volume of a star body with radial function
$\rho$ is $\tfrac14\int_{S^3}\rho^4$. The radial function of a polar
dual is the reciprocal of the support function of the body, so by
\cref{prop:polar-form},
$\rho_{V_c(W)}(\theta)=1/h_{\conv(W)}(\theta)$ with
$h_{\conv(W)}(\theta)=\max_{w}\langle\theta,w\rangle$. Since every
$w$ is a unit vector, that maximum is
$\max_w\cos\angle(\theta,w)=\cos\delta(\theta)$, which is positive
exactly when $\delta(\theta)<\pi/2$; where it fails the cell contains
the ray through $\theta$ and is unbounded. Substituting
$\rho=1/\cos\delta$ gives \eqref{eq:radial}.
\end{proof}

So the problem is to show that a set of $m\le24$ points on $S^3$ with
pairwise angular distance at least $60^\circ$ cannot be so evenly spread
that the mean of $\sec^4\delta$ drops below $32/(2\pi^2)$. Written this
way, a crude covering estimate already gives most of what is wanted.

\begin{thm}[Covering bound]\label{thm:covering-bound}
Let $W$ be a contact configuration with $|W|=m$ and bounded cell, and
let $r_m$ be the unique solution in $(0,\pi/2)$ of
\begin{equation}\label{eq:rm}
  2r-\sin2r=\frac{2\pi}{m} .
\end{equation}
Then
\begin{equation}\label{eq:covering-bound}
  \vol\bigl(V_c(W)\bigr)\;\ge\;\frac{\pi m}{3}\,\tan^3 r_m .
\end{equation}
\end{thm}

\begin{proof}
A geodesic cap of angular radius $r$ in $S^3$ has surface measure
\[
  C(r)=\int_0^r\bigl(4\pi\sin^2t\bigr)\,dt=\pi\,(2r-\sin2r),
\]
the integrand being the area of the $2$-sphere of radius $\sin t$ that
the cap meets at angular distance $t$ from its centre; in particular
$C(\pi)=2\pi^2$, the whole measure. The set $\{\delta\le r\}$ is the
union of the $m$ caps of radius $r$ about the points of $W$, so
\[
  \sigma\bigl(\{\delta>r\}\bigr)\;\ge\;S(r):=2\pi^2-m\,\pi\,(2r-\sin2r),
\]
where $\sigma$ denotes surface measure. The function $2r-\sin2r$ has
derivative $2-2\cos2r\ge0$, vanishing only at $r=0$, and increases from
$0$ to $\pi$ on $[0,\pi/2]$; since $m\ge5$ forces $2\pi/m<\pi$, equation
\eqref{eq:rm} has exactly one root $r_m$ there, and $S>0$ on $[0,r_m)$,
$S(r_m)=0$.

Now $\sec^4\delta\ge1$, so by the layer-cake formula and the
substitution $s=\sec^4r$, under which
$ds=4\sec^4r\,\tan r\,dr=d(\sec^4r)$,
\begin{multline*}
  \int_{S^3}\sec^4\delta\,d\theta
  \;=\;2\pi^2+\int_1^{\infty}\sigma\bigl(\{\sec^4\delta>s\}\bigr)\,ds\\
  \;=\;2\pi^2+\int_0^{\pi/2}\sigma\bigl(\{\delta>r\}\bigr)\,d(\sec^4r).
\end{multline*}
The integrand is nonnegative, so restricting the range to $[0,r_m]$ and
replacing $\sigma(\{\delta>r\})$ by the smaller $S(r)$ only decreases
the value. Integrating by parts, with $S(0)=2\pi^2$, $S(r_m)=0$ and
$S'(r)=-4\pi m\sin^2r$,
\begin{multline*}
  \int_0^{r_m}S(r)\,d(\sec^4r)\\
  =\bigl[S(r)\sec^4r\bigr]_0^{r_m}+4\pi m\int_0^{r_m}\sin^2r\,\sec^4r\,dr\\
  =-2\pi^2+4\pi m\int_0^{r_m}\tan^2r\,\sec^2r\,dr,
\end{multline*}
and the last integral is $\tfrac13\tan^3r_m$. The two occurrences of
$2\pi^2$ cancel, leaving
$\int_{S^3}\sec^4\delta\ge\tfrac{4\pi m}{3}\tan^3r_m$, which with
\eqref{eq:radial} is \eqref{eq:covering-bound}.
\end{proof}

The right-hand side of \eqref{eq:covering-bound} decreases as $m$ grows,
since a larger $m$ both multiplies by more and shrinks $r_m$; the second
effect wins. The values for the whole admissible range
$5\le m\le24$, a bounded cell in $\mathbb{R}^4$ needing at least five
facets and \cref{prop:polar-form}(iii) allowing at most $24$, are these.

\begin{table}[htbp]
\centering
\small
\begin{tabular}{@{}rcc@{\hspace{2.6em}}rcc@{}}
\toprule
$m$ & $r_m$ (degrees) & bound & $m$ & $r_m$ (degrees) & bound\\
\midrule
 5 & $60.536972$ & $29.042559$ & 15 & $40.255989$ & $9.536151$\\
 6 & $56.406799$ & $21.440235$ & 16 & $39.340721$ & $9.227382$\\
 7 & $53.197354$ & $17.503703$ & 17 & $38.502427$ & $8.962031$\\
 8 & $50.601922$ & $15.119204$ & 18 & $37.730497$ & $8.731395$\\
 9 & $48.441431$ & $13.525786$ & 19 & $37.016289$ & $8.528959$\\
10 & $46.603050$ & $12.387172$ & 20 & $36.352679$ & $8.349750$\\
11 & $45.011524$ & $11.533082$ & 21 & $35.733735$ & $8.189905$\\
12 & $43.614387$ & $10.868477$ & 22 & $35.154469$ & $\mathbf{8.046376}$\\
13 & $42.373749$ & $10.336286$ & 23 & $34.610650$ & $7.916728$\\
14 & $41.261436$ & $9.900225$  & 24 & $34.098663$ & $7.798989$\\
\bottomrule
\end{tabular}
\caption{The covering bound \eqref{eq:covering-bound} over the whole
admissible range of $m$. It exceeds $8$ for every $m\le22$ and falls
below $8$ at $m=23$. Each $r_m$ is the root of \eqref{eq:rm}, computed
to the digits shown by bisection on an exactly evaluated monotone
function.}
\label{tab:covering-bound}
\end{table}

\begin{cor}[The local bound for at most $22$ contacts]\label{cor:m22}
Every contact configuration $W$ with $|W|\le22$ and bounded cell
satisfies $\vol(V_c(W))\ge8.0463\ldots>8$, with no hypothesis on how
many of its directions deviate from roots and no hypothesis on the size
of those deviations. Consequently
\[
  \vol(V_c)\;\ge\;\tfrac{\pi\cdot24}{3}\tan^3r_{24}=7.798989\ldots
\]
holds for every contact configuration whatever, and the local bound
$\vol(V_c)\ge8$ is left undecided only for configurations with exactly
$23$ or exactly $24$ contacts.
\end{cor}

\begin{proof}
Immediate from \cref{thm:covering-bound} and
\cref{tab:covering-bound}, using $|W|\le24$ from
\cref{prop:polar-form}(iii) and the monotonicity visible in the table,
which over the finite range $5\le m\le24$ is a comparison of twenty
numbers.
\end{proof}

\Cref{cor:m22} is \cref{thm:local-fewcontacts} of the introduction, and
it settles one point that the deviation picture leaves open. A
configuration in which several of the twenty-four contacts are absent and not merely deviated is not a perturbation of the root
configuration at all, and the argument of \cref{sec:cap-theorem} says
nothing about it. \Cref{cor:m22} disposes of every such configuration
with at most twenty-two contacts outright, whatever its directions.

\begin{prop}[The covering bound is the best area-based bound]\label{prop:area-optimal}
Let $W$ be a contact configuration with bounded cell and let
$\Omega_w=\{\theta\in S^3:\langle\theta,w\rangle\ge\langle\theta,w'\rangle
\ \text{for all }w'\in W\}$ be its spherical Voronoi cells, so that
$\delta(\theta)$ is the distance to $w$ on $\Omega_w$ and
$\sum_w\sigma(\Omega_w)=2\pi^2$. Write $\rho_w$ for the radius of the
cap of the same measure, $C(\rho_w)=\sigma(\Omega_w)$. Then
\begin{equation}\label{eq:percell}
  \vol\bigl(V_c(W)\bigr)\;\ge\;\frac{\pi}{3}\sum_{w\in W}\tan^3\rho_w ,
\end{equation}
and the function $\phi(s)=\tan^3\bigl(C^{-1}(s)\bigr)$ is strictly
convex, with $\phi'(s)=\tfrac{3}{4\pi}\sec^4\bigl(C^{-1}(s)\bigr)$.
Consequently the right-hand side of \eqref{eq:percell} is smallest, over
all ways of distributing the total measure $2\pi^2$ among $m$ cells,
exactly when the cells have equal measure, and its value there is
$\tfrac{\pi m}{3}\tan^3r_m$: the bound of \cref{thm:covering-bound}.
\end{prop}

\begin{proof}
On $\Omega_w$ the integrand of \eqref{eq:radial} is
$\sec^4 d(\theta,w)$, an increasing function of the distance to $w$, so
its sublevel sets are the caps about $w$. By the bathtub principle, among
all sets of measure $\sigma(\Omega_w)$ the integral is least on the
sublevel set of that measure, which is the cap of radius $\rho_w$; and
\[
  \int_{\{d(\theta,w)\le\rho\}}\sec^4d(\theta,w)\,d\theta
  =\int_0^{\rho}\sec^4t\cdot4\pi\sin^2t\,dt
  =4\pi\int_0^{\rho}\tan^2t\,\sec^2t\,dt
  =\frac{4\pi}{3}\tan^3\rho .
\]
Summing over $w$ and dividing by $4$ gives \eqref{eq:percell}. For the
convexity, $C'(\rho)=4\pi\sin^2\rho$, so
\[
  \phi'(s)=\frac{3\tan^2\rho\,\sec^2\rho}{4\pi\sin^2\rho}
  =\frac{3\sec^4\rho}{4\pi},\qquad \rho=C^{-1}(s),
\]
which is strictly increasing in $\rho$ and hence in $s$. Jensen's
inequality for a convex function then gives
$\sum_w\phi(\sigma(\Omega_w))\ge m\,\phi(2\pi^2/m)$, with equality only
when all the measures agree, and $\phi(2\pi^2/m)=\tan^3r_m$ by the
definition \eqref{eq:rm} of $r_m$.
\end{proof}

So \cref{thm:covering-bound} is not one estimate among several: it is
the strongest conclusion that the areas of the spherical Voronoi cells
can support, and at the root configuration those areas are all equal, so
no redistribution argument can improve it there. What the estimate gives
away is shape. It replaces each cell by the round set of the same
measure, and a round set is exactly the worst case for an integrand that
grows with the distance. At the root configuration each cell is a
spherical polytope with eight facets and circumradius $45^\circ$, while
the cap of the same measure has radius $34.0986\ldots^\circ$; that
difference is the whole of the remaining $0.2011$. Any improvement must
therefore begin by proving that a spherical Voronoi cell of a
$60^\circ$-separated set cannot be round, which is a statement about the
shape of the cells and not about their areas.

\begin{remark}[The shape is not local to a cell]\label{rem:no-local-cell}
The natural way to get at that shape is one cell at a time, and it fails
in both of the decompositions the configuration provides.

On $\Omega_w$ the integrand of \eqref{eq:radial} is $\sec^4d(\theta,w)$,
so the Voronoi decomposition localises the identity exactly, and
$\vol(V_c(W))\ge8$ would follow from
\begin{equation}\label{eq:local-cell}
  \int_{\Omega_w}\sec^4 d(\theta,w)\,d\theta
  \;\ge\;\frac{16}{\pi^2}\,\sigma(\Omega_w)
\end{equation}
at every $w$, that is from the mean of $\sec^4$ over each cell being at
least $16/\pi^2=1.621138\ldots$. At the root configuration every cell
gives exactly that, so \eqref{eq:local-cell} is sharp. It is also false.
A direction may carry nine contacts at $60^\circ$: writing
$w'=\tfrac12 w+\tfrac{\sqrt3}{2}n$ with $n$ a unit vector orthogonal to
$w$, the condition $\langle w',w''\rangle\le\tfrac12$ becomes
$\langle n',n''\rangle\le\tfrac13$, and nine unit vectors of
$\mathbb{R}^3$ with pairwise inner products at most $\tfrac13$ exist, by
the case of nine points in Sch\"utte and van der Waerden's determination
of the optimal spherical codes of small size \cite{SW51}. The cell of
$w$ is then smaller than any cell of the root configuration, and the
mean of $\sec^4$ over it is $1.5732\ldots$; moving the nine contacts out
to $61^\circ$, so that every inner product is strictly below
$\tfrac12$, still gives $1.5954\ldots$. Enlarging the configuration only
shrinks $\Omega_w$ and lowers the mean, so every completion violates
\eqref{eq:local-cell} there.

The Delaunay decomposition does not localise at all. Inside a Delaunay
cell the nearest contact direction need not be a vertex of that cell, so
$\int_T\sec^4\delta$ is not a functional of $T$; on a $22$-point
configuration there are directions whose nearest contact lies outside
the cell containing them and beats every vertex of that cell by more
than $0.1$ in cosine. Even the functional it would have to be, the
distance to the cell's own vertices, fails the corresponding inequality.
For the regular spherical simplex of edge $60^\circ$, parametrised by
the barycentric coordinates $\lambda$ of the cone it spans, the
identities $|x|^2=\tfrac12(1+q)$ with $q=\sum_i\lambda_i^2$ and
$\langle x,v_i\rangle=\tfrac12(1+\lambda_i)$ turn the mean of $\sec^4$
into
\begin{equation}\label{eq:regular-simplex-mean}
  \frac{4\,\mathbb{E}\bigl[(1+M)^{-4}\bigr]}
       {\mathbb{E}\bigl[(1+q)^{-2}\bigr]}
  =\frac{4/5}{0.518254\ldots}=1.543643\ldots ,
  \qquad
  \mathbb{E}\bigl[(1+M)^{-4}\bigr]=\tfrac15 ,
\end{equation}
$M$ being the largest coordinate and both expectations taken over a
uniformly distributed $\lambda$. The value falls short of $16/\pi^2$ by
$0.077495\ldots$, which is $4.78\%$. That simplex has circumradius
$37.7612\ldots^\circ$, so \cref{prop:extendable} below does not exclude
it from a saturated configuration, and it does occur as a Delaunay cell
of a contact configuration; at edge $62^\circ$, where every inner
product is strictly below $\tfrac12$, the mean is still
$1.5955\ldots$. It occurs nowhere in the root system, whose $60^\circ$
graph has clique number three, and that is the point of the example: the
locally densest arrangement is not one the extremal configuration
contains.

So the deficit that \cref{prop:area-optimal} shows is needed is not
additive over cells. Cells worse than round exist in both
decompositions, and what makes the root configuration extremal is the way its cells fit together, not the shape of any one of them. All
of this is checked by \texttt{local\_cell\_obstruction.py}, whose
integrals over $S^3$ are Monte Carlo and are labelled as such; the
margins above are percentages, not last digits.
\end{remark}

\begin{remark}[What the estimate leaves]\label{rem:what-is-left}
Two contact counts survive \cref{cor:m22}, and they are of different
kinds. At $m=24$ the configuration is a full contact shell and the
question is the one \cref{conj:multidir-new} asks, whether two or more
of the twenty-four directions can leave the roots; the estimate returns
$7.7989\ldots$ there, a shortfall of $0.2011$. At $m=23$ the
configuration need not be a deviation of the root system at all, so
\cref{conj:multidir-new} does not reach it and \cref{cor:m22} stops one
short; the estimate returns $7.9167\ldots$, a shortfall of $0.0833$.
\Cref{sec:m24} disposes of the first count and reduces the second to
the configurations with bounded cell, and
\cref{sec:strict-inequality,sec:certificate} prove the bound for those.
\end{remark}

\begin{prop}[Every cell reaches out to $37.7612\ldots^\circ$]\label{prop:covering-radius}
Let $W$ be a contact configuration with bounded cell. Then every vertex
of every spherical Voronoi cell $\Omega_w$ is at angular distance at
least $R_0=\arccos\sqrt{5/8}=37.7612\ldots^\circ$ from $w$, and in
particular the covering radius of $W$ on $S^3$ is at least $R_0$.
\end{prop}

\begin{proof}
A vertex $V$ of $\Omega_w$ is a point of $S^3$ at which at least four
cells meet, so there are at least four points of $W$, one of them $w$,
at the common distance $R=d(V,w)$ from $V$ and none nearer. Let
$u_1,\dots,u_4$ be unit vectors in the tangent space at $V$ pointing
along the geodesics from $V$ to four of them. Four unit vectors in
$\mathbb{R}^3$ cannot have all six pairwise inner products below
$-\tfrac13$, since
$0\le\bigl|\sum_iu_i\bigr|^2=4+2\sum_{i<j}\langle u_i,u_j\rangle$ forces
$\sum_{i<j}\langle u_i,u_j\rangle\ge-2$ and there are six terms. Take a
pair with $\langle u_i,u_j\rangle\ge-\tfrac13$. The two corresponding
points of $W$ lie at distance $R$ from $V$ in directions at angle
$\alpha$ with $\cos\alpha\ge-\tfrac13$, so by the spherical law of
cosines their separation $\gamma$ satisfies
\[
  \cos\gamma=\cos^2R+\sin^2R\,\cos\alpha
  \;\ge\;\cos^2R-\tfrac13\sin^2R .
\]
Packing validity gives $\gamma\ge60^\circ$, hence
$\cos\gamma\le\tfrac12$, so $\cos^2R-\tfrac13(1-\cos^2R)\le\tfrac12$
and $\cos^2R\le\tfrac58$. The farthest point of $S^3$ from $W$ is a
vertex of the cell containing it, so the covering radius is at least
$R_0$ as well.
\end{proof}

This is a genuine improvement on what \cref{thm:covering-bound} extracts
along the way, since that argument gives up at
$r_{24}=34.0986\ldots^\circ$. It does not help. Knowing that the
survival function is positive out to $R_0$ bounds it below only by the
measure of a cap of radius $R_0-r$ about one deep hole, and feeding that
into the integral changes the constant in the fifth decimal:
$7.79899$ becomes $7.79902$, and even allowing one such hole per cell
and treating them as disjoint gives only $7.79966$. The missing
$0.140$ is a statement about how much of the sphere lies far from $W$,
not about how far the farthest point is.

The overlaps the estimate throws away can be put back exactly, and it is
worth doing, because the answer says how much room the method has left.

\begin{lem}[Three caps do not meet]\label{lem:no-triples}
Three directions of a contact configuration have circumradius at least
$\arccos\sqrt{2/3}=35.2644\ldots^\circ$ on their span. Consequently, if
$r<\arccos\sqrt{2/3}$, no point of $S^3$ lies within $r$ of three
contact directions at once.
\end{lem}

\begin{proof}
Let $G$ be the Gram matrix of the three directions, with unit diagonal
and off-diagonal entries at most $\tfrac12$, and let $z$ be the point
of their span equidistant from all three, at angular distance $R$.
Writing $z=\sum a_iw_i$, the conditions $\langle z,w_i\rangle=\cos R$
and $|z|=1$ give $Ga=(\cos R)\mathbf1$ and
$1=a^{\top}Ga=(\cos R)^2\,\mathbf1^{\top}G^{-1}\mathbf1$. Now
$\mathbf1^{\top}G^{-1}\mathbf1=\max_{x\ne0}(\mathbf1^{\top}x)^2/(x^{\top}Gx)$,
and $x=\mathbf1$ already gives
\[
  \mathbf1^{\top}G^{-1}\mathbf1\;\ge\;
  \frac{9}{3+2\sum_{i<j}g_{ij}}\;\ge\;\frac{9}{3+3}\;=\;\frac32 ,
\]
so $\cos^2R\le\tfrac23$. A point within $r$ of all three would give a
triple of circumradius at most $r$.
\end{proof}

\begin{prop}[Putting the overlaps back]\label{prop:second-order}
Let $W$ be a contact configuration with $|W|=m$ and bounded cell, write
$\gamma_{ij}$ for the angle between $w_i$ and $w_j$, and let
\begin{equation}\label{eq:lens}
  \Lambda(r,\gamma)=4\pi\int_{\gamma/2}^{r}
    \sin^2t\,\Bigl(1-\tan\tfrac{\gamma}{2}\,\cot t\Bigr)\,dt
  \qquad (r\ge\tfrac{\gamma}{2}),
\end{equation}
the measure of the intersection of two caps of radius $r$ whose centres
are $\gamma$ apart, and $\Lambda(r,\gamma)=0$ for $r<\gamma/2$. Then for
every $R\le r_*:=\arccos\sqrt{2/3}=35.2644\ldots^\circ$,
\begin{equation}\label{eq:second-order}
  \vol\bigl(V_c(W)\bigr)\;\ge\;
  \frac14\Bigl[2\pi^2+\int_0^{R}\bigl(2\pi^2-mC(r)\bigr)\,d(\sec^4r)\Bigr]
  \;+\;\frac14\sum_{i<j}\int_0^{R}\Lambda(r,\gamma_{ij})\,d(\sec^4r) .
\end{equation}
At $R=r_m$ the first bracket is the covering bound $\tfrac{\pi m}{3}\tan^3r_m$
of \cref{thm:covering-bound}, and $r_m<r_*$ for every $m\ge23$. Below
$r_*$ the right-hand side is not an estimate at all: by
\cref{lem:no-triples} the caps of radius $r$ meet only in pairs there, so
$\sigma(\{\delta\le r\})=mC(r)-\sum_{i<j}\Lambda(r,\gamma_{ij})$
exactly, and \eqref{eq:second-order} is the true value of the integral
truncated at $R$, that is the volume of $V_c(W)$ inside the ball of
radius $\sec R$; the largest admissible $R$ is $r_*$, where the ball has
radius $\sqrt{3/2}$.
\end{prop}

\begin{proof}
For $\Lambda$, split the lens by the hyperplane bisecting $w_i$ and
$w_j$; the half on the $w_i$ side is the part of the cap about $w_j$
lying beyond that hyperplane, which is at angular distance $\gamma/2$
from $w_j$. Parametrise the cap about $w_j$ by the distance $t$ and a
direction on the $2$-sphere it meets at distance $t$, of measure
$4\pi\sin^2t\,dt$; the coordinate of that direction along the bisector
normal is uniform on $[-1,1]$, and the half-space condition reads
$s\ge\cot t\,\tan(\gamma/2)$, which is vacuous for $t\le\gamma/2$ and
otherwise excludes a fraction $\tfrac12(1+\cot t\tan\tfrac\gamma2)$.
Doubling gives \eqref{eq:lens}. The rest is the layer-cake step in the
proof of \cref{thm:covering-bound}, run with the exact survival function
on $[0,R]$ and with the integrand discarded above $R$: writing
$\sec^4\min(\delta,R)=1+\int_0^R\mathbf 1[\delta>r]\,d(\sec^4r)$ and
integrating over $S^3$ gives the bracket, and the survival function is
$2\pi^2-mC(r)+\sum_{i<j}\Lambda(r,\gamma_{ij})$ on $[0,R]$ by
\cref{lem:no-triples}.
\end{proof}

At the configuration obtained by deleting one root, where $m=23$ and the
$88$ close pairs are all at exactly $60^\circ$, the right-hand side of
\eqref{eq:second-order} is
\begin{equation}\label{eq:second-order-value}
  7.997885\ldots\ \text{at } R=r_{23},
  \qquad
  8.034340\ldots\ \text{at } R=r_* ,
\end{equation}
against $7.916728\ldots$ for the plain estimate and $\tfrac{25}3$ for
the true volume. The first of these is the number one obtains by
stopping the integration where the covering estimate stops, at the
radius $r_{23}=34.6106\ldots^\circ$ at which the caps would exactly
exhaust the sphere if they were disjoint, and it falls short of $8$ by
$0.002115\ldots$. The second carries the same exact identity to its
natural limit $r_*$, the smallest circumradius of three contact
directions, and it passes $8$ with $0.034340\ldots$ to spare. The
pairwise estimate is therefore not blocked at the deletion
configuration; what it asks of an arbitrary configuration of $23$
directions is worked out in \cref{prop:pair-budget}, and
\cref{sec:strict-inequality} takes the matter further. At $m=24$ the
same computation gives $7.858738\ldots$ at $R=r_{24}$ and
$7.907070\ldots$ at $R=r_*$, against the true value $8$; there the
truncation discards the deep holes of the root system, which sit at
$45^\circ$.

\begin{remark}[What the constant is and is not]\label{rem:covering-scope}
Two things should be said about $7.798989\ldots$ so that it is not
mistaken for more than it is.

As a density statement it is weak. A per-cell bound $v$ gives packing
density at most $\vol(B^4)/v$, so $7.798989\ldots$ gives $0.632749\ldots$,
against $\pi^2/16=0.616850\ldots$ for $D_4$; the linear programming
bound recalled in \cref{sec:comparison-methods} is stronger as a density
bound.
What the covering estimate has that a density bound does not is that it
is a statement about one cell at a time, which is what
\cref{conj:multidir-new} asks for, and that is why it settles the range
$m\le22$ outright.

The estimate itself is crude, and deliberately so: the only fact about
the configuration it uses is that $m$ caps cannot cover more than $m$
times a cap's worth of the sphere. \Cref{prop:second-order} says how
much of that crudeness is repairable from pair data alone, and the
answer is almost all of it at $m=23$ and rather little at $m=24$. What
neither version reaches is the tail: both stop at $r_m$, while the
survival function stays positive out to the covering radius, $45^\circ$
at the root configuration. Going further means a lower bound on the
covering radius better than $r_m$, and \cref{prop:covering-radius} shows
what that is worth, which is the fifth decimal. The rest of the tail is
a statement about how much of the sphere lies far from $W$, and it is
the same statement as the rigidity of the kissing arrangement.

Everything in this subsection, the cap-area formula, the radial
identity, the layer-cake step, the closed form, the whole of
\cref{tab:covering-bound} including its monotonicity, and the absence of
any violation of \eqref{eq:covering-bound} at the root configuration, at
its subsets and at random configurations, is checked by
\texttt{covering\_bound.py} in the supplement, and \cref{lem:no-triples}
together with \eqref{eq:lens} and the two values above by
\texttt{second\_order\_estimate.py}.
\end{remark}

\subsection{Twenty-four contacts, and what is left}\label{sec:m24}

At the top of the range the problem changes character completely, because
a contact configuration of the maximal size is rigid.

\begin{thm}[Twenty-four contacts]\label{thm:m24}
Let $W$ be a contact configuration with $|W|=24$. Then $W$ is the image
of the $24$ normalised roots of $D_4$ under an orthogonal map, and
consequently $V_c(W)$ is a regular $24$-cell and $\vol(V_c(W))=8$.
\end{thm}

The theorem was first proved by de Laat, Leijenhorst and de Muinck Keizer
\cite[Thm.~5.3]{LLM24}, by a semidefinite computation. It is proved here in
full. The proof has two halves. The first is an inequality: a contact
configuration has at most $24$ elements, and one with exactly $24$ has all
its pairwise inner products in $\{-1,-\tfrac12,0,\tfrac12\}$. That is
\cref{thm:twentyfour-points} below, and it is a consequence of an explicit
positive kernel on the two-element subsets of $S^3$. The kernel was found by
the authors of \cite{LLM24} and deposited with their paper
\cite{LLM24data}; what is proved here is that any kernel with its
properties bounds the size of a contact configuration
(\cref{lem:las2}), that a kernel built from positive definite blocks in the
way theirs is built is positive (\cref{lem:gram}), that the deposited data
have every property required, which is a finite computation carried out in
\cref{sec:certificate-checked}, and that at the value $24$ the bound
forces the inner products. No statement of \cite{LLM24} is used; what is
taken from that work is the certificate itself, which is data, and the
construction of the kernels in which it is expressed, which is a
definition. The second half is classical, and is the following lemma.

\begin{lem}[From the inner products to the root system]\label{lem:root-lattice}
Let $W$ be a set of $24$ unit vectors in $\mathbb{R}^4$ whose pairwise
inner products lie in $\{-1,-\tfrac12,0,\tfrac12\}$. Then $\sqrt2\,W$ is
the set of roots of a copy of the $D_4$ root lattice, and $W$ is the
image of the normalised roots of $D_4$ under an orthogonal map.
\end{lem}

\begin{proof}
The vectors of $\sqrt2\,W$ have norm $2$ and integer inner products, so
they generate an integral lattice $L\subset\mathbb{R}^4$ of rank at most
$4$, and they are among the vectors of norm $2$ of $L$. The vectors of
norm $2$ of an integral lattice form a root system, simply laced, since
two of them have inner product in $\{0,\pm1,\pm2\}$ and the reflection
in one of them preserves the lattice
\cite[Ch.~4]{CS1988}, and the lattice they generate is a root
lattice, an orthogonal sum of lattices of the types $A_n$, $D_n$ and
$E_n$ by the classification of simply laced root systems
\cite[Ch.~2]{Hum90}; it has a basis of roots, its simple roots. Since
$\sqrt2\,W$ generates $L$, the lattice $L$ is such a root lattice, of
rank $r\le4$, and $|\sqrt2\,W|=24$ is at most the number of roots of
$L$. The root lattices of rank at most $4$ and their numbers of roots
are $A_1$, $A_2$, $A_3$, $A_4$ with $2$, $6$, $12$, $20$, $D_4$ with
$24$, and the orthogonal sums $A_1^2$, $A_1A_2$, $A_1^3$, $A_1A_3$,
$A_2^2$, $A_1^2A_2$, $A_1^4$ with $4$, $8$, $6$, $14$, $12$, $10$, $8$.
Only $D_4$ reaches $24$, so $L$ is isometric to the $D_4$ lattice and
$\sqrt2\,W$ is its full set of roots. The $D_4$ root system is unique up
to isometry, which gives the last statement.

The finite content of this argument is checked, without appeal to the
classification, by \texttt{root\_lattices\_rank4.py} and again by Lean's
kernel in \texttt{D4RootLattices.lean}, which also exhibits in every
lattice with $24$ roots four of them with the Cartan matrix of $D_4$ as
Gram matrix, so that the lattice contains a copy of the $D_4$ lattice of
its own determinant and is that lattice:
a lattice generated by norm-$2$ vectors has a basis of such vectors,
whose Gram matrix has $2$ on the diagonal and $-1$, $0$ or $1$ off it;
the script enumerates every positive definite matrix of that shape of
order $1$ to $4$, there are $1$, $3$, $23$ and $393$ of them, counts the
integer solutions of $x^{\mathsf T}Gx=2$ for each, exactly and within a
rigorous box, and finds the largest counts $2$, $6$, $12$ and $24$, the
count $24$ occurring only at determinant $4$ with each root having $8$
neighbours at inner product $1$, $6$ at $0$, $8$ at $-1$ and one at
$-2$, which is $D_4$; the next count at rank $4$ is $20$, at determinant
$5$, which is $A_4$.
\end{proof}

\subsubsection{The bound of the second level, and its equality case}
\label{sec:las2}

For $t\ge0$ let $\mathcal I_t$ be the set of subsets of $S^3$ of size at
most $t$ whose pairwise inner products are at most $\tfrac12$; the empty
set belongs to $\mathcal I_t$, and a contact configuration is an element
of $\mathcal I_t$ for $t$ large enough. Every subset of an element of
$\mathcal I_t$ is again one. For a symmetric real function $K$ on
$\mathcal I_2\times\mathcal I_2$ and a set $Q\in\mathcal I_4$ put
\begin{equation}\label{eq:A2}
  A_2K(Q)\;=\;\sum_{\substack{J_1,J_2\in\mathcal I_2\\ J_1\cup J_2=Q}}K(J_1,J_2),
\end{equation}
the sum running over ordered pairs, so that $A_2K(\emptyset)=K(\emptyset,\emptyset)$.
Call $K$ \emph{positive} if for every finite family $J_1,\dots,J_m$ of
elements of $\mathcal I_2$ the matrix $\bigl(K(J_i,J_j)\bigr)_{i,j}$ is
positive semidefinite.%
\footnote{The hierarchy of \cite{dLV15} works with continuous positive
  kernels on the compact space $\mathcal I_2$, because it is stated for
  an optimisation problem whose feasible set has to be closed. For the
  bound on one given code only the finite Gram matrices above enter, so
  continuity is not needed and is not assumed here; the certificate
  kernel is in fact a polynomial in inner products and is continuous.}
This is the second level of the hierarchy that de
Laat and Vallentin \cite{dLV15} set up for packing problems in discrete
geometry, and the bound it gives is the following; we include the short
proof because the equality case, which is what the theorem needs, is
read off from it.%
\footnote{The bound is the generalisation to four points of the
  Delsarte linear programming bound and of the three-point bound of
  Bachoc and Vallentin \cite{BV08}: with $\mathcal I_1$ in place of
  $\mathcal I_2$ the same argument gives $|C|\le K(\emptyset,\emptyset)$
  for a positive kernel on the one-element subsets of $S^3$, which is
  the Delsarte bound in the form of \cite{DGS77}, and the equality
  argument is the one that identifies the extremal codes there, for
  instance in \cite{CK07}. The only thing that changes with the level
  is the set of subsets over which the pairs are grouped.}

\begin{lem}[The bound of the second level]\label{lem:las2}
Let $K$ be positive, and suppose that $A_2K(\{x\})\le-1$ for every
$x\in S^3$ and that $A_2K(Q)\le0$ for every $Q\in\mathcal I_4$ with
$|Q|\ge2$. Then every contact configuration $C$ satisfies
\[
  |C|\;\le\;K(\emptyset,\emptyset).
\]
If $|C|=K(\emptyset,\emptyset)$, then $A_2K(\{x\})=-1$ for every $x\in C$
and $A_2K(Q)=0$ for every $Q\subseteq C$ with $2\le|Q|\le4$.
\end{lem}

\begin{proof}
The set $C$ is finite, since its points are at pairwise angular distance
at least $\pi/3$ on a compact sphere, so we may argue with finite
sums. Let $\mathcal J$ be the family of subsets of $C$ of size at most
$2$; it lies in $\mathcal I_2$. Positivity of $K$ on this family, with every
coefficient equal to $1$, gives
\[
  0\;\le\;\sum_{J_1,J_2\in\mathcal J}K(J_1,J_2)
  \;=\;\sum_{\substack{Q\subseteq C\\ |Q|\le4}}A_2K(Q),
\]
the second expression being the first with the pairs grouped by their
union $Q$, which is a subset of $C$ of size at most $4$ and hence an
element of $\mathcal I_4$. The term $Q=\emptyset$ is $K(\emptyset,\emptyset)$.
The $|C|$ terms with $|Q|=1$ are each at most $-1$, and every term with
$|Q|\ge2$ is at most $0$. Hence $0\le K(\emptyset,\emptyset)-|C|$, which is
the bound. If $|C|=K(\emptyset,\emptyset)$ the chain
\[
  0\;\le\;\sum_{Q}A_2K(Q)
  \;=\;K(\emptyset,\emptyset)+\sum_{|Q|=1}A_2K(Q)+\sum_{|Q|\ge2}A_2K(Q)
  \;\le\;K(\emptyset,\emptyset)-|C|+0\;=\;0
\]
collapses, so the two sums take the values $-|C|$ and $0$ that bound
them; since each singleton term is at most $-1$ and each remaining term at
most $0$, every singleton term equals $-1$ and every remaining term
equals $0$.%
\footnote{This is the whole of the argument by which an exact bound
  identifies a configuration, and it is worth seeing how little it
  uses. It does not use that $C$ is maximal in any sense other than
  $|C|=K(\emptyset,\emptyset)$, it does not use that the kernel is
  optimal, and it does not use the structure of $K$ beyond positivity
  and the two families of inequalities. What it needs, and what a
  bound of $24.10$ can never give, is that the number
  $K(\emptyset,\emptyset)$ is attained by a code; that is why the
  certificate has to have objective value exactly $24$ and not merely
  below $25$, and why its rounding to an exact rational point with
  objective exactly $24$ is the decisive step of its construction.}
\end{proof}

The kernels to which the lemma is applied have a particular form, and
their positivity is a matter of that form alone.

\begin{lem}[Gram kernels are positive]\label{lem:gram}
Let $\Lambda$ be a finite set, let $H$ be a real or complex inner product
space, and for each $\lambda\in\Lambda$ let
$\psi_{\lambda,1},\dots,\psi_{\lambda,n_\lambda}$ be maps from
$\mathcal I_2$ to $H$. Put
\[
  Z_\lambda(J_1,J_2)_{ab}\;=\;\bigl\langle\psi_{\lambda,a}(J_1),\psi_{\lambda,b}(J_2)\bigr\rangle ,
\]
the inner product being conjugate linear in its first argument when $H$
is complex, and let $X_\lambda$ be a real symmetric positive semidefinite
matrix of order $n_\lambda$ for each $\lambda$. Then the function
\[
  K(J_1,J_2)\;=\;\sum_{\lambda\in\Lambda}\ \sum_{a,b=1}^{n_\lambda}(X_\lambda)_{ab}\,Z_\lambda(J_1,J_2)_{ab}
\]
satisfies $K(J_2,J_1)=\overline{K(J_1,J_2)}$, and its real part is a
positive kernel on $\mathcal I_2$; when $K$ is real valued, $K$ itself is
one.%
\footnote{Nothing about the maps $\psi_{\lambda,a}$ is used except that
  they take values in an inner product space, so the positivity of the
  certificate kernel does not depend on the representation theory that
  produces the zonal matrices, only on the zonal matrices being Gram
  matrices, which is their definition \eqref{eq:zonal}. The
  representation theory is what makes the entries of $Z_\lambda$
  polynomials in the inner products and what makes every invariant
  positive kernel expressible in this form; the second of these
  properties, the completeness of the family, is what an optimiser
  needs and is not needed here at all, since only one kernel is ever
  used and it is given.}
\end{lem}

\begin{proof}
Write $X_\lambda=R_\lambda^{\mathsf T}R_\lambda$ with $R_\lambda$ real, which
is possible for a real symmetric positive semidefinite matrix. For
$J_1,\dots,J_m\in\mathcal I_2$ and real numbers $c_1,\dots,c_m$,
\begin{align*}
  \sum_{i,j}c_ic_j\,K(J_i,J_j)
  &=\sum_{\lambda}\sum_{c}\Bigl\langle\,\sum_{i,a}c_i(R_\lambda)_{ca}\psi_{\lambda,a}(J_i),\ 
  \sum_{j,b}c_j(R_\lambda)_{cb}\psi_{\lambda,b}(J_j)\Bigr\rangle\\
  &=\sum_{\lambda,c}\Bigl\|\sum_{i,a}c_i(R_\lambda)_{ca}\psi_{\lambda,a}(J_i)\Bigr\|^2\;\ge\;0 ,
\end{align*}
using $(X_\lambda)_{ab}=\sum_c(R_\lambda)_{ca}(R_\lambda)_{cb}$ and the
sesquilinearity of the inner product. The symmetry
$K(J_2,J_1)=\overline{K(J_1,J_2)}$ follows from
$\langle v,w\rangle=\overline{\langle w,v\rangle}$ and the symmetry of
$X_\lambda$. The displayed quantity is real and nonnegative, so it equals
its real part, which is $\sum_{i,j}c_ic_j\operatorname{Re}K(J_i,J_j)$; that
is positivity of $\operatorname{Re}K$, and of $K$ when the two agree.
\end{proof}

\subsubsection{The certificate kernel}
\label{sec:certificate-kernel}

The kernel that settles twenty-four points is of the form in
\cref{lem:gram}, with the index set $\Lambda$ the sixty signatures
$\lambda=(\lambda_1,\lambda_2)$, $\lambda_1\ge\lambda_2\ge0$,
$\lambda_1+\lambda_2\le14$, of irreducible representations of
$\mathrm O(4)$, and with the maps $\psi_{\lambda,a}$ the functions that
\cite[\S2--3]{LLM24} attach to a subset $J$ of $S^3$ of size at most
$2$: an irreducible representation $\rho_\lambda$ of
$\mathrm{GL}(2,\mathbb C)$ on a space $W_\lambda$ carrying a Hermitian inner
product, a fixed vector $w_k\in W_\lambda$ for each admissible index
$a=(i,j,k)$, a section $s(J)\in\mathrm O(4)$ carrying a reference
configuration $\epsilon$ onto $J$, a fixed $2\times4$ complex matrix
$\omega$, and a scalar $\xi_{\lambda,a}(J)$ which is a polynomial in the
inner product of the two points of $J$ times powers of
$\sqrt{2\pm2\langle x,y\rangle}$, and then
\begin{equation}\label{eq:psi}
  \psi_{\lambda,a}(J)\colon\ \gamma\ \longmapsto\ \xi_{\lambda,a}(J)\,
  \rho_\lambda\bigl(\omega\gamma\,s(J)\,\epsilon\bigr)\,w_k ,
  \qquad \gamma\in\mathrm O(4),
\end{equation}
an element of the space $H=L^2(\mathrm O(4),W_\lambda)$ of square
integrable $W_\lambda$-valued functions on $\mathrm O(4)$ with the inner
product $\langle f_1,f_2\rangle=\int_{\mathrm O(4)}\langle f_1(\gamma),f_2(\gamma)\rangle\,d\gamma$
against Haar probability measure. The \emph{zonal matrices} are the Gram
matrices
\begin{equation}\label{eq:zonal}
\begin{aligned}
  Z_\lambda(J_1,J_2)_{ab}&=\bigl\langle\psi_{\lambda,a}(J_1),\psi_{\lambda,b}(J_2)\bigr\rangle\\
  &=\xi_{\lambda,a}(J_1)\,\xi_{\lambda,b}(J_2)\int_{\mathrm O(4)}
  \bigl\langle\rho_\lambda(\omega\gamma s(J_1)\epsilon)w_{k_1},\ \rho_\lambda(\omega\gamma s(J_2)\epsilon)w_{k_2}\bigr\rangle\,d\gamma ,
\end{aligned}
\end{equation}
which are invariant under the diagonal action of $\mathrm O(4)$ on
$(J_1,J_2)$ because Haar measure is, and whose entries are polynomials in
the inner products among the points of $J_1\cup J_2$; the scalars
$\xi$ are chosen so that this is so, and the $\sqrt{2\pm2u}$ they carry
are what absorb the square roots that the section $s(J)$ introduces. This
is the construction of \cite{LLM24}, which extends the one of
\cite{dLV15} and of the equiangular lines computations that preceded it,
and it is recomputed from scratch in \cref{sec:certificate-checked}. The
\emph{certificate kernel} is
\begin{equation}\label{eq:Kcert}
  K(J_1,J_2)\;=\;\sum_{\lambda\in\Lambda}\ \sum_{a,b}(X_\lambda)_{ab}\,Z_\lambda(J_1,J_2)_{ab},
\end{equation}
where the $X_\lambda$ are the sixty rational matrices of the deposited
data \cite{LLM24data}. Its values are rational polynomials in the inner
products, so it is real valued, and it is $\mathrm O(4)$-invariant. For
$Q\in\mathcal I_4$ the quantity $A_2K(Q)$ of \eqref{eq:A2} is therefore a
polynomial in the inner products among the points of $Q$; call it
$p_1$, $p_2(u)$, $p_3(u_1,u_2,u_3)$ and $p_4(u_1,\dots,u_6)$ according
to $|Q|=1,2,3,4$, the $u$ being the inner products of the pairs of $Q$ in a
fixed order. Since a set of at most four points of $\mathbb R^4$ is
determined up to an orthogonal map by its Gram matrix, these polynomials
determine $A_2K$ on all of $\mathcal I_4$.%
\footnote{Four points of $\mathbb R^4$ with a given Gram matrix $G$ are
  the columns of a $4\times4$ matrix $P$ with $P^{\mathsf T}P=G$, and any
  two such matrices differ by an orthogonal factor on the left, since
  $P_2=P_1U$ with $U$ orthogonal whenever $P_1^{\mathsf T}P_1=P_2^{\mathsf
  T}P_2$; for fewer than four points the same holds after adding zero
  columns. This is why the constraints of the second level are
  polynomial identities in one, three and six inner products and why
  the domain of each is described by the positivity of a Gram matrix.
  It is also why the level stops at four: five points of $\mathbb R^4$
  have Gram matrices constrained by a determinant condition, and the
  third level would need it.}

\begin{prop}[What the verification establishes]\label{prop:verified}
The following statements about the deposited data hold, each by an exact
computation described in \cref{sec:certificate-checked}.
\begin{enumerate}[label=\textup{(\roman*)},leftmargin=2.6em,itemsep=0.12em]
\item Each of the sixty matrices $X_\lambda$, and each of the sixty-seven
      further matrices $M$ that enter the sums of squares of
      \textup{(iii)}, is a real symmetric positive definite matrix with
      rational entries.
\item $K(\emptyset,\emptyset)=24$.
\item There are polynomials $\sigma_1,\dots,\sigma_4$, each a sum of terms
      $w\cdot\langle M,\,m\,m^{\mathsf T}\rangle$ with $M$ one of the
      matrices of \textup{(i)}, $m$ a vector of monomials in the inner
      products of $Q$, and $w$ a weight which is a nonnegative constant, or
      $(u+1)(\tfrac12-u)$ for an inner product $u$ of $Q$, or the
      determinant of the Gram matrix of $Q$, or an elementary symmetric
      function of the family of the $(u+1)(\tfrac12-u)$ over the inner
      products of $Q$, or, for $|Q|=4$, an elementary symmetric function
      of the four $3\times3$ principal minors of the Gram matrix of $Q$,
      such that the four identities
      \[
        p_1+1+\sigma_1=0,\qquad p_2+\sigma_2=0,\qquad p_3+\sigma_3=0,\qquad p_4+\sigma_4=0
      \]
      hold identically as polynomials in the inner products.
\item The polynomial $\sigma_2$ has degree $16$, and its zeros in
      $[-1,\tfrac12]$ are exactly $-1$, $-\tfrac12$, $0$ and $\tfrac12$, of
      multiplicities $1$, $2$, $2$ and $1$.
\end{enumerate}
\end{prop}

Every weight in (iii) is nonnegative at every $Q\in\mathcal I_4$: the
inner products of $Q$ lie in $[-1,\tfrac12]$, so each $(u+1)(\tfrac12-u)$
is nonnegative; the Gram matrix of $Q$ is positive semidefinite, so its
determinant and its principal minors are nonnegative; and an elementary
symmetric function of nonnegative quantities is nonnegative. Each
$\langle M,mm^{\mathsf T}\rangle=m^{\mathsf T}Mm$ is nonnegative because
$M$ is positive definite. Hence $\sigma_k(Q)\ge0$ on $\mathcal I_4$ for
$k=1,\dots,4$, and (iii) says
\begin{equation}\label{eq:constraints-hold}
  A_2K(\{x\})=-1-\sigma_1\le-1,\qquad A_2K(Q)=-\sigma_{|Q|}(Q)\le0\quad(2\le|Q|\le4),
\end{equation}
for every $x\in S^3$ and every $Q\in\mathcal I_4$ with $|Q|\ge2$.

\begin{figure}[tb]
\centering
\includegraphics[width=\textwidth]{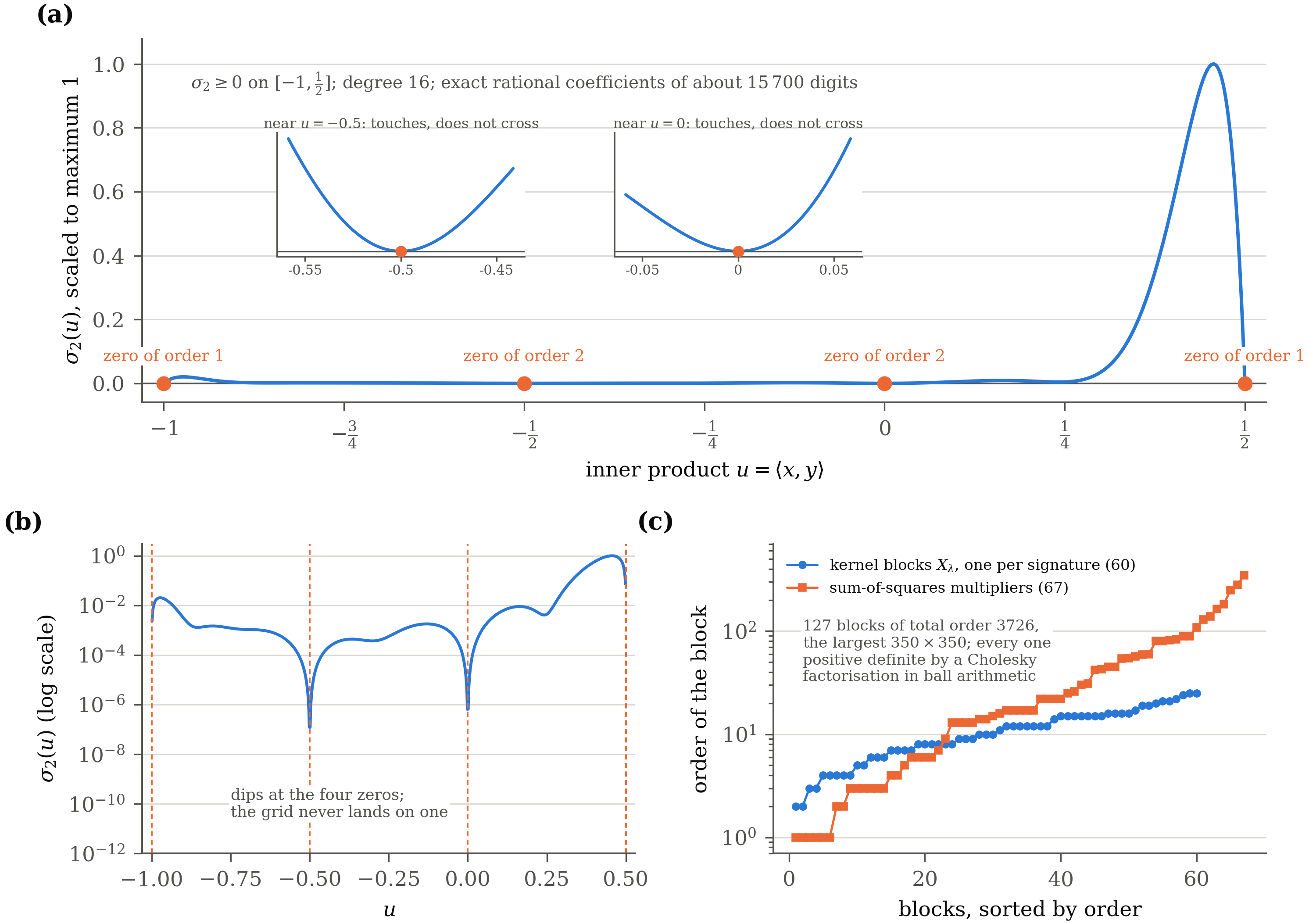}
\caption{The certificate, as \cref{prop:verified} reads it. (a) The
two-point polynomial $\sigma_2$ on $[-1,\tfrac12]$, evaluated exactly
from the deposited rational coefficients and scaled by its maximum: it is
nonnegative, it vanishes to first order at $-1$ and $\tfrac12$ and to
second order at $-\tfrac12$ and $0$, and nowhere else; the insets show the
two double zeros, where the curve touches the axis without crossing it.
(b) The same on a logarithmic scale. (c) The orders of the $127$
positive definite blocks, the sixty kernel blocks $X_\lambda$ and the
sixty-seven sum-of-squares multipliers; the largest is $350\times350$.}
\label{fig:sigma2}
\end{figure}

\begin{thm}[Twenty-four points]\label{thm:twentyfour-points}
Let $C\subset S^3$ be a set whose pairwise inner products are at most
$\tfrac12$. Then $|C|\le24$. If $|C|=24$, then every pairwise inner
product of $C$ lies in $\{-1,-\tfrac12,0,\tfrac12\}$.
\end{thm}

\begin{proof}
By \cref{prop:verified}(i) and \cref{lem:gram}, the certificate kernel
\eqref{eq:Kcert} is positive; it is real valued, so no real part has to
be taken. By \eqref{eq:constraints-hold} it satisfies the hypotheses of
\cref{lem:las2}, and by \cref{prop:verified}(ii) that lemma gives
$|C|\le K(\emptyset,\emptyset)=24$. If $|C|=24$ the equality case of the
lemma gives $A_2K(\{x,y\})=0$ for every pair of distinct points
$x,y\in C$, that is $\sigma_2(\langle x,y\rangle)=0$ by
\eqref{eq:constraints-hold}, and \cref{prop:verified}(iv) places
$\langle x,y\rangle$ in $\{-1,-\tfrac12,0,\tfrac12\}$.
\end{proof}

\begin{figure}[tb]
\centering
\includegraphics[width=\textwidth]{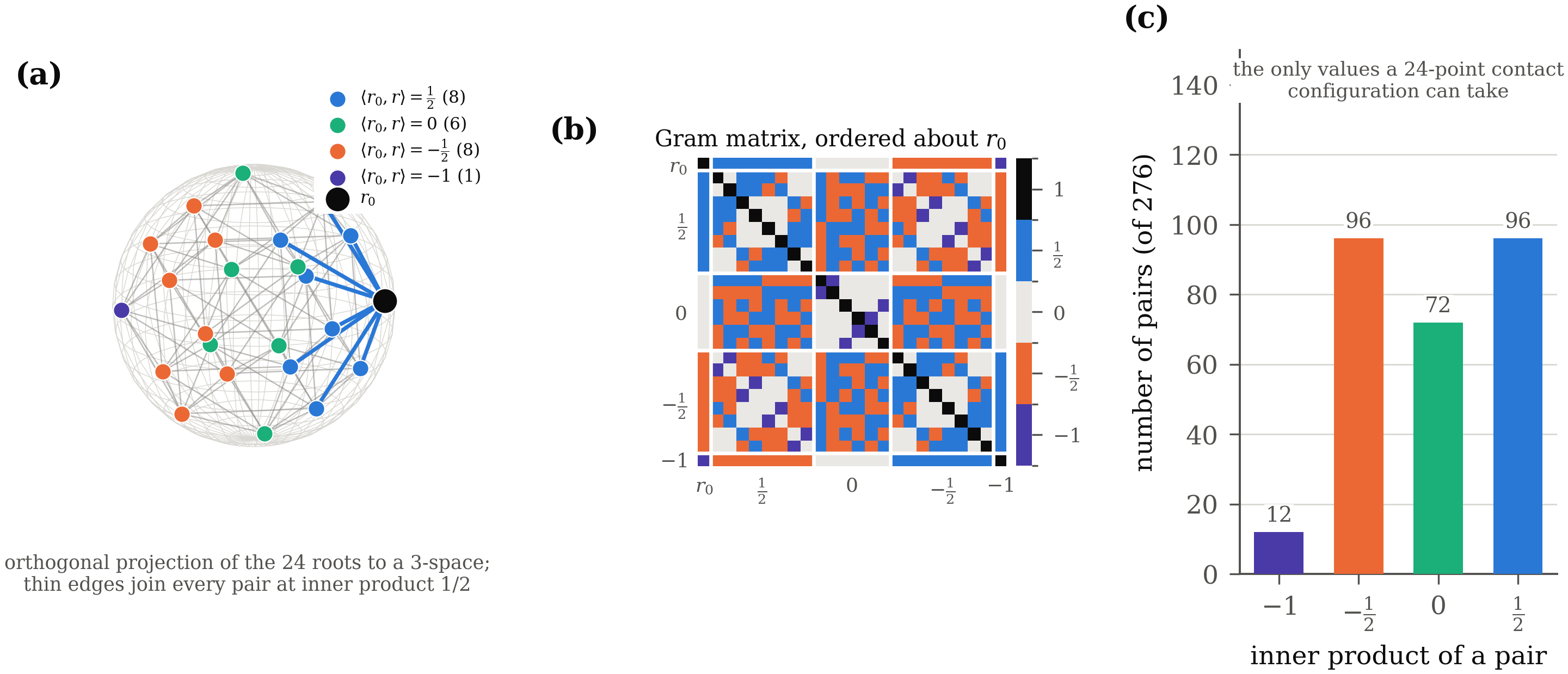}
\caption{What \cref{thm:twentyfour-points} forces. (a) The twenty-four
normalised roots of $D_4$, projected orthogonally to a generic $3$-space
so that no two coincide, with every pair at inner product $\tfrac12$
joined; one root $r_0$ is singled out and the others are coloured by
their inner product with it, eight at $\tfrac12$, six at $0$, eight at
$-\tfrac12$ and the antipode at $-1$. (b) The Gram matrix of the
twenty-four directions, rows and columns ordered by these classes; it
takes no value outside $\{1,\tfrac12,0,-\tfrac12,-1\}$. (c) The $276$
pairs by inner product: $96$, $72$, $96$ and $12$. The theorem says that
every contact configuration of size $24$ has a Gram matrix with the
values of (b), and \cref{lem:root-lattice} that it then is (a) up to an
orthogonal map.}
\label{fig:d4contacts}
\end{figure}

The first assertion is Musin's theorem \cite{Mus08} that the kissing
number in four dimensions is $24$, obtained here a second time from a
different certificate. It is the second assertion that the present paper
needs, and it is what the earlier bound $24.10$ of Bachoc and Vallentin
\cite{BV08}, at the first level of the same hierarchy, and Musin's
argument, which does not identify the configuration, do not give: it is
the equality case, and the equality case is available only because the
bound is exactly $24$ and not merely below $25$.

\begin{proof}[Proof of \cref{thm:m24}]
A contact configuration is a subset of $S^3$ whose pairwise inner
products are at most $\tfrac12$, that is a spherical code of minimal
angle at least $\pi/3$. By \cref{thm:twentyfour-points} the pairwise
inner products of its $24$ points lie in $\{-1,-\tfrac12,0,\tfrac12\}$, and
\cref{lem:root-lattice} makes it a copy of the root system. The cell of
the root configuration is the $24$-cell of volume $8$ by
\cref{lem:refcell-vol}, and volume is invariant under orthogonal maps.
\end{proof}

\subsubsection{The certificate, checked}
\label{sec:certificate-checked}

What remains is \cref{prop:verified}, which is a statement about a finite
set of rational numbers, and we are exact about where those numbers come
from and what was done with them. The certificate is the discovery of
de Laat, Leijenhorst and de Muinck Keizer \cite{LLM24}: they computed the
second level of the hierarchy of \cite{dLV15} for spherical codes in
$S^3$ with maximal inner product $\tfrac12$, reduced by the symmetry of
the problem to the form \eqref{eq:Kcert} with truncation degrees $14$
and $16$, solved it to forty digits in $256$-bit floating-point
arithmetic by the clustered low-rank solver of \cite{LdL24}, about two
weeks on eight cores with $128$ gigabytes of memory, rounded the
numerical optimum to an exact rational feasible point with objective
value exactly $24$, and deposited that point, with their own
verification code, at 4TU.ResearchData \cite{LLM24data}, under the MIT
licence. The archive is \texttt{LasserreSphericalCodes.zip}, MD5
\texttt{02acd5270f7b3fa799abdeb5291706fd}, $234$ megabytes in $331$ files
under \texttt{proofs/4\_24}; recording the checksum fixes which bytes the
statements below refer to, and the licence permits the files to be
redistributed, which the supplementary archive of this paper
\cite{D4Zenodo} does, with the licence notice, so that every input to
\cref{prop:verified} travels with the paper.%
\footnote{The point of redistributing rather than pointing is that a
  proof which depends on a file should not depend on a server. The MIT
  licence asks only that its notice accompany the copy, and it does.
  The checksum is recorded so that the copy in the supplementary
  archive, the copy at 4TU.ResearchData and the copy any reader
  downloads can be compared byte for byte; the verification scripts
  check it before they read anything.} Their paper is, at the time
of writing, a preprint; nothing below depends on it.

We have verified \cref{prop:verified} from the deposited data in code of
our own that shares nothing with theirs
(\texttt{llm24\_certificate\_check.py}, Python with the exact rational
arithmetic of FLINT and the ball arithmetic of Arb \cite{Joh17}, the
programs of the directory \texttt{zonal}, and
\texttt{D4InnerProducts.lean}). The verification has seven steps, and
we describe them in the order in which they are run. The data
are read and their format checked: $127$ positive semidefinite blocks
of total dimension $3726$, the largest $350\times350$, sixty of them
indexed by the signatures $\lambda$ of $O(4)$ with $|\lambda|\le14$ and
each of those with exactly as many rows as the paper's description of
the admissible index tuples prescribes, the rest the multipliers of the
sums of squares for the constraints on two, three and four points;
every one of the $127$ blocks is positive definite, by a Cholesky
factorisation in ball arithmetic at $256$ bits in which every pivot is
a ball inside the positive reals (the least about $1.4\times10^{-15}$),
and the $81$ blocks of size at most $16$ also by an exact $LDL^{\mathsf
T}$ over the rationals; each of the $125$ prefactors of the
sum-of-squares terms is a nonnegative multiple of one of the weights
that describe the domain, so that every such term is nonnegative on it;
the objective, the entry of the block $\lambda=(0,0)$ that is the
value of the bound, is exactly $24$; and the two-point polynomial
$p_2$, computed exactly from the data, has degree $16$ and vanishes on
$[-1,\tfrac12]$ at $-1,-\tfrac12,0,\tfrac12$ and nowhere else, with
multiplicities $1,2,2,1$, by a Sturm sequence over the rationals, the
last step being checked by Lean's kernel as well.

The two steps that the memory had put out of reach are done here as well.
Step 3 is the construction of the zonal matrices $Z_\lambda$, the polynomials
in the six inner products of four points through which the kernel is
expressed; step 5 is the check, which uses them, that the four constraint
polynomials $A_2K(Q)+\text{(sum of squares)}$ are identically equal to their
right-hand sides. The authors' implementation of step 3 takes about three days
and $128$ gigabytes. That memory is not intrinsic to it. What step 3 wants,
for each signature $\lambda$ and each admissible pair of indices $k_1,k_2$, is
\[
  P(S)_{k_1,k_2}\;=\;\int_{O(4)}\rho_{0,k_1}(\omega\gamma\varepsilon)\,
  \rho_{0,k_2}(\omega\gamma S)\,d\gamma,
\]
and the integral is taken monomial by monomial in the sixteen entries of
$\gamma$. Their code expands the product of the two representation
polynomials first and holds the expansion in a table indexed by the exponent
matrix of $\gamma$, with a polynomial in $S$ against each index; that table is
the $128$ gigabytes. It need never be formed. The answer is homogeneous of
degree $|\lambda|\le14$ in the seven free entries of $S$, so it has at most
$\binom{20}{6}=38\,760$ coefficients whatever $\lambda$ is, and each of the
two factors of the integrand is a product of powers of single entries of a
$2\times2$ matrix; so their monomials can be walked in triples, each triple
added straight into the coefficient of $P(S)$ it belongs to, with nothing held
but the factors themselves. Two facts keep the walk short. The integral
vanishes unless every row sum and every column sum of the exponent matrix is
even, which leaves $7.85\times10^{9}$ of the $2.39\times10^{11}$ triples; and
it is unchanged by permuting rows or columns, so its values can be stored
against the canonical form of the matrix, of which at most two and a half
million occur in any one of the $490$ entries. Carried out this way, in exact
rational arithmetic, the whole of step 3 takes about two hours on two cores
and never needs more than $400$ megabytes.

\begin{figure}[tb]
\centering
\includegraphics[width=\textwidth]{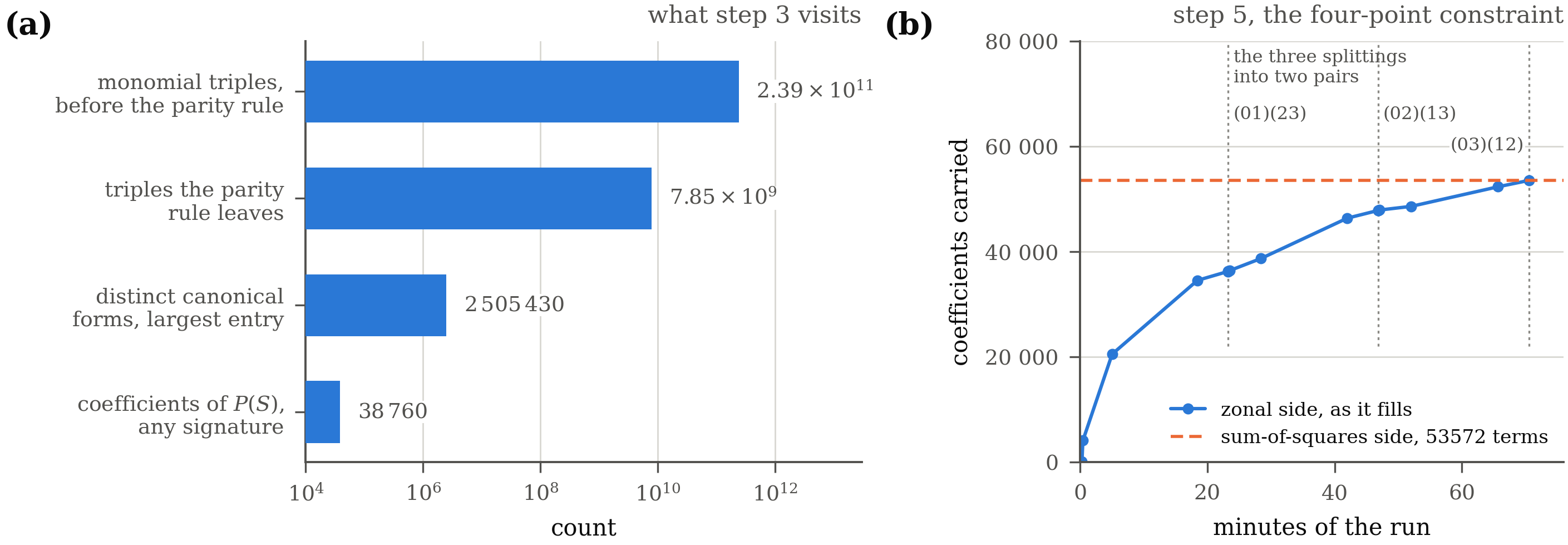}
\caption{The two steps of the verification of the certificate that its
    authors' own code \cite{LLM24data} leaves out, as they were carried out
    here. (a) Why step 3 fits in a few
    hundred megabytes. The integral over $O(4)$ is taken one monomial triple
    at a time; the parity rule on the exponent matrix discards all but about a
    thirtieth of the triples, the invariance under permuting rows and columns
    collapses those onto at most $2\,505\,430$ canonical forms in the largest
    of the $490$ entries, and the answer never has more than
    $\binom{20}{6}=38\,760$ coefficients whatever the signature. The first two
    counts are the counters in the kernel; the third is the size of its memo
    table. (b) Step 5 for the four-point constraint, the only one large enough
    to be difficult. The curve is the zonal side of the constraint filling up,
    one signature at a time, through the three ways of splitting four points
    into two pairs; it stops at $53\,572$ coefficients, which is where the
    sum-of-squares side stops as well. The two sides sit on the same monomials
    and each coefficient of one is the exact negative of the other, so the
    difference is zero at every one of the $53\,572$ places.}
\label{fig:verify}
\end{figure}

Step 5 then follows. The four constraint polynomials are assembled from the
zonal matrices and from the certificate's own blocks and factored
sum-of-squares data, the sum \eqref{eq:A2} being taken over ordered
pairs, which gives it $3$, $9$, $12$ and $6$ terms for $|Q|=1,2,3,4$,
an unordered pair $J_1\ne J_2$ entering with the factor $2$ that its two
orders supply and a pair $J_1=J_2$ once; and all four are identically
zero: the one-point
constraint, the two-point constraint in one variable, the three-point
constraint in three variables and the four-point constraint in six, that last
one assembled out of pieces carrying up to seventy-five thousand coefficients
of several thousand digits. The two sides of that last constraint end on the
same $53\,572$ monomials, and on each one the coefficient from the
sum-of-squares blocks is the exact negative of the coefficient from the zonal
matrices, both about $15\,700$ digits long, so the constraint is zero term by
term and not by any cancellation of rounding (\cref{fig:verify}).
None of this shares code with the authors' Julia package; the inputs are the
deposited data and the description of the construction in their paper. Three
checks that do not depend on the outcome fix the construction itself. The
integration routine reproduces the seven values of the test distributed with
their own integration code, and satisfies seventy further identities that
follow from the rows and the columns of an orthogonal matrix being unit
vectors. For the signatures $\lambda=(k,0)$ and two single points at inner
product $u$ the matrix $Z_\lambda$ has to be a multiple of the zonal spherical
harmonic kernel of $S^3$; it is, exactly, with ratio $8^k/(k+1)^2$ for every
$k\le14$, the denominator being the dimension of the space of harmonics of
degree $k$ on $S^3$. And $Z_\lambda$ has to be a positive semidefinite kernel
on the subsets of size at most two of any set of points, which nothing in the
computation forces; on random four-point sets the least eigenvalue is at the
precision of the arithmetic for every signature, including those with
$\lambda_2>0$, which the comparison with the Gegenbauer polynomials does not
reach.

The seven steps establish \cref{prop:verified} item by item. Step 2 is
(i): every block is positive definite, by a Cholesky factorisation in
ball arithmetic whose pivots are balls inside the positive reals, and the
small ones also by an exact $LDL^{\mathsf T}$ over the rationals. Step 6 is
(ii). Steps 3, 4 and 5 are (iii): the zonal matrices are built from the
definition \eqref{eq:zonal}, the prefactors are checked to be
nonnegative multiples of the domain weights, and the four identities are
checked to hold with every coefficient exactly zero. Step 7 is (iv), by a
Sturm sequence over the rationals, and the last statement is checked by
Lean's kernel as well. Step 1 is the reading of the data. What the
computation therefore establishes is exactly the list of
\cref{prop:verified}; \cref{lem:las2,lem:gram} are what turn that list
into \cref{thm:twentyfour-points}, and they are proved above. Nothing in
the argument is read from \cite{LLM24} rather than run or proved here: the
authors of the certificate proved the same bound and the same equality
case in their Sections~2 to~5, and a reader who wants a second account
of \cref{lem:las2} will find it there, but the present proof does not
cite it.

Everything after \cref{thm:twentyfour-points}, the passage from the four
inner products to the root system, is \cref{lem:root-lattice}, proved
above by a classical argument and checked by an exact enumeration, in
Python and in Lean. Every
statement below that invokes \cref{thm:m24} therefore rests on
\cref{prop:verified}, that is on a finite exact computation reproducible
from the supplementary archive, and on nothing outside this paper. All
the other results in this paper are independent even of that.

\begin{cor}[The multi-direction conjecture]\label{cor:conj-resolved}
\Cref{conj:multidir-new} holds. Indeed a configuration at the
all-contact corner with a full active set is a contact configuration of
size $24$, so by \cref{thm:m24} it is a copy of the root system; if the
$24-m$ undeviated directions are the roots $\base_i$, $i\notin D$, the
orthogonal map fixing them is the identity, so every tilt $\tilt_j$
vanishes and $\vol(V_c)=8$ with equality.
\end{cor}

\subsection{A volume identity exact at the root system}\label{sec:m24-volume}

The truncated estimates of \cref{sec:covering-bound} and
\cref{sec:strict-inequality} cut the cell at a radius below $\sqrt2$.
At $m=24$ no such truncation can be sharp: the $24$-cell has
circumradius $\sqrt2$ (\cref{lem:refcell-vol}), so a ball of smaller
radius removes volume from it, and the lower bound falls strictly below
the true value $8$ at the root configuration. This subsection truncates
at $\sqrt2$ itself and keeps three terms of inclusion and exclusion. The
resulting lower bound is exact at the root configuration, and it
converts the statement $\vol(V_c)\ge8$ for $24$ contacts into one
inequality between the pair and triple angles of the contact
directions, of the same form as the inequality proved for $23$ contacts
in \cref{sec:certificate}. \Cref{rem:m24-three-point} records what the
three-point relaxation of that inequality does and does not deliver.

Throughout, $B=\{x\in\mathbb{R}^4:|x|\le\sqrt2\}$, and for a unit vector
$w$ we write $K(w)=\{x\in B:\langle x,w\rangle>1\}$. By orthogonal
invariance the numbers
\begin{align*}
  c&=\vol K(w),\\
  \omega_2(u)&=\vol\bigl(K(w)\cap K(w')\bigr),\\
  \omega_3(u,v,t)&=\vol\bigl(K(w)\cap K(w')\cap K(w'')\bigr)
\end{align*}
depend only on the inner products $u=\langle w,w'\rangle$ and, for
$\omega_3$, on
\[
  (u,v,t)=\bigl(\langle w,w'\rangle,\ \langle w,w''\rangle,\
  \langle w',w''\rangle\bigr).
\]

\begin{lem}[Third-order bound]\label{lem:bonferroni3}
Let $w_1,\dots,w_m$ be unit vectors in $\mathbb{R}^4$, write
$u_{ij}=\langle w_i,w_j\rangle$ and
$V=\{x:\langle x,w_i\rangle\le1\ \text{for all } i\}$. Then
\begin{equation}\label{eq:bonferroni3}
  \vol(V)\;\ge\;2\pi^2-mc+\sum_{i<j}\omega_2(u_{ij})
  -\sum_{i<j<k}\omega_3(u_{ij},u_{ik},u_{jk}),
\end{equation}
with equality if and only if $V\setminus B$ and the set of points of $B$
lying in at least four of the sets $K(w_i)$ both have measure zero.
\end{lem}

\begin{proof}
First $\vol(V)\ge\vol(V\cap B)=\vol(B)-\vol\bigl(\bigcup_iK(w_i)\bigr)$,
since $B\setminus V=\bigcup_iK(w_i)$, and $\vol(B)=\tfrac12\pi^2(\sqrt2)^4
=2\pi^2$. For a point lying in exactly $n\ge1$ of the sets, the identity
$\sum_{j=0}^{3}(-1)^j\binom nj=-\binom{n-1}{3}$ gives
\[
  \binom n1-\binom n2+\binom n3\;=\;1+\binom{n-1}3\;\ge\;1,
\]
with equality exactly when $n\le3$. So the indicator of the union is
bounded above, pointwise, by the sum of the indicators of the sets, less
the sum over pairs, plus the sum over triples, with equality off the
points lying in four or more of the sets. Integration over $B$ gives
\eqref{eq:bonferroni3}, and the two inequalities used are equalities
exactly under the stated conditions.
\end{proof}

\begin{lem}[The cap and the pair term]\label{lem:pair-volume}
The cap volume is $c=\tfrac12\pi^2-\tfrac43\pi$. For $-1\le u\le0$ one has
$\omega_2(u)=0$, and $\omega_3(u,v,t)=0$ unless $u$, $v$ and $t$ are all
positive. For $0<u<1$, with $\gamma=\arccos u$,
\begin{align}\label{eq:pair-volume}
  \omega_2(u)&=2\pi\int_{\gamma/2}^{\pi/4}\Bigl(1-\tfrac12\sec^2\theta\Bigr)^2d\theta
  =2\pi\bigl(F(\tfrac\pi4)-F(\tfrac\gamma2)\bigr),\\
  F(\theta)&=\theta-\tfrac34\tan\theta+\tfrac1{12}\tan^3\theta .\nonumber
\end{align}
In particular $\omega_2(\tfrac12)=\pi(9\pi+26\sqrt3-72)/54=0.0760762770\ldots$.
\end{lem}

\begin{proof}
Write $x=(y,z)$ with $y=\langle x,w\rangle$ and $z$ in the orthogonal
complement of $w$. Then $K(w)$ is the set $1<y\le\sqrt2$,
$|z|^2\le2-y^2$, and
\begin{align*}
  c&=\int_1^{\sqrt2}\tfrac43\pi(2-y^2)^{3/2}\,dy
  =\tfrac{16}3\pi\int_{\pi/4}^{\pi/2}\cos^4\phi\,d\phi\\
  &=\tfrac{16}3\pi\Bigl(\tfrac{3\pi}{32}-\tfrac14\Bigr),
\end{align*}
by the substitution $y=\sqrt2\sin\phi$. If $x\in K(w)\cap K(w')$ then
$\langle x,w+w'\rangle>2$ and $|w+w'|^2=2+2u$, so
$|x|^2>2/(1+u)$; for $u\le0$ this exceeds $2$ and the intersection is
empty, and the same argument applied to any two of three directions
gives the statement on $\omega_3$. For $0<u<1$ decompose $x=y+z$ with
$y$ in the plane $P$ spanned by $w,w'$ and $z\perp P$. The conditions are
$\langle y,w\rangle>1$, $\langle y,w'\rangle>1$ and $|z|^2\le2-|y|^2$,
so $\omega_2(u)$ is the integral of $\pi(2-|y|^2)_+$ over the region of
$P$ cut out by the first two. In polar coordinates on $P$ with $w$ at
angle $0$ and $w'$ at angle $\gamma$, the ray at angle
$\theta\in(\gamma-\tfrac\pi2,\tfrac\pi2)$ meets that region in
$\rho>\max(\sec\theta,\sec(\gamma-\theta))$, and the reflection
$\theta\mapsto\gamma-\theta$ exchanges the two terms, so it suffices to
take $\theta\ge\gamma/2$, where the maximum is $\sec\theta$, and to
double. Since
\[
  \int_{\sec\theta}^{\sqrt2}\pi(2-\rho^2)\rho\,d\rho
  =\pi\Bigl(1-\tfrac12\sec^2\theta\Bigr)^2
\]
for $\sec\theta\le\sqrt2$, while for $\theta>\pi/4$ the ray does not
meet the region inside $B$, \eqref{eq:pair-volume} follows, and $\gamma/2<\pi/4$
because $u>0$. Differentiation gives
$F'(\theta)=1-\sec^2\theta+\tfrac14\sec^4\theta=(1-\tfrac12\sec^2\theta)^2$.
At $\gamma=\pi/3$, $F(\tfrac\pi4)=\tfrac\pi4-\tfrac23$ and
$F(\tfrac\pi6)=\tfrac\pi6-\tfrac{13}{18\sqrt3}$.
\end{proof}

\begin{prop}[Exactness at the root system]\label{prop:m24-exact}
For the $24$ normalised roots of $D_4$, \eqref{eq:bonferroni3} holds with
equality, and
\begin{equation}\label{eq:m24-identity}
  8\;=\;32\pi-10\pi^2+96\,\omega_2(\tfrac12)-96\,\omega_3(\tfrac12,\tfrac12,\tfrac12).
\end{equation}
Consequently
$\omega_3(\tfrac12,\tfrac12,\tfrac12)=\tfrac1{16}\pi^2-\pi+\tfrac{13\sqrt3}{27}\pi-\tfrac1{12}
=0.0118567031\ldots$.
\end{prop}

\begin{proof}
The cell is the $24$-cell, of volume $8$, whose vertices all have norm
$\sqrt2$ (\cref{lem:refcell-vol}); so it lies in $B$. By
\cref{lem:pair-volume} a point of $K(w_i)\cap K(w_j)$ forces
$u_{ij}>0$, and the only positive inner product between distinct
normalised roots is $\tfrac12$. Suppose a point of $B$ lies in four of
the sets. The four corresponding roots, scaled to norm $\sqrt2$, then
have Gram matrix $I+J$, with $2$ on the diagonal and $1$ off it, whose
determinant is $5$. They are therefore linearly independent and span a
sublattice $L$ of full rank of the $D_4$ lattice, whose determinant is
$4$; but $\det L=4\,[D_4:L]^2$, and $5$ is not of that form. So every
point of $B$ lies in at most three of the sets, and
\cref{lem:bonferroni3} gives equality. In the sum over pairs only the
pairs at inner product $\tfrac12$ contribute; each root has eight such
neighbours ($e_1+e_2$ has $e_1\pm e_3$, $e_1\pm e_4$, $e_2\pm e_3$,
$e_2\pm e_4$), giving $96$ pairs. In the sum over triples only triangles
of roots pairwise at $\tfrac12$ contribute; each such pair lies in
exactly three of them (for $e_1+e_2$ and $e_1+e_3$ the third root is
$e_1+e_4$, $e_1-e_4$ or $e_2+e_3$), giving $96\cdot3/3=96$ triangles.
With $2\pi^2-24c=32\pi-10\pi^2$ from \cref{lem:pair-volume} this is
\eqref{eq:m24-identity}, and the value of $\omega_3$ follows from the
value of $\omega_2(\tfrac12)$.
\end{proof}

The value of $\omega_3(\tfrac12,\tfrac12,\tfrac12)$ obtained this way
agrees with a direct quadrature of the triple intersection to
$1\times10^{-7}$, and the cap, pair and triple volumes agree with
Monte Carlo estimates in $\mathbb{R}^4$ within their standard errors
(\texttt{m24\_exact\_reduction.py}, \cref{app:code-index}).

\begin{cor}[The $24$-contact case as an inequality between angles]\label{cor:m24-reduction}
For a contact configuration $W=\{w_1,\dots,w_{24}\}$ put
\[
  E(W)=\sum_{i<j}\omega_2(u_{ij})-\sum_{i<j<k}\omega_3(u_{ij},u_{ik},u_{jk}).
\]
Then $\vol(V_c(W))\ge32\pi-10\pi^2+E(W)$. The root configuration has
$E=\tau$, where
\[
  \tau=8-32\pi+10\pi^2=6.1650790960\ldots,
\]
and every $W$ with $E(W)\ge\tau$ has $\vol(V_c(W))\ge8$.
\end{cor}

\begin{proof}
\Cref{lem:bonferroni3} with $m=24$, and \cref{prop:m24-exact}.
\end{proof}

\begin{remark}[The three-point relaxation of \cref{cor:m24-reduction}]\label{rem:m24-three-point}
The inequality $E(W)\ge\tau$ for all contact configurations of $24$
directions would give the volume statement of \cref{thm:m24} at the three-point
level, without the four-point kernel of \cref{thm:twentyfour-points},
and it has the form treated in \cref{sec:certificate}.
Splitting each pair over the $22$ triples that contain it,
$E(W)=\sum_{i<j<k}h(u_{ij},u_{ik},u_{jk})$ with
$h(u,v,t)=\tfrac1{22}(\omega_2(u)+\omega_2(v)+\omega_2(t))-\omega_3(u,v,t)$.
In the notation of \cref{sec:certificate}, a function
$f=\sum_kf_kG_k$ with $f_k\ge0$ for $k\ge1$ and positive semidefinite
matrices $F_k$ satisfying
\begin{align*}
  h(u,v,t)\;\ge\;&\tfrac1{22}\bigl(f(u)+f(v)+f(t)\bigr)+F(u,v,t)\\
  &+\tfrac1{22}\bigl(F(1,u,u)+F(1,v,v)+F(1,t,t)\bigr)
\end{align*}
on every admissible triple give $E(W)\ge\tfrac12(576f_0-24f(1))-4F(1,1,1)$
for every $W$ of size $24$. Because the root configuration attains $\tau$,
a certificate reaching $\tau$ must be sharp, and exact rounding needs it
to be attained at a finite degree.

The dual of that problem is a probability measure $p$ on admissible
triples subject to the two-point conditions
$1+\tfrac{23}{3}\,\mathbb{E}_p[P_k(u)+P_k(v)+P_k(t)]\ge0$ for
$1\le k\le d$ and, writing
$\Sigma_k=S_k(1,u,u)+S_k(1,v,v)+S_k(1,t,t)$, the three-point conditions
\[
  24\,S_k(1,1,1)+552\,\mathbb{E}_p[\Sigma_k]
  +12144\,\mathbb{E}_p[S_k]\succeq0
\]
for $0\le k\le d$, in the matrices $S_k$ of \cite{BV08}. Every such $p$
satisfies $2024\,\mathbb{E}_p[h]\ge V_d$, where $V_d$ is the value of the
best degree-$d$ certificate, so a single feasible $p$ with
$2024\,\mathbb{E}_p[h]<\tau$ rules out every degree-$d$ certificate at
once. Restricting $p$ to a finite set of triples only shrinks the dual
feasible set, so such a $p$ is a proof and not an estimate.

At degree $6$ there is one. On a rational grid of $670$ admissible
triples a measure supported on $632$ of them, with denominator
$10^{10}$, has value $5.9826875$, least two-point slack
$1.09\times10^{-2}$, and all seven moment matrices positive definite,
every pivot of an exact rational $LDL^\top$ being at least
$1.01\times10^{-2}$. The rational arithmetic is exact; only $\omega_3$
is numerical, and it enters through the value alone, where the
quadrature error is below $10^{-12}$ against a margin of $0.18$. Since
$5.9826875<\tau=6.1650791$, no degree-$6$ certificate of this form
exists, with no appeal to the conditioning of any solver.

At degrees $8$ through $20$, on a grid of $1264$ triples, the dual
minimum is $6.0395$, $6.0700$, $6.0879$, $6.1060$, $6.1256$, $6.1416$
and $6.1516$, short of $\tau$ by $0.126$, $0.095$, $0.077$, $0.059$,
$0.039$, $0.023$ and $0.014$. These are double-precision values, but the
direction of the discretisation is favourable: restricting the measure
to a grid only shrinks the dual feasible set and so only raises the
minimum, and refining the grid by a factor of $2.6$ at degrees $6$ and
$8$ moves the value by less than $0.005$. At degree $20$ the shortfall
is still $1.4\times10^{-2}$, where a certificate that is to be rounded
exactly would need it to be zero.

Four further computations say that the shortfall is structural rather
than numerical. First, the same two-point and three-point conditions,
read as conditions on a code of $N$ points instead of $24$, remain
feasible up to $N=26.61$ at degree $4$, $N=24.77$ at degree $6$ and
$N=24.24$ at degree $8$: the relaxation does not know that $24$ is the
largest number of contacts available. Since $E(W)\ge\tau$ holds with
equality at the root configuration and, by \cref{thm:m24}, only there,
a certificate for it is a certificate of uniqueness, and a relaxation
that does not see maximality does not see uniqueness. Second, the
minimising $p$ is not a perturbation of the root configuration: it
carries about a third of its pair mass at $u=\tfrac12$ and spreads the
rest continuously over $(-1,0)$, where the root configuration has atoms
at $-1$, $-\tfrac12$ and $0$ alone. Requiring the two-point slacks at
$k=2$ and $k=4$, which vanish exactly at the root configuration, to
exceed $\varepsilon$ raises the degree-$8$ minimum only from $6.0395$ to
$6.0489$ at $\varepsilon=0.1$, and is infeasible at $\varepsilon=0.3$;
so the shortfall is not sharpness at the root configuration, and
treating a neighbourhood of it separately does not remove it. Third,
the largest $s$ for which the moment matrices admit $\succeq sI$ falls
from $0.75$ at degree $4$ to $1.7\times10^{-2}$ at degree $6$ and to
about $10^{-4}$ at degrees $8$ and above: from degree $8$ on the dual
has no usable interior. That is a property of the feasible set, not of
the arithmetic, and it is what makes interior-point solvers degrade
there. Fourth, the root configuration never becomes a local minimum of
the relaxation. Minimising the directional derivative of the objective
over the directions the degree-$d$ conditions leave open at the root
configuration gives $-1.5705\times10^{-4}$ at degree $6$ and
$-1.0951\times10^{-4}$, $-8.7210\times10^{-5}$,
$-7.8417\times10^{-5}$ and $-7.6223\times10^{-5}$ at degrees $8$,
$10$, $12$ and $14$, each computed on the support of the direction with
the closed form of $\omega_3$ rather than a table, and each agreeing to
four figures with the value the table gives. Fitting $L+Cq^d$ to the
five returns $L=7.25\times10^{-5}$ with $q=0.654$ and residuals below
$9\times10^{-7}$, and dropping any one of the five moves $L$ only
within $[7.03,7.43]\times10^{-5}$. The improvement available at the
root configuration therefore does not tend to zero with the degree,
where a certificate reaching $\tau$ would need it to be zero. The sweep
is \texttt{m24\_descent\_sweep.py}.

That the three-point level is the wrong level is not a conclusion of
ours alone. Cohn, de Laat and Leijenhorst \cite{CdLL24} settle a family
of spherical codes by three-point bounds with exact rounding, and record
that the one member of it they cannot settle is the $24$ points of
$\Rt$ in $\mathbb{R}^4$, for which, in their words, three-point bounds
do not seem to be sharp. Two of the three are authors of \cite{LLM24},
and the four-point computation there is what they went to instead.

Extended precision therefore does not decide this question in the
affirmative, and \cref{thm:m24} is proved from the four-point kernel
of \cref{thm:twentyfour-points}. What \cref{cor:m24-reduction} adds is an
exact reduction of its volume consequence to a single inequality
between pair and triple angles; what the computations above add is the
location of the obstruction. Closing the $24$-contact case at the
three-point level calls not for a longer computation at higher
precision, but for a method that certifies the uniqueness of the root
configuration, which is what the second level of the hierarchy
supplies and the first does not.
\end{remark}

\begin{remark}[How far the elementary constraints reach]\label{rem:m24-remains}
They reach a definite distance and stop, and we record where. The strongest of them is a second moment. For $24$ unit vectors
in $\mathbb{R}^4$ the Gram matrix $G=(u_{ij})$ is positive semidefinite
of rank at most $4$ with trace $24$, so Cauchy--Schwarz applied to its
nonzero eigenvalues gives $\|G\|_F^2\ge24^2/4=144$, that is
\begin{equation}\label{eq:second-moment}
  \sum_{i<j}u_{ij}^2\;\ge\;60,
\end{equation}
with equality exactly when the configuration is a tight frame. The root
system attains it, its frame operator being $6I$: the $96$ pairs at
$\tfrac12$ contribute $24$, and the $96$ at $-\tfrac12$ and the $12$ at
$-1$ contribute $36$.

Write $P$ for the set of pairs with $u_{ij}>0$, the only ones on which
$\omega_2$ and $\omega_3$ are supported. Since $u_{ij}\le\tfrac12$
there, each such pair contributes at most $\tfrac14$ to the second
moment, so a bound $\sum_Pu_{ij}^2\ge24$ would force $|P|\ge96$, exactly
the count of the root system. That bound does not follow from
\eqref{eq:second-moment} alone: it needs the pairs outside $P$ to carry
at most $36$, which is the root system's own value and is not available
in advance.

Even granting it the route stops short. The ratio $\omega_2(u)/u^2$
increases on $(0,\tfrac12]$, from $0.0487$ at $u=0.05$ to $0.3043$ at
$u=\tfrac12$, so for a fixed second moment $\sum_P\omega_2(u_{ij})$ is
smallest when the mass is spread over as many pairs as possible at as
small a value as possible. A budget $\sum_Pu_{ij}^2=24$ spread over all
$276$ pairs, at $u=0.2949$, gives $\sum_P\omega_2=5.1357$, short of
$\tau=6.1651$ by $1.0294$, with the triple term still to be subtracted;
the same budget concentrated on $96$ pairs at $\tfrac12$ gives the root
system's $7.3033$. The second moment fixes how much inner-product mass a
configuration carries but not where that mass sits, and the difference
between the two placements exceeds the entire margin to be proved.

What is missing here and in \cref{rem:m24-three-point} is the same
thing: a reason for the positive inner products to concentrate at $\tfrac12$ instead of spreading. That reason is the uniqueness of the root
system among $24$-point configurations with $u_{ij}\le\tfrac12$, and a
treatment of $m=24$ has to supply it. The equality case of the
second-level bound, \cref{thm:twentyfour-points}, supplies it, which is
why the argument goes through that kernel.
\end{remark}

The extendability of a configuration turns out to be exactly a statement
about the reach of its own cell, which is worth recording because it
converts the remaining case into a question about one polytope.

\begin{prop}[Extendability and the circumradius]\label{prop:extendable}
Let $W$ be a contact configuration with bounded cell and write
\[
  g(W)=\min_{\theta\in S^3}\ \max_{w\in W}\ \langle\theta,w\rangle
      =\cos\bigl(\text{covering radius of }W\bigr).
\]
The following are equivalent:
\begin{enumerate}[label=\textup{(\roman*)},leftmargin=2.2em]
\item some $v\in S^3$ has $\langle v,w\rangle\le\tfrac12$ for all $w\in
W$, that is $W\cup\{v\}$ is again a contact configuration;
\item $V_c(W)$ contains a point of norm $2$, that is its circumradius is
at least $2$;
\item $g(W)\le\tfrac12$, that is the covering radius of $W$ on $S^3$ is
at least $60^\circ$.
\end{enumerate}
We call $W$ \emph{saturated} when these fail, so that $g(W)>\tfrac12$.
\end{prop}

\begin{proof}
For a unit vector $v$, the condition $\langle v,w\rangle\le\tfrac12$ for
all $w$ says exactly $\langle 2v,w\rangle\le1$ for all $w$, that is
$2v\in V_c(W)$; and $|2v|=2$. So (i) and (ii) are the same statement,
since $V_c(W)$ is convex and contains the origin, so it meets the sphere
of radius $2$ if and only if its circumradius is at least $2$. For (iii),
the circumradius is $\max_\theta\rho(\theta)=\max_\theta1/\cos\delta(\theta)$
by \cref{lem:radial-form}, which is $1/g(W)$; and $1/g(W)\ge2$ is
$g(W)\le\tfrac12$.
\end{proof}

\begin{cor}[The reduction to twenty-three contacts]\label{cor:remaining}
Every contact configuration with bounded cell satisfies $\vol(V_c)\ge8$, with equality only for the root
configuration, except possibly for a configuration of exactly $23$
directions whose cell has circumradius less than $2$.
\end{cor}

\begin{proof}
By \cref{cor:m22} the bound holds when $|W|\le22$, with strict
inequality. By \cref{thm:m24} it holds with equality when $|W|=24$. Let
$|W|=23$. If the circumradius of $V_c(W)$ is at least $2$ then by
\cref{prop:extendable} there is a $v$ with $W'=W\cup\{v\}$ a contact
configuration of size $24$, and $V_c(W')\subseteq V_c(W)$, so
$\vol(V_c(W))\ge\vol(V_c(W'))=8$ by \cref{thm:m24}; the inequality is
strict because $V_c(W)$ contains a point of norm $2$ while the $24$-cell
has circumradius $\sqrt2$. The only case not covered is $|W|=23$ with
circumradius below $2$.
\end{proof}

The remaining case is therefore a single, sharply posed statement, and
it is no longer about deviations at all.

\begin{thm}[The $23$-point case]\label{thm:m23}
Every contact configuration of $23$ directions with bounded cell has
$\vol(V_c)>8$. In particular a contact configuration of $23$ directions
whose cell has circumradius less than $2$, if one exists, has cell
volume above $8$.
\end{thm}

\begin{figure}[tb]
\centering
\includegraphics[width=\textwidth]{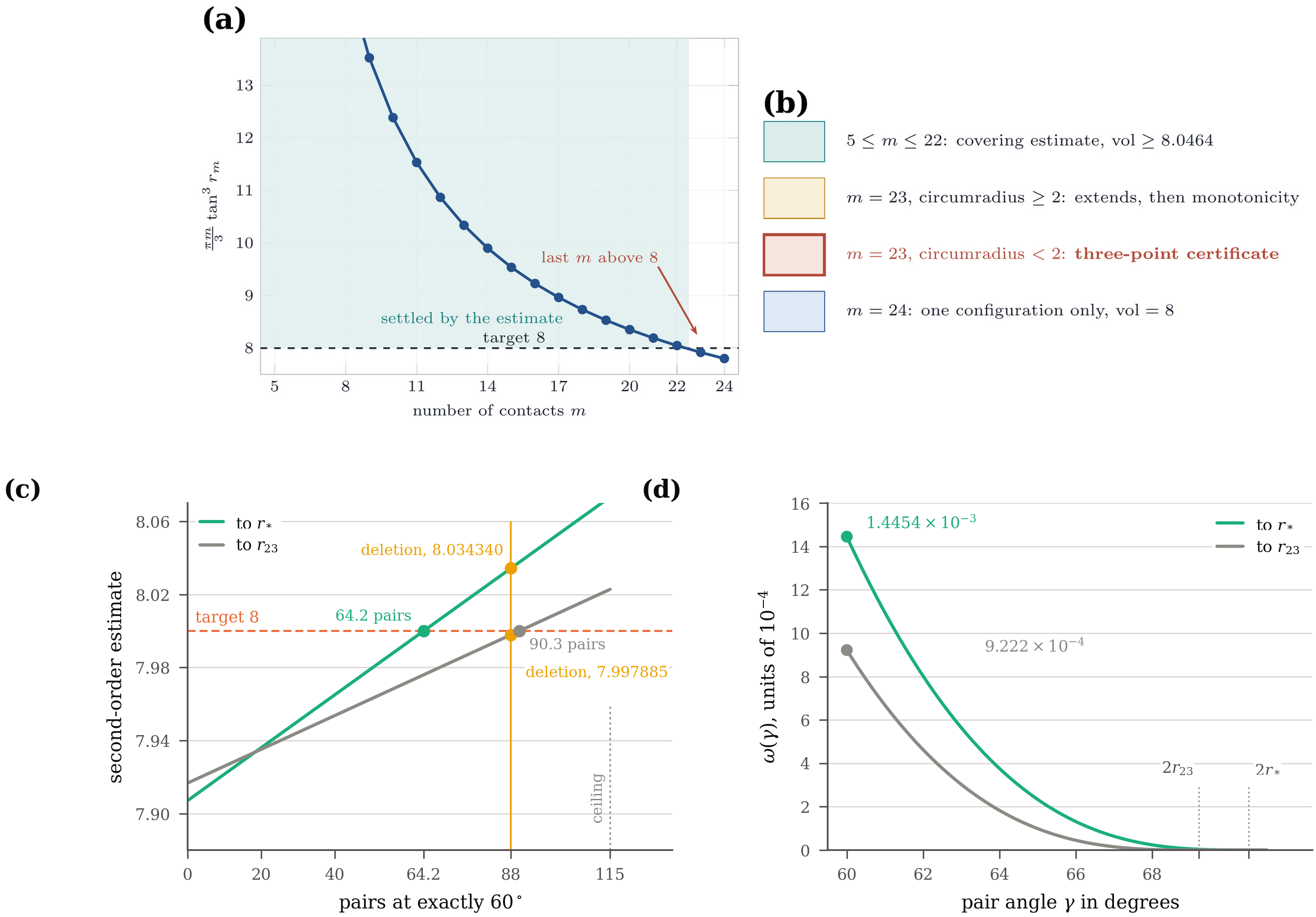}
\caption{(a) The covering bound \eqref{eq:covering-bound} at every
admissible contact count $m$, against the target $8$. It clears the
target for all $m\le22$ and falls below at $m=23$. (b) How each part of the range of contact counts is settled: the
estimate below $22$ (\cref{cor:m22}), extension followed by monotonicity
at $23$ when the circumradius reaches $2$ (\cref{prop:extendable}), and
the classification at $24$ (\cref{thm:m24}). The third row is settled by
the pair inequality of \cref{thm:strict-reduction} and its certificate,
\cref{thm:certificate}. (c) The pairwise estimate at $m=23$, for the two integration limits
$r_{23}$ and $r_*$: the right-hand side of \eqref{eq:second-order} as a
function of the number of pairs at exactly $60^\circ$, with the bracket
alone at $0$ pairs, the deletion configuration at $88$ on both lines,
the crossings of the target at $90.3$ and at $64.2$ pairs, and the
elementary ceiling of \cref{prop:pair-budget} at $115$. (d) The per-pair weight $\omega(\gamma)$ for each integration
limit: $0.000922\ldots$ and $0.001445\ldots$ at $60^\circ$, vanishing at
$2r_{23}=69.2213\ldots^\circ$ and at $2r_*=70.5288\ldots^\circ$.}
\label{fig:counting}
\end{figure}

The proof occupies the rest of this section and is of a different kind
from everything before it. \Cref{sec:m23-explicit} to
\cref{sec:slack-continuation} approach the statement through the
configuration: what a saturated set of $23$ directions would have to
look like, how much of a root system it can contain, and how the codes
of $23$ points behave as their minimal angle is pushed to $60^\circ$;
none of that is the proof, and all of it shapes the one that
\cref{sec:strict-inequality} and \cref{sec:certificate} give. There the
statement is turned into an inequality between the pair angles of any
$23$ contact directions (\cref{thm:strict-reduction}), and the
inequality is established by a semidefinite certificate that reads the
triples of directions, verified in exact and interval arithmetic
(\cref{thm:certificate}). Two remarks first on the size of what is
asked. The configuration obtained by deleting
one root from $D_4$ has cell the $24$-cell together with the pyramid of
\cref{lem:Q-structure}, of volume $\tfrac{25}3=8.333\ldots$, and its
circumradius is exactly $2$: the apex $2\base_0$ is the point that the
deleted root would occupy. So that configuration sits exactly on the
boundary of \cref{prop:extendable}(ii), and is extendable by the single
direction $\base_0$. And the covering estimate of
\cref{thm:covering-bound} gives $7.916728\ldots$ at $m=23$, short of $8$
by $0.083272\ldots$, against a true value of $\tfrac{25}3$ at the
configuration just described; so the estimate has to be improved by
about a fifth of its slack there, not by all of it, which is a different
situation from $m=24$, where the true value is $8$ and no slack exists
at all.
\Cref{fig:counting}(a) plots the range, and \cref{fig:counting}(b)
records what settles each part of it.

\subsection{The saturated case, written out}\label{sec:m23-explicit}

Everything \cref{thm:m23} asserts can be written down explicitly, and it
is worth doing so, because what it concerns is a single polytope with
twenty-three facets and no deviation parameters at all. Fix a contact
configuration $W=\{w_1,\dots,w_{23}\}$ with bounded cell and
circumradius below $2$, and write $s_{ij}=\langle w_i,w_j\rangle$. The
same statement takes three exact forms.

\emph{Radially.} By \cref{lem:radial-form} the quantity to bound is
$\tfrac14\int_{S^3}\sec^4\delta$, so the target is
\begin{equation}\label{eq:m23-target}
  \int_{S^3}\sec^4\bigl(\delta(\theta)\bigr)\,d\theta\;\ge\;32 .
\end{equation}
Saturation says $\delta<60^\circ$ everywhere, and
\cref{prop:covering-radius} says $\max\delta\ge37.7612\ldots^\circ$, so
the distance function is confined to a band of width about $22^\circ$.

\emph{On the boundary.} By \cref{lem:surface-form} the same target is
$\sum_{i}\vol_3(F_i)\ge32$, an average of $\tfrac{32}{23}=1.391304\ldots$
per facet. That number is larger than the $\tfrac43$ of a facet of the
root configuration, and larger than the
$16\sqrt2-\tfrac{64}3=1.294084\ldots$ of the antiprism facet of
\cref{rem:facet-local-obstruction}. So a counterexample would be a
configuration whose facets, each containing the ball of radius
$1/\sqrt3$ by \cref{lem:facet-inball}, average less than
$\tfrac{32}{23}$; and every facet shape that occurs at twenty-four
contacts already averages less than that. This is the exact sense in
which the remaining case is not a perturbation of the root
configuration: it must have larger facets than the root configuration
has, not smaller.

\emph{At the vertices.} If $I$ indexes four active directions at a
vertex $z_I$ of the cell, then $\langle z_I,w_i\rangle=1$ for $i\in I$,
so writing $z_I=\sum_{i\in I}a_iw_i$ gives $G_Ia=\mathbf1$ and
\begin{equation}\label{eq:vertex-norm}
  |z_I|^2=\mathbf1^{\top}G_I^{-1}\mathbf1 ,
  \qquad G_I=(s_{ij})_{i,j\in I}.
\end{equation}
Circumradius below $2$ therefore says
$\mathbf1^{\top}G_I^{-1}\mathbf1<4$ at every vertex. Two cases of
\eqref{eq:vertex-norm} are worth recording. Four pairwise orthogonal
active directions give $G_I=\mathrm{Id}$ and $|z_I|^2=4$ exactly, so a
saturated configuration has no vertex at which four active directions
are mutually orthogonal. And four directions with all pairwise inner
products equal to $g$ give $|z_I|^2=4/(1+3g)$, which is at least $4$
precisely when $g\le0$; so at every vertex of a saturated configuration
the active directions must be positively correlated on average.

\begin{prop}[What the pairwise estimate would need]\label{prop:pair-budget}
For $m=23$ write
\[
  \omega(\gamma)=\frac14\int_0^{r_*}\Lambda(r,\gamma)\,d(\sec^4r)
\]
for the contribution a single pair at angle $\gamma$ makes to the
right-hand side of \eqref{eq:second-order} at $R=r_*$. Then $\omega$ is
strictly decreasing on $[60^\circ,2r_*]$ and vanishes from
$2r_*=70.5288\ldots^\circ$ on, with
$\omega(60^\circ)=0.001445408\ldots$, and a contact configuration of
$23$ directions with bounded cell satisfies $\vol(V_c)\ge8$ as soon as
\begin{equation}\label{eq:pair-budget}
  \sum_{i<j}\omega(\gamma_{ij})\;\ge\;0.092855570\ldots ,
\end{equation}
which holds in particular once it has $65$ pairs at exactly $60^\circ$.
The deletion configuration has $88$, and no configuration of $23$
directions has more than $115$.
\end{prop}

\begin{proof}
That $\omega$ is decreasing and eventually zero is immediate from
\eqref{eq:lens}, since $\Lambda(r,\gamma)$ is decreasing in $\gamma$ and
vanishes for $r\le\gamma/2$. Inequality \eqref{eq:pair-budget} is
\eqref{eq:second-order} at $R=r_*$ with the value $7.907144430\ldots$ of
the bracket subtracted from $8$. Numerically
$7.907144430+65\times0.001445409=8.001096\ldots$, while $64$ pairs
give $7.999651\ldots$, which is the sense in which $65$ suffice and
$64$ do not; the crossing is at $64.24$ pairs.

For the last assertion, fix $w\in W$ and let $S$ be the set of $w_i$
with $\langle w,w_i\rangle=\tfrac12$. Write
$w_i=\tfrac12w+\tfrac{\sqrt3}2u_i$ with $u_i$ a unit vector of
$w^{\perp}$. Then
$\langle w_i,w_j\rangle=\tfrac14+\tfrac34\langle u_i,u_j\rangle$, so the
contact condition $\langle w_i,w_j\rangle\le\tfrac12$ reads
$\langle u_i,u_j\rangle\le\tfrac13$: the $u_i$ are unit vectors of
$\mathbb{R}^3$ at pairwise angle at least
$\arccos\tfrac13=70.5287\ldots^\circ$, and the caps of angular radius
$35.2643\ldots^\circ$ about them are disjoint on $S^2$. Each such cap
has area $2\pi(1-\cos35.2643\ldots^\circ)=1.152986\ldots$ against the
$4\pi$ of $S^2$, so $|S|\le10$. Every vertex of the graph of pairs at
$60^\circ$ therefore has degree at most $10$, and a graph on $23$
vertices of maximum degree $10$ has at most
$\lfloor 230/2\rfloor=115$ edges.
\end{proof}

\begin{remark}[The one route that is not blocked]\label{rem:pair-budget}
With the integration stopped at $r_{23}$, as in
\eqref{eq:second-order-value}, the estimate misses $8$ at the deletion
configuration by three pairs, and that is the number often quoted as
closing off the pairwise route. Carried to $r_*$ it does not miss: the
deletion has $88$ pairs against the $64.24$ that \eqref{eq:pair-budget}
asks for, a surplus of $0.034$, and the weight falls away quickly once a
pair opens up, with $\omega(62^\circ)/\omega(60^\circ)=0.55504\ldots$
and $\omega(66^\circ)/\omega(60^\circ)=0.09077\ldots$. So what would
settle the local bound at twenty-three contacts, for every
configuration and not only for the saturated ones, is a lower bound on
the weighted number of pairs at or near $60^\circ$ that every contact
configuration of $23$ directions satisfies. \Cref{prop:two-point-barrier}
says how much of it the pair angles supply on their own, and
\cref{sec:strict-inequality} sets the statement out in full.
\Cref{fig:counting}(c) and (d) set the numbers down, for both
integration limits. The script is \texttt{pair\_budget.py} for the
limit $r_{23}$ and \texttt{three\_point\_reduction.py} for $r_*$.
\end{remark}

Two further facts about the pair angles are worth setting down, because
together they say how far a bound of this kind can be pushed and where it
stops. The first is a lower bound on the total overlap of the contact caps
that a saturated configuration is forced to have.

\begin{prop}[Covering multiplicity]\label{prop:cov-mult}
Let $W$ be a saturated contact configuration with $|W|=m$. Then
\begin{equation}\label{eq:cov-mult}
  \sum_{i<j}\Lambda\bigl(\tfrac\pi3,\gamma_{ij}\bigr)
  \;\ge\;\frac12\left(\frac{\bigl(m\,C(\pi/3)\bigr)^2}{2\pi^2}
    -m\,C(\pi/3)\right),
\end{equation}
where $C(r)=\pi(2r-\sin2r)$ is the measure of a cap of radius $r$. At
$m=23$ the right-hand side is $155.172054\ldots$, against
$162.566121\ldots$ at the deletion configuration, so the inequality is
within five per cent of binding there.
\end{prop}

\begin{proof}
Write $N(\theta)$ for the number of $w\in W$ with
$\langle\theta,w\rangle>\tfrac12$. Then $\int_{S^3}N=m\,C(\pi/3)$, and
$\int_{S^3}N(N-1)$ counts ordered pairs of caps containing $\theta$, so
$\int_{S^3}N^2=m\,C(\pi/3)+2\sum_{i<j}\Lambda(\tfrac\pi3,\gamma_{ij})$.
Cauchy--Schwarz gives $\int N^2\ge(\int N)^2/2\pi^2$, and rearranging is
\eqref{eq:cov-mult}.
\end{proof}

The second fact is negative, and it is the sharpest thing we can say
about the route. Everything the estimate of \cref{prop:second-order} sees
of a configuration is the multiset of its pair angles. One can ask how
far that multiset alone can carry it, by relaxing the multiset to an
arbitrary measure subject to every inequality of the same character that
a genuine configuration satisfies.

\begin{prop}[What the pair angles alone can give]\label{prop:two-point-barrier}
Let $\mu$ range over the nonnegative measures on $[60^\circ,180^\circ]$
of total mass $\binom{23}2=253$ subject to
\[
  \int\Lambda(t,\gamma)\,d\mu\;\ge\;23\,C(t)-2\pi^2
  \quad\text{for all } t,
  \qquad
  \int G_k(\cos\gamma)\,d\mu\;\ge\;-\tfrac{23}2
  \quad\text{for }k\ge1,
\]
where $G_k$ are the Gegenbauer polynomials of $S^3$ normalised by
$G_k(1)=1$. The pair-angle multiset of every contact configuration of
$23$ directions gives such a $\mu$. The minimum of $\int\omega\,d\mu$
over this family, for the weight $\omega$ of \cref{prop:pair-budget},
is $0.07380\ldots$, four fifths of the $0.092855\ldots$ that
\eqref{eq:pair-budget} requires; for the weight integrated only to
$r_{23}$ it is $0.04157\ldots$ against $0.083272\ldots$, one half.
\end{prop}

\begin{proof}
Feasibility is Bonferroni's inequality
$\sigma(\bigcup_iC_i(t))\ge\sum_i\sigma(C_i(t))-\sum_{i<j}\Lambda(t,\gamma_{ij})$
together with $\sigma(\bigcup_iC_i(t))\le2\pi^2$ for the first family,
and the positive definiteness of $G_k$ on $S^3$, which gives
$\sum_{i,j}G_k(\langle w_i,w_j\rangle)\ge0$ and hence
$\sum_{i<j}G_k(\cos\gamma_{ij})\ge-\tfrac{23}2$, for the second. The
values are the optima of finite linear programmes: in
\texttt{covering\_multiplicity.py} on a grid of $601$ angles with $24$
radii between $r_{23}$ and $45^\circ$ and $12$ degrees, for the weight
integrated to $r_{23}$, and in \texttt{three\_point\_reduction.py} on a
grid of $1201$ angles with $60$ radii between $r_{23}$ and $90^\circ$
and $24$ degrees, for both weights. Adding radii or degrees does not
raise either optimum, and the same measure attains both: $103.7$ pairs
at $62.3^\circ$, $90.0$ at $103.0^\circ$, $39.2$ at $122.6^\circ$ and
$20.1$ at $154.1^\circ$.
\end{proof}

\begin{remark}[Where a proof would have to look]\label{rem:two-point-barrier}
\Cref{prop:two-point-barrier} is a stronger statement than the one
\eqref{eq:second-order-value} makes, and a more useful one. It is not
that the deletion configuration falls short, which it does not once the
integration is carried to $r_*$: it is that nothing which reads only
the pair angles gets past four fifths of the way, whatever constraints
of that type are imposed and however they are combined. The minimising
measure is instructive. It gives every direction nine neighbours at
$62.3^\circ$ and nearly eight at $103^\circ$, which is a perfectly good
local picture at any one direction, and the two constraints available
to a pair-angle relaxation, Bonferroni and positive definiteness, cannot
see that the twenty-three such pictures do not fit together on one
sphere. Saturation itself is invisible at that level too, since the
covering radius is not a function of the pair angles, which is why
\cref{prop:cov-mult} does not help either. A proof of the bound at
twenty-three contacts has to read more of the configuration than its
pair angles record; the classification of \cref{thm:m23-from-codes} is
one such thing, and the three-point form of
\cref{sec:strict-inequality} is another.
\end{remark}

One family of candidates can be ruled out completely, and it is the family
one would look at first. A configuration of a prime number of directions
with any symmetry at all transitive on them is a single orbit of a rotation
of that order, and in dimension four such an orbit carries a single
parameter.

\begin{prop}[No symmetric candidate]\label{prop:no-transitive}
No contact configuration of $23$ directions in $\mathbb{R}^4$ admits a
group of symmetries transitive on its directions.
\end{prop}

\begin{proof}
Let $W$ be such a configuration and $G\le O(4)$ transitive on it. Then
$23$ divides $|G|$, so $G$ has an element $g$ of order $23$ by Cauchy's
theorem. Since $\det g$ is a $23$rd root of unity in $\{\pm1\}$ it is $1$,
so $g\in SO(4)$ and $g$ is a rotation through angles $t_1=2\pi a/23$ and
$t_2=2\pi b/23$ in two orthogonal planes $P_1,P_2$, with $(a,b)\ne(0,0)$.
The cell is bounded, so $W$ spans $\mathbb{R}^4$ and $g$ cannot fix every
member of $W$; the orbit of any $w_0\in W$ under $\langle g\rangle$
therefore has $23$ elements and is all of $W$.

Write $w_0=p+q$ with $p\in P_1$, $q\in P_2$ and $c=|p|^2$, so $|q|^2=1-c$
and $c\in[0,1]$. Then $w_k=g^kw_0$ and
\[
  \langle w_0,w_k\rangle=\langle p,g^kp\rangle+\langle q,g^kq\rangle
  =c\cos(kt_1)+(1-c)\cos(kt_2),
\]
which does not depend on the directions of $p$ and $q$ inside their planes
and is affine in $c$. Transitivity makes every pairwise inner product one
of these, so the contact condition is the system of $22$ affine
inequalities
\[
  c\cos(kt_1)+(1-c)\cos(kt_2)\;\le\;\tfrac12,\qquad k=1,\dots,22,
\]
in the single variable $c\in[0,1]$. Each of the $528$ pairs $(a,b)$ is
therefore decided by a one-variable linear programme. None of the $528$ is
feasible. The quantity
\[
  \delta(a,b)=\min_{c\in[0,1]}\ \max_{1\le k\le22}
    \Bigl(c\cos(kt_1)+(1-c)\cos(kt_2)-\tfrac12\Bigr)
\]
is the least violation available to the type $(a,b)$, and it is positive
exactly when that type is infeasible. Its minimum over the $528$ types is
$0.0740002839\ldots$, attained at $(a,b)=(8,17)$. The maximum of affine
functions is convex, so each $\delta(a,b)$ is found exactly by evaluating
at the endpoints and at the crossings of pairs of the $22$ functions; the
computation is carried out at $60$ digits in
\texttt{saturation\_search.py}, and a margin of $0.074$ in an inner
product is not a rounding effect.
\end{proof}

\begin{remark}[What that removes]\label{rem:no-transitive}
\Cref{prop:no-transitive} is not the proof of \cref{thm:m23}, and it is
worth being exact about what it does. A saturated configuration, if one exists,
has no symmetry acting transitively on its directions: it is not a single
orbit of anything. That removes every construction of the kind one reaches
for first, and it is consistent with the rest of what is known about the
case, since the configurations that do exist at $23$ directions are the
deletions, whose symmetry group is the stabiliser of a root in the Weyl
group and has four orbits on the remaining directions, of sizes
$1,8,6,8$. The search of \texttt{saturation\_search.py} looks for the
rest: it minimises $\max_I|z_I|$ directly, from perturbed deletions, from
random starts and from partially completed root systems, and over $202$
starts it reached $12$ contact configurations, all with the inner-product
multiset of a deletion and all with $g$ exactly $\tfrac12$. That is
exploration, and we report it as such.
\end{remark}

The one configuration that comes close to the description is the
deletion of a root, and it is not merely on the wrong side of the
circumradius condition: it is rigid, so a saturated configuration cannot
be reached by bending it.

\begin{prop}[The deletion configurations are rigid]\label{prop:deletion-rigid}
Let $W=\Rt\setminus\{\base_0\}$ be the contact configuration obtained by
deleting one root. Every first-order motion of $W$ through contact
configurations is an infinitesimal rotation.
\end{prop}

\begin{proof}
Scale the roots to $a_i$ with $|a_i|^2=2$, so that the contact condition
reads $\langle a_i,a_j\rangle\le1$ and every quantity below is an
integer. Write $i\sim j$ when $\langle a_i,a_j\rangle=1$; deleting the
root $a_0$ leaves $23$ directions and a set $T$ of $88$ unordered pairs
$\{i,j\}$ with $i\sim j$. A first-order motion is a tuple
$d=(d_1,\dots,d_{23})$ with
\begin{equation}\label{eq:strut}
  \langle d_i,a_i\rangle=0\quad\text{for every }i,
  \qquad
  \langle d_i,a_j\rangle+\langle a_i,d_j\rangle\le0
  \quad\text{for every }\{i,j\}\in T,
\end{equation}
the first condition keeping each direction on its sphere and the second
keeping each tight pair from closing. The infinitesimal rotations
$d_i=Aa_i$ with $A$ antisymmetric satisfy both with equality and span a
space of dimension $6$.

Let $c_i=\langle a_i,a_0\rangle$ record the position of $a_i$ relative
to the deleted root; the values $-2,-1,0,1$ occur $1,8,6,8$ times. Put
\[
  y_{ij}=\begin{cases}
    1,&(c_i,c_j)\in\bigl\{(-1,-1),\,(-1,1)\bigr\},\\
    2,&(c_i,c_j)\in\bigl\{(-2,-1),\,(-1,0),\,(0,1)\bigr\},\\
    3,&(c_i,c_j)=(1,1),
  \end{cases}
  \qquad
  \mu_i=\begin{cases}-6,&c_i=-1,\\-8,&\text{otherwise},\end{cases}
\]
the first for $\{i,j\}\in T$ with the pair written so that $c_i\le c_j$,
the second for each $i$. Those six are the only types that occur among
the pairs of $T$, and a direct integer computation gives the equilibrium
relation
\begin{equation}\label{eq:stress}
  \sum_{j\,:\,j\sim i} y_{ij}\,a_j+\mu_i\,a_i=0
  \qquad\text{for every }i .
\end{equation}
Pair \eqref{eq:stress} with $d_i$ and sum over $i$. The diagonal terms
contribute $\mu_i\langle d_i,a_i\rangle=0$, and each unordered pair of
$T$ is met twice, contributing
$y_{ij}(\langle d_i,a_j\rangle+\langle a_i,d_j\rangle)$, so
\[
  \sum_{\{i,j\}\in T} y_{ij}
  \bigl(\langle d_i,a_j\rangle+\langle a_i,d_j\rangle\bigr)=0 .
\]
Every weight is positive and every bracket is nonpositive by
\eqref{eq:strut}, so every bracket vanishes: a first-order motion holds
all $88$ pairs at equality. The linear system consisting of those $88$
equalities together with the $23$ conditions
$\langle d_i,a_i\rangle=0$ has, over $\mathbb{Q}$, a solution space of
dimension exactly $6$, and the rotations already occupy $6$ dimensions
of it, so the two coincide.
\end{proof}

Since the deletion configuration is a tensegrity framework, $23$ bars
from the centre and $88$ struts, infinitesimal rigidity is rigidity in
the sense of Roth and Whiteley \cite[Thm.~5.7]{RW81}: no continuous
deformation through contact configurations moves it except by rotation.
A saturated $23$-point configuration therefore cannot be a small
perturbation of a deletion, which is the same conclusion the facet
average reached above, now from the contact graph, not from the volume. Both are checked by \texttt{rigidity23.py} in the supplement,
the stress and the dimension count in exact integer arithmetic.

The same argument, with nothing deleted, places the root system itself.

\begin{cor}[The root system is rigid]\label{cor:root-rigid}
Every first-order motion of $\Rt$ through contact configurations of
$24$ directions is an infinitesimal rotation, and $\Rt$ is rigid: a
contact configuration of $24$ directions sufficiently close to $\Rt$
is a rotation of it.
\end{cor}

The first half of this is not new. Cohn, Jiao, Kumar and Torquato
\cite[Prop.~3.3]{CJKT11} prove that the $D_4$ root system is
infinitesimally jammed, through the $A_2$ subsystems and not through a stress. We give the proof again because the stress, and not
the jamming alone, is what \cref{thm:local-uniqueness} needs, and
because the same stress handles the deletion of
\cref{prop:deletion-rigid}. Neither their argument nor ours carries a
radius as it stands; supplying one is what \cref{thm:local-uniqueness}
does.

\begin{proof}
The set $T$ is now the set of the $96$ pairs at $\langle a_i,a_j\rangle=1$,
the edges of the $24$-cell, and every direction lies in $8$ of them
(\cref{fig:pftwo}(a), drawn on page~\pageref{fig:pftwo}). The eight roots at $60^\circ$ from a root
$a_i$ sum to $4a_i$, so the constant stress, weight $1$ on every pair of
$T$ and $\mu_i=-4$, satisfies \eqref{eq:stress}; pairing it with a
first-order motion as before forces all $96$ brackets to vanish, and the
system of those $96$ equalities with the $24$ conditions
$\langle d_i,a_i\rangle=0$ has, over $\mathbb{Q}$, a solution space of
dimension $6$, which the rotations fill. Rigidity follows from
infinitesimal rigidity by \cite[Thm.~5.7]{RW81}, the pairs at inner
products $0$, $-1$ and $-2$ staying strictly slack under a small motion.
The count and the dimension are checked by \texttt{rigidity24.py}, and
the equilibrium relation, the sum of the eight neighbours of every root
being four times the root, by Lean's kernel in \texttt{D4Stress.lean}.
\end{proof}

\Cref{cor:root-rigid} is the local half of \cref{thm:m24}: the root
system is isolated among contact configurations of $24$ directions, and
that much is elementary. What the semidefinite kernel of
\cref{thm:twentyfour-points} supplies, and what no argument in a
neighbourhood of $\Rt$ can, is the global half, that there is no such
configuration at any distance from $\Rt$; it is the one ingredient of
the general theorems that rests on a certificate rather than on a
closed argument, and we have not found a way to derive it from anything
local.

The rigidity of \cref{cor:root-rigid} comes with no radius. Roth and
Whiteley produce a neighbourhood but not a measure of it, and a
statement of isolation with no number in it cannot be combined with
anything quantitative. It can be made explicit, and cheaply: the stress
argument that proves \cref{cor:root-rigid} proves rather more if the
quadratic terms are carried along instead of discarded, and what the
carrying needs is one constant, the smallest nonzero singular value of
the linear map that the argument linearises.

Write $\mathcal{T}=\{\tau\in(\mathbb{R}^4)^{24}:\langle\tau_i,a_i\rangle=0\}$ for the
tangent space, a subspace of dimension $72$, and
$\Lambda\colon\mathcal{T}\to\mathbb{R}^{96}$ for the \emph{rigidity operator}
\begin{equation}\label{eq:rigidity-operator}
  (\Lambda\tau)_{ij}
  =\tfrac1{\sqrt2}\bigl(\langle a_i,\tau_j\rangle+\langle a_j,\tau_i\rangle\bigr),
  \qquad \{i,j\}\in T,
\end{equation}
where $T$ is, as in the proof of \cref{cor:root-rigid}, the set of
the $96$ pairs at $\langle a_i,a_j\rangle=1$ and the factor normalises the
roots. The first-order motions of
\cref{cor:root-rigid} are the $\tau$ with $\Lambda\tau\le0$, and that
corollary says $\ker\Lambda$ is the $6$-dimensional space of
infinitesimal rotations.

\begin{lem}[Spectrum of the rigidity operator]\label{lem:rigidity-spectrum}
The singular values of $\Lambda$ are
$0,\ 1,\ \sqrt{5/2},\ \sqrt3,\ 2$, with multiplicities
$6$, $29$, $8$, $21$, $8$. In particular the smallest nonzero singular
value is exactly $1$, so
\begin{equation}\label{eq:sigma-one}
  \|\Lambda\tau\|_2\ \ge\ \|\tau\|_2
  \qquad\text{for every } \tau\in\mathcal{T} \text{ orthogonal to } \ker\Lambda .
\end{equation}
\end{lem}

\begin{proof}
What has to be established is a finite exact statement about one
integral matrix, and it is established by computation, carried out twice
on independent code paths: in integer arithmetic in Python and again
inside Lean's kernel. What follows sets out what is computed and why it
determines the spectrum.

Let $\Lambda'=\sqrt2\,\Lambda$, whose coefficients $\langle a_i,\cdot\rangle$
are integral, and let $P$ be the orthogonal projection onto $\mathcal{T}$, with
blocks $I-\tfrac12a_ia_i^{\mathsf T}$, which is rational. Then
$N=P\Lambda'^{\mathsf T}\Lambda'P$ has $4N$ integral, and the nonzero
eigenvalues of $N$ are twice the squares of the nonzero singular values
of $\Lambda$. A direct computation in integer arithmetic gives the
matrix identity
\[
  4N\,(4N-8I)\,(4N-20I)\,(4N-24I)\,(4N-32I)=0 ,
\]
so every eigenvalue of $N$ lies in $\{0,2,5,6,8\}$, and the rank of
$4N-cI$ over a field of positive characteristic, which bounds the rank
over $\mathbb{Q}$ from below, is $66$, $67$, $88$, $75$, $88$ for
$c=0,8,20,24,32$; the five multiplicities that result sum to $96$, so
each bound is met and the multiplicities are $30$, $29$, $8$, $21$, $8$.
The kernel of $N$ contains the $24$ radial directions, on which $P$
vanishes, together with $\ker\Lambda$, and $30=24+6$ accounts for it
exactly. Dividing by two and taking square roots gives the singular
values. The identity and the ranks are computed by
\texttt{rigidity\_spectrum.py} by fraction-free elimination over
$\mathbb{Z}$, with no floating point at any step; the $96$ tight pairs it
starts from, and the neighbour sum of \cref{cor:root-rigid}, are checked
again by Lean's kernel in \texttt{D4Stress.lean}.
\end{proof}

\begin{thm}[A radius for the rigidity]\label{thm:local-uniqueness}
Let $W=\{w_1,\dots,w_{24}\}\subset S^3$ be a contact configuration and
put
\[
  d(W)=\min_{R\in O(4)}\ \min_{\pi}\
  \Bigl(\sum_{i=1}^{24}\bigl|w_{\pi(i)}-Ra_i/\sqrt2\bigr|^2\Bigr)^{1/2},
\]
the minimum over orthogonal maps and over bijections $\pi$. Then either
$d(W)=0$, and $W$ is an orthogonal image of the normalised root system,
or
\[
  d(W)\ \ge\ \frac{2}{\sqrt{577}}\ =\ 0.0832611\ldots
\]
\end{thm}

\begin{proof}
Both minima are over compact sets and are attained; relabel and rotate
so that they are attained at $\pi=\mathrm{id}$ and $R=I$, write
$u_i=a_i/\sqrt2$, and put $\varepsilon_i=w_i-u_i$, $e_i=|\varepsilon_i|$
and $\|\varepsilon\|^2=\sum_ie_i^2=d(W)^2$.

Since $|w_i|=|u_i|=1$, expanding $|u_i+\varepsilon_i|^2=1$ gives
\begin{equation}\label{eq:radial-part}
  \rho_i:=\langle u_i,\varepsilon_i\rangle=-\tfrac12e_i^2 ,
\end{equation}
so the radial part of each displacement is quadratically small. Put
$\tau_i=\varepsilon_i-\rho_iu_i$, so that $\tau\in\mathcal{T}$ and
$|\tau_i|^2=e_i^2-\rho_i^2$.

Minimality in $R$ fixes the gauge: differentiating
$R\mapsto\sum_i|w_i-Ru_i|^2$ along $R=e^{sA}$ with $A$ antisymmetric and
setting $s=0$ gives $\sum_i\langle\varepsilon_i,Au_i\rangle=0$, and since
$\langle u_i,Au_i\rangle=0$ the same holds for $\tau$. So $\tau$ is
orthogonal to every infinitesimal rotation, hence to $\ker\Lambda$ by
\cref{cor:root-rigid}, and \eqref{eq:sigma-one} applies.

Now the contact condition, kept to second order. For a pair
$\{i,j\}\in T$, expanding $\langle w_i,w_j\rangle\le\tfrac12$ and
using $\langle u_i,\varepsilon_j\rangle=\langle u_i,\tau_j\rangle+\tfrac12\rho_j$,
\[
  g_{ij}:=(\Lambda\tau)_{ij}
  \ \le\ -\tfrac12(\rho_i+\rho_j)-\langle\varepsilon_i,\varepsilon_j\rangle
  \ \le\ \tfrac14\bigl(e_i^2+e_j^2\bigr)+e_ie_j
\]
by \eqref{eq:radial-part} and Cauchy--Schwarz. Every index lies in exactly
eight pairs of $T$, so
$\sum_{T}(e_i^2+e_j^2)=8\|\varepsilon\|^2$ and
$\sum_{T}e_ie_j\le\tfrac12\sum_{T}(e_i^2+e_j^2)=4\|\varepsilon\|^2$,
and therefore $\sum_{T}g_{ij}^{+}\le6\|\varepsilon\|^2$ for the
positive parts. The constant stress of \cref{cor:root-rigid} gives
$\sum_{T}g_{ij}=0$ exactly, so the negative parts sum to as much
as the positive ones and
\[
  \|\Lambda\tau\|_2\ \le\ \|\Lambda\tau\|_1
  \ =\ 2\sum_{T}g_{ij}^{+}\ \le\ 12\|\varepsilon\|^2 .
\]
With \eqref{eq:sigma-one} this reads $\|\tau\|_2\le12\|\varepsilon\|^2$.
Finally $\sum_i\rho_i^2=\tfrac14\sum_ie_i^4\le\tfrac14\|\varepsilon\|^4$
by \eqref{eq:radial-part}, so
\[
  \|\varepsilon\|^2=\|\tau\|_2^2+\sum_i\rho_i^2
  \ \le\ 144\|\varepsilon\|^4+\tfrac14\|\varepsilon\|^4
  \ =\ \tfrac{577}4\|\varepsilon\|^4 ,
\]
and $\|\varepsilon\|$ is either $0$ or at least $2/\sqrt{577}$.
\end{proof}

\begin{cor}[The twenty-four-contact case near the root
system]\label{cor:m24-local}
A contact configuration $W$ with $|W|=24$ and $d(W)<2/\sqrt{577}$ is an
orthogonal image of the normalised roots, so $V_c(W)$ is a regular
$24$-cell and $\vol(V_c(W))=8$. Nothing in this uses a certificate.
\end{cor}

\begin{proof}
\Cref{thm:local-uniqueness} gives $d(W)=0$, and the volume is then
\cref{lem:refcell-vol} together with the invariance of volume under
orthogonal maps.
\end{proof}

\begin{figure}[tb]
\centering
\includegraphics[width=0.82\textwidth]{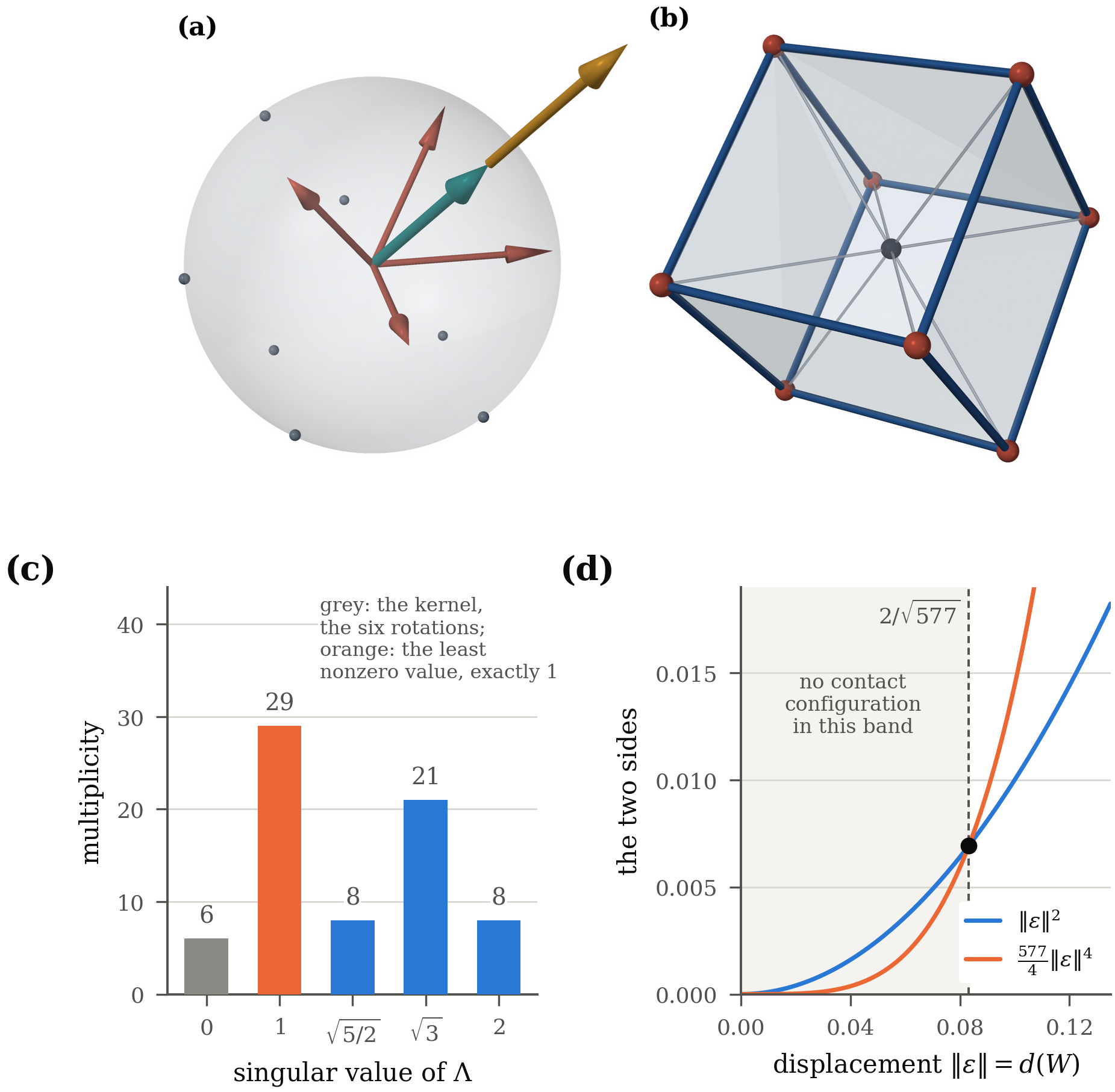}
\caption{(a) The neighbour sum in dimension three, where it fits in one picture: a contact direction of the face-centred cubic lattice (teal), the four contact directions at $60^\circ$ from it (clay), and their sum (ochre), which is twice the direction itself. (b) The same identity in dimension four, drawn where it lives. The eight roots at inner product $1$ from a root $\alpha$ fill the hyperplane $\langle\alpha,x\rangle=1$, a copy of $\mathbb{R}^3$, as the vertices of a cube centred at $\alpha/2$ (dark); eight times that centre is $4\alpha$. (c) The singular values of the rigidity operator $\Lambda$ on the $72$-dimensional tangent space, with multiplicities (\cref{lem:rigidity-spectrum}). The kernel is the six infinitesimal rotations, and the smallest nonzero value is the constant that turns the first-order argument into a radius. (d) The two sides of the inequality that closes \cref{thm:local-uniqueness}. Feasibility forces $\|\varepsilon\|^2\le\tfrac{577}{4}\|\varepsilon\|^4$, so the displacement is either zero or at least $2/\sqrt{577}$: the shaded band holds no contact configuration of $24$ directions.}
\label{fig:rig}
\end{figure}

We state plainly where the line between the certificate and the
closed arguments of this paper now falls. The global classification at
$24$ contacts is \cref{thm:twentyfour-points}, proved from the verified
kernel. What this section supplies is the quantitative statement that
the classification by itself does not give: an explicit isolation
radius about the normalised root system, $2/\sqrt{577}$ in the
root-sum-of-squares displacement and about $0.017$ per direction,
within which \cref{thm:m24} is proved by hand and the kernel is not
used at all. Beyond that radius the kernel is doing the work, and it is
doing it on the whole of the rest of a
$66$-dimensional quotient, $(S^3)^{24}$ modulo $O(4)$, which is by far
the larger part. The constant $577/4$ is not optimal and a sharper
reading of the same inequalities would enlarge the ball, but no argument
of this kind reaches a configuration that is nowhere near $\Rt$, which
is the whole of the remaining difficulty.
\Cref{fig:rig}(a) and (b) draw the neighbour sum that makes the stress
work, \cref{fig:rig}(c) the spectrum that supplies the constant, and
\cref{fig:rig}(d) the two sides of the final inequality.

Beyond that, what we have is negative and numerical. Searches started
from each of the twenty-four deletions, from greedy random packings and
from annealed random starts produced no contact configuration of $23$
directions with $g(W)>\tfrac12$; the largest value found was $\tfrac12$
itself, attained exactly at the deletions, and the random starts reached
no $23$-point contact configuration at all. That is exploration and not
evidence of nonexistence, and we report it as such; the script is
\texttt{extendability.py}. What a proof has to supply is clear
from \cref{prop:area-optimal}: the covering estimate already uses the
areas of the spherical Voronoi cells to the full, so the missing
$0.083272\ldots$ has to come from their shape, and that is what the
truncated volume of \cref{sec:strict-inequality} reads off, pair by
pair, from the lenses in which neighbouring caps overlap.

\subsection{How much of a root system a configuration can keep}\label{sec:root-meet}

\Cref{prop:deletion-rigid} says that a saturated configuration of $23$
directions is not a bent deletion. The following makes that quantitative
and turns it into a statement with no smallness in it: such a
configuration can share at most $20$ of its directions with any root
system, so at least three of them are new. The argument is elementary and
uses only the coordinates.

Take $\Rt=\{(e_i\pm e_j)/\sqrt2:1\le i<j\le4\}$ and let $S\subseteq\Rt$.
Write
\begin{equation}\label{eq:addable-set}
  A(S)=\bigl\{v\in S^3:\langle v,r\rangle\le\tfrac12\ \text{for every }
  r\in S\bigr\}
\end{equation}
for the directions that may be added to $S$, so that a contact
configuration $W$ with $W\cap\Rt=S$ satisfies $W\subseteq S\cup A(S)$.
Call $\{i,j\}$ the \emph{support} of the four roots
$(\pm e_i\pm e_j)/\sqrt2$, and call a support \emph{filled} for $S$ when
all four of its roots lie in $S$. The six supports pair off into three
couples of complementary index pairs,
\begin{equation}\label{eq:couples}
  \{1,2\}\,\big|\,\{3,4\},
  \qquad
  \{1,3\}\,\big|\,\{2,4\},
  \qquad
  \{1,4\}\,\big|\,\{2,3\},
\end{equation}
and these three couples are what carry the whole argument. Throughout,
$t=1/\sqrt2$ and $v=(x_1,x_2,x_3,x_4)$.

\begin{lem}[A filled support]\label{lem:filled-support}
If $\{i,j\}$ is filled for $S$, then every $v\in A(S)$ satisfies
$|x_i|+|x_j|\le t$, and hence $x_i^2+x_j^2\le\tfrac12$.
\end{lem}

\begin{proof}
The four conditions $\langle v,(\pm e_i\pm e_j)/\sqrt2\rangle\le\tfrac12$
read $\pm x_i\pm x_j\le t$ over the four sign patterns, and the pattern
matching the signs of $x_i$ and $x_j$ gives $|x_i|+|x_j|\le t$. Then
$x_i^2+x_j^2\le(|x_i|+|x_j|)^2\le t^2=\tfrac12$.
\end{proof}

\begin{lem}[A filled couple]\label{lem:filled-couple}
If both supports of one of the couples \eqref{eq:couples} are filled for
$S$, then $A(S)=\Rt\setminus S$.
\end{lem}

\begin{proof}
Let $\{i,j\}$ and $\{k,l\}$ be the couple and let $v\in A(S)$. By
\cref{lem:filled-support}, $x_i^2+x_j^2\le\tfrac12$ and
$x_k^2+x_l^2\le\tfrac12$, and the two sides add to $|v|^2=1$, so both are
equalities. Equality in $x_i^2+x_j^2\le(|x_i|+|x_j|)^2$ forces
$x_ix_j=0$, and equality in $(|x_i|+|x_j|)^2\le t^2$ forces
$|x_i|+|x_j|=t$; so one of $x_i,x_j$ is zero and the other is $\pm t$,
and the same holds for $x_k,x_l$. Therefore $v=\pm te_p\pm te_q$ with
$p\in\{i,j\}$ and $q\in\{k,l\}$, which is a root. A root in $A(S)$ cannot
lie in $S$, since $\langle v,v\rangle=1>\tfrac12$; and every root outside
$S$ does lie in $A(S)$, because distinct roots have inner product at most
$\tfrac12$.
\end{proof}

\begin{prop}[Twenty-two roots]\label{prop:meet22}
If a contact configuration $W$ satisfies $|W\cap\Rt|\ge22$ for some root
system $\Rt$, then $W\subseteq\Rt$.
\end{prop}

\begin{proof}
Put $S=W\cap\Rt$. At most two roots are missing, and they lie in at most
two of the six supports, so at least one of the three couples has both
supports filled. By \cref{lem:filled-couple}, $A(S)=\Rt\setminus S$, and
$W\subseteq S\cup A(S)=\Rt$.
\end{proof}

At $21$ the conclusion changes, and exactly one family of exceptions
appears.

\begin{prop}[Twenty-one roots]\label{prop:meet21}
If a contact configuration $W$ satisfies $|W\cap\Rt|=21$ for some root
system $\Rt$, then $W\subseteq\Rt$ or $|W|\le22$.
\end{prop}

\begin{proof}
Let $S=W\cap\Rt$ and let $\rho_1,\rho_2,\rho_3$ be the missing roots. If
their supports occupy at most two of the six index pairs, some couple is
filled, so $A(S)=\{\rho_1,\rho_2,\rho_3\}$ by \cref{lem:filled-couple}
and $W\subseteq\Rt$. So assume the three supports are distinct and meet
each couple of \eqref{eq:couples} once. Choosing one index pair from each
couple gives eight possibilities, and inspection of the list shows that
four of them are \emph{stars}, three pairs sharing an index, as with
$\{1,2\},\{1,3\},\{1,4\}$, and four are \emph{triangles}, the three pairs
inside a triple of indices, as with $\{1,2\},\{1,3\},\{2,3\}$. Permuting
the coordinates reduces to those two.

Let $v\in A(S)$ and suppose $v$ is not a root. For each $m$ the set
$S\cup\{\rho_m\}$ misses only two roots, so $A(S\cup\{\rho_m\})$ consists
of roots by \cref{lem:filled-couple}, and therefore does not contain $v$.
Since $v\in A(S)$, the only condition that can fail at $S\cup\{\rho_m\}$
is the one at $\rho_m$ itself, so
\begin{equation}\label{eq:inside-sixty}
  \langle v,\rho_m\rangle>\tfrac12,\qquad m=1,2,3 .
\end{equation}

Suppose first that the supports form the star $\{1,2\},\{1,3\},\{1,4\}$.
Sign changes in the coordinates $2,3,4$ are symmetries of $\Rt$, so we
may take $\rho_m=(\varepsilon_me_1+e_m)/\sqrt2$ with
$\varepsilon_m\in\{\pm1\}$ for $m=2,3,4$. The supports $\{2,3\}$,
$\{2,4\}$ and $\{3,4\}$ are filled, so $|x_m|\le t$ for $m=2,3,4$ by
\cref{lem:filled-support}, and \eqref{eq:inside-sixty} reads
$\varepsilon_mx_1+x_m>t$. If $\varepsilon_m=1$ this gives
$x_1>t-x_m\ge0$, and if $\varepsilon_m=-1$ it gives $x_1<x_m-t\le0$; the
two cannot both occur, so the three signs agree, and after a sign change
in the first coordinate $\rho_m=(e_1+e_m)/\sqrt2$ for $m=2,3,4$. Those
three roots are pairwise at $60^\circ$.

Suppose instead that the supports form the triangle
$\{1,2\},\{1,3\},\{2,3\}$, and write
$\rho_1=(ae_1+be_2)/\sqrt2$, $\rho_2=(ce_1+de_3)/\sqrt2$,
$\rho_3=(fe_2+ge_3)/\sqrt2$ with signs $a,\dots,g$. The supports
$\{1,4\},\{2,4\},\{3,4\}$ are filled, so $|x_m|\le t$ for $m=1,2,3$.
Adding the three inequalities \eqref{eq:inside-sixty} gives
\[
  (a+c)\,x_1+(b+f)\,x_2+(d+g)\,x_3>3t ,
\]
and the left side is at most $2tk$, where $k$ counts the nonzero terms
among $a+c$, $b+f$, $d+g$. Hence $k=3$, that is $a=c$, $b=f$ and $d=g$,
and $\rho_1,\rho_2,\rho_3$ are again pairwise at $60^\circ$. Changing the
signs of the coordinates $1,2,3$ as needed makes the triple
$\{(e_1+e_2),(e_1+e_3),(e_2+e_3)\}/\sqrt2$, and the orthogonal map with
matrix $\tfrac12\bigl(\begin{smallmatrix}1&1&1&1\\1&1&-1&-1\\
1&-1&1&-1\\1&-1&-1&1\end{smallmatrix}\bigr)$, which permutes $\Rt$,
carries it to $\{(e_1+e_2),(e_1+e_3),(e_1-e_4)\}/\sqrt2$; one more sign
change in the fourth coordinate returns the star.

So in every case we may take $\rho_m=(e_1+e_m)/\sqrt2$ for $m=2,3,4$, and
it remains to read off $A(S)$ for that $S$. The supports $\{2,3\}$,
$\{2,4\}$, $\{3,4\}$ are filled, so $|x_2|+|x_3|$, $|x_2|+|x_4|$ and
$|x_3|+|x_4|$ are all at most $t$. The maximum of $p^2+q^2+r^2$ over
$\{p,q,r\ge0:\ p+q\le t,\ p+r\le t,\ q+r\le t\}$ is attained at a vertex
of that polytope, and its vertices are the origin, the three points
$te_k$ and $\tfrac t2(1,1,1)$, with values $0$, $t^2$ and $\tfrac34t^2$.
Hence
\begin{equation}\label{eq:doorway-coords}
  x_2^2+x_3^2+x_4^2\le t^2=\tfrac12,
  \qquad\text{so}\qquad
  x_1^2\ge\tfrac12 .
\end{equation}
Each support $\{1,m\}$ keeps three of its four roots, all but
$(e_1+e_m)/\sqrt2$, and these give $x_1-x_m\le t$ and $|x_m|-x_1\le t$.
If $x_1\le-t$ the second forces $x_m=0$ for $m=2,3,4$, hence $v=-e_1$,
which the first contradicts at $m=2$. So $x_1\ge t$, and then
$x_m\ge x_1-t\ge0$ for $m=2,3,4$: on $A(S)$ the first coordinate is at
least $t$ and the other three are nonnegative. Consequently, for
$v,v'\in A(S)$,
\begin{equation}\label{eq:doorway-pair}
  \langle v,v'\rangle=x_1x_1'+x_2x_2'+x_3x_3'+x_4x_4'
  \;\ge\;x_1x_1'\;\ge\;t^2=\tfrac12 ,
\end{equation}
with equality only if $x_1=x_1'=t$. That makes
\eqref{eq:doorway-coords} an equality, so $(x_2,x_3,x_4)$ is one of the
vertices $te_k$, and $v\in\{\rho_1,\rho_2,\rho_3\}$; the same holds for
$v'$.

Now $W\setminus S$ lies in $A(S)$ and its elements are pairwise at
$60^\circ$ or more. By \eqref{eq:doorway-pair} any two of them are
missing roots, so either every element of $W\setminus S$ is a root, in
which case $W\subseteq\Rt$, or $W\setminus S$ has a single element and
$|W|=22$.
\end{proof}

\begin{cor}[A saturated configuration is far from every root system]\label{cor:meet20}
A contact configuration of $23$ directions whose cell has circumradius
below $2$ meets every root system in at most $20$ directions.
\end{cor}

\begin{proof}
Let $|W|=23$ and suppose $|W\cap\Rt|\ge21$. By \cref{prop:meet22} and
\cref{prop:meet21}, either $|W|\le22$, which is false, or $W\subseteq\Rt$
and $W$ is $\Rt$ with one root $\base_0$ removed. In the second case the
cell of $W$ contains the point $2\base_0$, by \cref{lem:Q-structure}, so
its circumradius is $2$.
\end{proof}

\begin{figure}[tb]
\centering
\includegraphics[width=\textwidth]{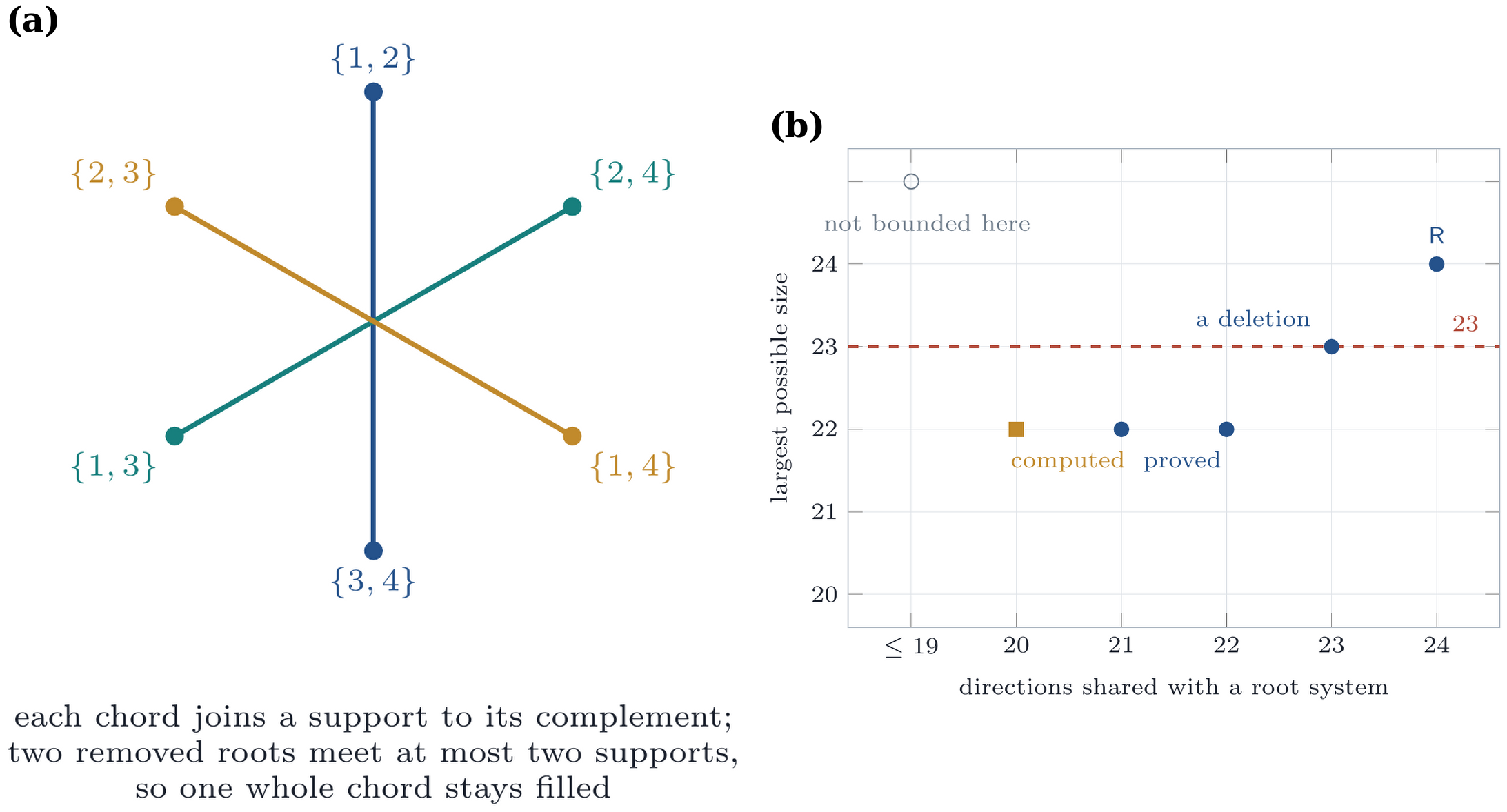}
\caption{(a) The six supports of the roots, joined in the three couples of
complementary index pairs. Removing two roots touches at most two
supports and leaves a whole couple filled, which is
\cref{lem:filled-couple}; removing three can meet each couple once, and
the eight ways of doing that are four stars and four triangles. (b) The largest size a contact configuration can have, against the
number of directions it shares with a root system. The value $23$ occurs
only inside the root system, where the circumradius is exactly $2$, so
the saturated case is pushed to twenty directions or fewer.}
\label{fig:rootmeet}
\end{figure}

\Cref{fig:rootmeet} collects the two ingredients. The bound of
\cref{cor:meet20} is the first statement we have that separates a
hypothetical saturated configuration from the root system by a fixed amount, not by an unspecified small one: three of its
twenty-three directions lie outside every copy of $\Rt$, so it is not
reachable from a deletion by exchanging one or two directions.

\begin{remark}[The next rung, and what it costs]\label{rem:meet-next}
The same reduction can be tried at twenty roots. Of the $10626$ ways to
remove four roots, only four orbits under the symmetry group of $\Rt$
leave a cell of circumradius above $2$, and in each of those the removed
set is a triple pairwise at $60^\circ$ together with one further root.
For each of the four we minimised the largest of the three pairwise
inner products over triples of directions in $A(S)$, requiring each
direction to stay an angle $\delta$ away from every root, so that the
removed roots themselves are excluded. At $\delta=0$ the value is
$\tfrac12$ in all four, and the minimiser is three of the removed roots;
at $\delta=2^\circ$ the smallest of the four values is already
$0.538256\ldots$, at $\delta=5^\circ$ it is $0.592985\ldots$, at
$\delta=10^\circ$ it is $0.676602\ldots$, and it increases from there. So the
doorways at twenty roots hold two directions but not three, which would
push \cref{cor:meet20} to nineteen. We record this as a computation and
not as a proof: it is a search over a family of local minima, not an
exhaustive argument, and the case analysis that would make it one is
four times the length of \cref{prop:meet21}. The script is
\texttt{root\_meet.py}, with \texttt{root\_deletions\_exact.py} carrying
out the same enumeration for one, two and three removed roots in exact
integer arithmetic as a check on \cref{prop:meet22} and
\cref{prop:meet21}.
\end{remark}

\subsection{The remaining case as a question about codes}\label{sec:m23-codes}

A contact configuration of $m$ directions is the same object as a
spherical code of $m$ points on $S^3$ with minimal angle at least
$60^\circ$, and the question is how much room a code of that size has.
Write $\theta(m)$ for the largest minimal angle available to $m$
points on $S^3$. Musin's theorem \cite{Mus08} gives
$\theta(25)<60^\circ$, and \cref{thm:m24} gives $\theta(24)=60^\circ$,
since a code of $24$ points at more than $60^\circ$ would be a contact
configuration of that size and hence a copy of the root system, whose
minimal angle is exactly $60^\circ$. The value at $m=23$ is what decides
the shape of the set of configurations that \cref{thm:m23} is about.

\begin{prop}[Slack, and its absence]\label{prop:code-slack}
If $\theta(m)>60^\circ$, the contact configurations of $m$ directions
form a set with nonempty interior in $(S^3)^m$. If $\theta(m)=60^\circ$,
that set has empty interior, and every contact configuration of $m$
directions attains $\theta(m)$.
\end{prop}

\begin{proof}
The contact configurations of size $m$ are exactly the $m$-point codes
of minimal angle at least $60^\circ$, and the minimal angle is a
continuous function on $(S^3)^m$. If $\theta(m)>60^\circ$, a code
attaining $\theta(m)$ has a neighbourhood on which the minimal angle
stays above $60^\circ$, and every configuration in it is a contact
configuration. If instead $\theta(m)=60^\circ$, then a contact
configuration $W$ has minimal angle exactly $60^\circ$, so some pair
$w_i,w_j$ is at $60^\circ$; rotating $w_i$ towards $w_j$ through any
positive angle, however small, gives a configuration at distance as
small as one likes from $W$ whose minimal angle is below $60^\circ$, and
which is therefore not a contact configuration. So $W$ is not an
interior point.
\end{proof}

At $m\le22$ the first alternative is the one that holds, and the
configurations at those sizes fill out an open set; that is why
\cref{cor:m22} has to hold uniformly over a family of positive
dimension, and it does. At $m=23$ every computation we know of, our own
in \cref{rem:code-search} and the published tables \cite{Slo}, returns
$\theta(23)=60^\circ$, and if that is the true value the second
alternative takes over: the configurations that \cref{thm:m23} is about
are extremal codes, not a family with interior. That is the structural
difference between $23$ and every size below it, and a set with empty
interior can of course turn out to have no saturated member at all. We
have no proof that $\theta(23)=60^\circ$; what \cref{sec:slack-continuation}
shows is how the codes of $23$ points behave as their minimal angle is
pushed towards $60^\circ$ from below.

\begin{thm}[The $23$-point case from a uniqueness statement]\label{thm:m23-from-codes}
Suppose every set of $23$ unit vectors in $\mathbb{R}^4$ with pairwise
inner products at most $\tfrac12$ is the image under an orthogonal map
of $\Rt\setminus\{\base\}$ for some root $\base$. Then no saturated
$23$-point contact configuration exists, and the case that
\cref{thm:m23} covers is empty.
\end{thm}

\begin{proof}
Let $W$ be a contact configuration with $|W|=23$, and take an orthogonal
$O$ and a root $\base$ with $W=O(\Rt\setminus\{\base\})$. Every $w\in W$
is $Or$ for a root $r\ne\base$, so
$\langle O\base,w\rangle=\langle\base,r\rangle\le\tfrac12$, and
$W\cup\{O\base\}$ is again a contact configuration. By
\cref{prop:extendable} the circumradius of $V_c(W)$ is then at least
$2$, so $W$ is not saturated, and \cref{cor:remaining} then needs
nothing beyond \cref{thm:m24}.
\end{proof}

\begin{remark}[Why a uniqueness statement and not a bound]\label{rem:m23-uniqueness}
The hypothesis of \cref{thm:m23-from-codes} is what
\cref{thm:m24} proves at $24$ points, asked at $23$, and one
difference between the two carries the whole of the difficulty. At $24$
the code has maximal size, so the semidefinite relaxation is tight
there and complementary slackness forces every inner product into
$\{-1,-\tfrac12,0,\tfrac12\}$, which is \cref{thm:twentyfour-points};
the classification is then combinatorial. At $23$ nothing is tight. No upper
bound on the size of a code in terms of its minimal angle can be of any
use here, whatever the method that produces it: such a bound is at least
$24$ at $60^\circ$, because a code of $24$ points at that angle exists,
so it cannot exclude $23$ points at the same angle. What is wanted is a
classification at one fixed size below the maximum, which is a different
kind of statement and, as far as we know, has not been attempted in this
dimension.
\end{remark}

\begin{remark}[What the search returns]\label{rem:code-search}
We looked for a $23$-point code of minimal angle above $60^\circ$, and
for a contact configuration of $23$ directions whose Gram matrix is
unlike that of a deletion, by Riesz continuation: minimise
$\sum_{i<j}|w_i-w_j|^{-s}$ over the sphere with $s$ running through
$4,16,64,256,1024$, then reduce the largest inner product directly.
Calibration comes first. Over $120$ starts the same procedure returns
$0.5000071\ldots$ at $m=24$, where the true value is $\tfrac12$;
$0.4980170\ldots$ at $m=22$, a minimal angle of $60.1311\ldots^\circ$,
against the $60.1398863^\circ$ of the published tables \cite{Slo}; and,
over $60$ starts, $0.5374065\ldots$ at $m=25$, where no contact
configuration exists. So the procedure locates the root system, sees
the slack at $22$ and sees the obstruction at $25$, and one can read how
close to $\tfrac12$ it gets when $\tfrac12$ is the answer. At $m=23$,
over $1500$ random starts, the smallest largest inner product it reached
was $0.5000078\ldots$, which is nearer to $\tfrac12$ than the calibration
run at $m=24$ managed; no start went below $\tfrac12$; and of the $32$
outcomes that came within $5\times10^{-3}$ of feasibility, all $32$ had
the inner-product multiset of a deletion, namely $-1$ eleven times,
$-\tfrac12$ eighty-eight times, $0$ sixty-six times and $\tfrac12$
eighty-eight times, which is an invariant of the orbit under $O(4)$.
The tables \cite{Slo} record $60.0000000^\circ$ at $23$ points, the same
value as at $24$, and $60.1398863^\circ$ at $22$. None of this proves
the hypothesis of \cref{thm:m23-from-codes}, and we do not offer it as a
proof; what it does is say which statement to aim at. The script is
\texttt{spherical\_code\_23.py}.

\end{remark}

One finite family can be searched completely instead of sampled, and it
is the natural place for a code that is not a subset of a root system to
hide. The $600$-cell has $120$ vertices, the unit icosians, and its
inner products are $0$, $\pm\tfrac12$, $\pm\cos36^\circ$,
$\pm\cos72^\circ$ and $\pm1$; two vertices are closer than $60^\circ$
exactly when they are joined by an edge, at $36^\circ$. Its subsets with
all inner products at most $\tfrac12$ are therefore the independent sets
of the edge graph, which is $12$-regular, and the largest of them have
$24$ elements and are the twenty-five inscribed copies of $\Rt$, five
through each vertex. The question is whether a $23$-element independent
set can avoid all twenty-five.

\begin{prop}[Codes among the vertices of the $600$-cell]\label{prop:cell600}
Let $I\subset S^3$ be the vertex set of the regular $600$-cell. A subset
of $I$ whose pairwise inner products are all at most $\tfrac12$ has at
most $24$ elements; the subsets with $24$ elements are the twenty-five
inscribed copies of $\Rt$; and every subset with $23$ elements is one of
those copies with a single vertex removed. In particular no saturated
contact configuration of $23$ directions lies in $I$.
\end{prop}

\begin{proof}
Write $\phi=\tfrac12(1+\sqrt5)$. Doubling the coordinates puts every
vertex of the $600$-cell in $\mathbb{Z}[\phi]^4$, with each coordinate
one of $0,\pm1,\pm2,\pm\phi,\pm(\phi-1)$, and $4\langle u,v\rangle$ is
then computed exactly in $\mathbb{Z}[\phi]$; the nine values that occur
are $0$, $\pm2$, $\pm4$, $\pm2\phi$ and $\pm(2\phi-2)$, of which only
$2\phi$, that is $\langle u,v\rangle=\cos36^\circ$, exceeds $2$. So the
sets in question are exactly the independent sets of the graph $\Gamma$
on $I$ whose edges join vertices at $36^\circ$, and every vertex of
$\Gamma$ has degree $12$, its vertex figure being an icosahedron. The
symmetry group of the $600$-cell is transitive on $I$, so it suffices to
enumerate the independent sets through one fixed vertex $v_0$, and this
is a finite search which we carry out to the end: a depth-first
enumeration over bitsets, with the standard pruning that an independent
set meets every clique of $\Gamma$ at most once, so that a greedy
partition of the remaining candidates into cliques bounds what can still
be added.

The search returns no independent set of size $25$ through $v_0$;
exactly five of size $24$, each of which has all its pairwise inner
products in $\{-1,-\tfrac12,0,\tfrac12\}$ and is therefore a copy of
$\Rt$ (by \cref{thm:m24}, or directly, since the $24$-point contact
configurations inside $I$ are checked one by one); and exactly $115$ of
size $23$, every one of which is contained in one of the five copies of
$\Rt$ through $v_0$. Since a copy of $\Rt$ through $v_0$ has $23$ other
vertices and $5\times23=115$, the sets of size $23$ through $v_0$ are
precisely those five copies with one of their other vertices removed,
and nothing else. Summing over the $120$ vertices, $115\times120/23=600
=25\times24$ independent sets of size $23$ in all, one for each pair
consisting of an inscribed copy of $\Rt$ and one of its vertices, and
$5\times120/24=25$ inscribed copies.

A set of $23$ vertices of one copy of $\Rt$ is a deletion, and is
extended by the vertex it lacks, so by \cref{prop:extendable} its cell
has circumradius $2$ and it is not saturated.

The enumeration is carried out three times, independently. Once in
Python, with the vertices in $\mathbb{Z}[\phi]$ and the graph and the
search as described; once in C, on the graph the Python script writes
out; and once in Lean~4, where the vertex list, the nine inner-product
values, the degree $12$, the twenty-five cells, their inner products and
the fact that five pass through each vertex are all checked by the
kernel from the $\mathbb{Z}[\phi]$ coordinates, and the three counts
$0$, $5$ and $115$ are established by the same depth-first search
compiled to native code and run by \texttt{native\_decide}. The scripts
are \texttt{cell600\_exact.py} and \texttt{cell600\_enum.c}, and the
Lean development is described in \cref{sec:reproducibility-new}.
\end{proof}

\begin{figure}[tb]
\centering
\includegraphics[width=0.86\textwidth]{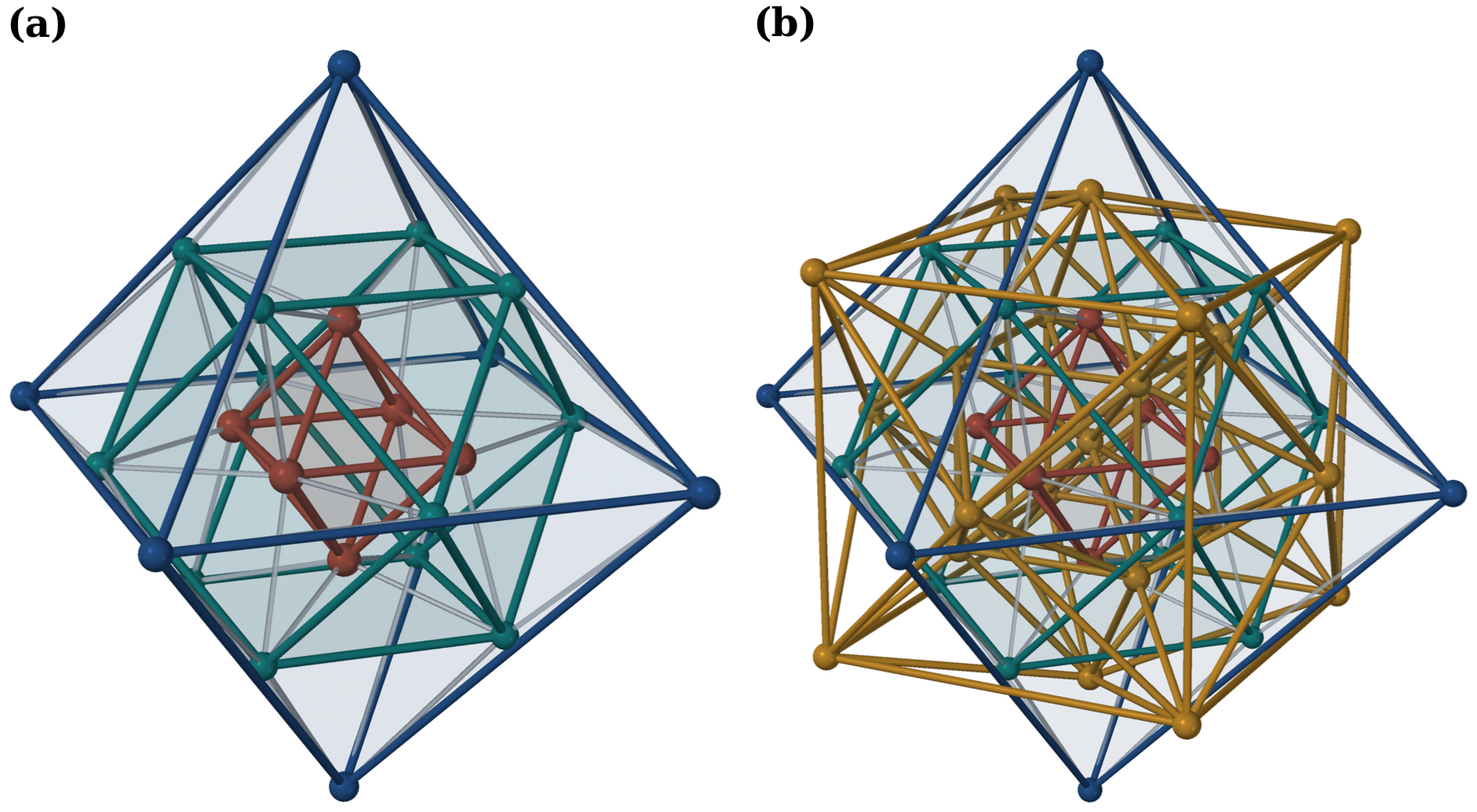}
\caption{(a) The $24$ contact directions of $D_4$ seen from a deep hole, radius proportional to angle: an octahedron at $45^\circ$ (clay), a cuboctahedron at $90^\circ$ (teal), an octahedron at $135^\circ$ (blue), and the $96$ pairs at $60^\circ$ as edges. (b) The same, with the dual $24$-cell added (ochre): the $48$ directions of the binary octahedral group, a cube, an octahedron and a cube about the hole, which is now itself a direction.}
\label{fig:structab}
\end{figure}

\begin{figure}[tb]
\centering
\includegraphics[width=0.86\textwidth]{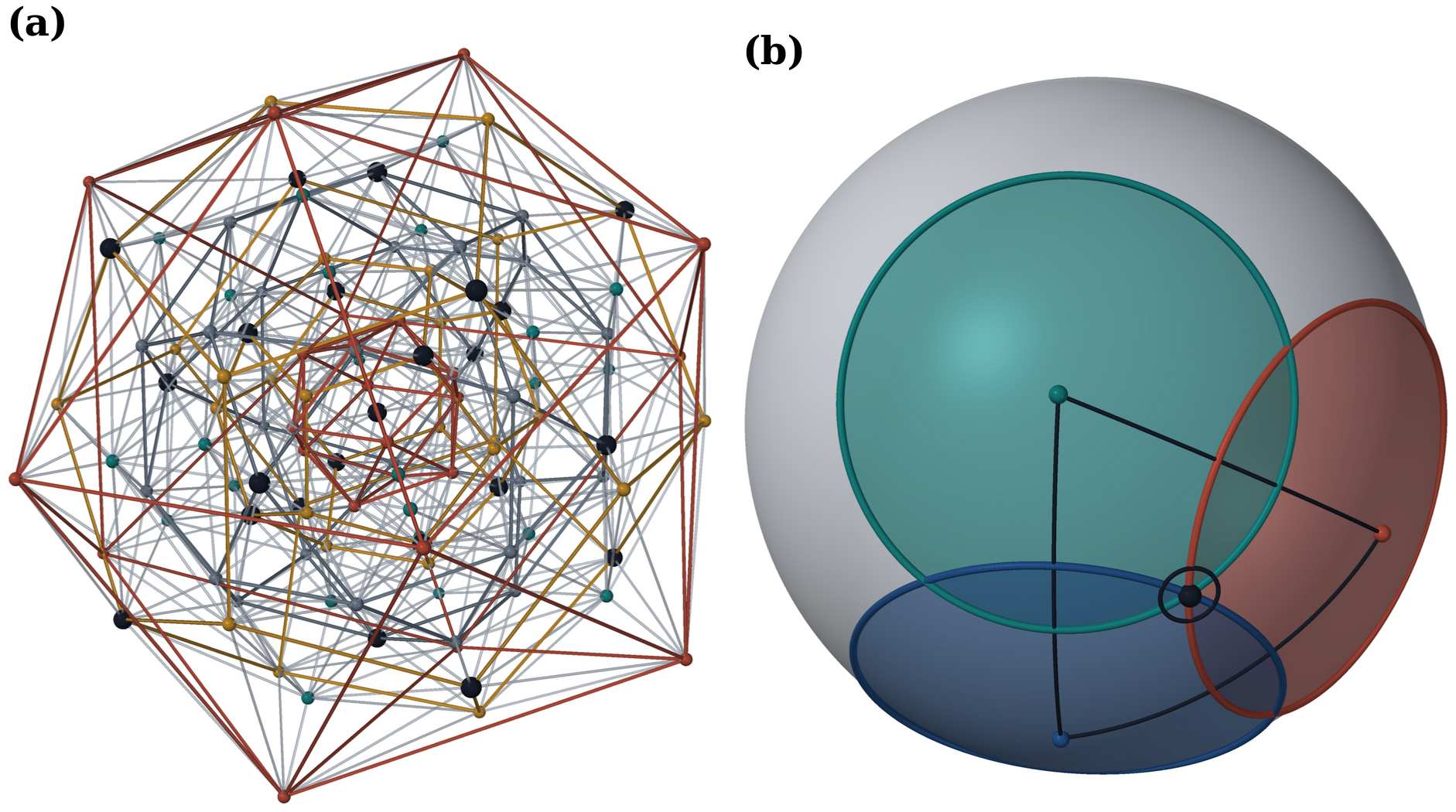}
\caption{(a) The $600$-cell seen from one vertex, its $119$ other vertices in shells at $36^\circ$ to $144^\circ$ with $708$ of the $720$ edges at $36^\circ$; the dark points are an inscribed $24$-cell whose far vertex is the one at infinity, that is, a deletion. (b) The cap lemma on a great $2$-sphere: three contact directions $60^\circ$ apart carry caps of radius $r_*=35.26^\circ$ that meet in exactly one point, the circumcentre of their triangle (ringed); any smaller radius separates them.}
\label{fig:structcd}
\end{figure}

\begin{remark}[What the $600$-cell does and does not decide]\label{rem:cell600-scope}
\Cref{prop:cell600} closes the family one reaches for first
(\cref{fig:structcd}(a), drawn on page~\pageref{fig:structcd}). The unit
icosians form the largest finite subgroup of the unit quaternions, the
binary icosahedral group of order $120$; they contain the twenty-four
Hurwitz units, which are $\Rt$, in twenty-five different ways; and a
subset of them with all inner products at most $\tfrac12$ is the first
thing to try when looking for a code built on the golden ratio rather
than on $\sqrt2$. The proposition says that no such subset of $23$
points escapes a copy of $\Rt$. It says nothing about a code whose coordinates lie outside
$\mathbb{Q}(\sqrt5)$, and nothing about the continuous families that
\cref{sec:slack-continuation} finds once the minimal angle is allowed to
drop below $60^\circ$, none of which has the inner-product spectrum of a
subset of the $600$-cell. The random sampling of maximal independent
sets that preceded the enumeration, two million of them drawn by greedy
growth along random orders of the vertices, had already found sizes
$10$ to $22$ and $24$ and never $23$; the enumeration replaces that
observation by a statement. The same enumeration on the $48$ unit
quaternions of the binary octahedral group, the vertices of two
$24$-cells in dual position (\cref{fig:structab}(b)), in
which two directions are closer than $60^\circ$ exactly when they lie
in different copies and are $45^\circ$ apart, returns the two copies of
$\Rt$ and their forty-eight deletions
and nothing else; that one is a second of computation, in
$\mathbb{Z}[\sqrt2]$, and is \texttt{octahedral48\_exact.py}.
\end{remark}

\subsection{The largest inradius as the minimal angle rises}\label{sec:slack-continuation}

The remaining case can be put as the value of one function at one
point, and the function can be computed away from that point. For
$\delta\ge0$ let
\[
  K(\delta)=\bigl\{W=(w_1,\dots,w_{23})\in(S^3)^{23} :
  \langle w_i,w_j\rangle\le\tfrac12+\delta\ \text{for all } i\ne j\bigr\},
\]
the configurations of $23$ directions whose minimal angle is at least
$\arccos(\tfrac12+\delta)$, and let
\[
  h(\delta)=\max_{W\in K(\delta)}\ g(W),
  \qquad
  g(W)=\min_{\theta\in S^3}\ \max_i\ \langle\theta,w_i\rangle ,
\]
so that $g$ is the function of \cref{prop:extendable}, the cosine of the
covering radius, and the inradius of $\conv(W)$ whenever that hull
contains the origin. The contact configurations are $K(0)$, and by
\cref{prop:extendable} a member of $K(0)$ is saturated exactly when
$g>\tfrac12$.

\begin{lem}[The limit from above]\label{lem:h-limit}
The function $h$ is nondecreasing on $[0,\infty)$, $h(0)\ge\tfrac12$,
and $h(\delta)\to h(0)$ as $\delta\downarrow0$. No saturated contact
configuration of $23$ directions exists if and only if $h(0)=\tfrac12$;
and this holds as soon as $h(\delta)\le\tfrac12+C\delta^{\kappa}$ on some
interval $(0,\delta_0)$, for any constants $C$ and $\kappa>0$.
\end{lem}

\begin{proof}
The map $(\theta,W)\mapsto\max_i\langle\theta,w_i\rangle$ is continuous,
and $|\max_i\langle\theta,w_i\rangle-\max_i\langle\theta,w_i'\rangle|
\le\max_i|w_i-w_i'|$, so $g$ is continuous on $(S^3)^{23}$, indeed
$1$-Lipschitz in the maximum norm. Each $K(\delta)$ is a closed subset
of the compact space $(S^3)^{23}$, so the maximum defining $h(\delta)$
is attained; the sets increase with $\delta$, so $h$ is nondecreasing;
and $K(0)$ contains the deletions of a root, at which $g=\tfrac12$, so
$h(0)\ge\tfrac12$. For the limit take $\delta_n\downarrow0$ and
$W_n\in K(\delta_n)$ with $g(W_n)=h(\delta_n)$. By compactness a
subsequence converges to some $W$, and $W\in K(\delta_m)$ for every $m$,
since $K(\delta_m)$ is closed and contains $W_n$ for all $n\ge m$; hence
$W\in\bigcap_mK(\delta_m)=K(0)$. By continuity $g(W)=\lim h(\delta_n)$,
and $g(W)\le h(0)$ because $W\in K(0)$, while $h(\delta_n)\ge h(0)$ by
monotonicity; so the limit is $h(0)$, along every such sequence. The
equivalence is the definition of saturation, and the last statement
follows by letting $\delta\downarrow0$ in the bound.
\end{proof}

So the hypothesis of \cref{thm:m23-from-codes} is the statement
$h(0)=\tfrac12$, and $h$ can be looked at
for $\delta>0$, where the constraint set has interior and a local
optimiser has room to move. We did this with the trust-region
sequential linear programme described in \cref{app:code-index}, which
maximises the smallest facet offset of $\conv(W)$ over tangent moves of
bounded size subject to the linearised contact constraints, accepting a
move only when the true inradius, read from a fresh convex hull,
improves without any constraint being violated. For each $\delta$ on a
decreasing schedule the starts are the twenty best configurations of
the previous level, restored to the new constraint, together with a
batch of fresh random starts; three passes were run, with different
seeds and with $40$, $60$ and $80$ fresh starts per level, and
\cref{tab:slack-continuation} records the largest inradius found at
each level over the three. Every number in the table is a lower bound
for $h(\delta)$ obtained by local search, and nothing in this
subsection is a proof.

\begin{table}[t]
\centering
\caption{The largest inradius $g$ found at each level of the
continuation, over three passes. The third column is the minimal angle
$\arccos(\tfrac12+\delta)$ imposed, the fifth the covering radius
$\arccos g$ of the best configuration, and the last two the number of
feasible endpoints across the three passes and how many of them are
deletions of a root. At $\delta=0$ every endpoint is a deletion and
$g=\tfrac12$ exactly.}
\label{tab:slack-continuation}
\small
\begin{tabular}{@{}rrrrrr@{}}
\toprule
$\delta$ & min.\ angle & $h(\delta)\ge$ & cov.\ radius & endpoints & deletions\\
\midrule
$0.05$ & $56.633^\circ$ & $0.676201$ & $47.453^\circ$ & $180$ & $0$\\
$0.03$ & $57.995^\circ$ & $0.695650$ & $45.921^\circ$ & $239$ & $0$\\
$0.02$ & $58.668^\circ$ & $0.686522$ & $46.645^\circ$ & $159$ & $0$\\
$0.015$ & $59.003^\circ$ & $0.673135$ & $47.691^\circ$ & $113$ & $0$\\
$0.01$ & $59.336^\circ$ & $0.657060$ & $48.924^\circ$ & $104$ & $0$\\
$0.009$ & $59.403^\circ$ & $0.653786$ & $49.172^\circ$ & $34$ & $0$\\
$0.008$ & $59.469^\circ$ & $0.508771$ & $59.418^\circ$ & $2$ & $0$\\
$0.007$ & $59.536^\circ$ & $0.508968$ & $59.405^\circ$ & $11$ & $0$\\
$0.005$ & $59.669^\circ$ & $0.505913$ & $59.608^\circ$ & $15$ & $0$\\
$0.003$ & $59.801^\circ$ & $0.504484$ & $59.703^\circ$ & $22$ & $0$\\
$0.002$ & $59.868^\circ$ & $0.502914$ & $59.807^\circ$ & $26$ & $0$\\
$0.0015$ & $59.901^\circ$ & $0.502552$ & $59.831^\circ$ & $34$ & $0$\\
$0.001$ & $59.934^\circ$ & $0.501392$ & $59.908^\circ$ & $40$ & $0$\\
$0.0007$ & $59.954^\circ$ & $0.501123$ & $59.926^\circ$ & $42$ & $0$\\
$0.0005$ & $59.967^\circ$ & $0.500878$ & $59.942^\circ$ & $46$ & $0$\\
$0.0003$ & $59.980^\circ$ & $0.500480$ & $59.968^\circ$ & $52$ & $0$\\
$0.0002$ & $59.987^\circ$ & $0.500378$ & $59.975^\circ$ & $56$ & $0$\\
$0.0001$ & $59.993^\circ$ & $0.500202$ & $59.987^\circ$ & $56$ & $0$\\
$0.00005$ & $59.997^\circ$ & $0.500081$ & $59.995^\circ$ & $63$ & $0$\\
$0$ & $60.000^\circ$ & $0.500000$ & $60.000^\circ$ & $51$ & $51$\\
\bottomrule
\end{tabular}
\end{table}

\begin{figure}[tb]
\centering
\includegraphics[width=\textwidth]{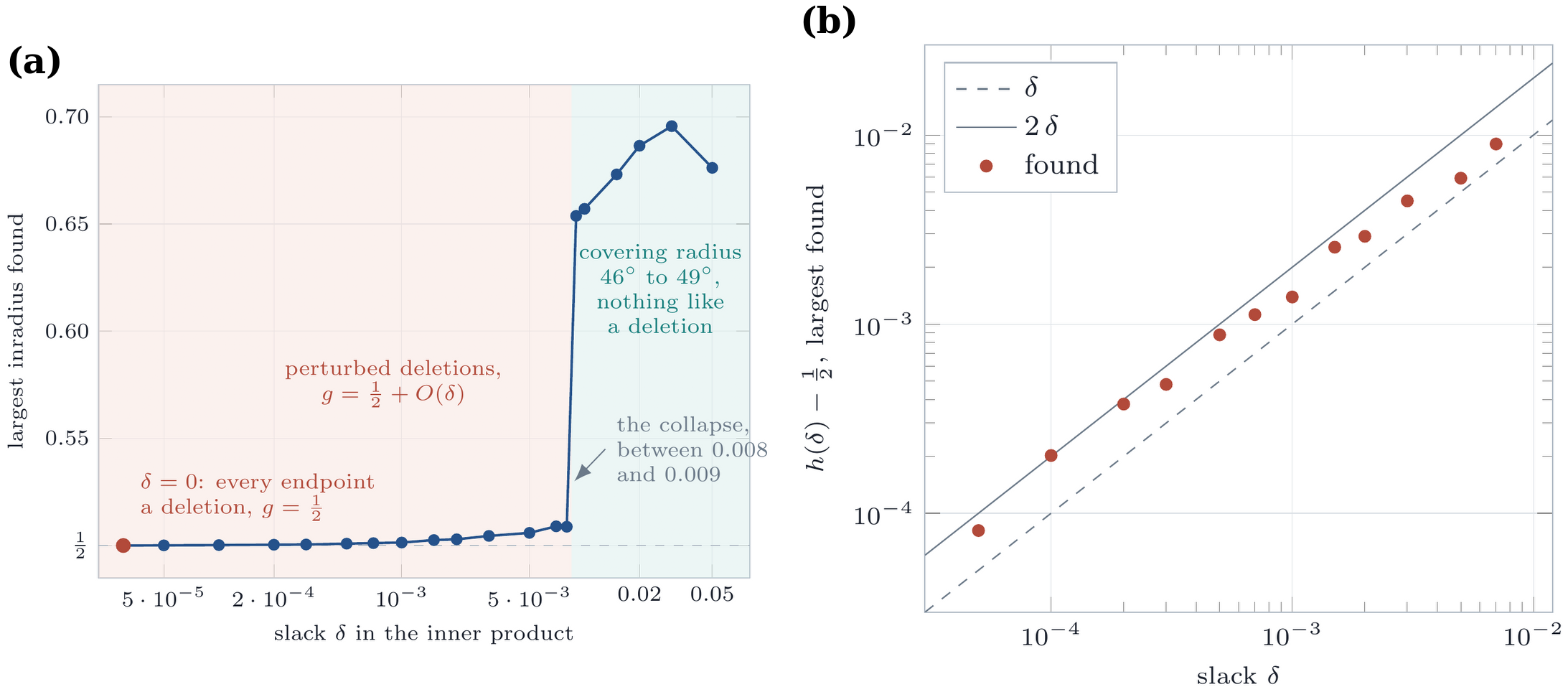}
\caption{(a) The continuation of \cref{tab:slack-continuation}: the largest
inradius found, against $\delta$ on a logarithmic axis. The plateau near
$\tfrac23$ runs down to $\delta=0.009$ and is carried by configurations
far from any deletion; from $\delta=0.008$ down it collapses to
$\tfrac12+O(\delta)$. (b) The same table below the collapse: $h(\delta)-\tfrac12$ against
$\delta$ on logarithmic axes, with the lines $\delta$ and $2\delta$ for
comparison. The approach to $\tfrac12$ is linear in $\delta$, which is
what a perturbed deletion gives.}
\label{fig:slack}
\end{figure}

Three things stand out. First, for $\delta\ge0.009$, that is for minimal
angles up to $59.40^\circ$, the best configurations have inradius
between $0.65$ and $0.70$: their covering radius is between $46^\circ$
and $49^\circ$, so twenty-three caps of radius $49^\circ$ about their
directions cover the sphere, where the deletion of a root needs caps of
radius $60^\circ$ and the root system itself, with twenty-four
directions, needs $45^\circ$. These configurations are not perturbed
deletions in any sense. Their inner-product multisets have no pair at
$-1$ at all, where a deletion has eleven antipodal pairs, and none
within $10^{-3}$ of $-\tfrac12$ or of $0$; between three and
fifty-three of their pairs sit at the imposed bound $\tfrac12+\delta$,
and fewer than three per cent of their inner products lie within
$10^{-3}$ of a value taken by the $600$-cell, so none lies in the
$600$-cell either. Second, between $\delta=0.009$ and $\delta=0.008$ they
vanish. At every level from $0.008$ down the best inradius is
$\tfrac12+O(\delta)$, with $h(\delta)-\tfrac12$ between $1.1\delta$ and
$2\delta$ across the whole range down to $\delta=5\times10^{-5}$, and
the configurations attaining it are deletions of a root with their
$88$ tight pairs opened by amounts of order $\delta$: every one of them
keeps the eleven antipodal pairs of a deletion, and as $\delta$
decreases their inner-product multisets settle onto the deletion's,
which is reached to $10^{-3}$ by $\delta=3\times10^{-4}$. The transition is
not gradual: a configuration of the first kind, restored to the
constraint at $\delta=0.008$ by the feasibility programme, does not
reach it, because its minimal angle sits at a local maximum below
$59.47^\circ$ from which no local move raises it. Third, at $\delta=0$
every feasible endpoint of every pass is a deletion of a root, with the
inner-product multiset $-1$ eleven times, $-\tfrac12$ eighty-eight
times, $0$ sixty-six times, $\tfrac12$ eighty-eight times, and
$g=\tfrac12$ to fifteen digits. \Cref{fig:slack} draws the table.

\smallskip\noindent\textit{Configurations with a symmetry.}
A configuration invariant under a nontrivial element of $O(4)$ is a
union of orbits of that element, and inside such a class the search has
far fewer parameters and can be run from enough starts to cover it.
\Cref{prop:no-transitive} disposes of the transitive case. The script
\texttt{symmetric\_search.py} treats every rotation of order $n\le12$,
through angles $2\pi a/n$ and $2\pi b/n$ in two orthogonal planes with
$\gcd(a,b,n)=1$, together with the two improper involutions, a
reflection in a hyperplane and $\diag(1,-1,-1,-1)$; for each it lists
every way of writing $23$ as a sum of orbit sizes, the generic orbit
having $n$ points and a point of either plane an orbit of $n/\gcd(a,n)$
or $n/\gcd(b,n)$, with at most six points on either plane since they lie
on a great circle. That gives $73$ orbit structures, with between $5$
and $35$ parameters each, and in each the continuation of
\cref{tab:slack-continuation} is run from sixteen random starts on the
orbit representatives, through $\delta=0.05,0.02,0.01,0.005,0.002,0$.
Eleven of the seventy-three structures reached a feasible endpoint at
$\delta=0$, sixty-nine endpoints in all, and every one of them is a
deletion of a root: the structures are those with a rotation of order
$2$, $3$ or $4$ fixing a plane pointwise, a reflection with five fixed
points, and $\diag(1,-1,-1,-1)$, which are symmetries that a deletion
has, its symmetry group being the stabiliser of the deleted root in the
symmetry group of $\Rt$. The remaining sixty-two structures never reached
the constraint at $\delta=0$ from any start. At $\delta=0.01$ the best
symmetric configuration has inradius $0.6599$, in line with
\cref{tab:slack-continuation}. No symmetric configuration crossed into
$K(0)$ except the deletions.
\begin{remark}[What this says about the remaining case]\label{rem:slack-continuation}
The picture is the one the classification of \cref{thm:m23-from-codes}
would present if it holds, and it says something about where a proof of
it would have to work. By \cref{lem:h-limit} it would be enough to show
$h(\delta)\le\tfrac12+C\delta$ on some interval $(0,\delta_0)$, and the
data put the constant near $2$ for small $\delta$; but the same data
put $h(0.009)$ above $0.65$, so no such bound holds with a uniform
constant beyond $\delta_0\approx0.008$, and an argument would have to
work inside a window of about half a degree in the minimal angle. That
is the quantitative form of what \cref{rem:m23-uniqueness} said: the
classification at $23$ points is a statement about codes of minimal
angle exactly $60^\circ$, and codes at $59.3^\circ$ already behave in a
completely different way. It also says what a counterexample, if there
is one, would have to look like: a saturated configuration is a member
of $K(0)$ with $g>\tfrac12$, hence a member of every $K(\delta)$ with
$h(\delta)\ge g$, and it would have to have escaped a search that, at
every level below $0.009$, found nothing but deletions. We found no
such configuration, in three passes of the continuation and in $300$
direct maximisations of $g$ over $K(0)$ from the three kinds of start
described in \cref{app:code-index}: $94$ of those reached a feasible
endpoint, every one of them a deletion with $g=\tfrac12$, and the $100$
random starts, brought down through the same schedule of slacks, never
reached the constraint at $\delta=0$ at all, which is the collapse seen
from the other side. That is exploration, and it is reported as such.
\end{remark}

\subsection{The remaining case as a strict inequality}\label{sec:strict-inequality}

Everything in this section so far has approached \cref{thm:m23} through
the configuration: what a saturated set of $23$ directions would have
to look like, and why none has been found. This subsection approaches
it through the cell, and arrives at a statement of a different
character: an inequality between the pair angles of the configuration
which, if it holds for every contact configuration of $23$ directions,
proves \cref{thm:m23} outright, which the one configuration we can write
down satisfies with room to spare, and which \cref{sec:certificate}
then establishes. The tools are the ones already in
hand, the covering bound of \cref{thm:covering-bound} and the exact
overlaps of \cref{prop:second-order}; what is new is the observation of
where the integration ought to stop, and what stopping there buys.

\begin{lem}[Contact directions in a cap]\label{lem:cap-count}
If $k$ contact directions lie in a closed cap of angular radius $r$ on
$S^3$, then $\sin^2r\ge\dfrac{k-1}{2k}$. In particular two contact
directions need a cap of radius at least $30^\circ$, three need
$r_*=\arcsin(1/\sqrt3)=35.2644\ldots^\circ$, and four need
$r_4=\arcsin\sqrt{3/8}=37.7612\ldots^\circ$.
\end{lem}

\begin{proof}
A point $\theta$ of the cap about $c$ satisfies
$|\theta-(\cos r)c|^2=1-2\cos r\,\langle\theta,c\rangle+\cos^2r\le\sin^2r$,
so the cap lies in the Euclidean ball of radius $\sin r$ about
$(\cos r)c$. For points $x_1,\dots,x_k$ of that ball with centroid $q$,
\[
\begin{aligned}
  \sum_{i<j}|x_i-x_j|^2 &= k\sum_i|x_i-q|^2\\
  &\le k\sum_i|x_i-(\cos r)c|^2\ \le\ k^2\sin^2r ,
\end{aligned}
\]
the first step being the standard identity and the second the fact
that the centroid minimises the sum of squared distances. Contact
directions have $|x_i-x_j|^2=2-2\langle x_i,x_j\rangle\ge1$, so
$\binom k2\le k^2\sin^2r$.
\end{proof}

The case $k=3$ is drawn in \cref{fig:structcd}(b), and the
lens in which the caps of two contact directions overlap, whose measure
is the $\Lambda(r,\gamma)$ of \eqref{eq:three-point-survival} below, in
\cref{fig:pfone}(a). The cases $k=3$ and $k=4$ are \cref{lem:no-triples} and
\cref{prop:covering-radius} again, since
$\arccos\sqrt{2/3}=\arcsin\sqrt{1/3}$ and
$\arccos\sqrt{5/8}=\arcsin\sqrt{3/8}$; the lemma puts them on one
footing and is the reason the two radii $r_*$ and $r_4$ govern what
follows.

\begin{prop}[The cell inside a ball]\label{prop:truncated}
Let $W$ be a contact configuration of $m$ directions with bounded cell,
let $U(r)$ be the measure of the set of points of $S^3$ farther than $r$
from every direction of $W$, and let $R\ge1$. Then
\begin{equation}\label{eq:truncated}
  \vol\bigl(V_c(W)\cap B(R)\bigr)
  =\frac14\Bigl[2\pi^2+\int_0^{\arccos(1/R)}U(r)\,d(\sec^4r)\Bigr]
  \;\le\;\vol\bigl(V_c(W)\bigr),
\end{equation}
and for $r<r_4$
\begin{equation}\label{eq:three-point-survival}
  U(r)=2\pi^2-mC(r)+\sum_{i<j}\Lambda(r,\gamma_{ij})
       -\sum_{i<j<k}\Lambda_3(r;\gamma_{ij},\gamma_{ik},\gamma_{jk}),
\end{equation}
where $\Lambda_3$ is the measure of the intersection of three caps of
radius $r$ about directions at the given pairwise angles, and the last
sum is absent for $r<r_*$. Consequently, with
\[
\begin{gathered}
  A_4=\frac14\Bigl[2\pi^2+\int_0^{r_4}\bigl(2\pi^2-23\,C(r)\bigr)\,d(\sec^4r)\Bigr]
  =7.647558\ldots,\\
  \omega_4(\gamma)=\frac14\int_0^{r_4}\Lambda(r,\gamma)\,d(\sec^4r),
\end{gathered}
\]
and $\omega_3$ the same integral of $\Lambda_3$, every contact
configuration of $23$ directions with bounded cell satisfies
\begin{equation}\label{eq:three-point-bound}
  \vol(V_c)\;\ge\;\vol\bigl(V_c\cap B(\sqrt{8/5})\bigr)
  \;=\;A_4+\sum_{i<j}\omega_4(\gamma_{ij})
  -\sum_{i<j<k}\omega_3(\gamma_{ij},\gamma_{ik},\gamma_{jk}) .
\end{equation}
\end{prop}

\begin{proof}
The radial function of $V_c(W)$ in the direction $\theta$ is
$\sec\delta(\theta)$ by \cref{lem:radial-form}, so that of
$V_c(W)\cap B(R)$ is $\min(\sec\delta,R)$ and the volume is
$\tfrac14\int_{S^3}\min(\sec\delta,R)^4$. Writing $\rho=\arccos(1/R)$,
\[
  \min(\sec\delta,R)^4=\sec^4\min(\delta,\rho)
  =1+\int_0^{\rho}\mathbf 1[\delta(\theta)>r]\,d(\sec^4r),
\]
and integrating over $\theta$ gives \eqref{eq:truncated}, the inequality
because $V_c(W)\cap B(R)\subseteq V_c(W)$. For
\eqref{eq:three-point-survival}, $U(r)$ is $2\pi^2$ less the measure of
the union of the closed caps of radius $r$ about the directions of $W$,
and by \cref{lem:cap-count} no point lies in four of those caps when
$r<r_4$, nor in three when $r<r_*$; the indicator of a union of sets no
four of which meet is the sum of the indicators of the sets, less that
of the pairwise intersections, plus that of the triple intersections,
exactly. Integrating \eqref{eq:three-point-survival} against
$d(\sec^4r)$ over $[0,r_4]$ gives \eqref{eq:three-point-bound}, and
$\sec r_4=\sqrt{8/5}$.
\end{proof}

\begin{figure}[tb]
\centering
\includegraphics[width=0.86\textwidth]{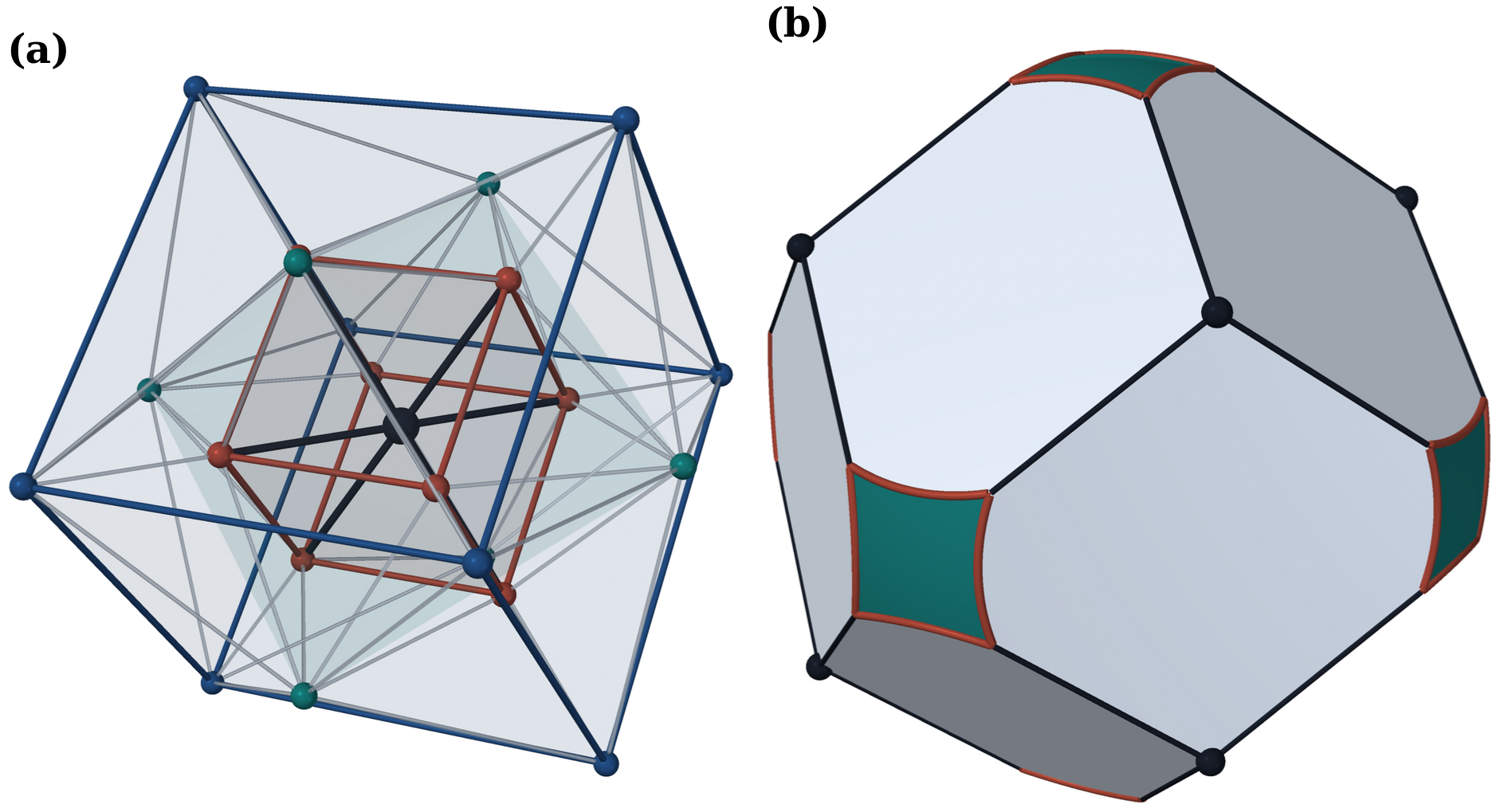}
\caption{(a) The contact graph of the root system from one root (dark, at the centre), the other twenty-three at radius proportional to angle: a cube at $60^\circ$ (clay), an octahedron at $90^\circ$ (teal), a cube at $120^\circ$ (blue), and the $88$ tight pairs among them; the eight at the centre are drawn dark, and they sum to four times it. (b) The truncated cell in dimension three: the Voronoi cell of the face-centred cubic lattice at inradius $1$, a rhombic dodecahedron, cut by the ball of radius $(3/2)^{1/2}$. The ball passes through the eight three-valent vertices and replaces the six four-valent corners by spherical squares (teal).}
\label{fig:pftwo}
\end{figure}

The three-dimensional analogue of all this can be drawn, and
\cref{fig:pfone}(b) and \cref{fig:pftwo}(b) draw it: the Voronoi cell of the face-centred
cubic lattice, whose contact directions are the twelve vertices of a
cuboctahedron, cut by the ball of radius $\sqrt{3/2}$, and the caps of
radius $r_*$ about those twelve directions, which overlap in pairs and
meet three at a time only at the eight triangle circumcentres; the
uncovered measure there is exactly the two-term inclusion and
exclusion, as it is on $S^3$ below $r_*$.

\begin{figure}[tb]
\centering
\includegraphics[width=0.86\textwidth]{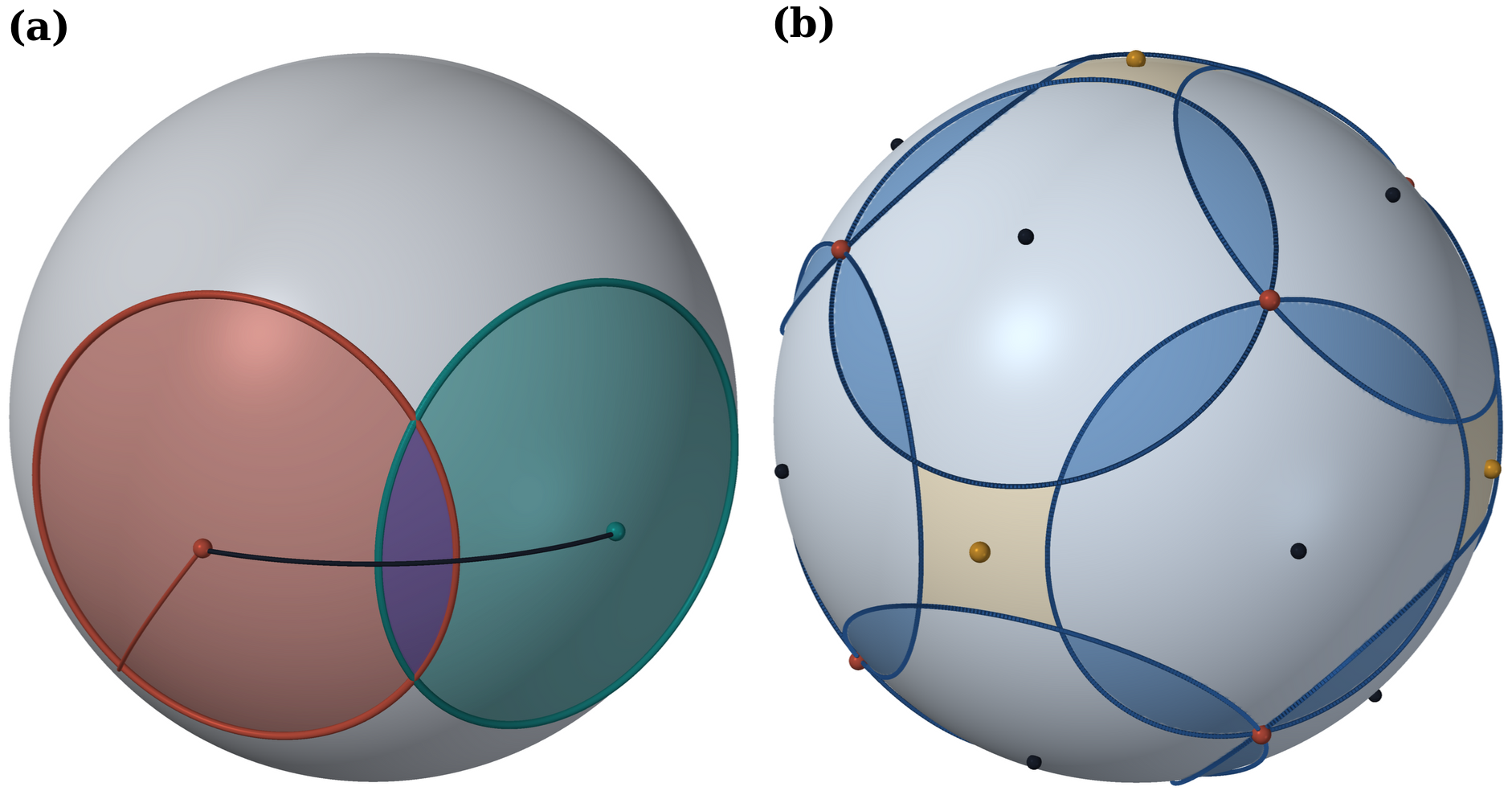}
\caption{(a) Two contact directions at angle $\gamma=60^\circ$ and their caps of radius $r=r_*$. The lens in which the caps overlap (violet) has measure $\Lambda(r,\gamma)$, and $\omega(\gamma)$ integrates it against $d(\sec^4r)$; the dark arc is $\gamma$ and the clay arc is $r$. (b) The caps of radius $r_*$ about the twelve contact directions of the same lattice, the sphere shaded by how many caps cover it: one (pale blue), two in the lenses, three at the eight triangle circumcentres (clay), and none in the six square holes (ochre).}
\label{fig:pfone}
\end{figure}

At the deletion of a root the three-point bound
\eqref{eq:three-point-bound} is $8.140848\ldots$: the $88$ tight pairs
contribute $88\,\omega_4(60^\circ)=88\times0.00565947\ldots=0.498033\ldots$,
the $84$ triangles of tight pairs subtract
$84\,\omega_3(60^\circ,60^\circ,60^\circ)=84\times5.65\times10^{-5}=0.00474\ldots$,
and no other pair or triple contributes, since $\Lambda(r,\gamma)=0$
for $\gamma\ge2r_4=75.52^\circ$. A direct quadrature of
$\min(\sec\delta,\sqrt{8/5})^4$ over $S^3$ returns $8.14068$, and the
same quadrature at the full root system, with $24$ directions, returns
$7.96841$: there the truncation discards the deep holes at $45^\circ$
and the bound falls below the true value $8$, which says that at
$24$ contacts the three-point information below $r_4$ is not enough,
as it need not be. At $23$ contacts it is enough at the one
configuration known, with $0.14$ to spare, and the pair-only bound
\eqref{eq:second-order} at $R=r_*$ is enough there too, with $0.034$.
That is the content of the following statement, which we regard as the
right form of the remaining case.

\begin{thm}[The remaining case as an inequality between pair angles]\label{thm:strict-reduction}
Suppose that every contact configuration of $23$ directions with
bounded cell satisfies
\begin{equation}\label{eq:strict-pairs}
  \sum_{i<j}\omega(\gamma_{ij})\;\ge\;8-A_*\;=\;0.092855570\ldots ,
\end{equation}
where
\[
  A_*=\frac14\Bigl[2\pi^2+\int_0^{r_*}\bigl(2\pi^2-23\,C(r)\bigr)\,d(\sec^4r)\Bigr]
  =7.907144430\ldots
\]
and $\omega$ is the weight of \cref{prop:pair-budget}; equivalently, that
the cell of every such configuration has volume at least $8$ inside the
ball of radius $\sqrt{3/2}$. Then every such configuration has
$\vol(V_c)\ge8$, with strict inequality wherever \eqref{eq:strict-pairs}
is strict, and \cref{thm:m23} follows. The same conclusion follows if
every such configuration satisfies the three-point inequality
$\sum_{i<j}\omega_4(\gamma_{ij})-\sum_{i<j<k}\omega_3\ge8-A_4=0.352442\ldots$.
\end{thm}

\begin{proof}
By \cref{prop:second-order} at $R=r_*$, or by \eqref{eq:three-point-bound},
the hypothesis gives $\vol(V_c)\ge\vol(V_c\cap B(\sqrt{3/2}))\ge8$ for
every contact configuration of $23$ directions with bounded cell,
strictly when the pair sum exceeds $8-A_*$. A configuration whose cell
has circumradius below $2$ has bounded cell, so \cref{thm:m23} is the
bound just proved once \eqref{eq:strict-pairs} is known to hold
strictly, which is what \cref{thm:certificate} supplies.
\end{proof}

\begin{remark}[What kind of statement this is]\label{rem:strict-kind}
Three things distinguish \eqref{eq:strict-pairs} from the classification
of \cref{thm:m23-from-codes}, and they are the reason for stating it.

First, it is strict where the classification is not. The classification
asks for a property of codes of minimal angle exactly $60^\circ$, and
the deletion of a root sits on the boundary of that property, with
$g=\tfrac12$; \cref{rem:m23-uniqueness} explains why no relaxation can
see a boundary. Inequality \eqref{eq:strict-pairs} holds at the
deletion with $0.034$ to spare, and at every configuration reached in
\cref{sec:slack-continuation} with more: the least value of
$\vol(V_c\cap B(\sqrt{3/2}))$ found over the configurations of
\cref{tab:slack-continuation} and over a direct minimisation of that
quantity on $K(\delta)$ is $8.028$, at $\delta=0.05$, and the least
value of $\vol(V_c\cap B(\sqrt{8/5}))$ is $8.129$
(\texttt{truncated\_volume.py}). A statement with slack is one that a
relaxation can certify, and \cref{sec:certificate} does.

Second, it needs no classification and no saturation. It is the
$m=23$ case of the $24$-cell conjecture itself, in the form
$\vol(V_c)\ge8$ for every configuration of $23$ contacts, saturated or
not; what \cref{sec:m24} contributes is only that this case is the
last one, the others being settled by \cref{cor:m22} and by
\cref{thm:m24}.

Third, it is the kind of statement a semidefinite relaxation can
certify, and the two-point relaxation already delivers four fifths of
it. \Cref{prop:two-point-barrier} minimises the left-hand side of
\eqref{eq:strict-pairs} over every pair-angle distribution that
satisfies the Bonferroni inequality at every radius and is positive
definite in every degree, and finds $0.07380\ldots$ against
$0.09285\ldots$. The measure that attains the minimum gives each
direction nine neighbours at $62.3^\circ$ and nearly eight at
$103^\circ$, a local picture that is consistent at any one direction
and cannot be consistent at all twenty-three at once; telling the two
apart is exactly what three-point information is for. The
semidefinite bounds of Bachoc and Vallentin \cite{BV08}, which carry
the distribution of triples of directions and which de Laat,
Leijenhorst and de Muinck Keizer \cite{LLM24} sharpened, at the second
level, to an exact classification at $24$ points
(\cref{thm:twentyfour-points}), are the instrument built for that, and
the objective \eqref{eq:strict-pairs} is linear in the pair
distribution, as such bounds require. A relaxation can fall short of a
true inequality; what makes this one a good candidate is that the
target it has to reach is fixed, $0.092855\ldots$, that the pair
angles alone reach $0.0738$, and that the gap is not a boundary
phenomenon but a quantity with room on either side of it. The next
subsection carries the computation out, and the relaxation does reach
the target.
\end{remark}

\subsection{The inequality between pair angles holds}\label{sec:certificate}

The inequality \eqref{eq:strict-pairs} is a lower bound for a linear
functional of the pair distribution of a spherical code of fixed size
and fixed minimal angle, and the method built for statements of that
shape is the semidefinite programming bound of Bachoc and Vallentin
\cite{BV08}, which reads the distribution of triples of directions. We
set it up in the form we need, prove the one lemma that turns a
feasible point of the programme into an inequality for every
configuration, and then describe the certificate and its verification,
which is carried out in exact rational arithmetic where the objects are
rational and in interval arithmetic where they are not. Throughout,
$N=23$, $u_{ij}=\langle w_i,w_j\rangle=\cos\gamma_{ij}$, and $\omega$ is
the weight of \cref{prop:pair-budget}, regarded as the function of the
inner product $u=\cos\gamma$: increasing on $[\tfrac13,\tfrac12]$ from
$0$ to $0.00144540\ldots$ and identically $0$ below $\tfrac13$, since
$\Lambda(r,\gamma)=0$ for $\gamma\ge2r_*$.

\subsubsection*{Positive definite kernels on the three-sphere}

For $k\ge0$ let $G_k$ be the Gegenbauer polynomial of $S^3$ normalised
by $G_k(1)=1$, that is $G_k(\cos\theta)=\sin((k+1)\theta)/((k+1)\sin\theta)$,
the Chebyshev polynomial of the second kind divided by $k+1$. By
Schoenberg's theorem \cite{Sch42} every finite $W\subset S^3$ satisfies
$\sum_{x,y\in W}G_k(\langle x,y\rangle)\ge0$, which is the linear
programming bound of Delsarte, Goethals and Seidel \cite{DGS77} in its
primal form. The three-point matrices are built from the Legendre
polynomials $P_k$, which play the same role on $S^2$. Fix a degree $d$
and, for $0\le k\le d$, let $Y_k(u,v,t)$ be the square matrix of order
$d-k+1$ with entries
\begin{equation}\label{eq:Yk}
\begin{aligned}
  Y_k(u,v,t)_{ij}&=T_i(u)\,T_j(v)\,
  \bigl((1-u^2)(1-v^2)\bigr)^{k/2}\\
  &\quad\times P_k\Bigl(\frac{t-uv}{\sqrt{(1-u^2)(1-v^2)}}\Bigr),
\end{aligned}
\end{equation}
for $0\le i,j\le d-k$,
where $T_i$ is the Chebyshev polynomial of the first kind; since $P_k$
has the parity of $k$, each entry is a polynomial in $u,v,t$. Bachoc
and Vallentin write $u^iv^j$ in place of $T_i(u)T_j(v)$; the two
families differ by a fixed congruence of the matrices and are
interchangeable in everything below, and the Chebyshev basis is the
better conditioned one on $[-1,\tfrac12]$. Let $S_k$ be the average of
$Y_k$ over the six permutations of $(u,v,t)$.

\begin{lem}[Positivity of the three-point matrices]\label{lem:three-point-psd}
For every finite $W\subset S^3$ and every $k$, the matrix
$\sum_{(x,y,z)\in W^3}S_k(\langle x,y\rangle,\langle x,z\rangle,\langle y,z\rangle)$
is positive semidefinite.
\end{lem}

\begin{proof}
Fix $e\in S^3$ and write $x^\perp=x-\langle e,x\rangle e$ for $x\in S^3$.
With $u=\langle e,x\rangle$, $v=\langle e,y\rangle$, $t=\langle x,y\rangle$
one has $|x^\perp|^2=1-u^2$ and $\langle x^\perp,y^\perp\rangle=t-uv$, so
the scalar factor of \eqref{eq:Yk} is
$Q_k(x,y)=|x^\perp|^k|y^\perp|^kP_k(\langle x^\perp,y^\perp\rangle/(|x^\perp||y^\perp|))$,
read as $0$ for $k\ge1$ and as $1$ for $k=0$ when one of the projections
vanishes. The kernel $P_k(\langle\xi,\eta\rangle)$ is positive definite
on the unit sphere of $e^\perp$, a copy of $S^2$, by \cite{Sch42}, and
multiplying a positive definite kernel by $a(x)a(y)$ keeps it positive
definite; so $Q_k$ is positive definite on $S^3$, and for any functions
$a_0,\dots,a_{d-k}$ the matrix with entries
$\sum_{x,y\in W}a_i(x)a_j(y)Q_k(x,y)$ is positive semidefinite, because
$c^{\mathsf T}Mc=\sum_{x,y}b(x)Q_k(x,y)b(y)$ with $b=\sum_ic_ia_i$. Taking
$a_i(x)=T_i(\langle e,x\rangle)$ gives
$\sum_{x,y\in W}Y_k(\langle e,x\rangle,\langle e,y\rangle,\langle x,y\rangle)\succeq0$
for every $e$; summing over $e\in W$ gives the same for the sum over
$W^3$ of $Y_k$, and that sum is unchanged by permuting the three
arguments, so it equals the sum of $S_k$.
\end{proof}

\subsubsection*{From a feasible point to an inequality}

\begin{lem}[The certificate lemma]\label{lem:certificate}
Let $f=\sum_{k=0}^df_kG_k$ with $f_k\ge0$ for $k\ge1$, let
$F(u,v,t)=\sum_{k=0}^d\langle F_k,S_k(u,v,t)\rangle$ with each $F_k$
positive semidefinite of order $d-k+1$, and suppose that
\begin{equation}\label{eq:C}
  \omega(u)+\omega(v)+\omega(t)-f(u)-f(v)-f(t)-F(u,v,t)
  -\frac{F(1,u,u)+F(1,v,v)+F(1,t,t)}{N-2}\;\ge\;0
\end{equation}
for every admissible triple, that is every $(u,v,t)\in[-1,\tfrac12]^3$
with $1+2uvt-u^2-v^2-t^2\ge0$. Then every contact configuration of $N=23$
directions satisfies
\begin{equation}\label{eq:certificate-bound}
  \sum_{i<j}\omega(u_{ij})\;\ge\;B:=\frac N2\bigl(Nf_0-f(1)\bigr)-\frac{N\,F(1,1,1)}{6(N-2)} .
\end{equation}
\end{lem}

\begin{proof}
The inner products of three distinct directions of a contact
configuration form an admissible triple, since their Gram matrix is
positive semidefinite with unit diagonal and its determinant is
$1+2uvt-u^2-v^2-t^2$. Sum \eqref{eq:C} over the $\binom N3$ triples
$i<j<k$. Each pair lies in $N-2$ triples, so the $\omega$ terms give
$(N-2)\sum_{i<j}\omega(u_{ij})$ and the $f$ terms give
$(N-2)\sum_{i<j}f(u_{ij})=\tfrac{N-2}2\bigl(\sum_{i,j}f(u_{ij})-Nf(1)\bigr)
\ge\tfrac{N-2}2\bigl(N^2f_0-Nf(1)\bigr)$ by Schoenberg's theorem, since
$G_0=1$. The triple sum of $F$ and the pair sum of $F(1,u,u)$ are read
off from the sum of $F$ over all $N^3$ ordered triples of $W$, which is
nonnegative by \cref{lem:three-point-psd}: the $N$ triples with all
three entries equal contribute $NF(1,1,1)$; the $3N(N-1)$ with exactly
two equal contribute $F(1,u_{ij},u_{ij})$ each, by the symmetry of $F$,
that is $6\sum_{i<j}F(1,u_{ij},u_{ij})$; and the $N(N-1)(N-2)$ with all
entries distinct contribute $6\sum_{i<j<k}F(u_{ij},u_{ik},u_{jk})$.
Hence
$\sum_{i<j<k}F\ge-\sum_{i<j}F(1,u_{ij},u_{ij})-\tfrac N6F(1,1,1)$, while
the last term of \eqref{eq:C} summed over triples is exactly
$\sum_{i<j}F(1,u_{ij},u_{ij})$. Putting the three parts together,
$(N-2)\sum_{i<j}\omega(u_{ij})\ge\tfrac{N-2}2(N^2f_0-Nf(1))-\tfrac N6F(1,1,1)$,
which is \eqref{eq:certificate-bound}.
\end{proof}

The lemma uses the minimal angle only through the domain of
\eqref{eq:C}, and the number $23$ only through $N$. Its content is that
the pair sum on the left of \eqref{eq:strict-pairs} is bounded below by
the value of a semidefinite programme in the variables $f_k$ and $F_k$,
and that a feasible point of that programme with value above
$8-A_*$ is a proof of \cref{thm:strict-reduction}'s hypothesis. Such a
point exists at degree $8$.

\begin{thm}[The pair inequality]\label{thm:certificate}
There are numbers $f_0,\dots,f_8$ and symmetric matrices
$F_0,\dots,F_8$, of orders $9,8,\dots,1$, listed in the supplementary
material, with $f_k>0$ for $k\ge1$, every $F_k$ positive definite,
satisfying \eqref{eq:C} on the admissible domain, and with
\[
  B=0.0929000002\ldots\;>\;8-A_*=0.0928555702\ldots
\]
Consequently every contact configuration of $23$ directions satisfies
\eqref{eq:strict-pairs} with $0.0000444$ to spare, every such
configuration with bounded cell has
$\vol(V_c)\ge\vol(V_c\cap B(\sqrt{3/2}))=A_*+\sum_{i<j}\omega(\gamma_{ij})\ge A_*+B>8$,
and \cref{thm:m23} holds.
\end{thm}

The proof is the verification of the three properties of the
certificate, and we describe it in the detail that a reader who wants
to repeat it would need; the scripts named are in \cref{app:code-index}.

\begin{figure}[tb]
\centering
\includegraphics[width=\textwidth]{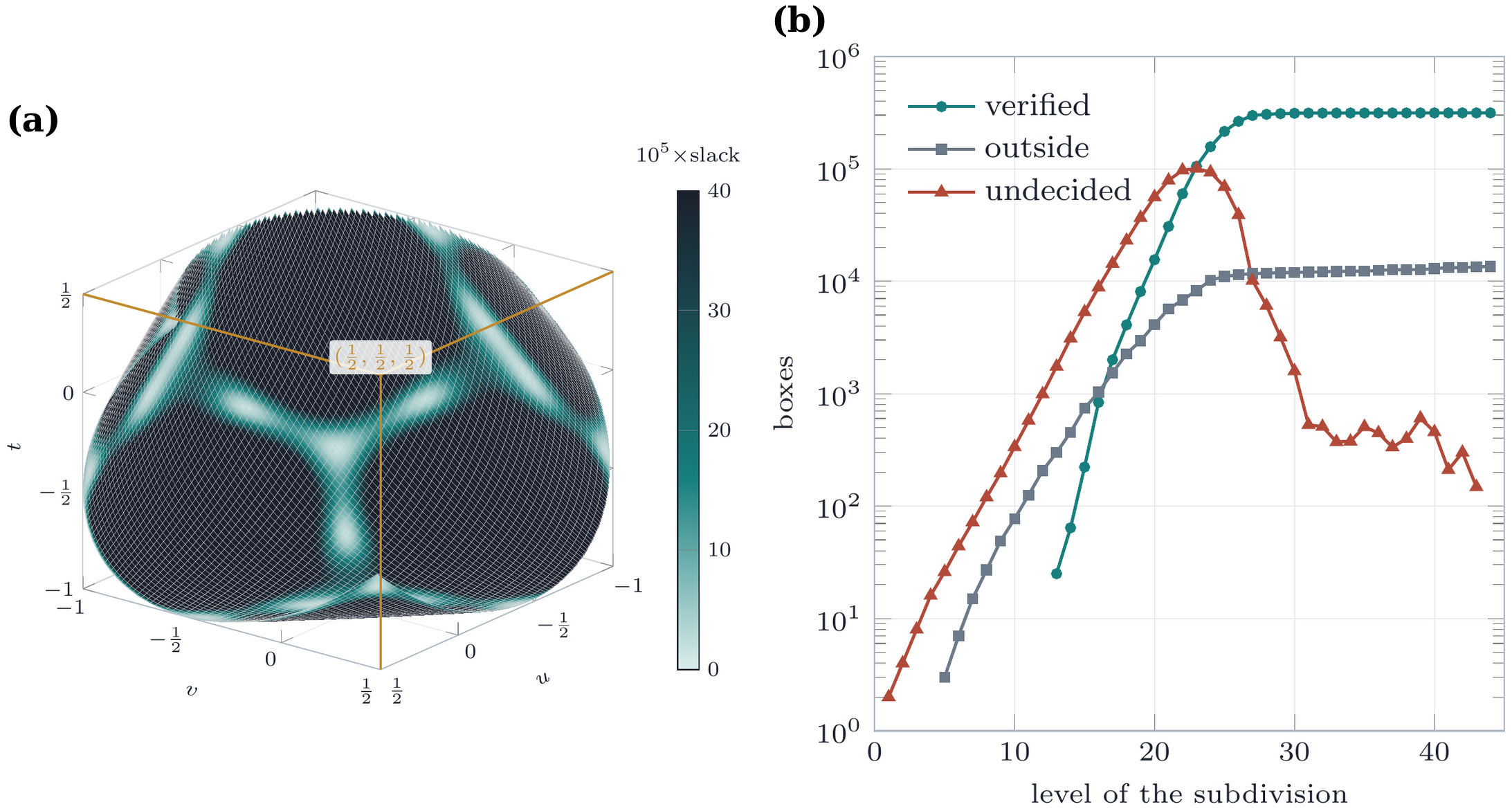}
\caption{(a) The domain of the certificate condition: the admissible region
of \cref{lem:certificate} in three dimensions, where the triples
$(u,v,t)$ of pairwise inner products of three contact directions fill
the part of the cube $[-1,\tfrac12]^3$ enclosed by the surface of
coplanar triples. The surface is coloured by the slack of \eqref{eq:C},
which is least along the curves where two of the three directions are
nearly antipodal. (b) The branch and bound of the proof of \cref{thm:certificate},
level by level: the boxes live, the boxes discarded on the enclosure,
and the width reached.}
\label{fig:dom3d}
\end{figure}

\begin{proof}
\emph{The certificate.} The numbers are double-precision floating-point
numbers, which are dyadic rationals, and they are taken as the exact
rationals they denote; nothing is rounded. All computations below in
which only these rationals and the integers enter are done in exact
rational arithmetic.

\emph{Positivity.} The nine matrices are symmetric and positive
definite: an exact $LDL^{\mathsf T}$ elimination over $\mathbb Q$ returns
positive pivots throughout, and $f_1,\dots,f_8$ are positive. The
bound $B$ is computed exactly from \eqref{eq:certificate-bound}; here
$F(1,1,1)$ is the sum of the entries of $F_0$, because every entry of
$Y_k(1,1,1)$ vanishes for $k\ge1$ and $T_i(1)=1$. These two steps, with
the certificate written into the source as exact rationals, are also
carried out by the Lean kernel (\texttt{lean/D4Certificate.lean}).

\emph{The constants.} $A_*$ is elementary. With $S(r)=\sec^4r$ one has
$\int r\,dS=r\sec^4r-\tan r-\tfrac13\tan^3r$ and
$\int\sin2r\,dS=\tfrac83\tan^3r$, and $\sec^2r_*=\tfrac32$,
$\tan r_*=1/\sqrt2$, so that
\begin{equation}\label{eq:A-star-closed}
  A_*=\frac{9\pi^2}8-\frac{207\pi}8\,r_*+\frac{253\pi}{12\sqrt2},
  \qquad r_*=\arctan\frac1{\sqrt2},
\end{equation}
and $8-A_*$ lies in $[0.092855570293,0.092855570294]$ by interval
arithmetic at forty digits. The weight $\omega$ is elementary too. The
lens measure of \cref{prop:second-order} integrates to
\[
  \Lambda(r,\gamma)=4\pi\Bigl[\frac{r-\gamma/2}2-\frac{\sin2r-\sin\gamma}4
  -\frac{\tan(\gamma/2)}2\bigl(\sin^2r-\sin^2(\gamma/2)\bigr)\Bigr],
\]
and $\omega(\gamma)=\int_{\gamma/2}^{r_*}\Lambda(r,\gamma)\tan r\sec^4r\,dr$
is a combination of the four primitives
$\int\tan r\sec^4r\,dr=\tfrac14\sec^4r$,
$\int r\tan r\sec^4r\,dr=\tfrac14r\sec^4r-\tfrac14(\tan r+\tfrac13\tan^3r)$,
$\int\sin2r\tan r\sec^4r\,dr=\tfrac23\tan^3r$ and
$\int\sin^2r\tan r\sec^4r\,dr=\tfrac14\tan^4r$, with coefficients
depending on $\gamma$ alone. The closed form agrees with the quadrature
of \cref{prop:pair-budget} to fifteen digits at every angle tested, and
it is differentiated symbolically to give $d\omega/du$ and $d^2\omega/du^2$
as functions of $\gamma$, with $du=-\sin\gamma\,d\gamma$. On
$(60^\circ,2r_*)$ the function $\omega$ is analytic, and it is decreasing
in $\gamma$, since $\Lambda(r,\gamma)$ is and the range of integration
shrinks; so $\omega$ is increasing in $u$ on $[\tfrac13,\tfrac12]$.

\emph{Tables.} At $12001$ angles $\gamma_i$ equally spaced from $2r_*$ down
to $60^\circ$, interval arithmetic gives an upper bound $u_i$ of
$\cos\gamma_i$, a lower bound of $\omega(\gamma_i)$ and an enclosure of
$\omega'(u)$ at $u=\cos\gamma_i$; consecutive $u_i$ differ by at most
$1.44\times10^{-5}$, and $|\omega''(u)|\le0.289712$ on the whole range,
by interval evaluation of the second derivative on $3000$ subintervals.
For $u\ge u_i$ the monotonicity gives $\omega(u)\ge\omega(\gamma_i)$, and
the mean value theorem gives $|\omega'(u)-\omega'(\cos\gamma_i)|\le0.289712\,(u-\cos\gamma_i)$.

\emph{The inequality \eqref{eq:C}.} The polynomial
$P(u,v,t)=f(u)+f(v)+f(t)+F(u,v,t)+(F(1,u,u)+F(1,v,v)+F(1,t,t))/21$ is
expanded exactly, from the recurrences for $T_i$, $P_k$ and $G_k$, into
$449$ monomials of total degree at most $16$ with rational coefficients;
it is symmetric in $u,v,t$, so \eqref{eq:C} need only be verified on
the ordered part $-1\le u\le v\le t\le\tfrac12$ of the domain, and it
reads $Q:=1000(\omega(u)+\omega(v)+\omega(t))-P\ge0$ after the scaling
by $1000$ under which the certificate was computed. The verification is
a branch and bound over boxes, in floating-point interval arithmetic
with every operation rounded outward. On a box with centre $c$ and
half-widths $\rho_1,\rho_2,\rho_3$ whose $u$, $v$ and $t$ ranges lie
inside $(\tfrac13,\tfrac12]$, Taylor's theorem gives
\[
  Q\;\ge\;Q(c)-\sum_v\bigl|\partial_vQ(c)\bigr|\rho_v
  -\frac12\sum_{v,w}\Bigl(\sup_{\rm box}\bigl|\partial_v\partial_wQ\bigr|\Bigr)\rho_v\rho_w ,
\]
where $Q(c)$ and $\partial_vQ(c)$ are bounded through the tables and the
exact polynomial evaluated at the point $c$, and the second derivatives
through the bound on $\omega''$ and the interval evaluation of the
second derivatives of $P$ over the box. A variable whose range is not
inside $(\tfrac13,\tfrac12]$ has its $\omega$ term replaced by $0$, which
is a lower bound, with derivatives $0$. A box on which the right-hand
side is nonnegative is verified; a box on which $1+2uvt-u^2-v^2-t^2$ is
negative throughout is outside the domain and is discarded; any other
box is bisected along the direction that contributes most to the
Taylor remainder. Starting from $[-1,\tfrac12]^3$, the process ends
after $44$ levels with $313\,780$ boxes verified and $13\,517$
discarded, in under ten minutes (\cref{fig:dom3d}(a)), and no box reaching the minimum width
$10^{-6}$ undecided. The same branch and bound is then carried out a
second time, independently and in exact arithmetic, by a Lean program
(\texttt{lean/certificate}): it expands $P$ from the certificate on its
own, multiplies it by $315$ so that every coefficient becomes an integer
over a power of two, and runs the subdivision with the same Taylor form
in dyadic arithmetic, integers and exponents only, with no rounding at
any point, taking from outside Lean nothing but the tables of bounds for
$\omega$, $\omega'$ and $\omega''$; it processes $421\,881$ boxes,
fewer than the Python run because its nested evaluation of the
polynomials encloses more tightly than a sum of monomials, subdivides
down to the width $2^{-20}$ where it has to, and returns true, which
Lean records as a theorem under \texttt{native\_decide}. Hence \eqref{eq:C} holds on the
admissible domain,
\cref{lem:certificate} gives \eqref{eq:strict-pairs} for every contact
configuration of $23$ directions, and \cref{thm:strict-reduction}
gives the rest.
\end{proof}

\begin{figure}[tb]
\centering
\includegraphics[width=\textwidth]{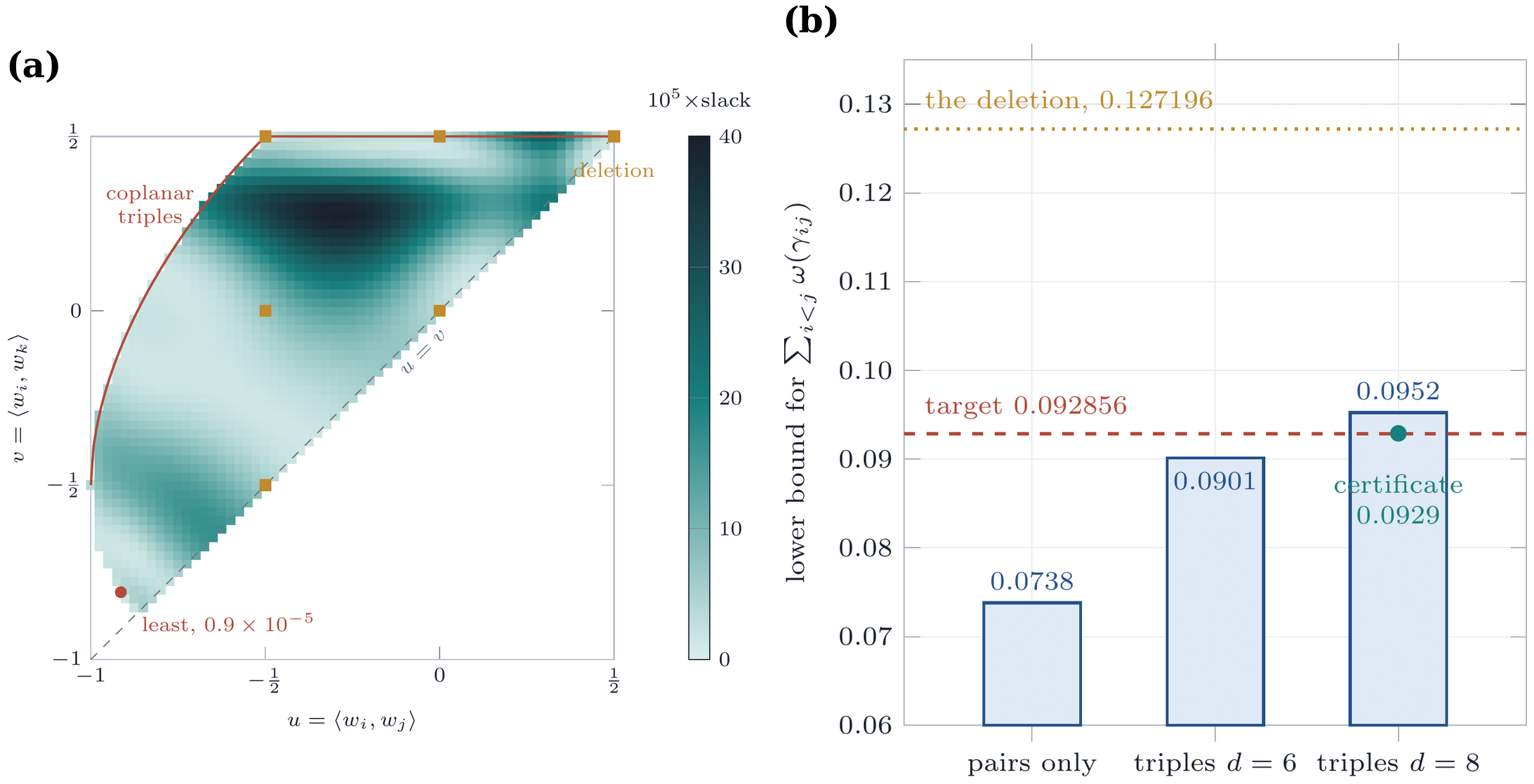}
\caption{(a) The certificate of \cref{thm:certificate}: the slack of the
condition \eqref{eq:C}, the amount by which its left-hand side exceeds
$0$, on the slice of the admissible domain in which one pair of the
triple is at $60^\circ$. It is positive throughout, least near a
coplanar triple with one nearly antipodal pair, and largest, about
$4\times10^{-4}$, in the interior. The six squares are the triples that
occur in a deletion of a root. (b) The values of the relaxations against the target $8-A_*$: the
pair angles alone (\cref{prop:two-point-barrier}), the three-point
relaxation at degrees $6$ and $8$, and the certificate, at which the
bound is fixed at $0.0929$ and the least slack is maximised.}
\label{fig:certpair}
\end{figure}

\begin{remark}[How the certificate was found, and how tight it is]\label{rem:certificate-found}
The programme of \cref{lem:certificate}, with \eqref{eq:C} imposed on a
sample of admissible triples and $B$ maximised
(\cref{fig:certpair}), is solved by an
interior-point method (\texttt{three\_point\_sdp.py}); the sample is
refined four times by adding the triples at which the current solution
violates \eqref{eq:C} most, on a check over about $1.15$ million further
triples, and the value reported is reduced by $\binom{23}3/21$ times
the largest violation that remains, so that it is a bound over the
whole check. At degree $6$ the value is $0.09011$, that is $97.0$ per
cent of $8-A_*$; at degree $8$ it is $0.09523$, or $102.6$ per cent. The
pair angles alone, in \cref{prop:two-point-barrier}, gave $79.5$ per
cent; the triples supply the rest. The certificate of
\cref{thm:certificate} is then obtained by fixing $B$ at $0.0929$ and
maximising instead the least slack of \eqref{eq:C} over the sample, with
$F_k-10^{-7}I$ and $f_k-10^{-7}$ constrained to be positive
semidefinite so that the solution survives being read as exact
rationals: the least slack is $1.35\times10^{-5}$ on the final sample
of $21\,648$ triples and $1.30\times10^{-5}$ on the check, against
$\omega(60^\circ)=0.00145$. Two features of the solution are worth
recording. The two-point part is almost inert, $f_0=0.000927$ and
$f_1,\dots,f_8$ below $3\times10^{-8}$, so it is the matrices $F_k$ that
carry the argument. And the inequality is far from tight at the
deletion of a root, where the pair sum is $0.1272$ against
$B=0.0929$; the relaxation had room, which is the situation
\cref{rem:strict-kind} asked for, and it is why a relaxation could
succeed here where, by \cref{rem:m23-uniqueness}, none could have
classified the codes. Finally, the checker itself was tested on a
certificate it should reject: adding $10^{-5}$ to $f_0$, which raises
$P$ by $3\times10^{-5}$ and so exceeds the slack, produces a certificate
for which the same programme exhibits the admissible triple
$(-\tfrac1{16},-\tfrac1{16},-\tfrac1{16})$ at which \eqref{eq:C} fails,
in under two minutes.
\end{remark}

\begin{remark}[What is verified where]\label{rem:certificate-trust}
The proof of \cref{thm:certificate} rests on three kinds of
computation, and the reader should know which is which. The exact
rational steps, positivity of the matrices and the value of $B$, are
verified by the Lean kernel from the same numbers, and the expansion of
$P$ is an exact rational computation in Python whose result is checked
against an independent floating-point evaluation of the same
expressions at random triples, agreeing to $10^{-13}$. The constants
$A_*$, $\omega$, $\omega'$ and $\omega''$ are evaluated from closed forms,
\eqref{eq:A-star-closed} among them, in the interval arithmetic of
\texttt{mpmath}. The branch and bound is run twice: in Python, in an
interval arithmetic on IEEE doubles with outward rounding by one unit
in the last place after every operation, and in Lean, in exact dyadic
arithmetic, where the expansion of $P$, its derivatives, every
enclosure and every bisection are integer computations and the only
data taken from outside are the tables of bounds for $\omega$ and its
two derivatives, so that no floating-point operation enters the Lean
run at all. The Lean theorem states that the program returns true and
is settled by \texttt{native\_decide}, which runs the compiled program
and therefore trusts the compiler as well as the kernel; what remains
outside any proof assistant is the evaluation of the closed forms at
the table points and the mathematics that makes the program a proof,
the Taylor form and the monotonicity of $\omega$, which are written
out above. The Taylor form is applied only on boxes inside the open
interval where $\omega$ is analytic.
\end{remark}
\section{The deviation domain and its fundamental triangle}\label{sec:deviation-triangle}

The identities of \cref{sec:deviation-domain} hold for any number of
deviating directions. From here to \cref{sec:cap-theorem} we return to
the case of one, which is where the exact structure lives, and set up the
domain that case is played out on.

\subsection{The two remaining parameter spaces}\label{sec:two-remaining}

\Cref{thm:freevol-new} makes clear that the domain covered by
\cref{thm:local-uncond} enlarges in two structurally different
directions:

\begin{itemize}[leftmargin=1.6em]
\item \emph{More directions at once} ($m\ge2$): the parameter space is
now a choice of $m$ pairwise-compatible directions on $S^3$, together
with an assignment of each to its own $3$-dimensional tangent
deviation. This is the subject of \cref{sec:hessian-technique} onward.
\item \emph{One direction, but off the root-aligned family}
($m=1$, $\dev$ not root-aligned): the parameter space is the full unit
$2$-sphere of deviation directions orthogonal to the undeviated root,
here reduced (\cref{sec:fundamental-domain-new}) to a triangle by two
exact symmetries. This is the subject of
\cref{sec:arc1-breakpoints} onward.
\end{itemize}

\subsection{The fundamental domain of deviation directions}\label{sec:fundamental-domain-new}

Fix a root $\alpha_0\in\Rt$ and let $u_0=\alpha_0/\sqrt2$ be the
corresponding contact direction. The tangent space orthogonal to $u_0$
is $3$-dimensional, and a deviation direction is any unit vector
$\dev$ in it; the deviating cap direction is then
$u_1(\tilt,\dev)=\cos\tilt\,u_0+\sin\tilt\,\dev$.

\begin{lem}[Two exact symmetries]\label{lem:two-symmetries-new}
Two elements of $\Stab(u_0)\le W(D_4)$ act on the space of deviation
directions $\dev$: a symmetry $S$ of order $2$ and a symmetry $H$
(a signed permutation realising a Hadamard-type transformation) also of
order $2$, generating a group of order $4$ whose fundamental domain on
the unit circle of directions orthogonal to two further fixed root
directions $v_1,w_1$ (themselves canonically determined by $u_0$) is a
single triangle with vertices $v_1$, $w_1$, and $v_2=Sw_1$.
\end{lem}

\begin{figure}[tb]
\centering
\includegraphics[width=0.40\textwidth]{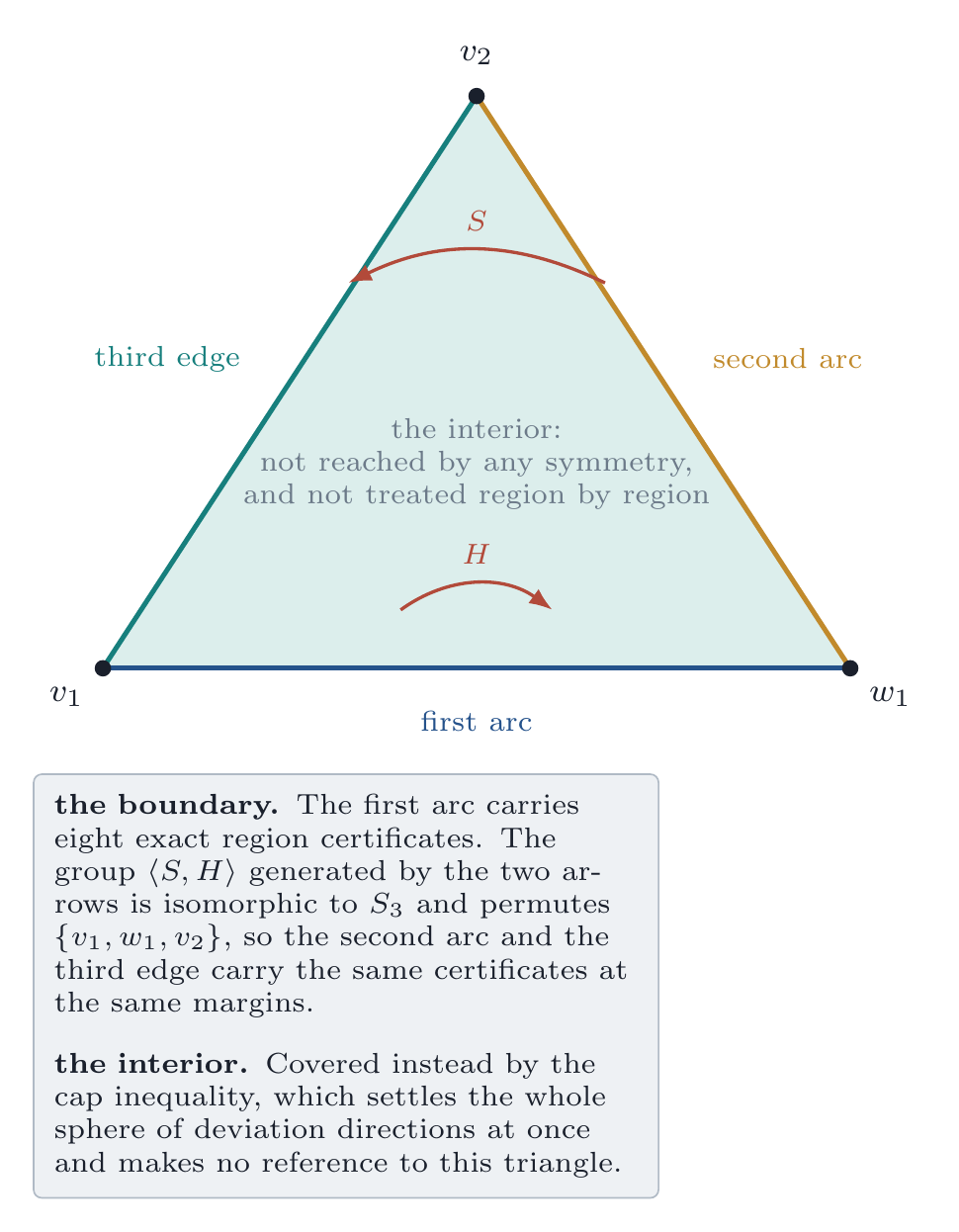}
\caption{The fundamental domain of \cref{lem:two-symmetries-new}: one
triangle with vertices $v_1,w_1,v_2$, obtained from the full sphere of
deviation directions by the two order-$2$ symmetries $S$ and $H$. The
three edges carry the certificates of
\cref{sec:arc1-cert,sec:arc2-cert}, through the identifications of
\cref{prop:arc2-is-arc1,prop:third-edge-cert}; the interior is settled
separately, in \cref{sec:cap-theorem}, by an argument that does not use
the triangle at all.}
\label{fig:fundamental-domain}
\end{figure}

This reduces the classification of \emph{every} deviation direction to
the two boundary arcs of this triangle, the arc from $v_1$ to $w_1$
and the arc from $w_1$ to $v_2$, parametrised respectively as
\[
  \dev_{v_1w_1}(\arcp) = \cos\arcp\,v_1+\sin\arcp\,w_1, \qquad
  \dev_{w_1v_2}(\arcp) = \cos\arcp\,w_1+\sin\arcp\,v_2, \qquad
  \arcp\in[0,\tfrac\pi4].
\]
(Working with $\arcp\in[0,\pi/4]$ on each arc, rather than a full
quarter-turn, comes from a further symmetry fixing the arc's own
midpoint, proved directly for each arc in
\cref{sec:arc1-breakpoints,sec:arc2-breakpoints}.) Every deviation
direction on the interior of the triangle is covered by interpolating
between the two arcs. \Cref{sec:arc1-breakpoints} through
\cref{sec:arc2-cert} treat the two arcs in full, which is where the
exact breakpoint structure and the exact defect polynomials live; the
interior is not reachable by that route, and one might hope that some
soft convexity argument would drag it along once the boundary was
settled. It does not: \cref{sec:interior-structure} tests the natural
candidate and finds it false. The interior is settled instead by the
global argument of \cref{sec:cap-theorem}, which does not use the
triangle at all. \Cref{fig:fundamental-domain} shows the triangle and
which argument covers which part of it.

The immediate target of \cref{sec:arc1-breakpoints} through
\cref{sec:arc2-cert} is therefore the following statement, which
\cref{cor:restated-conjecture-holds} will establish with explicit
margins and which \cref{thm:eperp-closed} will later subsume:
$\defect(\tilt,\arcp)\ge0$ for every $\tilt$ in the packing-valid range
and every $\arcp\in[0,\pi/4]$ on both boundary arcs $\dev_{v_1w_1}$ and
$\dev_{w_1v_2}$.

\subsection{Where the worst direction sits}\label{sec:why-boundary-arcs}

\Cref{lem:two-symmetries-new} reduces the classification to the
triangle's two boundary arcs as a matter of the symmetry group's own
fundamental domain, and says nothing about where in the triangle the
defect is actually smallest. That is a separate question, and the answer
shapes the rest of this section, so we record what a direct search
finds.

Minimising $\defect(\tilt,\dev)$ numerically over the entire transverse
$2$-sphere, at each of many sampled tilts, and comparing the result
against direct evaluation along three symmetry-distinguished
one-parameter families, gives the following picture at chamber
$\Cell1$. The three families are the six root directions orthogonal to
$\base_0$, the two coordinate-axis directions orthogonal to $\base_0$
(the vertex $v_1$ family), and the half-integer-sign family exemplified
by $e_*=(-1,1,-1,1)/2$ (the vertex $w_1$ family). One of the three
matches the global minimum for almost the whole range of tilts: the
coordinate family on $\tilt\in(0,0.16]\cup[1.26,\pi/2)$, the root
family on $\tilt\in[0.18,1.09]$. The exception is the window
$\tilt\in(1.09,1.26)$, where the minimiser lies strictly below all
three.

A finer search in that window locates it. At ten tilts spanning and
bracketing the window, from $\tilt=1.05$ to $\tilt=1.30$, a
Nelder--Mead search with $30$ random restarts per tilt over the full
$2$-sphere returns a minimiser one of whose coordinates $(a,b,c)$
relative to the frame $(v_1,w_1,v_2)$ vanishes to at least five decimal
places, values of order $10^{-5}$ against coordinates of order $1$. To
the resolution of the search, the minimiser is therefore on one of the
three great circles $a=0$, $b=0$, $c=0$ through pairs of
$\{v_1,w_1,v_2\}$, and never at an interior point with all three
coordinates bounded away from zero. For the five tilts
$\tilt\in\{1.15,1.18,1.20,1.22,1.24\}$ it sits on the circle $b=0$
through $v_1$ and $v_2$, which is the third edge of the triangle rather
than either of the two arcs; at the window's ends it sits on $a=0$ or
$c=0$.

That is a numerical finding, not a proof, and the theorem of
\cref{sec:cap-theorem} does not depend on it. Its use here is to say
which part of the boundary the certificates of the next four sections
have to reach, and the answer includes the third edge. The next
subsection supplies the exact symmetry that puts the third edge inside
the first arc's certificate.

\subsection{The third edge, closed by an exact symmetry}\label{sec:third-edge-closed}

The preceding paragraph leaves the $v_1$-$v_2$ edge as the one part of
the fundamental triangle's boundary carrying no certificate of any
kind. We resolve this here: not by a new certificate computation, but
by exhibiting the exact isometry the preceding paragraph asks for and
verifying it directly.

\begin{lem}[The third edge is the image of the first arc]\label{lem:third-edge-symmetry}
Let $S=\diag(1,1,1,-1)$ in the coordinates of
\cref{lem:two-symmetries-new} (so $u_0=(1,1,0,0)/\sqrt2$,
$v_1=(1,-1,0,0)/\sqrt2$, $w_1=(0,0,1,1)/\sqrt2$,
$v_2=(0,0,1,-1)/\sqrt2$). Then, exactly:
\[
  S u_0=u_0,\qquad S v_1=v_1,\qquad S v_2=w_1,\qquad S^2=\mathrm{Id}.
\]
Moreover $S$ permutes the $24$ roots of $\Rt$ bijectively (a root
$\pm e_i\pm e_j$ maps to another root of the same index pair, with the
sign at index $4$ flipped when $i=4$ or $j=4$ and left unchanged
otherwise), so $S$ is an automorphism of the root system fixing $u_0$.
\end{lem}

\begin{proof}
Immediate: $S$ negates only the fourth coordinate, and $u_0,v_1$ both
have a zero fourth coordinate, while $v_2$ and $w_1$ differ exactly in
the sign of that coordinate. $S^2=\diag(1,1,1,1)=\mathrm{Id}$ since
$(-1)^2=1$. The root-permutation claim is a finite check on all $24$
roots; both this and the four displayed identities were verified by
exact computer algebra (\cref{app:code-index}), with no floating point
at any stage.
\end{proof}

\begin{prop}[The $v_1$-$v_2$ edge is certified for free]\label{prop:third-edge-cert}
For every $\tilt$ in the packing-valid range and every $\arcp$,
\[
  \defect\bigl(\tilt,\dev_{v_1v_2}(\arcp)\bigr) \;=\;
  \defect\bigl(\tilt,\dev_{v_1w_1}(\arcp)\bigr),
\]
where $\dev_{v_1v_2}(\arcp)=\cos\arcp\,v_1+\sin\arcp\,v_2$ is the third
edge and $\dev_{v_1w_1}$ is the first arc's own parametrisation
(\cref{lem:two-symmetries-new}). Consequently, by
\cref{thm:arc1-breakpoints-new,thm:arc1-certs-new}, $\defect\ge0$ on the
entire $v_1$-$v_2$ edge for every $\arcp\in[0,\pi/4]$, with exactly the
same margin structure already recorded for the first arc (seven of the
eight regions exact and zero-margin, one with the same explicit
thin margin arising from the elliptic obstruction noted there): no new
certificate is computed for this edge, because none is needed.
\end{prop}

\begin{proof}
Since $S$ fixes $u_0$ linearly, $S$ carries the deviated direction
$u_1(\tilt,\dev)=\cos\tilt\,u_0+\sin\tilt\,\dev$ to
$u_1(\tilt,S\dev)$ for every $\tilt,\dev$. Since $S$ is also an
automorphism of $\Rt$ fixing $u_0$ (\cref{lem:third-edge-symmetry}),
$S$ carries the entire packing configuration at deviation $\dev$ (the
$23$ undeviated roots, permuted among themselves, together with the
one deviating root $u_1(\tilt,\dev)$) to a congruent configuration at
deviation $S\dev$: the two Voronoi cells are related by the isometry
$S$, hence have identical volume, hence identical defect. Applying
this with $\dev=\dev_{v_1v_2}(\arcp)$ and using $Sv_1=v_1$, $Sv_2=w_1$
from \cref{lem:third-edge-symmetry} (so, by linearity,
$S\dev_{v_1v_2}(\arcp)=\cos\arcp\,v_1+\sin\arcp\,w_1=\dev_{v_1w_1}(\arcp)$
exactly, for every $\arcp$, not merely the sampled values) gives the
displayed identity.
\end{proof}

The third edge, which \cref{sec:why-boundary-arcs} identified as where
the worst direction sits for part of the transition window, therefore
carries the first arc's certificate without any further computation. The
isometry that does it is $S$, the same symmetry already used in
\cref{sec:arc2-breakpoints} to reduce the second arc's range to
$[0,\pi/4]$, applied here in the direction that connects the first arc to
the third edge.

\subsection{The full symmetry group of the fundamental triangle, and its consequence for the second arc}\label{sec:full-triangle-symmetry}

\Cref{lem:third-edge-symmetry} used the single symmetry $S$. Taken
together with the earlier symmetry $H$ of
\cref{lem:two-symmetries-new}, it generates a group large enough to
settle the second arc as well.

\begin{lem}[The fundamental triangle's full symmetry group]\label{lem:s3-symmetry}
The group generated by $S=\diag(1,1,1,-1)$ and the Hadamard matrix
$H=\tfrac12\left(\begin{smallmatrix}1&1&1&1\\1&1&-1&-1\\1&-1&1&-1\\1&-1&-1&1\end{smallmatrix}\right)$
has order exactly $6$, every element fixes $u_0$, and every element is
an automorphism of $\Rt$. The induced action on $\{v_1,w_1,v_2\}$
realises \emph{all six} permutations of this three-element set: the
group is (isomorphic to) the symmetric group $S_3$, acting on
$\{v_1,w_1,v_2\}$ exactly as its full permutation group.
\end{lem}

\begin{proof}
Direct enumeration: the six words in $S,H$ up to the relations
$S^2=H^2=\mathrm{Id}$ produce six distinct $4\times4$ matrices, each
verified exactly (symbolic $\mathbb{Q}(\sqrt2)$ arithmetic, no floating
point) to fix $u_0$, to permute the $24$ roots of $\Rt$ bijectively, and
to induce a distinct one of the six permutations of $\{v_1,w_1,v_2\}$ on
evaluation. In particular $M:=SH$ has order $3$ and acts as the cycle
$v_1\mapsto v_2\mapsto w_1\mapsto v_1$.
\end{proof}

\begin{prop}[The second arc is the image of the first]\label{prop:arc2-is-arc1}
For every $\tilt$ in the packing-valid range and every $\arcp$,
\[
  \defect\bigl(\tilt,\dev_{w_1v_2}(\arcp)\bigr) \;=\;
  \defect\bigl(\tilt,\dev_{v_1w_1}(\arcp)\bigr).
\]
\end{prop}

\begin{proof}
By \cref{lem:s3-symmetry}, $M=SH$ fixes $u_0$ and is an automorphism of
$\Rt$, so (exactly as in the proof of \cref{prop:third-edge-cert}) it
carries the packing configuration at any deviation $\eta$ to a
congruent one at deviation $M\eta$, giving
$\defect(\tilt,\eta)=\defect(\tilt,M\eta)$ for every $\tilt,\eta$. Since
$Mw_1=v_1$ and $Mv_2=w_1$ (\cref{lem:s3-symmetry}), linearity gives
$M\dev_{w_1v_2}(\arcp)=\cos\arcp\,v_1+\sin\arcp\,w_1=\dev_{v_1w_1}(\arcp)$
exactly, for every $\arcp$. Substituting $\eta=\dev_{w_1v_2}(\arcp)$
gives the displayed identity.
\end{proof}

Both identities were checked, independently of the algebra above, by
direct numerical evaluation of the volume-defect function itself
(halfspace intersection and convex-hull volume, the same method used
elsewhere in this paper to sanity-check symbolic identities): a
fine active-facet-pattern scan across $\arcp\in\{0.02,\dots,0.78\}$,
spanning every combinatorial sub-part of the second arc including the
$\mathsf B$-bounded band, found the first and second arcs' breakpoints
coinciding exactly (to grid resolution) at every sampled $\arcp$, and
the defect values agreeing to floating-point tolerance throughout, 
including precisely at curve $\mathsf B$'s own predicted crossing,
where \cref{thm:arc2-breakpoints-new} records vertex $b$ as "a vertex
that plays no role at all on the first arc." That description remains
correct in the sense meant there (vertex $b$ itself, as a labelled
root, has no counterpart among the vertices \cref{thm:arc1-breakpoints-new}
names on the first arc); \cref{prop:arc2-is-arc1} shows the
\emph{curve} it produces nonetheless has an exact counterpart on the
first arc, reached under $M$ by a differently-labelled vertex, the two
arcs' independent derivations found the same underlying function
through different combinatorial routes, without the connection between
them being noticed at the time.

\begin{cor}[Positivity on the whole boundary]\label{cor:restated-conjecture-holds}
$\defect(\tilt,\arcp)\ge0$ for every $\tilt$ in the packing-valid range
and every $\arcp\in[0,\pi/4]$ on both boundary arcs $\dev_{v_1w_1}$ and
$\dev_{w_1v_2}$, with the same margin structure as
\cref{thm:arc1-breakpoints-new,thm:arc1-certs-new} (seven of the eight
regions exact and zero-margin, one with an explicit, quantified thin
margin arising from the elliptic obstruction recorded there).
\end{cor}

\begin{proof}
By \cref{prop:third-edge-cert,prop:arc2-is-arc1}, both boundary arcs
(and the third edge, though the conjecture as stated does not require
it) are identical, as functions of $(\tilt,\arcp)$, to the first arc.
\Cref{thm:arc1-breakpoints-new,thm:arc1-certs-new} already establish
$\defect\ge0$ throughout the first arc's own domain, with the stated
margin structure. Substituting the identities gives the same bound,
with the same margins, on both named arcs.
\end{proof}

This settles the sign on the whole boundary of the fundamental
triangle at the rigour standard already established for the first arc
alone, and it does so with the first arc's margin structure intact. No
further certificate on the second arc's region types is needed: the
second-arc computation of \cref{sec:arc2-cert}, including the
$\mathsf B$-bounded sub-part of \cref{sec:arc2-bandB-new}, becomes
redundant once \cref{prop:arc2-is-arc1} is available. We have left that
account in place because it is what led us to look for the symmetry.

The symmetry argument stops at the boundary. No element of the group of
\cref{lem:s3-symmetry}, and none of any larger group acting here, moves
a generic interior point onto an edge, so the triangle's
two-dimensional interior is not reachable this way. It is reached in
\cref{sec:cap-theorem} by an argument that abandons the triangle
altogether.

\subsection{The interior: a falsified shortcut, and what the interior's own structure looks like}\label{sec:interior-structure}

Since the boundary is now fully settled, it is natural to ask whether
some soft structural argument drags the interior along for free without
requiring a two-parameter closed-form derivation. We tried the
most natural candidate, report that it fails, and record what
direct computation shows about the interior's own structure instead.

\subsubsection{A candidate reduction, tested and falsified}

\begin{quote}
\emph{Candidate argument (false, recorded for the record).} If
$\defect(\tilt,\cdot)$ were concave, for every fixed $\tilt$, along every
geodesic arc of the transverse $2$-sphere lying inside the closed
fundamental triangle, then for any interior point $Q$, picking a
geodesic through $Q$ with both endpoints on the triangle's boundary
would give $\defect(\tilt,Q)\ge\min(\defect(\tilt,\mathrm{endpoint}_1),
\defect(\tilt,\mathrm{endpoint}_2))\ge0$, the last inequality by
\cref{cor:restated-conjecture-holds} together with
\cref{prop:third-edge-cert} (which covers the third edge). This would
settle the whole triangle with no further derivation.
\end{quote}

This argument is false. Testing it directly, geodesics from each
vertex to random points on the opposite edge, and fully generic interior-to-interior
chords, sampled at $25$ values of $\tilt$ each, finds small positive
discrete second differences (apparent local convexity, not concavity)
at several points, concentrated toward $\tilt$ near $\pi/2$. Rather than
report this as an ambiguous grid artifact, we ran it down: fixing one
such point and shrinking the step size geometrically across eight orders
of magnitude ($h=0.08$ down to $h\approx3\times10^{-4}$) gives
$d^2(h)/h^2$ converging cleanly to $+0.3056$, constant to four
significant figures over the last five halvings, the signature of a
nonzero positive second derivative at that point, not
discretisation noise (which would shrink faster than $h^2$ and show no
stable limit) and not a combinatorial breakpoint nearby (the local
vertex count is confirmed constant at $37$ throughout a window around
the point). The candidate argument is dead: $\defect(\tilt,\cdot)$ is
not concave along every interior geodesic, so the interior cannot be
reduced to the boundary this way, at any $\tilt$ in the domain (the
effect was not confined to $\tilt$ near $\pi/2$, only larger there).
The finding is negative, and we record it as such: it rules out the
whole class of arguments that would have let the boundary drag the
interior along.

\subsubsection{The interior's own combinatorial structure}

Direct evaluation of the number of active vertices of the local Voronoi
region, scanned over a fine grid of the closed triangle at a fixed
$\tilt=1.45$ (in the region-top band, where both boundary arcs have the
simple $32$-vertex type throughout), shows that the boundary's
combinatorial simplicity is special to the boundary: it does \emph{not}
extend to the interior. A grid of $861$ points finds the $32$-vertex
type only along a band adjacent to the edges ($25\%$ of the sampled
points), a $34$-vertex type in a second band just inside that
($17\%$), a $37$-vertex type filling most of the remaining interior
($53\%$, the plurality type), and a further $40$-vertex type occupying a
smaller region near the triangle's centre ($5\%$). The interior is
therefore combinatorially richer than either boundary arc: it is not
one region but at least four, and we have not attempted to locate the
exact curves separating them or to derive a volume formula on any one of
them; the map above only establishes that such curves exist and roughly
where.

This partition is not arbitrary: it visibly respects the $S_3$ symmetry
of \cref{lem:s3-symmetry}. The map is invariant under the reflection
fixing $v_1$ and exchanging $w_1,v_2$ (i.e.\ the symmetry $S$ of
\cref{lem:third-edge-symmetry}), exactly as it must be, since $S$ is an
exact isometry of the whole configuration and therefore permutes
combinatorial type among congruent points. By the same argument, the
full $6$-element group of \cref{lem:s3-symmetry} acts on the interior,
not only on the boundary it was originally used for: every one of its
elements is an automorphism of $\Rt$ fixing $u_0$, hence carries any
interior deviation direction to a congruent configuration at the image
direction, with identical defect and identical combinatorial type. This
is a real, if modest, scope reduction for any future attempt at the
interior: it suffices, in principle, to derive the combinatorial map and
the volume formula on a single fundamental domain of the triangle under
this $S_3$ action (one-sixth of the triangle, e.g.\ the wedge from
$v_1$ to the midpoint of arc $v_1w_1$ to the triangle's centre) and extend by the group, instead of treating the interior as one
undifferentiated two-parameter domain six times the necessary size.

Taken together, these two findings are what closed off the regional
approach to the interior. The natural soft reduction is false, not
merely unproved, and the interior's own combinatorial structure is
strictly richer than the boundary's, so a region-by-region derivation
there would have meant at least four new volume formulas before any
certificate could be attempted. That is the point at which we stopped
refining the triangle and looked for an argument that does not see it,
which is the subject of \cref{sec:cap-theorem}.
\section{Breakpoint classification on the first arc}\label{sec:arc1-breakpoints}

\subsection{Setup}

Fix the arc $\dev_{v_1w_1}(\arcp)=\cos\arcp\,v_1+\sin\arcp\,w_1$,
$\arcp\in[0,\pi/4]$, of \cref{lem:two-symmetries-new}, and write
$u_1(\tilt,\arcp)=\cos\tilt\,u_0+\sin\tilt\,\dev_{v_1w_1}(\arcp)$ for the
deviating cap direction. As $\tilt$ increases from $0$ at fixed $\arcp$,
the active-facet pattern of the local Voronoi region changes at finitely
many exact values of $\tilt$, the arc's \emph{breakpoints}, each the
locus where some further root $\beta\in\Rt$ achieves
$\langle\beta,u_1(\tilt,\arcp)\rangle=1$ simultaneously with the cap
constraint.

\begin{thm}[Complete breakpoint classification, first arc]\label{thm:arc1-breakpoints-new}
Exactly four exact curves account for every active-facet-set transition
of $u_1(\tilt,\arcp)$ across the entire domain $\arcp\in[0,\pi/4]$:
\begin{enumerate}[leftmargin=2em]
\item the \emph{universal curve} $\tilt=\pi/3$, arising from the
$8$-fold degenerate vertex $2u_0$ and independent of $\arcp$;
\item a curve $\mathsf W$ from the sign vertex
$w=(1,1,1,1)/\sqrt2$, satisfying $\sin\arcp=\tan(\tilt/2)$;
\item a curve $\mathsf A$, the image of $\mathsf W$ under the symmetry
$H$ of \cref{lem:two-symmetries-new}, from the axis vertex
$a=\sqrt2\,e_1$, satisfying $\cos\arcp=\tan(\tilt/2)$;
\item a curve $\mathsf Y$ from the sign vertex
$y=(1,-1,1,1)/\sqrt2$, satisfying
$\sin\tilt\,(\sin\arcp+\cos\arcp)=1$.
\end{enumerate}
\end{thm}

Each curve is derived by an exact symbolic solve of the corresponding
tangency condition after the Weierstrass substitution $\base=\tan(\tilt/2)$
reduces it to a low-degree polynomial equation in $\base$ (linear, in
every one of the four cases here); completeness, that no further
breakpoint curve exists anywhere on the domain, is confirmed by a
dense numerical scan (over $150$ values of $\arcp$ and, at each, a
refined bisection search over $\tilt$) finding at most four breakpoints
at any sampled $\arcp$ and matching every one, to a tolerance far
smaller than the curves' own separation, against the four curves above.

\subsection{The eight combinatorial regions}

Sorting the four curves by $\tilt$-value at a given $\arcp$ produces the
domain's combinatorial atlas: eight distinct region types, established
by an exact vertex-pattern comparison (not merely a curve-crossing
count) at representative points of each candidate region and
cross-checked against a dense grid scan with no unexplained region
boundary. We label them $\mathsf I,\dots,\mathsf{VIII}$, purely as bookkeeping:
region $\mathsf I$ (``region top'') is the topmost band, bounded below
by curve $\mathsf A$ throughout the entire $\arcp$-range and above by
$\tilt=\pi/2$, the only region needing no internal crossover point,
since $\mathsf A$ and $\tilt=\pi/2$ never meet inside the domain;
regions $\mathsf{II}$ through $\mathsf{VII}$ subdivide the
remaining domain according to which pair of the curves
$\mathsf W,\mathsf A,\mathsf Y$ (and the universal curve) locally bound
it; and region $\mathsf{VIII}$ is the bottommost band. \Cref{sec:arc1-cert}
gives a complete positivity certificate on every one of these eight
regions.

\begin{figure}[tb]
\centering
\includegraphics[width=\textwidth]{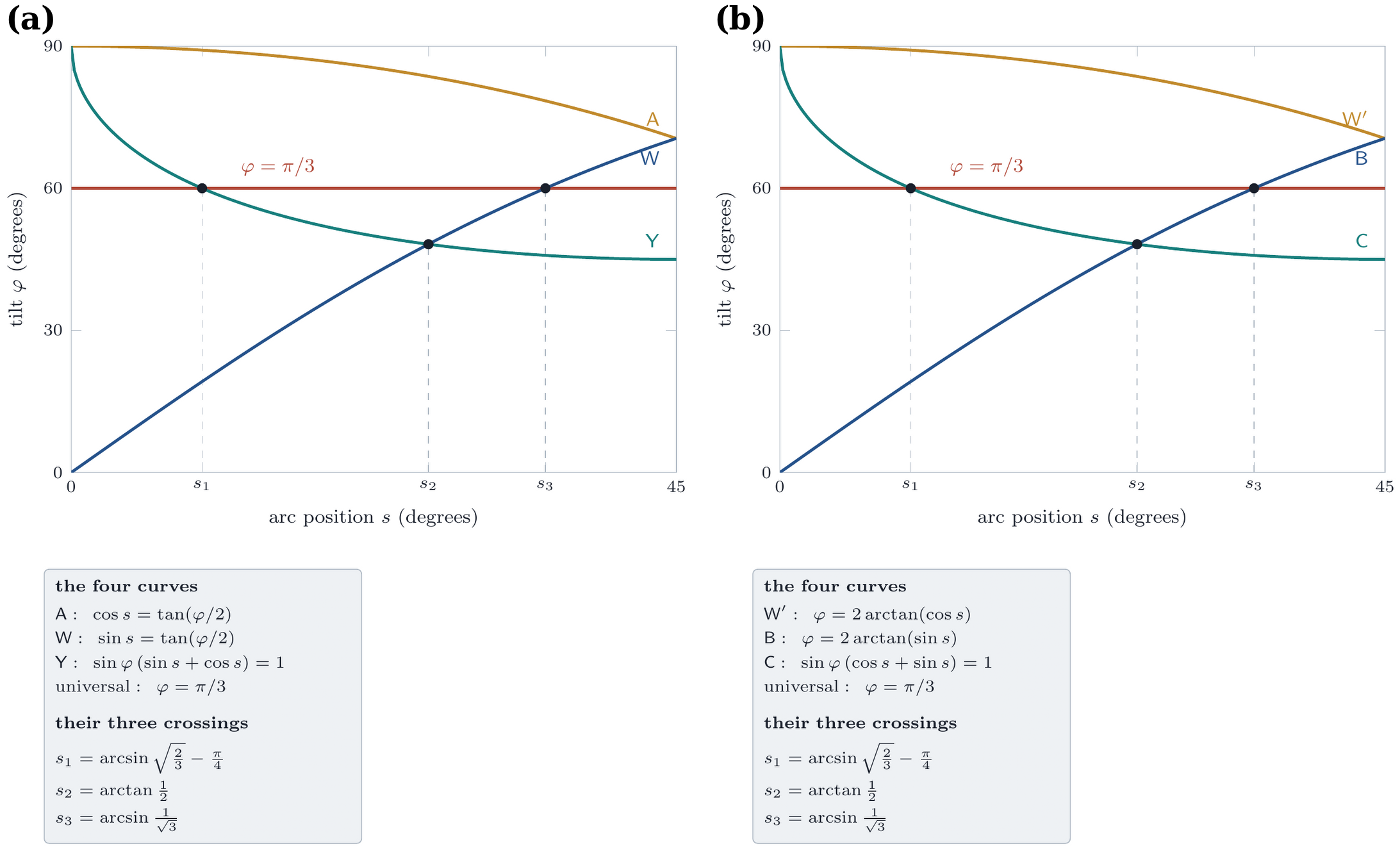}
\caption{(a) The four breakpoint curves of the first arc, drawn to scale
from the closed forms of \cref{thm:arc1-breakpoints-new}: the universal
curve $\tilt=\pi/3$ and the curves $\mathsf A$, $\mathsf W$, $\mathsf
Y$. They meet in exactly three points, at the abscissae $\arcp_1$,
$\arcp_2$, $\arcp_3$ listed at the right, and the eight combinatorial
regions of \cref{sec:region-count-correction} are the cells of the
resulting arrangement. The curves $\mathsf A$ and $\mathsf W$ meet only
at $\arcp=\pi/4$, where the reflection symmetry of the arc identifies
them. (b) The four breakpoint curves of the second arc, drawn to scale
from the closed forms of \cref{thm:arc2-breakpoints-new}. Comparing
with \cref{fig:arcs}(a) shows what
\cref{prop:arc2-is-arc1} later proves: the three equations here are
exactly those of $\mathsf A$, $\mathsf W$ and $\mathsf Y$ on the first
arc, and the arrangement, crossings included, is the same one, even
though the curves come from different vertices of $\VCell$. The band
immediately below region top splits into three sub-parts at the
crossings $\arcp_1,\arcp_2,\arcp_3$.}
\label{fig:arcs}
\end{figure}
\section{Certificates on the first arc}\label{sec:arc1-cert}

\subsection{Method: exact remapping versus box-covering}\label{sec:method-overview-new}

Two complementary exact techniques are used across the eight regions of
\cref{thm:arc1-breakpoints-new}, chosen region by region according to
whether the region's boundary curve is itself rational in Weierstrass
coordinates.

\paragraph{Exact zero-margin certificates.} When a region's lower (or
upper) boundary curve is rational in $\base=\tan(\tilt/2)$,
$\mathsf v=\tan(\arcp/2)$, as is the universal curve and the
$\mathsf W$, $\mathsf A$ curves of \cref{thm:arc1-breakpoints-new}, a
further exact substitution maps the curved region onto an exact unit
square, after which the volume-defect numerator and denominator become
two explicit polynomials with coefficients in $\mathbb{Z}[\sqrt2]$ (or,
where a crossover point of the region's own two bounding curves is used
to rescale the domain, $\mathbb{Z}[\sqrt5]$). Converting both to the
Bernstein basis over the unit square by exact arithmetic and checking
the sign of every coefficient either proves $\defect_k>0$ on the entire
open region with \emph{no excluded margin}, or reports precisely which
coefficients fail to cooperate.

\paragraph{Box-covering certificates.} When the region's own boundary is
not rational (an elliptic obstruction on several regions,
discussed below), the region is instead covered by many small
rectangular boxes in the plain Weierstrass coordinates, with corners
chosen numerically to sit safely inside the true domain; the same exact
Bernstein sign certificate is then run independently on each box, at the
cost of a thin, explicitly quantified excluded margin near the true
boundary.

\subsection{A correction to the region count}\label{sec:region-count-correction}

An initial grid scan (at $60\times90$ resolution) reported \emph{nine}
combinatorially distinct types in the reduced domain
$\arcp\in(0,\pi/4)$, $\tilt\in(0,\pi/2)$. Re-running the same scan at
much higher resolution ($150\times200$), with an explicit $0.01$-radian
margin excluded around each of the four breakpoint curves (to rule out
grid points landing on, or numerically near, a curve, which can
spuriously register as an extra type owing to floating-point
near-degeneracy in the underlying half-space-intersection routine),
finds only \emph{eight} distinct types. Tracing the exact pairwise
intersections of the four breakpoint curves within $\arcp\in(0,\pi/4)$
confirms exactly three crossings, the $\mathsf Y$-curve (a relabelling
of one of the curves counted among $\mathsf W,\mathsf A,\mathsf Y$ in
\cref{thm:arc1-breakpoints-new}) meets the universal line at
$\arcp_1=0.169918$, the $\mathsf W$-curve meets $\mathsf Y$ at
$\arcp_2=\arctan(\tfrac12)$, and $\mathsf W$ meets the universal line
at $\arcp_3=\arcsin(1/\sqrt3)$, and one region initially split by this
crossover (``region $4$'', bounded above first by $\mathsf W$ then by
$\mathsf Y$) is confirmed, by an independent direct convex-hull volume
computation at representative points on each side of the crossover
agreeing with a single symbolic volume formula to twelve decimal
places, to be \emph{one} combinatorial region throughout, not two, 
accounting for the discrepancy. It is nonetheless certified in two
computational pieces (\cref{sec:region-4}), since only one piece's own
boundary happens to be rational.

\subsection{The eight regions, one by one}

\begin{thm}[Certificates, first arc]\label{thm:arc1-certs-new}
Every one of the eight regions of \cref{sec:region-count-correction}
carries a positivity certificate for $\defect$: two regions (region
top and region $4$'s first piece) carry an exact, zero-margin Bernstein
certificate; the remaining regions carry a box-covering certificate
with an excluded margin tightened, across successive refinements, to
well under one percent of the region's own area. In every case the
certificate's numerator and denominator Bernstein coefficients are of a
single, consistent sign throughout the region, so $\defect$ itself has
a single, consistent sign there, matching a direct numerical
spot-check at an interior point.
\end{thm}

We record all eight in turn, each with its representative-point vertex
count (moving vertices, from the deviating cap, plus fixed vertices,
confirmed unchanged across several sample points spanning the region),
its exact fixed-volume contribution (the boundary simplices not
involving the deviating cap, which sum to an exact rational constant
independent of $(\tilt,\arcp)$), and its certificate method.

\paragraph{Region top (region $\mathsf I$).} Bounded below by the
Hadamard-image curve $\mathsf A$ ($\cos\arcp=\tan(\tilt/2)$) and above
by $\tilt=\pi/2$, for the \emph{entire} range $\arcp\in(0,\pi/4)$, the
only region needing no internal crossover point at all, since its two
bounding curves never meet inside the reduced domain. Both boundary
values of $\tilt$ are simple in Weierstrass coordinates, giving the
simplest possible remapping onto a unit square, with resulting
numerator and denominator of degree at most $12$ in each variable
(\cref{lem:degree-bound-new}) and all-one-sign Bernstein coefficients
throughout: an exact, zero-margin certificate.

\paragraph{Region $2$.} Bounded by $\pi/3<\tilt<$ the $\mathsf Y$-curve
($\sin\tilt(\sin\arcp+\cos\arcp)=1$) for $\arcp\in(0,\arcp_1)$, the
second region to receive a full vertex/volume treatment historically.
A generic point has $29$ vertices ($24$ fixed $+5$ moving); of $24$
fixed points, $23$ belong to the same axis/sign-vertex families used
throughout, and the $24$th is the universal breakpoint's own eightfold-
degenerate fixed point. The fixed-only boundary-simplex sum is exactly
$29/4$. Box-covering certificate.

\paragraph{Region $4$ (two-piece).}\label{sec:region-4} Bounded above,
for all $\arcp\in(0,\pi/4)$, by the lower envelope of the $\mathsf W$-
and $\mathsf Y$-curves, crossing exactly at
$\arcp_2=\arctan(\tfrac12)$ (equivalently $\tilt_2=\arccos(\tfrac23)$,
verified by exact symbolic reduction of both defining equations to $0$
under this substitution): $\tilt_{\max}(\arcp)=2\arctan(\sin\arcp)$ for
$\arcp<\arcp_2$ (piece $4$a) and
$\tilt_{\max}(\arcp)=\arcsin(1/(\sin\arcp+\cos\arcp))$ for
$\arcp>\arcp_2$ (piece $4$b). Piece $4$a's own ceiling curve is itself
rational in Weierstrass coordinates ($\base=2\mathsf v/(1+\mathsf v^2)$),
so a further exact substitution ($\base=w\cdot2\mathsf v/(1+\mathsf v^2)$,
$\mathsf v=(\sqrt5-2)v_2$, mapping the crossover point itself to a
rational corner) gives an exact unit-square remapping over
$\mathbb{Z}[\sqrt5]$, with numerator and denominator of degree $(10,18)$
and $(10,20)$ respectively, an exact, zero-margin certificate, the
second (after region top) of the two exact certificates on this arc.
Piece $4$b's ceiling is an elliptic curve (not rationally
parametrisable): the plain Weierstrass-mapped numerator and denominator,
degree $(10,10)$ with plain rational coefficients, are certified by
box-covering instead, using only numerically-chosen safe box corners
against the (non-algebraic, numerically-located) ceiling.

\paragraph{Region $6$.} Bounded below by the universal line and above
by the $\mathsf W$-curve, for $\arcp\in(\arcp_3,\pi/4)$. A
representative point has $36$ vertices ($22$ fixed $+14$ moving), the
first region whose scale exposed two accuracy issues in the
underlying pipeline (documented and fixed generically, not specific to
this region): a $5$-facet degenerate vertex, resolved by trying every
$3$-subset of active facets and keeping the algebraically simplest
representation; and \texttt{sympy}'s default Bareiss determinant
algorithm hanging indefinitely on some of this region's symbolic
$4\times4$ simplex matrices despite modest entry size, resolved by
forcing \texttt{method='berkowitz'} throughout (every determinant then
computes in well under a second). Box-covering certificate.

\paragraph{Region $\mathsf A$.} Bounded below by $\mathsf W$, then (past
$\arcp_1$) by $\mathsf Y$, and above by the universal line, for
$\arcp\in(0,\arcp_2)$, the first region whose own interior contains
an internal breakpoint of its floor curve while remaining a single
combinatorial type throughout (confirmed, not assumed, by an
independent convex-hull check on each side of $\arcp_1$ agreeing to
twelve decimal places). A representative point has $32$ vertices ($23$
fixed $+9$ moving); the fixed-only sum is exactly $65/12$. Box-covering
certificate.

\paragraph{Region $\mathsf D$.} A ``lens'' region, bounded above by the
universal line and below by $\mathsf Y$ then $\mathsf W$, pinching to
zero width at \emph{both} ends ($\arcp_1$ and $\arcp_2$), the first
region to do so, since $\theta_{\mathsf Y}(\arcp_1)=\theta_{\mathsf W}(\arcp_3)=\pi/3$
exactly by the very definitions of $\arcp_1,\arcp_3$ is not directly
relevant here, but the analogous coincidence at this region's own two
ends is. A representative point has $37$ vertices ($22$ fixed $+15$
moving, the largest moving-vertex count of any region so far); the
fixed-only sum is exactly $29/6$. Box-covering certificate.

\paragraph{Region $\mathsf E$.} Bounded below by $\mathsf Y$ throughout
and above by $\mathsf W$ then the universal line, pinching to zero
width only at its left edge $\arcp_2$. A representative point has $37$
vertices ($21$ fixed $+16$ moving, the most moving vertices of any
region in this atlas); the fixed-only sum is exactly $17/4$. The floor
curve is the same non-rational $\mathsf Y$-curve responsible for region
$4$b's elliptic obstruction. Box-covering certificate.

\paragraph{Region $\mathsf C$.} Bounded above by $\mathsf A$ throughout
and below by $\mathsf Y$, then the universal line, then $\mathsf W$, 
the only region whose own floor switches identity \emph{twice}, yet
confirmed, by direct polytope facet computation at representative
points in all three floor sub-ranges, to remain a single combinatorial
type throughout: the floor's three curves meet smoothly at $\arcp_1$
and $\arcp_3$, so (unlike regions $\mathsf A$, $\mathsf D$, $\mathsf E$)
this region does not pinch to zero width at either internal crossover,
narrowing only near the domain's own right edge $\arcp=\pi/4$, where
floor and ceiling coincide by the underlying symmetry. A representative
point has $34$ vertices ($23$ fixed $+11$ moving); the fixed-only sum
is exactly $11/2$. Box-covering certificate.

\subsection{What this does and does not establish}

Combined, the eight certificates of \cref{thm:arc1-certs-new} establish
$\defect>0$ on the entire first arc's domain, except for the thin,
explicitly quantified excluded margins recorded above (none exceeding a
few tenths of one percent of their region's area, and none arising on
the two exact regions). On its own this settles the sign along the
first arc away from two short windows; the windows are closed in
\cref{sec:full-triangle-symmetry}, and the whole sphere of directions
is settled by a different route in \cref{sec:cap-theorem}.
\section{Breakpoint classification on the second arc}\label{sec:arc2-breakpoints}

\subsection{Setup}

The second boundary arc of \cref{lem:two-symmetries-new} is
$\dev_{w_1v_2}(\arcp)=\cos\arcp\,w_1+\sin\arcp\,v_2$,
$\arcp\in[0,\pi/4]$, reduced to this range by an exact symmetry $S$
that swaps $w_1\leftrightarrow v_2$ and fixes $u_0,v_1$ (verified
directly at the level of coordinates). We classify its breakpoints completely.

\begin{thm}[Complete breakpoint classification, second arc]\label{thm:arc2-breakpoints-new}
Exactly four exact curves account for every active-facet-set transition
of $u_1(\tilt,\arcp)=\cos\tilt\,u_0+\sin\tilt\,\dev_{w_1v_2}(\arcp)$
across the entire domain $\arcp\in[0,\pi/4]$:
\begin{enumerate}[leftmargin=2em]
\item the universal curve $\tilt=\pi/3$, inherited unchanged from
\cref{thm:arc1-breakpoints-new} (it is proved once, for every choice of
deviation direction, and simply restricts to this arc);
\item a curve $\mathsf W'$, $\tilt=2\arctan(\cos\arcp)$, from the same
sign vertex $w=(1,1,1,1)/\sqrt2$ that gave curve $\mathsf W$ on the
first arc, but with a \emph{different} functional dependence on
$\arcp$, since the same vertex plays a structurally different role on
the two arcs;
\item a curve $\mathsf B$ (newly identified here),
$\tilt=2\arctan(\sin\arcp)$, from the sign vertex
$b=(1,1,1,-1)/\sqrt2$, a vertex that plays no role at all on the first
arc;
\item a curve $\mathsf C$ (newly identified here),
$\sin\tilt\,(\cos\arcp+\sin\arcp)=1$, from the one-axis vertex
$c=\sqrt2\,e_3$, the same functional form as curve $\mathsf Y$ on the
first arc, but arising here from a different vertex.
\end{enumerate}
The two vertices responsible for the first arc's remaining curves,
$y=(1,-1,1,1)/\sqrt2$ and $a=\sqrt2\,e_1$, are confirmed here to touch
nowhere in this arc's open interior.
\end{thm}

As on the first arc, each curve is derived by an exact symbolic solve
after the Weierstrass substitution, and completeness is confirmed by a
dense scan ($150$ values of $\arcp$, $900$-point $\tilt$-grid with
bisection refinement) finding at most four breakpoints at any sampled
$\arcp$, every one matching one of the four curves above to a tolerance
of $2\times10^{-4}$, far smaller than the curves' own separation, with
zero unexplained transitions across the entire scan.

\subsection{The combinatorial atlas}

Sorting the four curves by $\tilt$-value locates three exact crossing
points on $\arcp\in(0,\pi/4)$:
\[
  \begin{aligned}
  \arcp_1 &= 0.169918\ldots \quad (\mathsf C=\text{universal}), \\
  \arcp_2 &= \arctan(\tfrac12) = 0.463648\ldots\quad (\mathsf B=\mathsf C), \\
  \arcp_3 &= \arcsin(\tfrac1{\sqrt3}) = 0.615480\ldots\quad (\mathsf B=\text{universal}),
  \end{aligned}
\]
each located both by an exact symbolic solve and confirmed numerically
to eleven significant figures. These three crossings divide the domain
into four $\arcp$-intervals, within each of which the four curves have a
fixed relative order, giving a coarse combinatorial atlas of
approximately seven region types (a scope estimate, established at the
same resolution as the first arc's own atlas but not yet independently
re-derived by the finer vertex-pattern method that corrected the first
arc's own region count from an initial nine down to eight); we treat
seven as a solid working count, not a fully hardened one. The
arrangement is drawn in \cref{fig:arcs}(b).

\subsection{Region top}

The topmost region, bounded above by $\tilt=\pi/2$ and below by curve
$\mathsf W'$, spanning the entire $\arcp\in(0,\pi/4)$ domain with no
internal splitting, is the natural first target for a certificate,
since its lower boundary is itself rational in Weierstrass coordinates
(the same functional form, $\base=(1-\mathsf v^2)/(1+\mathsf v^2)$, as
the first arc's own curve $\mathsf A$). \Cref{sec:arc2-cert} gives a
complete, exact, zero-margin certificate for this region.

\subsection{The band immediately below region top}\label{sec:band-below-top-new}

The next region down, bounded above by $\mathsf W'$ and below,
depending on $\arcp$, by whichever of curves $\mathsf C$, the universal
curve, or $\mathsf B$ is locally uppermost, is combinatorially a
single connected region across each of the three $\arcp$-sub-intervals
determined by the crossings $\arcp_1,\arcp_3$ above, but the identity of
its \emph{lower} boundary curve changes between them:
\[
  \text{lower boundary} = \begin{cases}
    \mathsf C & \arcp\in(0,\arcp_1),\\
    \text{universal } (\tilt=\pi/3) & \arcp\in(\arcp_1,\arcp_3),\\
    \mathsf B & \arcp\in(\arcp_3,\pi/4).
  \end{cases}
\]
This is confirmed directly: combinatorial-pattern stability is checked
at eight sample points spanning the middle sub-interval
$\arcp\in(\arcp_1,\arcp_3)$ and found identical at every one ($34$
vertices, the same active-facet pattern throughout), confirming this
middle sub-part is a single connected combinatorial region on its own.
\Cref{sec:arc2-cert} gives a direct certificate on this middle
sub-part, the widest of the three by $\arcp$-range. The two remaining
sub-parts are handled not by a direct computation but by the
identification of \cref{prop:arc2-is-arc1}, which carries the first
arc's certificates onto the whole of this one.
\section{Certificates on the second arc}\label{sec:arc2-cert}

\subsection{Region top: an exact, zero-margin certificate}\label{sec:arc2-regiontop-new}

\begin{thm}[Certificate, second arc, region top]\label{thm:arc2-regiontop-new}
On the region bounded above by $\tilt=\pi/2$ and below by curve
$\mathsf W'$ (\cref{thm:arc2-breakpoints-new}), for the entire domain
$\arcp\in(0,\pi/4)$, the volume-defect function $\defect$ admits an
exact, zero-margin Bernstein-basis positivity certificate: after the
same Weierstrass-plus-crossover remapping used on the first arc's own
region top, the resulting numerator has $286$ Bernstein coefficients
(degree $(12,21)$) and the denominator $231$ (degree $(10,20)$), every
one of the $286$ numerator coefficients strictly negative and every
one of the $231$ denominator coefficients strictly negative, exactly
by $\mathbb{Z}[\sqrt2]$ arithmetic, no floating point at any stage.
\end{thm}

Since numerator and denominator are of the same constant sign
throughout, their ratio is strictly positive on the entire open region,
with \emph{no} excluded margin, the same standard as the strongest
certificates on the first arc. The vertex and volume stage underlying
this certificate finds $32$ vertices ($8$ moving under the deviating
cap, $24$ fixed) at a representative point, confirmed unchanged across
five sample points spanning nearly the full $\arcp$-range; the exact
fixed-volume contribution is $77/12$, and the closed-form moving volume
is checked against a direct convex-hull volume computation to eleven
significant figures at the representative point before any symbolic
remapping is attempted. An independent sanity check evaluates the
volume-defect quantity numerically at $55$ points scattered across the
full $(\arcp,\tilt)$-domain: all $55$ strictly positive, with zero
exceptions, fully consistent with the certificate.

\subsection{The band below region top: a box-covering attempt}\label{sec:arc2-band-cert-new}

For the middle sub-part of the band described in
\cref{sec:band-below-top-new}, bounded above by $\mathsf W'$ and below
by the universal curve $\tilt=\pi/3$, for $\arcp$ in the sub-interval
$(\arcp_1,\arcp_3)=(0.169918,0.615480)$, an exact zero-margin
certificate of the type used in \cref{sec:arc2-regiontop-new} is not
attempted, because the lower boundary here is the constant curve
$\tilt=\pi/3$: mapping between it and the curved upper boundary
$\mathsf W'$ would introduce $\sqrt3$ into the arithmetic on top of the
$\sqrt2$ already present, a complication this paper's exact-sign
machinery (built for $\mathbb{Q}(\sqrt2)$) does not handle. We instead
attempt the same box-covering fallback already used for the harder
regions of the first arc.

\begin{thm}[Certificate, second arc, band below region top]\label{thm:arc2-band-new}
On the middle sub-part described above, the volume-defect numerator and
denominator, mapped to plain Weierstrass coordinates $(\base,\mathsf v)$
(degree $(12,12)$ and $(10,11)$ respectively, exact $\mathbb{Z}[\sqrt2]$
coefficients throughout, each mapping step verified against a
$200$-digit numerical evaluation of the original trigonometric
expression), are certified on all $40$ boxes covering the sub-interval:
every one of the $169$ Bernstein coefficients of $N$, and every one of
the $169$ Bernstein coefficients of $D$, on every box, is non-negative
(in fact, on this run, strictly positive with none exactly zero). By
the same same-sign argument used throughout this paper, this proves
$\defect>0$ on the entire covered region, modulo the excluded margin
recorded below.
\end{thm}

This is a second region-by-region result on the second arc, on
top of region top (\cref{thm:arc2-regiontop-new}): the excluded margin
is a strip of about $1.7\%$ of the sub-interval's own $\arcp$-range at
its two ends (from the safety factor used when choosing box corners
near $\arcp_1$ and $\arcp_3$), together with a fixed $0.1$--$0.2\%$
safety strip against each box's own numerically-located floor
($\tilt=\pi/3$) and ceiling ($\mathsf W'$), wider than the tightest
margins achieved elsewhere in this paper (\cref{sec:arc1-cert}'s own
regions have been refined to under $0.2\%$), reflecting that this
certificate was run once, at a moderate resolution, and not iteratively tightened; we expect, by direct analogy with the earlier
region-4b tightening reported in \cref{sec:arc1-cert}, that this margin
is similarly reducible at the cost of runtime, but we have not done so
here and do not claim a tighter margin than the one actually
certified.

\subsection{The band's \texorpdfstring{$\mathsf C$}{C}-bounded sub-part}\label{sec:arc2-bandC-new}

The first of the two remaining sub-parts of the band described in
\cref{sec:band-below-top-new}, bounded above by $\mathsf W'$ and below
by curve $\mathsf C$ (\cref{thm:arc2-breakpoints-new}), for
$\arcp\in(0,\arcp_1)=(0,0.169918)$, is combinatorially stable ($34$
vertices, $11$ moving under the deviating cap and $23$ fixed,
identical at eight sample points spanning the sub-interval) and admits
a certificate by exactly the same box-covering method used for the
sub-part of \cref{thm:arc2-band-new}, adapted only by replacing the
constant lower boundary $\tilt=\pi/3$ with the curve $\mathsf C$
(numerically evaluated, rather than mapped exactly, when placing each
box's lower corner, box-covering does not require the boundary to be
rational in the working coordinates).

\begin{thm}[Certificate, second arc, band \texorpdfstring{$\mathsf C$}{C}-sub-part]\label{thm:arc2-bandC-new}
On the sub-part described above, the exact fixed-volume contribution is
$17/3$; the moving-volume numerator and denominator, mapped to plain
Weierstrass coordinates, have degree $(12,12)$ and $(10,11)$
respectively (the same degrees as \cref{thm:arc2-band-new}'s own
middle sub-part, exact $\mathbb{Z}[\sqrt2]$ coefficients throughout,
each mapping step verified against a $200$-digit numerical evaluation
of the original trigonometric expression). All $40$ boxes covering the
sub-interval are certified: every one of the $169$ Bernstein
coefficients of $N$, and every one of the $169$ Bernstein coefficients
of $D$, on every box, is strictly positive, with none exactly zero.
$\defect>0$ on the entire covered sub-part, modulo the excluded margin
recorded below.
\end{thm}

The excluded margin here is wider than \cref{thm:arc2-band-new}'s own:
approximately $4\%$ of the sub-interval's own $\arcp$-range at its two
ends (a wider safety factor was used near $\arcp=0$ and near
$\arcp_1$, where curve $\mathsf C$'s numerical evaluation is closer to
degenerate), together with the same $0.1$--$0.2\%$ safety strip against
each box's own floor and ceiling. As with \cref{thm:arc2-band-new}, we
expect this margin is reducible with a finer grid but have not pursued
that here, and report the certificate exactly at the resolution
actually run.

\subsection{The band's \texorpdfstring{$\mathsf B$}{B}-bounded sub-part}\label{sec:arc2-bandB-new}

One sub-part of the band, bounded above by $\mathsf W'$ and below by
curve $\mathsf B$ for $\arcp\in(\arcp_3,\pi/4)$, does not yield to the
same method. Its combinatorial data behave exactly as elsewhere: the
representative point has $34$ vertices, $11$ moving and $23$ fixed, and
the pattern is identical at eight sample points spanning the
sub-interval, with fixed volume contribution $6$. The moving-volume
step is where it stops. One of the $74$ moving boundary simplices has a
symbolic determinant whose exact simplification over the relevant
algebraic extension grows without settling, and we did not obtain a
closed form for the volume there. The certificate for this sub-part
therefore does not come from the direct computation.

It does not need to. \Cref{sec:full-triangle-symmetry} shows that the
second arc, this sub-part included, is the first arc's own defect
function in disguise, so the certificates of \cref{sec:arc1-cert}
already cover it. The three direct certificates recorded above remain
useful as an independent check of that identification on part of the
domain, computed by a route that shares nothing with it.

\subsection{Coverage of the boundary}\label{sec:arc2-coverage}

\begin{figure}[tb]
\centering
\includegraphics[width=\textwidth]{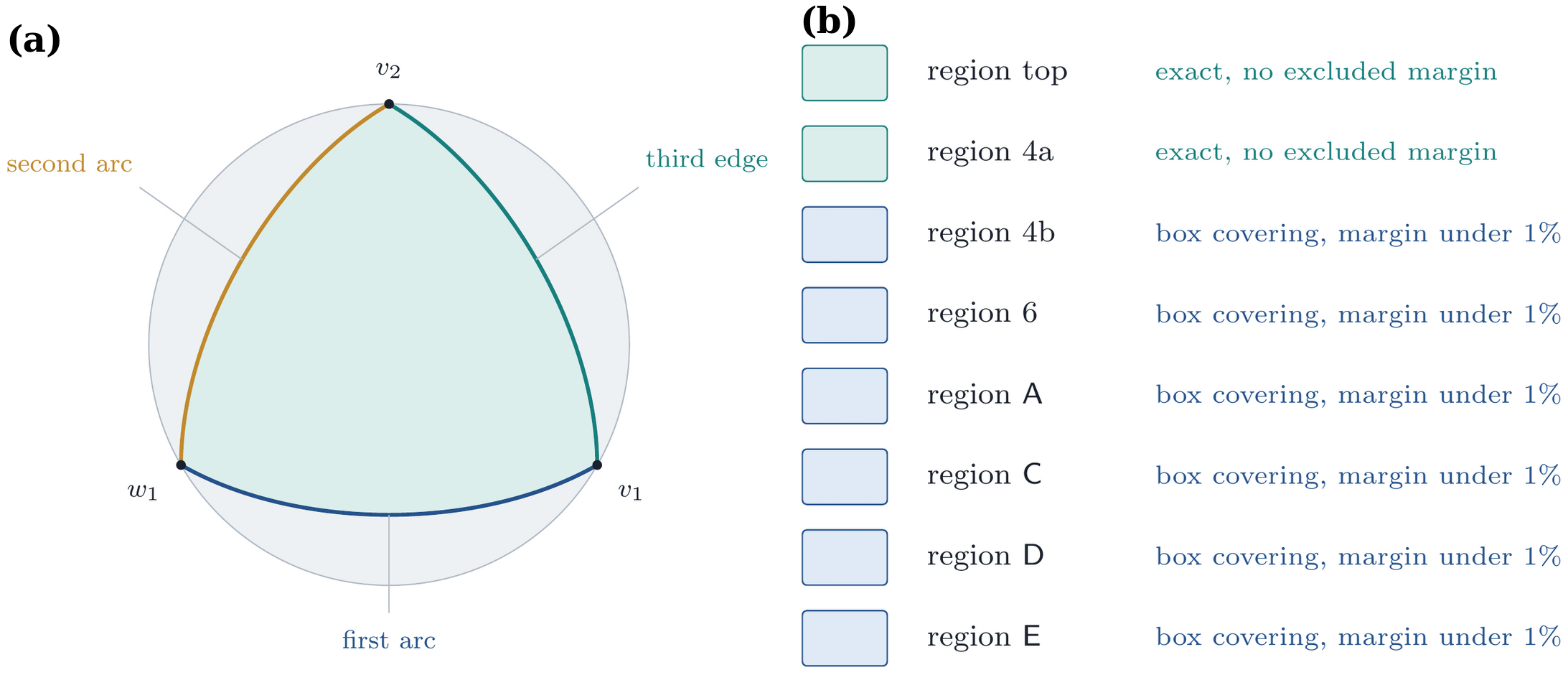}
\caption{(a) The transverse sphere of deviation directions and the
    spherical octant left by the two symmetries of
    \cref{lem:two-symmetries-new}. The vertices $v_1,w_1,v_2$ are
    pairwise orthogonal, so the octant's three edges are quarter great
    circles, and $\langle S,H\rangle\cong S_3$ fixes $u_0$ and permutes
    them; a certificate on the first arc is therefore a certificate on
    the whole boundary of the octant. (b) The eight combinatorial regions of the first arc and the kind
    of certificate each carries. Two are exact with no excluded margin;
    the other six are certified by box covering, with the excluded margin
    in each case a few tenths of one percent of the region's own area.
    Through the identifications of
    \cref{prop:arc2-is-arc1,prop:third-edge-cert}, these eight
    certificates cover the whole boundary of the octant.}
\label{fig:certmap}
\end{figure}

\Cref{fig:certmap}(a) shows the octant, and \cref{fig:certmap}(b) the
certificate carried by each region of the first arc.

Read directly, the region-by-region computations of this section settle
three of the second arc's units and leave the rest open. Read through
\cref{prop:arc2-is-arc1}, they are not needed at all: the second arc and
the third edge both carry the first arc's certificates, at the first
arc's margins. The interior of the octant is a separate matter, settled
by a different argument in \cref{sec:cap-theorem}.
\section{Positivity for every deviation direction}\label{sec:cap-theorem}

The certificates of \cref{sec:arc1-cert,sec:arc2-cert} are exact, but they
are local: each one governs a single combinatorial region of a single
one-dimensional arc. This section proves the positivity of
$\defect(\tilt,\dev)$ at every point of the full two-sphere of deviation
directions and at every tilt, by an argument of a different kind. The
defect is first rewritten as the deficiency of a cap cut from one fixed
polytope; that polytope is then enclosed in a cross-polytope, for which
the extremal cap can be determined outright.

\subsection{The defect as a cap deficiency}\label{sec:cap-reformulation}

Throughout, $\base_0=\alpha_0/\sqrt2$ is the undeviated contact direction,
$\base_1=\cos\tilt\,\base_0+\sin\tilt\,\dev$ the deviated one, and
\[
  Q \;=\; \bigcap_{i\ne0}\bigl\{x\in\mathbb{R}^4 :
          \langle x,\base_i\rangle\le1\bigr\}
\]
the intersection of the half-spaces of the twenty-three contact directions
that stay at roots. The local cell at the all-contact corner is
$Q\cap\{\langle x,\base_1\rangle\le1\}$, so
\begin{equation}\label{eq:defect-is-vol}
  \defect(\tilt,\dev)
  \;=\;\vol\bigl(Q\cap\{\langle x,\base_1\rangle\le1\}\bigr)-8 .
\end{equation}
Write $p=\langle x,\base_0\rangle$, and let $N$ denote the eight roots
making an angle of $60^\circ$ with $\alpha_0$, so that
$\langle\base_0,\base_i\rangle=\tfrac12$ for $i\in N$.

\begin{lem}[Structure of $Q$]\label{lem:Q-structure}
Let $F\subset\base_0^{\perp}$ be the octahedral facet of $\VCell$ at
$\base_0$, translated to the origin, so that the facet is $\base_0+F$.
Then $F$ is a regular octahedron whose six vertices are $\pm v_1,\pm
v_2,\pm v_3$, where $v_1,v_2,v_3$ are three pairwise orthogonal contact
directions orthogonal to $\base_0$; in particular
$(\base_0,v_1,v_2,v_3)$ is an orthonormal basis of $\mathbb{R}^4$,
$\vol_3(F)=\tfrac43$, and
\[
  F=\Bigl\{y\in\base_0^{\perp} : \textstyle\sum_{i}|\langle y,v_i\rangle|\le1\Bigr\} .
\]
Moreover $Q\cap\{p\ge1\}$ is the pyramid $\mathrm{Pyr}$ with apex
$2\base_0$ and base $\base_0+F$, of volume $\tfrac13$, while
$Q\cap\{p\le1\}=\VCell$. Hence $\vol(Q)=8+\tfrac13=\tfrac{25}{3}$.
\end{lem}

\begin{proof}
Six roots are orthogonal to $\alpha_0$, and they fall into three
orthogonal pairs $\pm\alpha_{i_1},\pm\alpha_{i_2},\pm\alpha_{i_3}$; put
$v_j=\alpha_{i_j}/\sqrt2$. Each $\base_0\pm v_j$ satisfies $\langle
\base_0\pm v_j,\base_i\rangle\le1$ for every root and meets six of those
constraints with equality, so the six points $\base_0\pm v_j$ are
vertices of $\VCell$ lying on the facet at $\base_0$; since a facet of the
$24$-cell is a regular octahedron with six vertices, these are all of
them. A regular octahedron with vertices $\pm v_j$ in an orthonormal frame
is $\{y:\sum_j|\langle y,v_j\rangle|\le1\}$ and has volume
$\tfrac{2^3}{3!}\cdot 2=\tfrac43$.

For $i\in N$ set $m_i=\base_i-\tfrac12\base_0$, so $m_i\perp\base_0$ and
$|m_i|=\tfrac{\sqrt3}{2}$. Writing $x=p\base_0+y$ with $y\perp\base_0$,
the constraint $\langle x,\base_i\rangle\le1$ reads $\langle
y,m_i\rangle\le1-\tfrac p2$. At $p=1$ these eight constraints cut out
exactly $F$, so at height $p$ they cut out $(2-p)F$; that is,
\begin{equation}\label{eq:cone-slice}
  \bigcap_{i\in N}\{\langle x,\base_i\rangle\le1\}
  \;=\;\bigl\{p\base_0+y \;:\; y\in(2-p)F\bigr\}.
\end{equation}
Let $1\le p\le2$ and $y\in(2-p)F\subseteq F$. The other fifteen
constraints are then slack: for the root $-\alpha_0$ one has $\langle
x,-\base_0\rangle=-p\le-1$; for a root $\alpha_j$ orthogonal to $\alpha_0$,
$\langle x,\base_j\rangle=\langle y,\pm v_k\rangle\le1$ because
$\max_{y\in F}\langle y,v_k\rangle=1$; and for a root $\alpha_j$ with
$\langle\base_0,\base_j\rangle=-\tfrac12$, necessarily $\base_j=-\base_i$
for some $i\in N$, whence $\langle x,\base_j\rangle=-\tfrac
p2+\langle y,-m_i\rangle\le-\tfrac12+\tfrac12=0$. So on $\{p\ge1\}$ the
set $Q$ agrees with \cref{eq:cone-slice}, which there is the cone over
$\base_0+F$ with apex $2\base_0$, of volume
$\int_1^2\vol_3\bigl((2-p)F\bigr)\,dp=\tfrac43\int_0^1\lambda^3d\lambda
=\tfrac13$. Finally $\VCell\subseteq\{p\le1\}$, and below that hyperplane
$Q$ and $\VCell$ are cut out by the same constraints apart from the one at
$\alpha_0$, which is inactive there.
\end{proof}

\begin{prop}[Cap form of the defect]\label{prop:cap-reformulation}
For a unit vector $u$ put $\capQ(u)=\{x\in Q:\langle x,u\rangle\ge1\}$.
Then
\[
  \defect(\tilt,\dev)\;=\;\tfrac13-\vol\bigl(\capQ(\base_1)\bigr).
\]
\end{prop}

\begin{proof}
By \cref{lem:Q-structure}, $\vol(Q)=\tfrac{25}{3}$, and
\cref{eq:defect-is-vol} gives $\defect=\vol(Q)-\vol(\capQ(\base_1))-8$.
\end{proof}

The deviation parameters have now left the body and entered the cutting
half-space: $Q$ is one fixed polytope with $25$ vertices, and the
question is how much of it a half-space at unit distance from the origin
can cut off. At $\base_1=\base_0$ the answer is $\tfrac13$, the pyramid.

\subsection{A cross-polytope around \texorpdfstring{$Q$}{Q}}\label{sec:crosspolytope}

Let $(\base_0,v_1,v_2,v_3)$ be the orthonormal frame of
\cref{lem:Q-structure} and write $x_0=\langle x,\base_0\rangle$,
$x_j=\langle x,v_j\rangle$ for the coordinates in that frame.

\begin{lem}[Enclosure]\label{lem:crosspolytope}
$Q\subseteq B$, where
\[
  B=\bigl\{x\in\mathbb{R}^4 : |x_0|+|x_1|+|x_2|+|x_3|\le2\bigr\}
\]
is the cross-polytope of radius $2$ in the frame
$(\base_0,v_1,v_2,v_3)$.
\end{lem}

\begin{proof}
Sixteen roots make an angle of $60^\circ$ or $120^\circ$ with $\alpha_0$,
and they come in eight antipodal pairs. If $\alpha_i$ is one of the eight
at $60^\circ$ then $\base_i=\tfrac12\base_0+\tfrac{\sqrt3}{2}n_i$ with
$n_i$ a unit vector orthogonal to $\base_0$; since the eight vectors $n_i$
are the outer normals of the eight faces of the regular octahedron $F$,
they must be $(\pm v_1\pm v_2\pm v_3)/\sqrt3$, one for each choice of
signs.
Hence
\[
  \langle x,\base_i\rangle
  =\tfrac12\bigl(x_0\pm x_1\pm x_2\pm x_3\bigr),
\]
and the eight constraints $\langle x,\base_i\rangle\le1$ say exactly that
$x_0\pm x_1\pm x_2\pm x_3\le2$ for all eight sign patterns, that is,
$x_0+|x_1|+|x_2|+|x_3|\le2$. The eight roots at $120^\circ$ are the
negatives of these, and give $-x_0+|x_1|+|x_2|+|x_3|\le2$. All sixteen
roots are distinct from $\alpha_0$, so all sixteen constraints are among
those defining $Q$.
\end{proof}

\begin{figure}[tb]
\centering
\includegraphics[width=\textwidth]{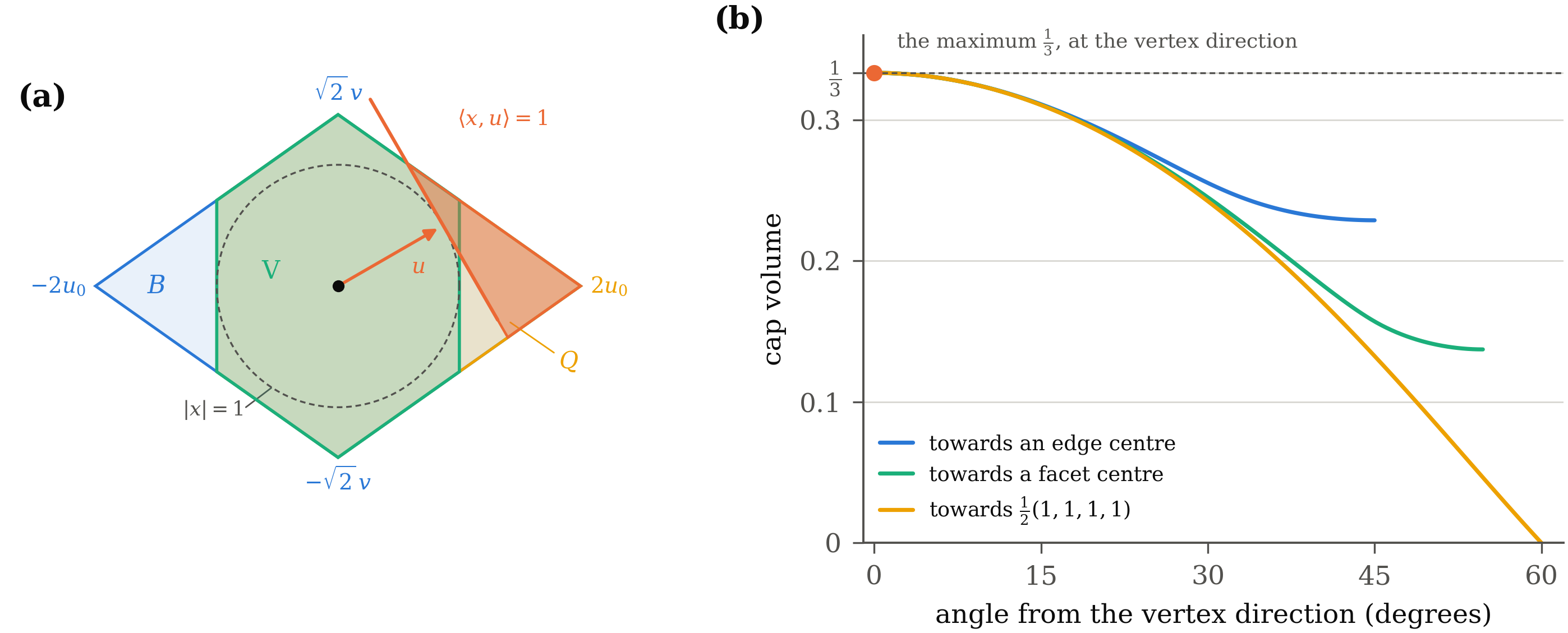}
\caption{(a) The section of the reference cell $\mathcal{V}$, of
    $Q$, and of the cross-polytope $B$ by the plane spanned by the base
    direction $\base_0$ and one transverse direction. Coordinates are
    exact: the section of $\mathcal{V}$ is a hexagon, that of $Q$ the
    same hexagon with its right edge replaced by the two edges meeting
    at the apex $2\base_0$, and that of $B$ a rhombus through the same
    apex and the same two vertices $\pm\sqrt2\,v$. The dashed circle is
    the unit sphere, on which every contact hyperplane is tangent, and
    the shaded triangle is the section of the cap cut by a half-space at
    unit distance whose normal is tilted $30^\circ$ from $\base_0$. (b) The volume of the cap cut from the cross-polytope by a
    half-space at unit distance, along three great circles leaving a
    vertex direction: towards an edge centre, towards a facet centre, and
    towards $\tfrac12(1,1,1,1)$, where the half-space meets $B$ only in
    its four positive vertices and the cap volume falls to zero. The
    maximum $\tfrac13$ of \cref{thm:cap-inequality} is attained at the
    vertex direction, marked.}
\label{fig:cap}
\end{figure}

The enclosure is efficient where it matters: $\vol(B)=\tfrac{32}{3}$
against $\vol(Q)=\tfrac{25}{3}$, and, more to the point, the half-space
$\{x_0\ge1\}$ cuts from $B$ the cone over the octahedron
$\{\sum_{j\ge1}|x_j|\le1\}$ with apex $2\base_0$, of volume exactly
$\tfrac13$. On the critical direction the enclosure is therefore not lossy
at all. \Cref{fig:cap}(a) shows the three bodies in section.

\subsection{The extremal cap of the cross-polytope}\label{sec:cap-inequality}

We now determine the largest volume that a half-space at unit distance
from the centre can cut from $B$. Write $B_1^4=\{\|x\|_1\le1\}$, so
$B=2B_1^4$, and let $e_0,e_1,e_2,e_3$ be the frame vectors.

\begin{lem}[Cap formula]\label{lem:cap-formula}
Let $u$ be a unit vector with coordinates $u_k$ in the frame, and set
$a_k=u_k^2$, so that $a_0+a_1+a_2+a_3=1$. Define
\[
  g(a)=\begin{cases} a\,(2\sqrt a-1)^4, & a\ge\tfrac14,\\[2pt] 0,& a\le\tfrac14 .\end{cases}
\]
Then
\[
  \vol\bigl(B\cap\{\langle x,u\rangle\ge1\}\bigr)
  \;=\;\tfrac13\,g[a_0,a_1,a_2,a_3],
\]
where $g[\,\cdot\,]$ denotes the third-order divided difference of $g$ at
the four nodes.%
\footnote{We use nothing about divided differences beyond the following,
  so we record it here instead of sending the reader elsewhere. For a
  function $h$ and pairwise distinct nodes $y_0,\dots,y_n$, the divided
  difference is defined by $h[y_0]=h(y_0)$ and
  \[
    h[y_0,\dots,y_n]
    =\frac{h[y_1,\dots,y_n]-h[y_0,\dots,y_{n-1}]}{y_n-y_0},
  \]
  and unwinding the recursion gives the closed form
  $h[y_0,\dots,y_n]=\sum_{i}h(y_i)\big/\prod_{j\ne i}(y_i-y_j)$, from
  which it is visibly a symmetric function of the nodes. Three
  consequences are used below. First, if $h$ is $C^n$ then
  $h[y_0,\dots,y_n]=h^{(n)}(\xi)/n!$ for some $\xi$ in the smallest
  interval containing the nodes; in particular the divided difference
  extends continuously to coincident nodes, with
  $h[y,\dots,y]=h^{(n)}(y)/n!$, so all the formulas below may be read at
  coincident nodes by continuity. Second, if $h^{(n-1)}$ is
  nondecreasing then $h[y_0,\dots,y_n]\ge0$ for all nodes, immediate by
  the mean-value form where $h$ is $C^n$ and in general by uniform
  approximation, which is the case we need, our $g$ being $C^3$ with
  $g'''$ nondecreasing but not $C^4$ at $a=\tfrac14$. Third, adding a
  node differentiates:
  $\partial h[y_0,\dots,y_n]/\partial y_i=h[y_0,\dots,y_n,y_i]$. These
  three facts turn a statement about the fourth derivative of one
  variable into a monotonicity statement in four.}
\end{lem}

\begin{proof}
Both sides are unchanged by a signed permutation of the coordinates, so we
may assume $u_k\ge0$ for all $k$. Split $B$ into the sixteen simplices
$S_\sigma=\conv\{0,2\sigma_0e_0,\dots,2\sigma_3e_3\}$,
$\sigma\in\{\pm1\}^4$, each of volume $\tfrac{2^4}{4!}=\tfrac23$. For a
$4$-simplex $\Delta$ with vertices $p_0,\dots,p_4$ at which a linear form
$f$ takes distinct values $f_0,\dots,f_4$, one has
\begin{equation}\label{eq:simplex-slice}
  \vol\bigl(\Delta\cap\{f\ge s\}\bigr)
  =\vol(\Delta)\sum_{i=0}^{4}\frac{(f_i-s)_+^4}{\prod_{j\ne i}(f_i-f_j)} .
\end{equation}%
\footnote{\Cref{eq:simplex-slice} is the only fact about simplices the
  proof uses, and it has a short derivation worth recording. Put
  $\Delta=\{\lambda\in\mathbb{R}^{5}:\lambda_i\ge0,\ \sum\lambda_i=1\}$
  in barycentric coordinates, so that $f=\sum_i f_i\lambda_i$ and the
  uniform probability measure on $\Delta$ pushes forward under $f$ to
  the distribution of $\sum_i f_i\lambda_i$ for $\lambda$ uniform on the
  simplex. The coordinates $\lambda_i$ are the normalised spacings of
  four independent uniform points on an interval, and the density of
  $\sum f_i\lambda_i$ is $4!$ times the third-order divided difference
  of $t\mapsto(t-x)_+^{3}/3!$ in the nodes $f_i$, which is the
  Curry--Schoenberg B-spline with those knots \cite{CSc66}. Integrating
  that density from $s$ upwards raises the truncated power by one degree
  and produces exactly the displayed sum. Equivalently, and with no
  probability at all: both sides of \cref{eq:simplex-slice} are, as
  functions of $s$, piecewise polynomials of degree $4$ that vanish for
  $s\ge\max_if_i$ and agree with $\vol(\Delta)$ for $s\le\min_if_i$,
  with matching derivatives at every breakpoint, so they agree
  identically. For $d$-simplices the exponent $4$ is replaced by $d$
  throughout.}
On $S_\sigma$ the form $\langle\cdot,u\rangle$ takes the value $0$ at the
origin and $2\sigma_ku_k$ at $2\sigma_ke_k$. With $s=1$ the vertex $0$
contributes nothing, and a vertex $2\sigma_ke_k$ contributes only when
$\sigma_k=+1$ and $2u_k>1$. Summing the contribution of the $k$-th vertex
over the eight sign choices of the remaining coordinates gives
\[
  \sum_{\sigma_l=\pm1,\;l\ne k}\ \prod_{l\ne k}\frac{1}{2u_k-2\sigma_lu_l}
  =\prod_{l\ne k}\Bigl(\frac{1}{2u_k-2u_l}+\frac{1}{2u_k+2u_l}\Bigr)
  =\prod_{l\ne k}\frac{u_k}{u_k^2-u_l^2},
\]
so that
\[
  \vol\bigl(B\cap\{\langle x,u\rangle\ge1\}\bigr)
  =\frac23\sum_{k}\frac{(2u_k-1)_+^4}{2u_k}\prod_{l\ne k}\frac{u_k}{u_k^2-u_l^2}
  =\frac13\sum_{k}\frac{(2u_k-1)_+^4\,u_k^2}{\prod_{l\ne k}(u_k^2-u_l^2)} .
\]
The $k$-th summand is $g(a_k)/\prod_{l\ne k}(a_k-a_l)$, and the sum is the
divided difference. Coincident nodes are handled by continuity, both sides
being continuous in $u$.
\end{proof}

The shape of the answer is not an accident. A divided difference of a
truncated power is the classical representation of a B-spline, and the
volume of a slice of a simplex is exactly such a quantity; the
identification goes back to Curry and Schoenberg \cite{CSc66}, and it is
what makes the next lemma available.

Two properties of $g$ drive everything that follows. Writing $t=\sqrt a$
for $a\ge\tfrac14$, repeated differentiation gives
\begin{equation}\label{eq:g-derivs}
  \begin{gathered}
  g'=(2t-1)^3(6t-1),\qquad
  g''=\frac{6(2t-1)^2(4t-1)}{t},\\
  g'''=96-\frac{60}{t}+\frac{3}{t^3},
  \end{gathered}
\end{equation}
and $g''''=\dfrac{3(20t^2-3)}{2t^5}$.

\begin{lem}[Monotonicity]\label{lem:monotone}
$g$ is of class $C^3$ on $[0,\infty)$, its third derivative is
nondecreasing, and consequently $g[a_0,a_1,a_2,a_3]$ is nondecreasing in
each node.
\end{lem}

\begin{proof}
All three expressions in \cref{eq:g-derivs} vanish at $t=\tfrac12$, so
$g,g',g'',g'''$ all vanish at $a=\tfrac14$ and match the identically zero
branch; hence $g\in C^3$. For $t\ge\tfrac12$ one has $20t^2-3\ge2>0$, so
$g''''>0$ there, while $g'''\equiv0$ on $[0,\tfrac14]$; therefore $g'''$ is
nondecreasing on $[0,\infty)$. A function whose $(n-1)$st derivative is
nondecreasing has all $n$th-order divided differences nonnegative, so
$g[y_0,\dots,y_4]\ge0$ for any five nodes.%
\footnote{The stronger conclusion one would like here is that
  $g[a_0,a_1,a_2,a_3]$ is Schur-convex in $(a_0,\dots,a_3)$, which would
  reduce \cref{thm:cap-inequality} to evaluating $F$ at the majorisation
  maximum $(1,0,0,0)$ in one step; see \cite[Ch.~3]{MOA11} for the
  criterion. It is not available. Schur-convexity of a symmetric
  divided difference of this form would follow from $g^{(5)}\ge0$, and
  $g^{(5)}(t^2)=\tfrac{45}{4}t^{-7}(1-4t^2)$ is negative at every
  $t>\tfrac12$, which is the whole interior of the support. The one
  place the fifth derivative is positive is the point mass at
  $a=\tfrac14$, where $g''''$ jumps from $0$ to $96$. Monotonicity in each node separately, which needs only
  $g''''\ge0$, is therefore the strongest tool of this kind that the
  function actually supports, and it is why the proof below splits into cases instead of reducing to a single extremal configuration.} Since
\[
  \frac{\partial}{\partial a_i}\,g[a_0,a_1,a_2,a_3]
  =g[a_0,a_1,a_2,a_3,a_i],
\]
the divided difference is nondecreasing in each node.
\end{proof}

\begin{thm}[Extremal cap of the cross-polytope]\label{thm:cap-inequality}
For every unit vector $u\in\mathbb{R}^4$,
\[
  \vol\bigl(B\cap\{\langle x,u\rangle\ge1\}\bigr)\;\le\;\tfrac13,
\]
with equality if and only if $u=\pm e_k$ for some $k$.
\end{thm}

\Cref{fig:cap}(b) shows the cap volume along three great circles
leaving a vertex direction, computed from the enumerated cell and
independent of the argument below.

\begin{proof}
By \cref{lem:cap-formula} the claim is that $F(a):=g[a_0,a_1,a_2,a_3]\le1$
on the simplex $a_k\ge0$, $\sum_ka_k=1$, with equality only at the
vertices. Let $S=\{k:a_k>\tfrac14\}$; since the nodes sum to $1$,
$|S|\le3$. We may relabel so that $a_0\ge a_1\ge a_2\ge a_3$.

\emph{Case $|S|=0$.} All four nodes lie in $[0,\tfrac14]$, where $g$
vanishes identically together with all its derivatives, so $F(a)=0$.

\emph{Case $|S|=1$.} Only $a_0$ exceeds $\tfrac14$, and again $g$ vanishes
to infinite order at the other three nodes, so the Newton form of the
divided difference collapses to a single term:
\[
  F(a)=\frac{g(a_0)}{(a_0-a_1)(a_0-a_2)(a_0-a_3)} .
\]
Write $c=\sqrt{a_0}\in(\tfrac12,1]$. The nodes $a_1,a_2,a_3$ range over
the polytope $P_c=\{\beta\in[0,\tfrac14]^3:\beta_1+\beta_2+\beta_3=1-c^2\}$,
and $\beta\mapsto\sum_l\log(a_0-\beta_l)$ is concave, so
$\prod_l(a_0-\beta_l)$ attains its minimum over $P_c$ at a vertex of
$P_c$, that is, at a point where at most one coordinate lies strictly
between $0$ and $\tfrac14$. Three ranges of $c$ occur.

If $c\ge\tfrac{\sqrt3}{2}$ the extreme point is $(1-c^2,0,0)$ and the
required inequality $g(a_0)\le\prod_l(a_0-\beta_l)$ reads
$(2c-1)^4\le c^2(2c^2-1)$, that is,
\[
  (1-c)\,\bigl(14c^3-18c^2+7c-1\bigr)\;\ge\;0 .
\]
The cubic has a single real root, at $c=0.7410\ldots<\tfrac{\sqrt3}{2}$,
and is positive beyond it, so the product is nonnegative on
$[\tfrac{\sqrt3}{2},1]$ and vanishes only at $c=1$.

If $\tfrac{1}{\sqrt2}\le c\le\tfrac{\sqrt3}{2}$ the extreme point is
$(\tfrac14,\tfrac34-c^2,0)$ and the inequality reads
$(2c-1)^4\le(c^2-\tfrac14)(2c^2-\tfrac34)$, that is,
$\tfrac{1}{16}(2c-1)\,q(c)\ge0$ with
\[
  q(c)=-112c^3+200c^2-102c+13 .
\]
Here $q'(c)=-2(168c^2-200c+51)$ has its roots at
$\tfrac{50\pm\sqrt{358}}{84}$, of which only
$\tfrac{50+\sqrt{358}}{84}=0.8205\ldots$ lies in
$[\tfrac1{\sqrt2},\tfrac{\sqrt3}{2}]$. So $q$ rises and then falls on
that interval and attains its minimum at an endpoint, where
$q(\tfrac1{\sqrt2})=113-79\sqrt2$ and $q(\tfrac{\sqrt3}{2})=163-93\sqrt3$.
Both are positive, since $113^2=12769>12482=2\cdot79^2$ and
$163^2=26569>25947=3\cdot93^2$.

If $\tfrac12\le c\le\tfrac{1}{\sqrt2}$ the extreme point is
$(\tfrac14,\tfrac14,\tfrac12-c^2)$, the required inequality becomes
$16(c-\tfrac12)^4c^2\le2(c-\tfrac12)^3(c+\tfrac12)^3$, that is,
$8(c-\tfrac12)c^2\le(c+\tfrac12)^3$, that is, $r(c)\ge0$, where
\[
  r(c)=8\Bigl((c+\tfrac12)^3-8(c-\tfrac12)c^2\Bigr)=-56c^3+44c^2+6c+1 .
\]
Again $r'(c)=-168c^2+88c+6$ has a single root,
$\tfrac{22+\sqrt{736}}{84}=0.5849\ldots$, in $[\tfrac12,\tfrac1{\sqrt2}]$,
so $r$ attains its minimum on that interval at an endpoint; $r(\tfrac12)=8$
and $r(\tfrac1{\sqrt2})=23-11\sqrt2>0$, the latter because
$23^2=529>242=2\cdot11^2$.

In all three ranges $F(a)\le1$, with equality only when $c=1$, that is
$a=(1,0,0,0)$ and $u=\pm e_0$.

\emph{Case $|S|\ge2$.} Here $a_0\ge a_1\ge\tfrac14$. By
\cref{lem:monotone}, $F$ is nondecreasing in each node, so on any box
$\prod_k[l_k,h_k]$ one has $F(a)\le F(h)$. The region
\[
  R=\Bigl\{a:\ a_0\ge a_1\ge a_2\ge a_3\ge0,\ \sum_ka_k=1,\ a_1\ge\tfrac14\Bigr\}
\]
is covered by $303$ boxes, obtained by bisecting a uniform initial grid of
$216$ cubes of side $\tfrac18$ and refining only where needed, and on each
box the value of $F$ at the upper corner is a rational number at most $1$;
the largest of the $303$ corner values is $0.99755\ldots$. The nodes at a corner
are taken at rational points $t_k$ with $t_k^2$ at least the corner
coordinate, and ties are broken upward, both of which only increase the
bound. The computation is exact rational arithmetic throughout, and is
reproduced by \texttt{cap\_inequality\_certificate.py} in the supplement.
Since $\sup_R F=0.76583\ldots$, attained at $(\tfrac34,\tfrac14,0,0)$, the
certificate has a margin of about a quarter and does not depend on the
fineness of the grid beyond what is stated.
\end{proof}

\subsection{The single-deviation theorem}\label{sec:single-dev-theorem}

\begin{thm}[Direction-of-deviation positivity]\label{thm:eperp-closed}
For every unit deviation direction $\dev$ orthogonal to $\base_0$ and
every tilt $\tilt\in[0,\pi]$,
\[
  \defect(\tilt,\dev)\;\ge\;0,
\]
with equality if and only if $\tilt=0$.
\end{thm}

\begin{proof}
By \cref{lem:crosspolytope}, $\capQ(\base_1)\subseteq
B\cap\{\langle x,\base_1\rangle\ge1\}$, so
\cref{thm:cap-inequality} gives $\vol(\capQ(\base_1))\le\tfrac13$, and
\cref{prop:cap-reformulation} turns this into $\defect\ge0$.

For the equality case, suppose $\defect(\tilt,\dev)=0$, so
$\vol(\capQ(\base_1))=\tfrac13$. Then the cap of $B$ also has volume
$\tfrac13$, so $\base_1=\pm e_k$ by \cref{thm:cap-inequality}. If
$\base_1=-\base_0$ then $\capQ(\base_1)=\{x\in Q:\langle
x,\base_0\rangle\le-1\}$, which is empty because $-\alpha_0$ is a root and
$Q$ therefore lies in $\{\langle x,\base_0\rangle\ge-1\}$; if
$\base_1=\pm v_j$ then $\pm v_j$ is itself one of the contact directions
defining $Q$, so again the cap is empty. In each of these cases
$\defect=\tfrac13$, not $0$. The remaining possibility is
$\base_1=\base_0$, that is $\tilt=0$, where the cap is the pyramid of
\cref{lem:Q-structure} and $\defect=0$.
\end{proof}

\Cref{thm:local-uncond} can now be proved in the form stated in the
introduction.

\begin{cor}[Single-deviation cell bound]\label{cor:single-dev-uncond}
Let $c$ be a centre of a unit-ball packing of $\mathbb{R}^4$ at which at
most one active neighbour direction fails to be a root direction of $\Rt$,
with nothing assumed about that neighbour's deviation direction or about
the size of its deviation. Then $\vol(V_c)\ge8$, with equality if and only
if the active neighbours of $c$ form a copy of $\sqrt2\,\Rt$ centred at
$c$.
\end{cor}

\begin{proof}
By \cref{sec:shell-new,sec:radial-new} the volume is minimised at the
all-contact corner, and by \cref{thm:freevol-new} the deficit there is
$\defect(\tilt,\dev)$ for the appropriate parameters, which is
nonnegative by \cref{thm:eperp-closed}. If the active set is a proper
subset of $\Rt\setminus\{\alpha_0\}$, the cell contains
$Q\cap\{\langle x,\base_1\rangle\le1\}$ and the same bound applies a
fortiori. Equality forces $\tilt=0$ and a full active set.
\end{proof}

Two features of the argument deserve comment. First, it is insensitive to
which combinatorial cell of \cref{sec:chamber-register} the parameters
fall in: all $176$ of them are covered at once, because the cap
inequality is a statement about a single polytope and makes no reference
to the combinatorial type of the perturbed cell. Second, the enclosure
$Q\subseteq B$ discards a seventh of the volume of $Q$ and yet loses
nothing at the critical direction, which is what makes the estimate
sharp; the same argument applied to the $24$-cell in place of $Q$, or to
a smaller enclosing body, would not have this property.

\subsection{A quantitative form near the base direction}\label{sec:quantitative-cone}

\Cref{thm:eperp-closed} gives a sign and no more. For small and moderate
tilts one can say how large the defect is, by replacing $Q$ with a cone in place of a cross-polytope. The eight roots at $60^\circ$ alone
cut out
\[
  K \;=\; \bigl\{x : x_0+|x_1|+|x_2|+|x_3|\le2\bigr\},
\]
a cone with apex $2\base_0$, and $Q\subseteq B\subseteq K$. The cap of a
cone is a union of simplices, so it has a closed form.

\begin{prop}[Cone bound]\label{cor:second-order}
Let $\dev$ be a unit vector orthogonal to $\base_0$, with coordinates
$\dev_j=\langle\dev,v_j\rangle$ in the frame of \cref{lem:Q-structure},
and write $c=\cos\tilt$, $\sigma=\sin\tilt$. If
$\sigma\max_j|\dev_j|<c$, then
\begin{equation}\label{eq:cone-bound}
  \defect(\tilt,\dev)\;\ge\;
  \frac13-\frac{c^2\,(2c-1)_+^4}{3\prod_{j=1}^{3}\bigl(c^2-\sigma^2\dev_j^2\bigr)} .
\end{equation}
The right-hand side equals $\tilt^2/3+O(\tilt^4)$ as $\tilt\to0$, and it
is strictly positive for every such $\dev$ if and only if $c>c_0$, where
$c_0=0.7410861\ldots$ is the unique real root of
$14c^3-18c^2+7c-1$; that is, for $\tilt<42.17598\ldots^\circ$.
\end{prop}

\begin{proof}
Every root at $60^\circ$ to $\alpha_0$ is one of the twenty-three, so
$Q\subseteq K$ by the first half of the proof of
\cref{lem:crosspolytope}, and $\vol(\capQ(\base_1))\le\vol(K\cap\{\langle
x,\base_1\rangle\ge1\})$.

Put $w=2\base_0-x$. Then $K$ becomes the cone
$C=\{w:w_0\ge|w_1|+|w_2|+|w_3|\}$ with apex at the origin, and the
half-space $\{\langle x,\base_1\rangle\ge1\}$ becomes $\{\langle
w,\base_1\rangle\le2c-1\}$, since $\langle
2\base_0,\base_1\rangle=2c$. If $2c-1\le0$ the intersection is empty,
because $\langle w,\base_1\rangle\ge w_0\,(c-\sigma\max_j|\dev_j|)\ge0$
on $C$ under the stated hypothesis, and then the cap is empty and
$\defect=\tfrac13$, which is \eqref{eq:cone-bound}. Assume $2c-1>0$.

Intersecting $C$ with an orthant $\{\sigma_jw_j\ge0\}$, $\sigma\in\{\pm
1\}^3$, gives the simplicial cone generated by $g_0=e_0$ and
$g_j=e_0+\sigma_je_j$, and $\det(g_0,g_1,g_2,g_3)=\sigma_1\sigma_2\sigma_3$.
Truncating it by $\{\langle w,\base_1\rangle\le h\}$, $h=2c-1$, gives the
simplex with vertices $0$ and $h\,g_i/\langle g_i,\base_1\rangle$, of
volume
\[
  \frac{h^4}{24}\cdot\frac{1}{\langle g_0,\base_1\rangle
  \prod_{j=1}^{3}\langle g_j,\base_1\rangle},
  \qquad
  \langle g_0,\base_1\rangle=c,\quad
  \langle g_j,\base_1\rangle=c+\sigma_j\sigma\dev_j .
\]
Summing over the eight sign patterns factorises:
\[
\begin{aligned}
  \vol\bigl(C\cap\{\langle w,\base_1\rangle\le h\}\bigr)
  &=\frac{h^4}{24c}\prod_{j=1}^{3}
   \Bigl(\frac{1}{c+\sigma\dev_j}+\frac{1}{c-\sigma\dev_j}\Bigr)\\
  &=\frac{h^4}{24c}\prod_{j=1}^{3}\frac{2c}{c^2-\sigma^2\dev_j^2}
  =\frac{c^2h^4}{3\prod_j(c^2-\sigma^2\dev_j^2)} .
\end{aligned}
\]
With $h=2c-1$ this is the subtracted term in \eqref{eq:cone-bound}, and
\cref{prop:cap-reformulation} converts the cap bound into the stated
bound on $\defect$.

For the expansion, $c=1-\tilt^2/2+O(\tilt^4)$ gives $2c-1=1-\tilt^2+O(\tilt^4)$
and $c^2=1-\tilt^2+O(\tilt^4)$, so the numerator $c^2(2c-1)^4$ is
$1-5\tilt^2+O(\tilt^4)$. Since $\sigma^2=\tilt^2+O(\tilt^4)$ and
$\sum_j\dev_j^2=1$,
\[
  \prod_{j=1}^{3}\bigl(c^2-\sigma^2\dev_j^2\bigr)
  =\prod_{j=1}^{3}\bigl(1-\tilt^2-\tilt^2\dev_j^2\bigr)+O(\tilt^4)
  =1-4\tilt^2+O(\tilt^4),
\]
so the quotient is
$\tfrac13(1-5\tilt^2)(1+4\tilt^2)+O(\tilt^4)=\tfrac13-\tfrac{\tilt^2}{3}+O(\tilt^4)$,
and the right-hand side of \eqref{eq:cone-bound} is
$\tilt^2/3+O(\tilt^4)$.

Finally, the bound is positive exactly when
$c^2(2c-1)^4<\prod_j(c^2-\sigma^2\dev_j^2)$. Writing $a_j=\dev_j^2$, the
constraint set is the simplex $\{a\ge0,\ \sum a_j=1\}$, and
$a\mapsto\sum_j\log(c^2-\sigma^2a_j)$ is a sum of concave functions,
hence concave, so its minimum over the simplex is at a vertex. The
worst case is therefore $\dev=\pm v_j$, where the product is
$(c^2-\sigma^2)c^4=(2c^2-1)c^4$, and the condition becomes
$(2c-1)^4<c^2(2c^2-1)$. Expanding,
\[
  c^2(2c^2-1)-(2c-1)^4=(1-c)\bigl(14c^3-18c^2+7c-1\bigr),
\]
which for $c\in(0,1)$ is positive exactly when $14c^3-18c^2+7c-1>0$. The
cubic has one real root $c_0=0.7410861\ldots$ and is negative below it and
positive above, so the bound is positive for all $\dev$ precisely when
$c>c_0$.
\end{proof}

\begin{figure}[tb]
\centering
\includegraphics[width=\textwidth]{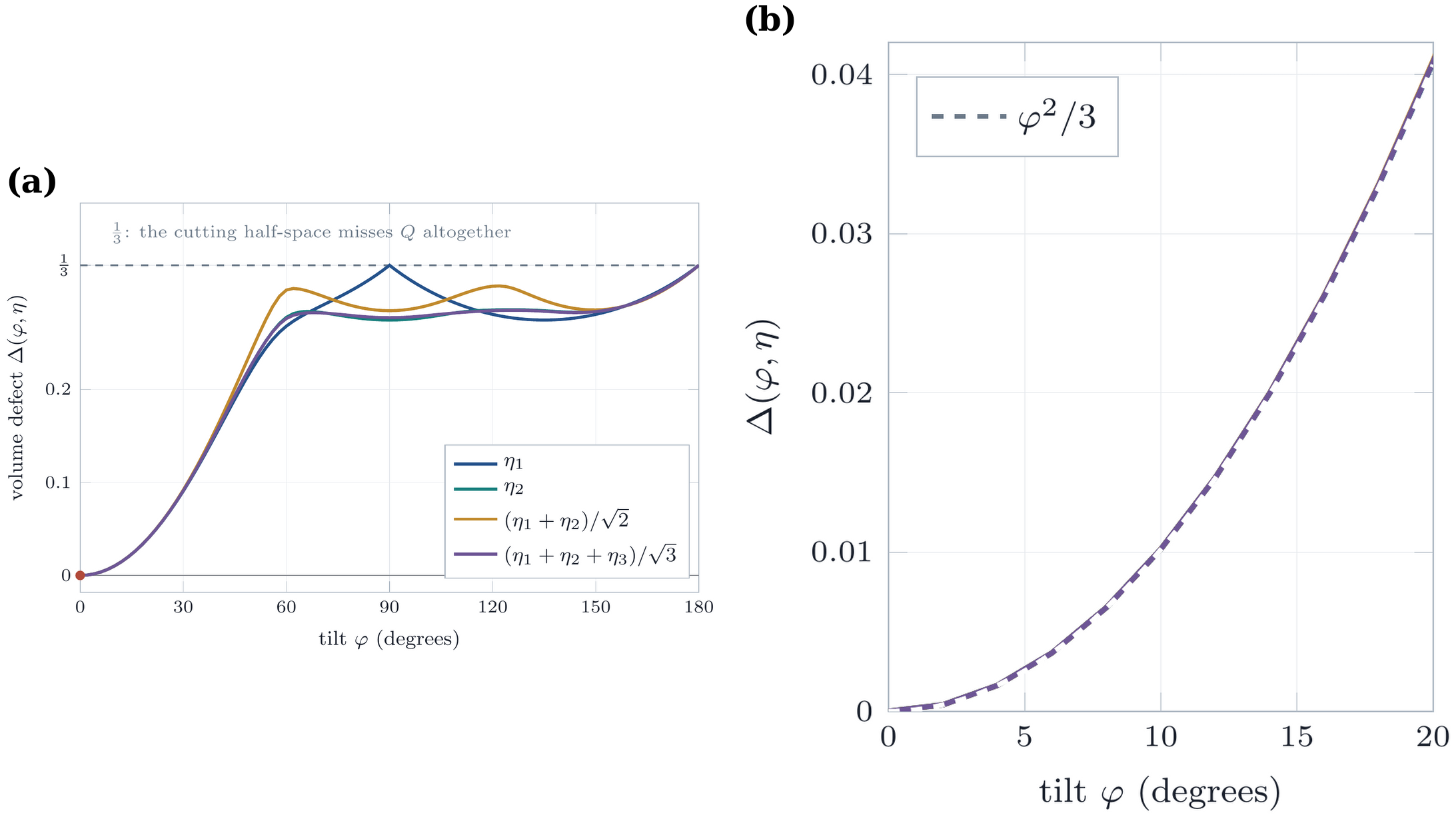}
\caption{(a) The volume defect $\defect(\varphi,\eta)$ over the whole
    range of tilts, for four deviation directions spanning the
    transverse sphere: two of the orthonormal frame directions and two
    diagonal combinations of them. The defect is nonnegative throughout
    and vanishes only at $\varphi=0$, as \cref{thm:eperp-closed}
    asserts, and it returns to the value $\tfrac13$ at $\varphi=\pi$,
    where the cutting half-space misses $Q$ entirely and the cell is
    exactly $Q$ itself. The intermediate maxima and minima are where the
    breakpoint structure of \cref{sec:arc1-breakpoints} shows itself.
    The curves are computed from directly enumerated cell volumes, by a
    route independent of the proof. (b) The same four curves at small tilt, against the second-order
    rate $\varphi^2/3$ of \cref{cor:second-order}. The four are
    indistinguishable there, and from each other, at plotting accuracy.}
\label{fig:defect}
\end{figure}

The threshold $c_0$ is the same cubic that governs the first case of
\cref{thm:cap-inequality}. That is not a coincidence: both come from
comparing $(2c-1)^4$ with $c^2(2c^2-1)$, once for the cone and once for
the divided difference at a single node. Beyond $42.17598\ldots^\circ$ the
cone is too generous an enclosure and \eqref{eq:cone-bound} goes
negative, which is exactly where the cross-polytope argument is needed.
\Cref{fig:defect}(a) plots the defect over the whole range and
\cref{fig:defect}(b) against that second-order rate.

\begin{remark}\label{rem:why-arcs-survive}
The certificates of \cref{sec:arc1-cert,sec:arc2-cert} are not made
redundant by \cref{thm:eperp-closed}. They locate the breakpoints of the
combinatorial structure along the two arcs exactly, give the defect as an
explicit rational function on each region, and quantify the margin there;
\cref{thm:eperp-closed} gives the sign and nothing more. The two
descriptions were also found in that order, the regional analysis having
suggested where to look for a global bound.
\end{remark}
\section{What the single-deviation case contains}\label{sec:deviation-scope-summary}

The single-deviation problem has been approached twice in this paper,
and the two approaches return different things. Only one of them proves
the theorem; only the other one gives the numbers.

The regional route occupies \cref{sec:arc1-breakpoints} through
\cref{sec:arc2-cert}. On the first boundary arc of the fundamental
deviation triangle, the breakpoint classification is complete: eight
combinatorial regions, each with its own defect polynomial in the
Weierstrass variables, seven certified exactly and with no margin, the
eighth with a thin but explicitly quantified margin arising from an
elliptic obstruction in the region's own boundary curve
(\cref{thm:arc1-breakpoints-new,thm:arc1-certs-new}). The second arc
was first attacked directly, three of its region-and-sub-part units
being certified and a fourth left incomplete, and then shown to be
unnecessary: the group generated by $S$ and $H$ is isomorphic to $S_3$,
fixes $\base_0$, realises every permutation of $\{v_1,w_1,v_2\}$, and
carries both the second arc and the third edge onto the first arc
exactly as functions of $(\tilt,\arcp)$
(\cref{lem:s3-symmetry,prop:third-edge-cert,prop:arc2-is-arc1}). The
whole boundary of the triangle therefore inherits the first arc's
certificate together with its margins
(\cref{cor:restated-conjecture-holds}).

That route cannot reach the interior. The symmetry group does not move
an interior point to an edge, the natural convexity reduction is false
(\cref{sec:interior-structure}), and the interior carries at least four
distinct combinatorial region types against the boundary's one, so a
regional treatment there would need four new volume formulas before any
certificate could be attempted.

The global route is \cref{sec:cap-theorem}. It discards the triangle,
the arcs, and the combinatorial cell register alike, replacing all of
them by a single question about one fixed polytope: how much can a
half-space at unit distance from the centre cut off? The answer,
$\tfrac13$, is exactly the volume of the pyramid that the undeviated
configuration already cuts off, which is why the resulting bound is
sharp and why equality pins down the $D_4$ configuration
(\cref{thm:eperp-closed}). The argument is insensitive to the tilt, to
the deviation direction, and to which of the $176$ cells the parameters
land in.

Neither route subsumes the other in practice. \Cref{thm:eperp-closed}
gives the sign everywhere and nothing else; it does not produce the
defect as a function, and it cannot say where along an arc the defect is
small or how the breakpoints sit. The regional certificates give exactly
that, region by region, and they are the reason we knew what the answer
had to look like before we could prove it. Between the two sits
\cref{cor:second-order}, which comes from the same reformulation as the
global route but replaces the cross-polytope by a cone: it gives an
explicit positive lower bound on the defect, of size $\tilt^2/3$ near
$\tilt=0$, on the range $\tilt<42.17598\ldots^\circ$, and says nothing
at all above it. A reader interested only in the theorem may read
\cref{sec:deviation-domain} and \cref{sec:cap-theorem} and skip what
lies between; a reader interested in the shape of the defect function
should not.
\section{The joint Hessian technique}\label{sec:hessian-technique}

\subsection{Simultaneous deviations}

\Cref{conj:multidir-new} concerns configurations in which $m\ge2$ active
directions $u_1,\dots,u_m$ deviate from their nearest roots
simultaneously. Restricting each deviation to be infinitesimal (a
second-order, local analogue of \stageone\ for the single-deviation
case), the joint behaviour is governed by a $3m\times3m$ Hessian, one
$3$-dimensional tangent block per deviating direction, with off-diagonal
cross-blocks determined by the Gram value between each pair of active
roots exactly as in \cref{def:hessian-new}.

\begin{definition}[Joint Hessian]\label{def:joint-hessian-new}
For an active set $A=\{k_1,\dots,k_m\}\subset\{1,\dots,24\}$ of pairwise
packing-valid root indices, the \emph{joint Hessian}
$H_A\in\mathbb{Q}^{3m\times3m}$ has diagonal blocks as in
\cref{def:hessian-new} and cross-block $(i,j)$ determined by the
Gram value $\langle\alpha_{k_i},\alpha_{k_j}\rangle/2\in\{-\tfrac12,0,\tfrac12\}$
via the same signed-adjacency rule.
\end{definition}

\Cref{conj:multidir-new} holds, to second order, at a configuration $A$
exactly when $H_A$ is positive semi-definite; a strict (not merely
second-order) resolution of the conjecture at that configuration would
require ruling out any \emph{finite}-angle failure as well, but the
Hessian is the natural first tool, and is what all computational work
on this conjecture to date, including in this paper, has used.%
\footnote{The distinction between second-order (Hessian) evidence and finite-angle proof matters here.  A positive-definite $H_A$ guarantees only
  that the volume-defect function $\defect_A$ increases away from the
  root-aligned configuration, it does not rule out a later decrease or
  a zero at some finite angle.  For this reason the Hessian technique, even
  when it succeeds completely (e.g., for chain active sets, \cref{lem:chains-new}),
  is a second-order certificate only.  The swap-configuration results of
  \cref{sec:swap-configs} are of a different character: they are exact,
  finite-angle proofs that $\Delta_A>0$ \emph{throughout} the packing-valid
  range, not just infinitesimally.}

\subsection{Small configurations: exact results}

\begin{lem}[Two-direction case]\label{lem:m2-exact-new}
For $m=2$ active directions with Gram value $\tfrac12$ (an adjacent
pair), the joint Hessian's minimum eigenvalue is exactly $\tfrac13$. For
the isolated $m=1$ case (a single deviating direction with no other
active directions considered jointly with it), $H_A=\tfrac23I_3$
exactly. This $m=1$ joint-Hessian construction is a different object
from the cell Hessian $H_k$ of \cref{def:hessian-new}: $H_k$ additionally
carries cross-terms from a cell's other active, non-deviating roots,
so the two need not (and, as the worked example of \cref{sec:hessian-def-new}
shows, generally do not) share the same numerical value.
\end{lem}

\subsection{A first step toward a first-principles derivation for \texorpdfstring{$m=2$}{m=2}}\label{sec:m2-firstprinciples-new}

The closed-form cross-term underlying \cref{lem:m2-exact-new} was obtained
by fitting the numerically observed cross-Hessian to a small family of natural
$W(D_4)$-covariant invariants and verifying the fit to within $10^{-6}$ at
several sample points and two independent adjacent pairs, not by a
derivation from the vertex structure of $\VCell$ itself, the way
\stageone's single-direction Hessian is derived from Cramer's rule
(\cref{sec:chamber-register}). We record here a first step toward closing
that gap, obtained during the preparation of this paper; it is a partial
structural result, not a completed derivation.

For an adjacent pair $\alpha,\beta$ (Gram value $\tfrac12$), we verified
directly (exact vertex enumeration over all $\binom{24}{4}$ candidate
facet-quadruples of $\VCell$) that their two corresponding facets, each an
octahedron with exactly $6$ vertices, share exactly $3$ of those
vertices, so the two facets meet along a triangular ridge. A simultaneous
perturbation of both facet normals therefore moves $9$ vertices in total
($6+6-3$), each of which the single-direction Cramer analysis of
\cref{sec:chamber-register} can in principle be applied to individually,
provided the correct governing set of active facets is identified at each
one under the joint perturbation.

That proviso is not automatic. Each of the $9$ affected vertices lies, in
the unperturbed cell, on strictly more than the generic $4$ active facets
of a $4$-dimensional vertex (specifically $6$, an artefact of the
$24$-cell's exceptional symmetry), so perturbing only one or two of those
facets can render one or more of the others redundant at the new vertex
position rather than leaving the vertex's governing facet set unchanged.
We checked this directly at the shared vertex $(\sqrt2,0,0,0)$ (active
facets $\{\alpha,\beta_1,\beta_2,\beta_3,\beta_4,\beta_5\}$ in the
unperturbed cell) under a perturbation of $\alpha$ alone, by testing every
$\binom53=10$ way of completing the perturbed facet to a determining
$4$-set using the vertex's other $5$ original active facets: exactly $4$
of the $10$ remain simultaneously feasible and agree on a single new
vertex position (with one previously-active facet becoming redundant
there), $1$ gives a singular (non-invertible) $4\times4$ system, and the
remaining $5$ give points violating some other facet's inequality. This
confirms the vertex deforms smoothly under this perturbation, and not combinatorially splitting into several nearby vertices, but it also
confirms that identifying the correct governing facet set is a real
case-by-case reduction, not a uniform rule read off in advance. Completing
this into a full first-principles derivation of \cref{lem:m2-exact-new}
requires repeating this reduction correctly at all $9$ affected vertices
and for a general (not merely one tested) perturbation direction, then
assembling the resulting exact volume-change formula and expanding it to
second order, a substantially larger undertaking than we complete here.
We record the partial result as a step taken, not a gap closed;
\cref{lem:m2-exact-new} itself continues to rest on the numerical
verification described above, not on this partial first-principles check.

\begin{lem}[Chains]\label{lem:chains-new}
For an induced path (a ``chain'') of $m$ pairwise-Gram-$\tfrac12$
directions with no other adjacency among them, the joint Hessian's
minimum eigenvalue, computed by Richardson extrapolation over three
step sizes for chain lengths $m=3,\dots,8$ (the true maximum induced
path length in the $24$-vertex, $8$-regular Gram-$\tfrac12$ adjacency
graph on $\Rt$, found by exhaustive search from every starting vertex),
decreases from approximately $0.276$ at $m=3$ to a floor plateauing
near $0.168$ for $m=6,7,8$, remaining strictly positive throughout.
\end{lem}

A geometric constraint underlies the plateau rather than a
continued collapse: at an interior chain vertex, the two tangent axes
pointing toward its left and right neighbours (within its own
$3$-dimensional orthogonal complement) are not independent, their
inner product is found to be exactly $-\tfrac13$ or $-1$ at every
interior vertex tested, i.e.\ partially or fully antipodal. This exact
algebraic constraint prevents the naive worst case (every pairwise
cross-term simultaneously at its individual extreme) from being jointly
realisable, and is a plausible mechanistic explanation for the observed
plateau, though we do not turn this observation into a proof.

\begin{figure}[tb]
\centering
\includegraphics[width=0.99\textwidth]{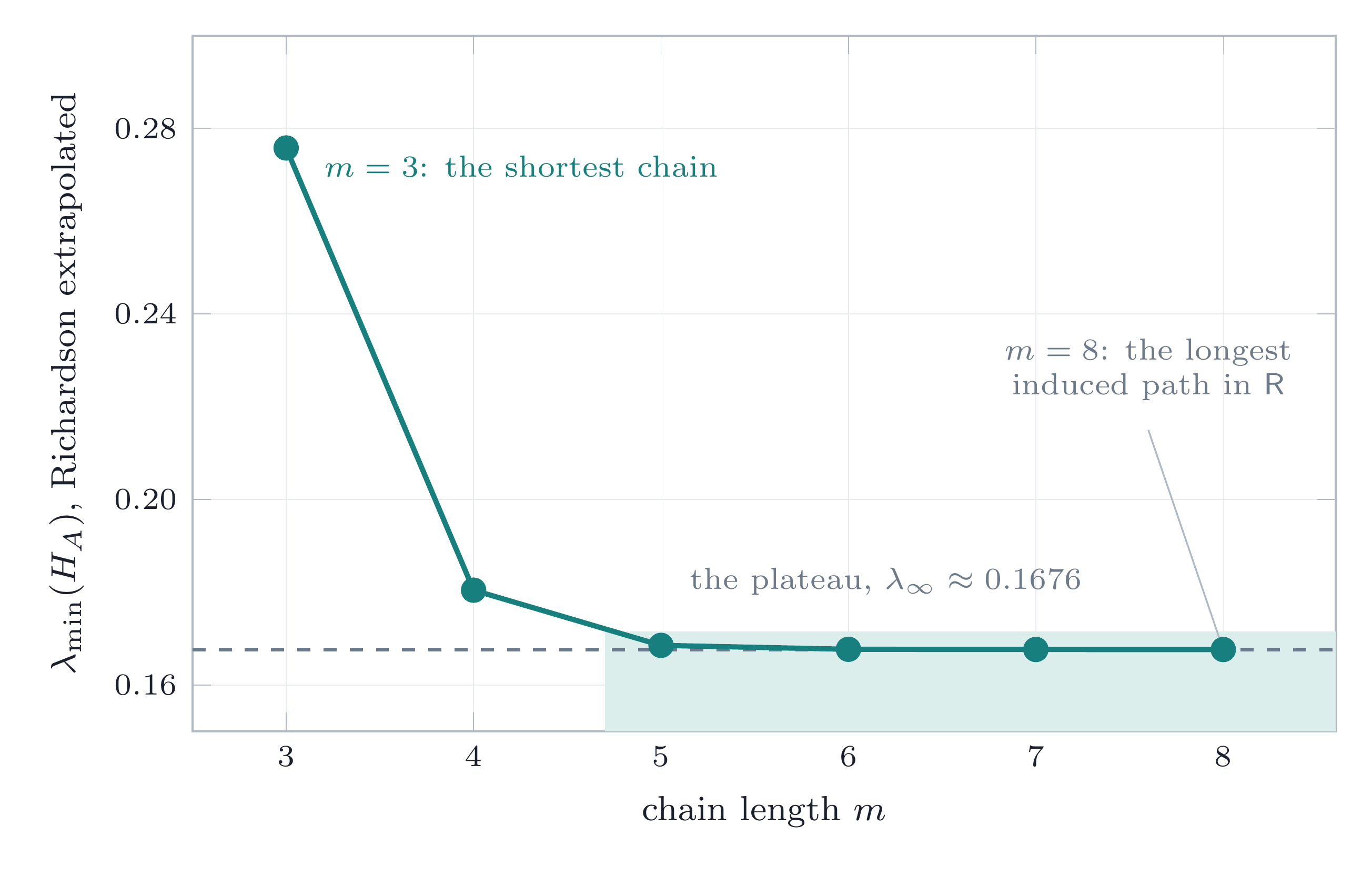}
\caption{Richardson-extrapolated minimum eigenvalue of the joint
Hessian $H_A$ for chain configurations of length $m=3,\dots,8$
(\cref{lem:chains-new}; full values in \cref{app:numerical-tables}).
The eigenvalue falls sharply from $m=3$ to $m=5$ and then plateaus near $0.1676$ instead of continuing to collapse, which is the
antipodal-axis mechanism described below. It stays well clear of zero
for every chain length up to $m=8$, the longest induced path in the
contact graph of $\Rt$, so the second-order method certifies positivity
on the whole chain family.}
\label{fig:chain-plateau}
\end{figure}

\subsection{Dense (non-chain) configurations: the floor does not plateau}\label{sec:dense-configs-new}

\begin{figure}[tb]
\centering
\includegraphics[width=0.99\textwidth]{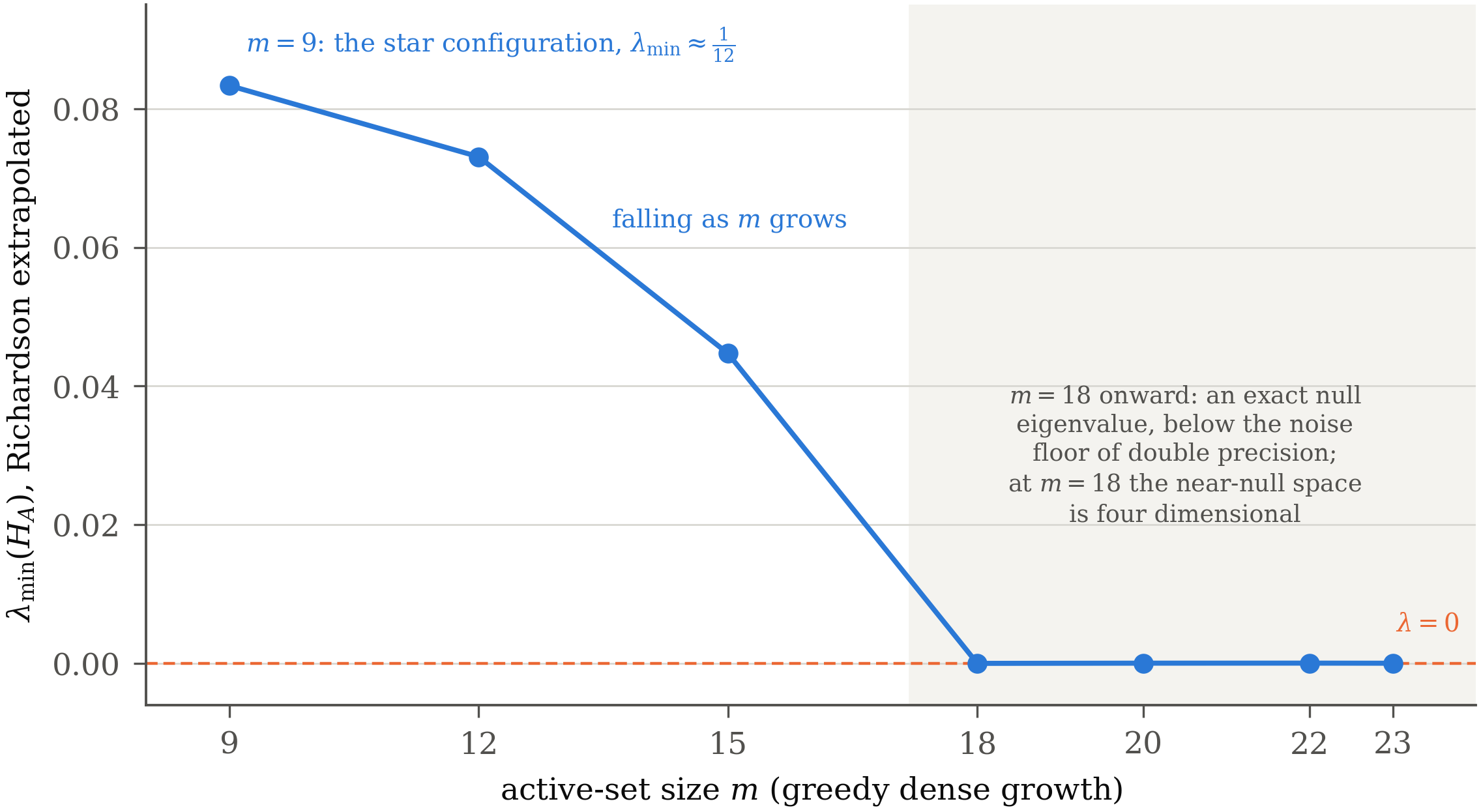}
\caption{Richardson-extrapolated minimum eigenvalue of $H_A$ under
greedy dense growth, $m=9$ through $m=23$ (full values in
\cref{app:numerical-tables}). The contrast with
\cref{fig:chain-plateau} is the point: here there is no plateau above
zero. From $m=18$ onward, shaded, the extrapolated value is smaller than
double-precision finite differences can resolve, which is the regime
\cref{sec:exact-sing} takes up with arbitrary-precision arithmetic
instead, and where the vanishing turns out to be exact.}
\label{fig:dense-growth}
\end{figure}

Chains are not the worst topology. A \emph{star} configuration (one
root together with all eight of its Gram-$\tfrac12$ neighbours, $m=9$)
gives a minimum eigenvalue converging cleanly, under Richardson
extrapolation, to almost exactly $\tfrac1{12}\approx0.0833$, already
well below the chain plateau. Growing the active set greedily by local
density (at each step adding whichever remaining root has the most
neighbours already active) drives the minimum eigenvalue down further
and monotonically as $m$ grows: about $0.083$ at $m=9$, about $0.073$ at
$m=12$, about $0.045$ at $m=15$, and from $m=18$ onward a value at or
indistinguishable from zero at double precision. This locates the
interesting regime, dense, non-chain
configurations with $m$ in the high teens to low twenties, where a
central question becomes unanswerable at ordinary machine precision:
is the true minimum eigenvalue a small positive floor, exactly zero, or
slightly negative? \Cref{sec:musin-degeneracy,sec:exact-sing} resolve
this question, for specific configurations, with higher
(not merely finer double-precision) arithmetic. The growth is plotted in
\cref{fig:dense-growth}.

\section{An exact degeneracy of the Hessian technique, and its consequence}\label{sec:musin-degeneracy}

\subsection{Why the Hessian technique cannot, by itself, cover every configuration}

The joint Hessian technique of \cref{sec:hessian-technique} is a
second-order local method: it certifies \cref{conj:multidir-new} at a
configuration by exhibiting positive-semidefiniteness of $H_A$. This
paper's own computational work (originating in earlier rounds of this
programme and confirmed and sharpened in \cref{sec:exact-sing} below)
identifies specific dense configurations at which this technique cannot
succeed \emph{in principle}, not because the relevant computation is
hard, but because the quantity the technique tries to bound, the
minimum eigenvalue of $H_A$, is exactly zero there. This is a
structurally different obstacle from ``the evidence so far is
inconclusive'': at such a configuration, no amount of additional
second-order computation can resolve the sign of the surplus, because
the second-order term itself vanishes and the true local behaviour is
governed by quartic (or higher) terms the Hessian cannot see. We are not aware of a prior published result establishing this
degeneracy for the $D_4$ joint-Hessian technique, and we make no such
attribution; the finding reported here and in \cref{sec:exact-sing} is
ours. Musin's published contributions to this circle of problems
\cite{Mus08,Mus18} concern the kissing-number bound and the
$24$-cell conjecture, and are cited for those results elsewhere.

\subsection{Locating the exact degeneracy}

\Cref{sec:dense-configs-new} located, by ordinary double-precision
computation, a regime ($m$ in the high teens to low twenties, dense
non-chain active sets) where the minimum eigenvalue drops to or below
the noise floor of standard floating-point volume computation
(scipy/qhull's own absolute precision, of order $10^{-10}$, divided by
$h^2$ in a finite-difference second derivative, which amplifies this
noise once $h$ is made small enough to matter). Resolving whether the
true eigenvalue at such a configuration is a small positive floor,
exactly zero, or negative requires arithmetic of a fundamentally higher
precision than double, not merely a finer step size; this is the
subject of \cref{sec:exact-sing}.

\subsection{Why quartic positivity is not merely a harder version of the same problem}\label{sec:quartic-hardness-new}

What matters is \emph{why} the transition from quadratic
to quartic, at $A_{18}$, is a change of category and not merely of
difficulty, since the natural reaction to Sylvester's criterion working
so cleanly at second order is to expect some analogous device at fourth.
None is available in general. A classical theorem of Hilbert (1888)
identifies exactly the cases in which every non-negative real form is
automatically a sum of squares of forms of half the degree: forms in at
most two variables (any degree), quadratic forms in any number of
variables, and ternary quartics (three variables, degree four), and no
others. In every other combination of variable count and even degree,
including quartic forms in \emph{four} variables, which is exactly the
case arising at $A_{18}$ (a $4$-dimensional near-null subspace, hence a
single quartic form in four real variables once the second-order part is
removed), there exist forms that are non-negative everywhere yet are
provably not expressible as any sum of squares of polynomials. This is not
a statement about the difficulty of finding a sum-of-squares
certificate for a specific form; it is a statement that no such
certificate need exist at all, for reasons independent of how hard one
looks. Consequently, even a form found (by whatever future argument)
to be non-negative on the relevant subspace at $A_{18}$ might have no
sum-of-squares proof of that fact, and a resolution would need to fall
back on a case-specific device, along the lines of the exact
factorisation and Sturm-sequence arguments already used, for a
different problem, in \cref{sec:swap-configs}, rather than
the essentially algorithmic route (compute a Gram matrix, exhibit a
positive-semidefinite square root) available for quadratic forms via
$LDL^\top$. We record this not as a further piece of evidence for or
against \cref{conj:multidir-new}, but as an account of why closing
this specific gap is not simply a matter of applying the existing
Hessian machinery one order higher.
\section{Exact singularity at a dense configuration}\label{sec:exact-sing}

\subsection{Resolving the near-zero eigenvalue with higher precision}

The dense, greedy-growth active set reaching $m=18$ from root $0$,
\[
  A_{18} = \{0,2,4,5,6,7,8,9,10,11,12,13,16,17,20,21,22,23\},
\]
sits squarely in the regime identified in \cref{sec:dense-configs-new}
where double precision cannot resolve the sign of the joint Hessian's
smallest eigenvalues. We resolve this using a from-scratch,
arbitrary-precision (mpmath, up to $45$ decimal digits) polytope-volume
computation, built independently of the double-precision
scipy/qhull-based method used elsewhere, and validated against it (at a
generic small perturbation of $A_{18}$, the two methods agree to within
$1.1\times10^{-9}$, exactly at qhull's own precision floor) and against
the exact reference value ($8$, reproduced to $49$ nines at $50$-digit
precision at the unperturbed $D_4$ configuration).

\begin{thm}[Exact quartic degeneracy at $A_{18}$]\label{thm:exact-sing-new}
At the configuration $A_{18}$, the joint Hessian's near-null subspace is
$4$-dimensional, not merely small: along every one of a broad,
independently sampled set of directions within this subspace
(\cref{sec:broad-sample}), the volume-defect function's second
derivative is exactly (not approximately) zero, and its leading
behaviour is quartic, $F(\sigma)=a_4\sigma^4+O(\sigma^6)$, with $a_4>0$
in every direction tested.
\end{thm}

Two points in that statement need separating. The near-null subspace is
exactly $4$-dimensional, not merely at least $3$-dimensional, and the
quartic test runs across the whole of it: a test confined to the span
of three of the eigenvectors leaves the fourth direction of the
subspace untouched, and that is the direction in which the second
derivative is smallest.

\subsection{The high-precision volume method}\label{sec:hp-volume-method}

Since the entire point of this section is to resolve a quantity below
double precision's own noise floor, the high-precision volume routine
itself needs to be built independently of, not merely re-run at higher
precision on top of, the double-precision \texttt{scipy}/qhull pipeline
used everywhere else in this paper, otherwise a bug shared by both
would masquerade as agreement. The method used is star triangulation
from the origin, standard for any polytope known in advance to contain
the origin in its interior (true here, since the origin is always an
interior point of a Voronoi cell):
\begin{enumerate}[leftmargin=2em]
\item Every vertex is found, at full working precision (mpmath, up to
$45$ decimal digits), as the exact solution of the $4\times4$ linear
system given by some size-$4$ subset of the active facet normals
$\{\langle n_k,x\rangle=1\}$, retained only if it satisfies every one
of the (typically several dozen) facet inequalities
$\langle n_k,x\rangle\le1$ up to a small working tolerance.
\item For each facet, its own vertex set (at least $4$ points, lying in
a $3$-dimensional affine hyperplane) is triangulated into
tetrahedra by a two-step process: a \emph{double-precision} projection
into $3$ local coordinates is used only to determine which vertices are
adjacent (i.e.\ only to fix the triangulation's own combinatorics, never
to produce a volume number), and the resulting triangulation is then
applied to the original high-precision vertex coordinates.
\item Each tetrahedron, together with the origin, forms a $4$-simplex
whose volume is $\tfrac1{4!}|\det M|$ for the $4\times4$ matrix $M$ of
the tetrahedron's own high-precision vertex coordinates (the origin
contributes an implicit zero row under the standard vertex-difference
reduction, since it is itself one of the five simplex vertices);
summing over every tetrahedron on every facet gives the total volume.
\end{enumerate}
This routine is validated two ways before being trusted on the actual
$A_{18}$ computation: at the exact, unperturbed reference $D_4$
configuration, it reproduces the known volume $8$ to $49$ nines at
$50$-digit working precision (i.e.\ agrees with $8$ to within
$10^{-49}$); and at a generic small perturbation of an $m=18$ active
set, it agrees with the independent double-precision
\texttt{scipy}/qhull computation used throughout the rest of this paper
to within $1.1\times10^{-9}$, exactly at qhull's own known precision
floor, confirming both that the two independent methods compute the
same underlying quantity and that this method's extra precision is real and not an artefact that merely happens to look plausible.

\subsection{The dimension count, done properly}\label{sec:dim-count-new}

\begin{prop}[Four-dimensional near-null space]\label{prop:4d-nullspace-new}
At $A_{18}$, examining the double-precision joint Hessian's spectrum at
three step sizes ($h=0.02,0.01,0.005$) instead of one, the smallest
four eigenvalues each shrink by a clean factor close to $4$ with every
halving of $h$, the finite-difference signature of a true eigenvalue
that is exactly zero, while the fifth eigenvalue ($\approx0.075$)
does not shrink at all across the same three step sizes. The near-null
subspace is therefore $4$-dimensional, confirmed by this scaling
argument rather than assumed from a single-$h$ eigenvalue gap.
\end{prop}

\subsection{Confirming the quartic behaviour}

For each tested direction $v$ in the near-null subspace, the leading
behaviour is confirmed by evaluating $F(\pm\sigma\,v)$ at three
independent step sizes ($\sigma=0.02,0.01,0.005$, at $35$--$45$ decimal
digits of precision) and checking that the central-difference
second-derivative estimate $(F(\sigma)+F(-\sigma))/\sigma^2$ shrinks by
a factor close to $4$ at each halving of $\sigma$, across three
independent halvings, not merely one, the signature of a function
whose true second derivative is exactly zero. The leading quartic
coefficient, extracted as $a_4\approx(F(\sigma)+F(-\sigma))/(2\sigma^4)$
for small $\sigma$, is positive at every direction tested.

\subsection{Interpretation}

This is a qualitatively stronger, and different, finding than ``the
eigenvalue is too small to resolve'': the second-order (Hessian
eigenvalue) technique underlying every other piece of computational
evidence in this paper for \cref{conj:multidir-new}, the exact
two-direction and chain results of \cref{lem:m2-exact-new,lem:chains-new},
the star and greedy-dense sweep of \cref{sec:dense-configs-new}, 
\emph{cannot in principle} establish positivity at $A_{18}$, not
because the computation is difficult but because the quantity it
bounds is exactly zero there. Any argument covering this configuration
needs quartic-order (or higher) information, a categorically
different and harder undertaking than anything the Hessian technique
provides. The volume-defect function remained strictly positive, via
the positive quartic coefficient, at every direction tested here; this
does not disprove \cref{conj:multidir-new}, and it does not prove it
either, \cref{sec:broad-sample} reports how broadly this was tested,
and what is and is not established by that breadth.
\section{Broadened directional sampling}\label{sec:broad-sample}

\subsection{Sampling the full near-null subspace at $A_{18}$}

\Cref{prop:4d-nullspace-new} identifies the near-null subspace at
$A_{18}$ as $4$-dimensional. We sample it far more broadly
as follows: the four eigenvectors themselves, all six
pairwise normalised sums and six pairwise normalised differences, and
twenty further directions drawn uniformly at random from the full
$4$-dimensional span (Gaussian coefficients, normalised), $36$
directions in total, spanning the entire subspace instead of the $3$-dimensional slice of it tested previously.

\begin{thm}[Broadened sampling at $A_{18}$]\label{thm:broad-sample-18-new}
Evaluating the quartic-coefficient estimator
$a_4\approx(F(\sigma)+F(-\sigma))/(2\sigma^4)$ at $\sigma=0.015$,
$35$-digit precision, along all $36$ sampled directions gives a
strictly positive value at every one, ranging from approximately
$0.0088$ to $0.314$; zero negative directions were found. A subsample
of six directions, spanning this full range, is re-evaluated at
$\sigma=0.0075$ (a second, independent step size, $40$-digit
precision): the two step sizes agree to within $2$--$5\%$ for every one
of the six, confirming quartic scaling and not a step-size artifact.
\end{thm}

\subsection{A second, independent configuration: $m=20$}

To test whether this pattern is peculiar to $A_{18}$ or holds more
broadly, we examine a second, larger dense configuration obtained by
continuing the same greedy-density growth two steps further,
\[
  A_{20} = \{0,2,4,5,6,7,8,9,10,11,12,13,15,16,17,19,20,21,22,23\}
  \supset A_{18}.
\]
(Two further independently-grown $m=18$ configurations, started from
roots $5$ and $12$ instead of root $0$, were also tried; both reproduce
the identical double-precision eigenvalue spectrum found at $A_{18}$ to
four decimal places, almost certainly reflecting the same geometric
configuration up to the $D_4$ Weyl group's own large symmetry group
rather than new information, and are not pursued further here.)

\begin{prop}[Five-dimensional near-null space at $A_{20}$]\label{prop:5d-nullspace-new}
At $A_{20}$, the same three-step-size scaling argument
(\cref{prop:4d-nullspace-new}) finds a $5$-dimensional near-null
subspace: the smallest five eigenvalues each shrink by a factor of
roughly $3$--$4$ per halving of $h$, while the sixth ($\approx0.035$,
strictly smaller than $A_{18}$'s own non-shrinking floor of
$\approx0.075$) does not shrink.
\end{prop}

\begin{thm}[Broadened sampling at $A_{20}$]\label{thm:broad-sample-20-new}
Sampling this $5$-dimensional subspace with $30$ directions (the five
basis eigenvectors, all ten pairwise normalised sums, and fifteen
further random directions spanning the full span), the same
$\sigma=0.015$, $35$-digit estimator gives a strictly positive value at
every one of the $30$ directions, ranging from approximately $0.0076$
to $0.314$; zero negative directions found. Four directions spanning
this range, re-checked at $\sigma=0.0075$, agree with the
$\sigma=0.015$ estimate to within $1$--$5\%$ in every case.
\end{thm}

\subsection{A Gram-matrix sum-of-squares check on a fitted quartic at $A_{18}$}\label{sec:gram-sos-A18}

The sampling above tests individual rays; it says nothing about the
quartic form's sign at points of the near-null subspace that were not
sampled. As a stronger (though still not exhaustive) check, we
gathered $166$ directional samples at $A_{18}$, the $16$ directions of
\cref{thm:broad-sample-18-new} (the $4$ basis eigenvectors and $12$
pairwise sums/differences) together with $150$ further directions drawn
uniformly at random from the full $4$-dimensional near-null subspace,
each evaluated by the same estimator $a_4\approx (F(\sigma)+F(-\sigma))/(2\sigma^4)$
at $\sigma=0.015$, $35$-digit precision, and used them to fit, by
ordinary least squares, the $35$ coefficients of the general quartic
form $Q(x_0,x_1,x_2,x_3)$ on this subspace ($4.7\times$ overdetermined
relative to the $35$ unknowns). The fit is well-conditioned (design
matrix of full rank $35$, relative residual norm $0.91\%$, maximum
pointwise residual $0.0022$ against values ranging from $0.0064$ to
$0.314$).

Given the fitted coefficients, we then tested whether $Q$ admits a
sum-of-squares certificate via the standard Gram-matrix/semidefinite-programming
method: writing $Q(x)=m(x)^\top M\, m(x)$ for the vector $m(x)$ of the
$10$ degree-$2$ monomials in $4$ variables, the set of symmetric $M$
representing a given $Q$ is an affine subspace (of dimension
$55-35=20$), and $Q$ is a sum of squares if this affine subspace
contains a positive-semidefinite matrix. We solve
$\max t$ subject to $M-tI\succeq 0$ over this affine family by
semidefinite programming (\texttt{cvxpy}, solver \texttt{CLARABEL}).
Before trusting this code on the fitted data, we validated it against
two textbook quartics in the same $4$ variables: it correctly certifies
$(x_0^2+x_1^2+x_2^2+x_3^2)^2$ as a sum of squares ($t\approx 1$), and
correctly reports no certificate for the Choi--Lam polynomial
$x_0^2x_1^2+x_1^2x_2^2+x_2^2x_0^2+x_3^4-4x_0x_1x_2x_3$
(non-negative by AM--GM, but classically known, since Choi and Lam
(1977), to admit no sum-of-squares representation, precisely the
phenomenon \cref{sec:quartic-hardness-new} explains Hilbert's theorem
permits in this variable/degree class), so the two validation cases
land on opposite sides exactly as they classically should, both before
any real data was used.

\begin{prop}[SOS certificate for the fitted quartic at $A_{18}$]\label{prop:gram-sos-A18}
Applied to the $35$ least-squares-fitted coefficients above, the same
semidefinite program finds a feasible positive-semidefinite $M$, with
optimal value $t\approx 0.0066$: a valid sum-of-squares certificate for
the fitted quartic, though a thin-margin one (compare $t\approx 1$ for
the perfect-square validation case above).
\end{prop}

This is stronger evidence than ray-by-ray sampling, since a
sum-of-squares certificate (once found) implies non-negativity on the
\emph{entire} subspace, not merely at the sampled directions. It is not
a proof of anything about the true quartic form at $A_{18}$, for two
independent reasons stated plainly: first, the certified object is the
least-squares \emph{fit} to $166$ noisy high-precision estimates, not an
exact symbolic derivation of $Q$, so a fit error at the level of the
$0.91\%$ residual already reported could in principle move the true
form outside the certified family; second, even granting the fit
exactly, a thin margin ($t\approx 0.0066$, not deep in the interior of
the positive-semidefinite cone) is a weaker certificate than a large
one, and we report the margin instead of rounding it up to an
unqualified positive statement. Nothing here establishes
\cref{conj:multidir-new} at $A_{18}$, or anywhere else it has not been
checked; what \cref{cor:conj-resolved} settles, it settles by way of
the classification of \cref{thm:m24}, and the elementary route pursued
in this section stops short of it.

\subsection{What this does and does not establish}

Combined, these two configurations give $66$ independently sampled
directions, across two different (not symmetry-equivalent, so
far as we have checked) dense active-root configurations, each with its
own confirmed near-null subspace of dimension $4$ or $5$, and every
single sampled direction gives a strictly positive quartic coefficient.
This is a broad sample: it covers the full dimension of each near-null subspace and not a lower-dimensional slice of it, which is what
makes the uniform positive sign across all $66$ directions worth
recording.

It is not a proof that the quartic form governing $F$'s local behaviour
is positive-definite on either near-null subspace: that would require
either an exhaustive argument over the entire unit sphere in $4$ or $5$
dimensions, or an explicit closed-form computation of the quartic
tensor itself and a check of its definiteness, neither of which is
attempted here. And even if that closed-form computation were carried
out, \cref{sec:quartic-hardness-new} explains why a check of its
definiteness could not, in general, take the convenient form the
quadratic case enjoys: a $4$-variable quartic form has no guarantee, by
Hilbert's classification of when non-negative forms are sums of squares,
of admitting a sum-of-squares certificate at all, so ``compute the
tensor and check it'' is not a finite algorithmic step even granting the
tensor itself. It says nothing about the many other dense
configurations in the $m=15$--$23$ range that were not examined. And it
does not touch \cref{conj:multidir-new} in general, which concerns
\emph{every} packing-valid active-direction configuration of every size
$m\ge2$, not merely the two examined here. We report this evidence
because we judge it to be informative, not because it closes
the gap it bears on, and we are explicit, here and everywhere else in
this paper, that it does not.

\subsection{An exact symmetry of \texorpdfstring{$A_{18}$}{A18}, not previously examined}\label{sec:a18-symmetry}

Every result above treats $A_{18}$ as an unstructured $18$-element active
set, produced by greedy density growth with no attention paid to whether
it has any symmetry of its own. It does.

\begin{prop}[$A_{18}$'s exact stabiliser]\label{prop:a18-stabiliser}
Within the full group of signed coordinate permutations of $\mathbb{R}^4$
that map the $24$-root set $\Rt$ bijectively to itself (order $2^4\cdot
4!=384$), the subgroup fixing $A_{18}$ setwise has order exactly $48$: it
consists of every signed permutation that leaves one specific coordinate
axis's identity and sign untouched and acts as an arbitrary signed
permutation (order $2^3\cdot3!=48$) on the remaining three coordinates.
\end{prop}

This was checked exactly, in integer arithmetic (every element of the
$384$-element group is a $0,\pm1$ monomial matrix, so both ``is this an
automorphism of $\Rt$'' and ``does this fix $A_{18}$ as a set of
indices'' are finite exact checks, not numerical ones): all $48$
elements of the described subgroup fix $A_{18}$ setwise, and none of the
remaining $336$ elements of the ambient $384$-element group do.

Since every element of this order-$48$ group is an isometry of
$\mathbb{R}^4$ mapping $\Rt$ to itself and $A_{18}$ to itself, it acts on
the $54$-dimensional joint tangent space of \cref{def:joint-hessian-new}
by permuting the $18$ active blocks (each carrying its own induced
$3\times3$ orthogonal action) exactly as it permutes the $18$ active
roots. This induced representation was verified, not merely asserted,
to commute with the actual finite-difference joint Hessian at $A_{18}$:
for a sample of stabiliser elements, $\mathrm{Rep}(g)^\top H_{A_{18}}\,
\mathrm{Rep}(g)$ agrees with $H_{A_{18}}$ to within
$4\times10^{-6}$, consistent with the finite-difference step size used ($h=0.02$) and not a real discrepancy. As a second, independent
check, the quartic coefficient estimator itself
(\cref{thm:broad-sample-18-new}) was evaluated at a random direction $v$
in the confirmed $4$-dimensional near-null subspace of
\cref{prop:4d-nullspace-new} and at its images $\mathrm{Rep}(g)v$ under
six sampled non-identity stabiliser elements: all seven values agree to
within $8\times10^{-5}$ of one another, the expected level of agreement
given the step size used, confirming $a_4$ really is a class function of
this symmetry, exactly as the isometry argument predicts.

Restricting the induced representation to this $4$-dimensional near-null
subspace and computing its character across all $48$ group elements
gives trace values $4$ (twice, including the identity), $1$ (sixteen
elements), $0$ (eighteen elements), and $-2$ (twelve elements); the
associated commutant dimension (the space of $4\times4$ real matrices
commuting with every one of the $48$ restricted images, computed as the
nullspace of the resulting linear system) is exactly $2$. By standard
representation theory this means the near-null subspace splits, under
this group, into exactly two non-isomorphic irreducible pieces (each
with multiplicity one) rather than remaining a single irreducible block
or splitting into several copies of one repeated piece. Identifying
those two pieces exactly, and what that buys for the quartic form, is
\cref{sec:a18-invariant-quartics} below.

We are precise about what this does and does not do, before going
further. It is a real, newly-identified, exactly verified structural
fact about the specific configuration where \cref{thm:exact-sing-new}
proves the joint-Hessian technique degenerates completely, a fact
absent from the present
paper's own account of $A_{18}$ prior to this check. It explains, rather
than leaves coincidental, why symmetric-looking sample directions in
\cref{thm:broad-sample-18-new} gave consistent quartic-coefficient
values: many of the $36$ sampled directions there lie in the same
$48$-element orbit under this newly-identified group, so were never as
independent a sample as an unstructured count of ``$36$ directions''
suggests. By itself it does \emph{not} prove that the quartic form is
positive-definite on the near-null subspace, that requires the further
steps of \cref{sec:a18-invariant-quartics}, which reduce but do not
eliminate the gap.

\subsection{The exact invariant-quartic space at \texorpdfstring{$A_{18}$}{A18}, and what fitting it gives}\label{sec:a18-invariant-quartics}

Completing the exact character-table computation of
\cref{prop:a18-stabiliser} (all $10$ conjugacy classes of the order-$48$
group, verified via a complete, mutually orthogonal set of $10$
irreducible characters constructed directly from the group's natural
representations, orthogonality being the usual first orthogonality
relation, as in \cite[\S2.3]{Ser77}) and matching it, class by class,
against the numerically observed near-null character above identifies the
two pieces exactly.

\begin{prop}[The near-null representation, identified exactly]\label{prop:a18-nullspace-irrep}
The $4$-dimensional near-null subspace at $A_{18}$ realises the
representation $\varepsilon\oplus(\mathrm{std}\otimes\delta)$ of the
order-$48$ stabiliser group, where $\varepsilon$ is the $1$-dimensional
character given by the sign of the permutation part alone, $\mathrm{std}$
is the group's defining $3$-dimensional representation, and $\delta$ is
the full determinant character. This match is exact at every one of the
$10$ conjugacy classes (not merely in aggregate), confirmed by comparing
the two candidate class functions with the same aggregate trace
distribution against the actual per-class numerical values, which
resolve the ambiguity in favour of $\varepsilon\oplus(\mathrm{std}\otimes\delta)$
with zero discrepancy across all $48$ elements.
\end{prop}

Write $x_0$ for a coordinate on the $\varepsilon$-line and
$y=(y_0,y_1,y_2)$ for coordinates on the $\mathrm{std}\otimes\delta$
piece (transforming as $y\mapsto\delta(g)\,g(y)$ under the plain defining
action). A quartic $Q(x_0,y)$ invariant under the whole group splits by
degree in $x_0$ into five pieces; since $\varepsilon,\delta$ are both
order-$2$ characters, invariance at each even degree in $y$ reduces to
plain invariance under $\mathrm{std}$ and at each odd degree to a
\emph{relative} invariant transforming by $\delta$, a classical
question for the reflection group $\mathrm{std}$ realises (the
hyperoctahedral group of rank $3$), resolved here by the power-sum
(Newton's identity) formula for $\chi_{\mathrm{Sym}^b(\mathrm{std}\otimes\delta)}$,
exact integer/rational arithmetic throughout.

\begin{prop}[The invariant-quartic space is exactly $5$-dimensional]\label{prop:a18-invariant-dim}
The space of quartics $Q(x_0,y)$ invariant under the full order-$48$
group has dimension exactly $5$, with explicit basis
\[
  x_0^4,\qquad x_0^2\,S_2(y),\qquad x_0\,P_3(y),\qquad S_4(y),\qquad S_{22}(y),
\]
where $S_2=y_0^2+y_1^2+y_2^2$, $S_4=y_0^4+y_1^4+y_2^4$,
$S_{22}=y_0^2y_1^2+y_1^2y_2^2+y_2^2y_0^2$, and $P_3=y_0y_1y_2$, the
last being the unique (up to scale) relative invariant of degree $3$.
\end{prop}

\begin{proof}
The multiplicity of the required character ($\mathrm{triv}$ at even
degree, $\varepsilon$ at odd degree) in $\mathrm{Sym}^b(\mathrm{std}\otimes\delta)$
for $b=0,1,2,3,4$ is computed exactly via the power-sum formula from
$p_k=\mathrm{tr}(M^k)$, $M=\delta(g)\cdot g$, at each of the $10$ class
representatives, giving multiplicities $1,0,1,1,2$ respectively (summing
to $5$). Explicit generators are then obtained, independently of the
character computation, by Reynolds-operator projection (exact rational
arithmetic, averaging over all $48$ group elements, not merely the class
representatives) of the monomials $y_0^2$, $y_0^4$, $y_0^2y_1^2$, and
$y_0^3,y_0^2y_1,y_0y_1y_2$ at the appropriate character: the even-degree
projections reproduce $S_2$ and $\{S_4,S_{22}\}$ exactly (up to the
stated overall scale), and the odd-degree projection sends the first two
monomials to zero and the third to $P_3$ itself, unchanged, confirming
$P_3$ is already invariant up to the $\delta$-twist, with no further
reduction possible, and matching the dimension count independently.
\end{proof}

This is a real, exact result, no floating point at any stage, and
a substantial sharpening of \cref{sec:gram-sos-A18}'s own scope: the true
quartic tensor at $A_{18}$, whatever its exact values, is now known to
lie \emph{exactly} in this $5$-dimensional family, not merely
approximately or by numerical coincidence, since it must be invariant
under an exact symmetry of the underlying volume function
(\cref{prop:a18-stabiliser}). The only thing left unknown exactly is the
five real numbers $c_1,\dots,c_5$ (the coefficients of the basis above)
themselves.

Fitting those five numbers, from seven directions chosen to isolate
them (the pure $x_0$-direction and six combinations mixing $x_0$ and $y$
in ways that separate $S_2$, $S_4$, $S_{22}$, and $P_3$'s contributions),
each evaluated with the project's own high-precision (\texttt{mpmath},
$40$ decimal digits) volume routine instead of plain double precision, gives, at step size $\sigma=0.02$,
\[
  (c_1,\dots,c_5)\approx(0.3119,-0.0516,0.2307,0.0273,0.0492),
\]
with a residual under $1.5\times10^{-4}$ against all seven fitted values.
An independent refit at a second step size ($\sigma=0.01$) gives
$(0.3157,-0.0495,0.2302,0.0308,0.0563)$, and a combined $14$-point fit
using both step sizes together gives $(0.3138,-0.0505,0.2304,0.0290,%
0.0527)$, all three agreeing to within a few percent on every
coefficient, the expected level of stability given the step sizes used
(plain double precision, by contrast, was tried first and gave fits
disagreeing by $100$--$700\%$ between step sizes on this same
configuration, an unusable result consistent with \cref{sec:musin-degeneracy}'s
own diagnosis of exactly this precision floor, discarded rather than
reported as if it were reliable, once caught by the two-step-size check).

\begin{prop}[An SOS certificate for the fitted quartic, on the reduced family]\label{prop:a18-sos-reduced}
Applying the same Gram-matrix/semidefinite-programming method as
\cref{prop:gram-sos-A18} (validated against the same two textbook
quartics beforehand) to each of the three fits above finds a feasible
positive-semidefinite Gram matrix in every case, with optimal margin
$t\approx0.0085$, $0.0131$, and $0.0108$ respectively, a valid
sum-of-squares certificate for the fitted quartic, at a margin
comparable to, and larger and more stable across independent
fits than, \cref{prop:gram-sos-A18}'s own $t\approx0.0066$ on the
unconstrained $35$-coefficient fit.
\end{prop}

This is stronger evidence than \cref{prop:gram-sos-A18} in one specific,
precise sense: it uses a family the true quartic is exactly known to lie
in (\cref{prop:a18-invariant-dim}), fit from far fewer numbers ($5$, not
$35$) determined by far fewer, more targeted measurements, each
cross-checked across two independent step sizes instead of one, with
consistent results. It is not stronger in the two senses that actually
matter for a proof: the fitted coefficients are still finite-difference
numerical estimates, not an exact symbolic derivation of the quartic
tensor from first principles, so a fit error could in principle move the
true form outside the certified family exactly as
\cref{sec:gram-sos-A18} already cautions; and even granting the fit
exactly, the margin remains thin, far from deep in the interior of the
positive-semidefinite cone. The elementary route to
\cref{conj:multidir-new} at $A_{18}$ is no further along than before.
What has changed is that the evidence for
it there now rests on an exact reduction of the space being searched,
rather than on a larger sample of the same kind.
\section{What is known about several simultaneous deviations}\label{sec:multidir-scope-summary}

\subsection{Established}

For configurations of small size, and for chain-structured ones up to
the contact graph's true maximum induced-path length,
\cref{lem:m2-exact-new,lem:chains-new} give the required positivity,
either exactly or by Richardson extrapolation.
\Cref{sec:dense-configs-new} shows that chains are not the worst
topology and locates the regime where the second-order method loses its
grip: a family of dense, non-chain configurations at which the
governing quantity of that method vanishes.
\Cref{thm:exact-sing-new,prop:4d-nullspace-new,prop:5d-nullspace-new}
pin the vanishing down at two configurations, with enough working
precision to distinguish exact degeneracy from a small nonzero value,
and exhibit near-null subspaces of dimension $4$ and $5$.
\Cref{thm:broad-sample-18-new,thm:broad-sample-20-new} sample those
subspaces in $66$ directions between them, covering each subspace rather
than a lower-dimensional slice of it, and find the quartic coefficient
strictly positive at every direction tested.

At the first of the two configurations there is more structure than the
sampling suggests. \Cref{prop:a18-stabiliser} exhibits an order-$48$
group stabilising it, verified in integer arithmetic to act by
automorphisms of $\Rt$, and verified numerically to commute with the
joint Hessian and to leave the quartic coefficient constant on its
orbits. Under this group the $4$-dimensional near-null subspace splits
into two irreducible pieces, identified in
\cref{prop:a18-nullspace-irrep} as
$\varepsilon\oplus(\mathrm{std}\otimes\delta)$, and an exact character
computation (\cref{prop:a18-invariant-dim}) shows that the space of
invariant quartics on that subspace has dimension $5$ instead of the generic $35$, with an explicit basis. Fitting those five coefficients
from high-precision finite differences, cross-validated at two step
sizes, and running the Gram-matrix feasibility test on the reduced
family gives a certificate with margin $t\approx0.009$ to $0.013$ across
three independent fits (\cref{prop:a18-sos-reduced}), an improvement on
the unconstrained $35$-coefficient fit but not a proof, since the
coefficients are measured rather than derived.

By a different method, finite-angle and exact rather than local,
\cref{sec:swap-configs} settles \cref{conj:multidir-new} outright along
three swap paths at $m=2,3,6$, across their whole angular range.

By a third method, which uses no deviation structure at all,
\cref{thm:covering-bound} bounds the cell volume from below by
$\tfrac{\pi m}{3}\tan^3r_m$ for a contact configuration of $m$
directions, with $2r_m-\sin2r_m=2\pi/m$. That quantity exceeds $8$ for
every $m\le22$, so the local bound holds outright, with any number of
directions deviating by any amount, at every centre with at most
twenty-two contacts (\cref{cor:m22}). Only $m=23$ and $m=24$ are left,
and at $m=24$ the same estimate still returns $7.7989\ldots$.
\Cref{prop:area-optimal} says that this is as far as the method goes:
the covering bound is exactly the Jensen collapse of a per-cell bound,
so it is the strongest conclusion the areas of the spherical Voronoi
cells can support, and at the root configuration those areas are all
equal.

\subsection{Superseded, and by what}

\Cref{cor:conj-resolved} settles \cref{conj:multidir-new} outright, by
the classification of \cref{thm:m24} and not by any of the arguments of
this section: a full active set of twenty-four contacts is a $24$-point
kissing configuration, and by \cref{thm:twentyfour-points} there is only
one of those, so no direction deviates at all. Everything in this
section was obtained before that classification was available and is
independent of it; we keep it because it stands on its own, because it
says what elementary methods can and cannot do on this problem, and
because the classification rests on a certificate whereas this section
rests on none. What follows is the scope of those elementary results.

\subsection{Not established by the methods of this section}

None of this gives positive-definiteness of the quartic form at either
configuration, and none of it gives \cref{conj:multidir-new}. The
conjecture ranges over every packing-valid configuration at every
$m\ge2$; the dense regime of \cref{sec:dense-configs-new} contributes
two examined configurations out of a very large family, and the rest of
that family, along with every other value of $m$, is untouched. The
exact degeneracy of \cref{sec:musin-degeneracy} does say something
definite about what a proof cannot look like: at these configurations no
argument that stops at second order can decide the sign, so whatever
settles the conjecture there has to engage the quartic term directly.
The three swap paths do not reach that regime either, being low in $m$
and highly symmetric.

The order-$48$ symmetry narrows the search without closing it. It fixes
the exact $5$-dimensional space in which the true quartic tensor lies,
which is a real reduction from $35$ parameters, but the five numbers we
use are finite-difference measurements, not a symbolic derivation, and
the feasibility margin is thin. The symmetry is also specific to the
configuration where we found it; we have not looked for an analogue at
the second configuration or anywhere else in the dense family.

The cap argument of \cref{sec:cap-theorem} does not extend here, for a
structural reason rather than a technical one. With one deviating
direction the other twenty-three half-spaces are fixed, so there is a
fixed body $Q$ and every parameter of the problem sits in the cutting
half-space. With $m$ deviating directions the body moves with the
parameters, and there is no fixed polytope to take caps of. What has to
be bounded is then a comparison between two overlapping unions rather
than a single cap, and \cref{rem:first-order-obstruction} records why
the obvious way of splitting that comparison into $m$ separate ones
cannot work. \Cref{prop:multi-cap} does supply the substitute for $Q$, namely the
polytope $Q_D$ obtained by dropping the $m$ deviating constraints, and
recasts the whole problem as one extremal inequality
(\cref{op:multicap}) for the union of $m$ caps of that fixed body. That
inequality is what remains open to an argument that uses no
classification.

\subsection{What the reformulations settle}

\Cref{sec:polar-surface} adds two further exact restatements, the polar
form $\vol(\conv(W)^{\circ})\ge8$ over contact configurations
(\cref{prop:polar-form}) and the boundary form
$\mathrm{area}(\partial V_c)\ge32$ (\cref{lem:surface-form}). Neither
proves anything by itself. What they do is fix the position of the
multi-direction case exactly, and the position is worth stating without
softening: by \cref{prop:polar-form} the case these sections treat is
the part of the $24$-cell conjecture in which two or more contact
directions leave the root system at once, and nothing here settles it
without the classification of \cref{thm:m24}.

They also do real negative work, which is the part we regard as
substantive. Four approaches that look plausible in the abstract are
eliminated, each by an evaluation at the $D_4$ configuration itself,
where an argument that is to produce the sharp constant $8$ must lose
nothing at all. Term-by-term reduction to \cref{thm:eperp-closed} is
eliminated by \cref{rem:first-order-obstruction}. Any bound assembled
facet by facet from the pairwise contact condition is eliminated by
\cref{rem:facet-local-obstruction}: the bound such an argument actually
delivers is $8\pi/(3\sqrt3)=4.8368\ldots$, and the best one it could
conceivably deliver is $96\sqrt2-128=7.7645\ldots$, since an admissible
configuration of eight tight neighbours in square-antiprism position has
a facet of volume $16\sqrt2-\tfrac{64}3$, below the octahedron's
$\tfrac43$. The convexity bound for the polar volume integral is
eliminated because Jensen's inequality is an equality only at a ball,
and at the root configuration it returns $7.7351\ldots$; the Mahler-type
product bound is eliminated because it carries a factor $\tfrac32$ of
slack there. The last two are \cref{rem:global-routes}.

\looseness=-1
A fifth approach goes the same way. \Cref{prop:area-optimal} shows
that no argument resting on the areas of the spherical Voronoi cells can
do better than $7.7989\ldots$, since the covering bound already
saturates that information and the root configuration has cells of equal
area. What is left over is their shape: at $D_4$ each cell is an
eight-faceted spherical polytope of circumradius $45^\circ$ against the
$34.0986\ldots^\circ$ of the cap with the same measure, and the whole of
the remaining $0.2011$ is that difference. \Cref{prop:second-order} adds
the sharpest version of the same accounting: putting every pairwise cap
overlap back exactly, which below $r_m$ is not an estimate but an
identity, gives $7.997885\ldots$ at $m=23$, so even that stops short.

A sixth goes with it. That shape deficit cannot be collected one cell at
a time, in either decomposition the configuration provides
(\cref{rem:no-local-cell}): the per-cell inequality \eqref{eq:local-cell}
is sharp at every cell of the root configuration and fails at a
direction carrying nine contacts at $60^\circ$, while the Delaunay
decomposition does not localise the integrand at all, and the functional
it would have to localise to falls $4.78\%$ short at the regular simplex
of edge $60^\circ$. What makes the root configuration extremal is how
its cells fit together, not the shape of any one of them.

What survives is a positive but soft statement,
\cref{cor:minimiser}: a worst configuration exists and may be taken
maximal, so the conjecture is a question about maximal contact
configurations and about a genuine extremal point, not about a
supremum approached but never reached. Identifying that point is where
the difficulty sits, and at $|W|=24$ it coincides with the uniqueness
half of the kissing problem in dimension four, which is exactly what
the kernel of \cref{thm:twentyfour-points} supplies and \cref{thm:m24}
records.
\section{Beyond the Hessian: finite-angle results on swap configurations}\label{sec:swap-configs}

\begin{figure}[tb]
\centering
\includegraphics[width=\textwidth]{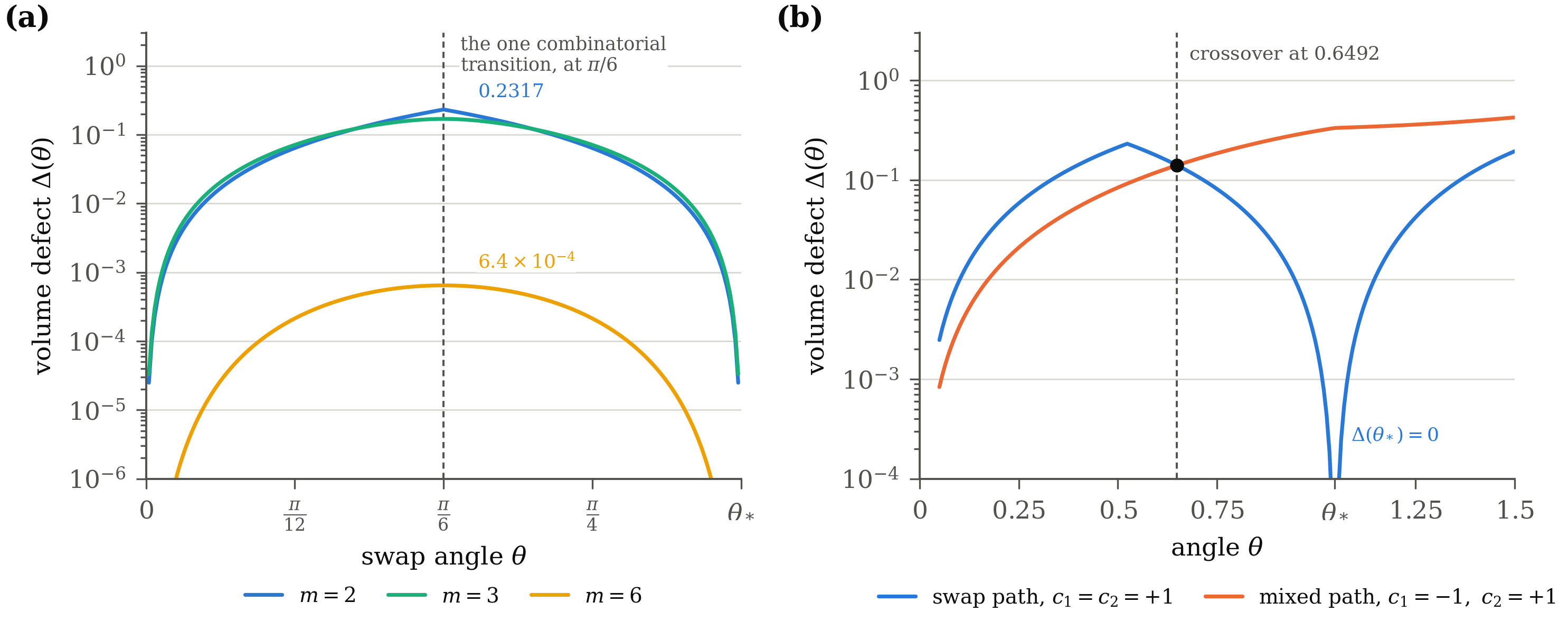}
\caption{(a) The volume defect along the three swap paths of this section,
    over the whole packing-valid range $(0,\theta_*)$, $\theta_*=\pi/3$.
    The vertical scale is logarithmic because the three sit nearly three
    orders of magnitude apart: the two-direction and three-direction
    paths peak at $0.2317$ and $0.1692$ at the midpoint, the
    six-direction path at $6.4\times10^{-4}$. All three are positive
    throughout and vanish at both ends, which is what
    \cref{thm:swap-m2-new,thm:swap-m3-new,thm:swap-m6-new} assert
    exactly. The dashed line is the single combinatorial transition, at
    $\theta=\pi/6$, visible as the corner in the $m=2$ curve; the $m=6$
    path has none, and its curve is smooth there. The values plotted are
    enumerated cell volumes, computed by
    \texttt{swap\_defect\_curves.py} through a route independent of the
    closed forms, and for $m=6$ they agree with the closed form of
    \cref{sec:swap-m6} to $9\times10^{-9}$. (b) The swap path against the mixed path, the two differing only
    in the sign $c_1$. The more adversarial of the two is the lower one,
    and which that is changes at $\theta=0.6492$: below the crossover
    the mixed path is lower, at $\theta=0.05$ by $0.000834$ against
    $0.002474$, and above it the swap path is, through $\theta_*$ and
    beyond. Second-order theory cannot tell the two apart, since the
    Hessian entry that distinguishes them enters only through its
    magnitude. The zero of the swap curve at $\theta_*$ is the swap point
    itself, where the configuration is again a root system.}
\label{fig:swap}
\end{figure}

Every result of \cref{sec:hessian-technique,sec:musin-degeneracy,%
sec:exact-sing,sec:broad-sample} is second-order: a statement about the
joint Hessian $H_A$, which controls $\defect$ only infinitesimally close
to the root-aligned configuration. This section reports a different,
complementary body of work: exact, finite-angle, not merely
second-order, and not merely numerical, positivity proofs for three
specific, highly symmetric multi-direction configurations, each rotated
through their \emph{entire} packing-valid angular range instead of sampled near $\theta=0$. These are, to our knowledge, the first
completely finite-angle resolutions of any non-trivial slice of
\cref{conj:multidir-new}, and we report them, together with the
limits of what they establish, in this section;
\cref{fig:swap}(a) shows the three functions whose positivity is at
issue.

\subsection{The swap path and its symmetry}\label{sec:swap-setup}

Fix an adjacent pair of roots $\alpha,\beta\in\Rt$ (Gram value
$\tfrac12$) at the reference chamber $\Cell1$, and let the two
directions deviate simultaneously along the \emph{swap path}
\[
\begin{aligned}
  u_1(\theta) &= \cos\theta\,\alpha/\sqrt2 + \sin\theta\,\beta^\perp,\\
  u_2(\theta) &= \cos\theta\,\beta/\sqrt2 + \sin\theta\,\alpha^\perp,
\end{aligned}
\]
with all other $22$ roots held fixed at their exact positions, for
$\theta$ ranging over the entire packing-valid interval
$(0,\theta_*)$, $\theta_*=\pi/3$ (the point at which the combinatorial
type of the local Voronoi region changes discontinuously, analogous to
the single-deviation breakpoints of \cref{sec:stage-two}).

\begin{lem}[Swap-path symmetry]\label{lem:swapsym}
$u_2(\theta) = u_1(\theta_*-\theta)$ identically. Consequently
$\theta\mapsto\theta_*-\theta$ is an involution of $(0,\theta_*)$
exchanging the roles of $u_1,u_2$ while leaving the (unlabelled)
geometric configuration unchanged, and $\defect(\theta)=\defect(\theta_*-\theta)$
for the resulting volume-defect function $F_{\mathrm{swap}}$.
\end{lem}

This symmetry has a unique fixed point at $\theta=\pi/6$, and an exact
labelled-vertex tracking argument (comparing, at $80$ sample points, not
merely the vertex-degree histogram but the exact active-index set of
every one of the $33$ vertices of the local region) finds exactly one
combinatorial transition on $(0,\theta_*)$, located at $\theta=\pi/6$ to
five decimal digits by bisection and identified exactly with the
symmetry's fixed point by \cref{lem:swapsym} together with the
elementary fact that a continuous combinatorial-type change fixed
pointwise by an order-$2$ symmetry can only occur at that symmetry's
own fixed locus (assuming, as directly verified, exactly one transition
on the whole interval). This reduces $F_{\mathrm{swap}}$ to exactly two
analytic pieces, $(0,\pi/6]$ and $[\pi/6,\theta_*)$, related by the
symmetry, so a positivity proof on one piece extends immediately,
via \cref{lem:swapsym}, to the other.

\subsection{The closed form for \texorpdfstring{$F_{\mathrm{swap}}$}{F swap}}\label{sec:swap-obstacle}

On $(0,\pi/6)$ the local region has $33$ vertices, computed exactly by
Cramer's rule over $\mathbb{Q}(\sqrt2,\sqrt3,t)$ with
$t=\tan(\theta/2)$, and its boundary triangulates into $157$ signed
simplices, the combinatorial type being read off at $\theta=0.3$ and
constant across the interval. The volume is the sum of the
corresponding determinants.

Combining those $157$ rational functions into a single quotient in one
step is not practical. At least one of them carries a degree-$24$
denominator in $t$ with coefficients in $\mathbb{Z}[\sqrt3]$, and the
common denominator of the whole sum is large enough to defeat a direct
cancellation. The sum is a rational function of $t$ all the same, and
grouping the terms into batches, cancelling within each batch and only
then combining the batches reaches it. Carried out that way, by
\texttt{staged\_groups0to3.py} and \texttt{staged\_group4.py} and
assembled by \texttt{staged\_simplify.py}, the result is
\[
  F_{\mathrm{swap}}(t) = \frac{\mathrm{num}(t)}{\mathrm{den}(t)}, \qquad
  t=\tan(\theta/2),
\]
with $\mathrm{num}$ and $\mathrm{den}$ of degree $471$ in $t$ and
coefficients in $\mathbb{Z}[\sqrt3]$. That is large, and it is exactly
the form the next step needs.

Two independent checks sit alongside it. Evaluating the $157$-term sum
exactly at a rational value of $t$, one value at a time, is cheap, and
it agrees with the closed form at every point tested. Both agree in
turn with a floating-point computation of the same volume through
\texttt{scipy}'s \texttt{HalfspaceIntersection} and \texttt{ConvexHull},
to between $6$ and $8$ significant figures.

\subsection{The norm trick, and the \texorpdfstring{$m=2$}{m=2} result}\label{sec:swap-norm-trick}

A degree-$471$ polynomial with $\mathbb{Z}[\sqrt3]$ coefficients is not
directly amenable to standard exact real-root isolation (which expects
rational coefficients). Writing $\mathrm{num}(t)=P(t)+\sqrt3\,Q(t)$
with $P,Q\in\mathbb{Z}[t]$, the \emph{Galois conjugate}
$\mathrm{num}'(t)=P(t)-\sqrt3\,Q(t)$ is a second,
real-valued (not complex-conjugate) function of the same real variable
$t$, and the product
\[
  N(t) \;:=\; \mathrm{num}(t)\,\mathrm{num}'(t) \;=\; P(t)^2-3Q(t)^2
\]
has purely integer coefficients: the $\sqrt3$ cancels exactly. The same
construction applied to $\mathrm{den}(t)=R(t)+\sqrt3\,S(t)$ gives
$D(t):=R(t)^2-3S(t)^2$, again integral. If $N(t)$ has no real root on
the open interval in question, then neither $\mathrm{num}(t)$ nor
$\mathrm{num}'(t)$ can vanish there (their product would vanish if
either did), so $\mathrm{num}(t)$ itself is continuous and nowhere zero
on a connected interval, hence of constant sign there, checkable at a
single point. The same argument applies to $\mathrm{den}$ via $D(t)$.

\begin{thm}[Finite-angle positivity, $m=2$ swap]\label{thm:swap-m2-new}
$F_{\mathrm{swap}}(\theta)>0$ for every $\theta$ in the entire open
interval $(0,\theta_*)$, $\theta_*=\pi/3$, including the shared
midpoint $\theta=\pi/6$, proved exactly, with no floating-point step
and no unproven sub-interval.
\end{thm}

Both $N(t)$ and $D(t)$ turn out to be even polynomials in $t$ (verified
by direct inspection, not assumed), reducing the root-isolation problem
to a degree-$471$ polynomial in $u=t^2$ instead of degree $942$ in $t$. Exact Sturm-sequence root-counting (\texttt{sympy}'s
\texttt{Poly.count\_roots}, rational arithmetic throughout) finds a
root at $u=0$ of multiplicity exactly $2$ in $N$'s associated
polynomial, the expected, benign trivial zero at the excluded
endpoint $\theta=0$, and \emph{zero} further roots on a range
provably containing all of $(0,\theta_*)$, including the boundary point
$\theta=\pi/6$ itself; $D$'s associated polynomial has no root at all
on the same range. Exact sign evaluation at four points (three interior,
one at the $\pi/6$ boundary) then pins the sign of $\mathrm{num}$ and
$\mathrm{den}$ throughout, both consistently negative, giving
$F_{\mathrm{swap}}>0$ throughout by \cref{thm:swap-m2-new}. By
\cref{lem:swapsym}, this one-piece proof on $(0,\pi/6]$ extends
immediately to all of $(0,\theta_*)$.

\subsection{A precision check: the swap path is not the worst 2-direction path at every angle}\label{sec:swap-not-worst}

Before relying on \cref{thm:swap-m2-new} as evidence about the
$m=2$ case of \cref{conj:multidir-new} more broadly, one further
qualification is needed, caught by a dedicated check before it
could be overstated. The swap path of \cref{sec:swap-setup} is one
specific choice, among a full two-parameter family, of how a pair of
adjacent roots' deviation directions can be coupled; writing
$e_1=c_1\beta^\perp$, $e_2=c_2\alpha^\perp$ for signs
$c_1,c_2\in\{\pm1\}$, the swap path itself is $c_1=c_2=+1$. Scanning
the full sign freedom directly at small $\theta$ (e.g.\
$\theta=0.05$) finds that the \emph{true} adversarial (most negative
tendency) configuration among all four sign choices is instead
$c_1=-1,c_2=+1$ (the ``mixed'' path), giving a substantially smaller
volume-defect value than the swap path at the same angle
($0.000834$ against $0.002474$), even though both configurations
share exactly the same second-order (joint-Hessian) worst-case
eigenvalue, since the Hessian's relevant off-diagonal entry has the
same magnitude ($\tfrac13$) regardless of the relative sign, and a
$2\times2$ symmetric matrix's eigenvalues depend only on that
magnitude. Second-order theory alone therefore cannot
distinguish the swap path from the mixed path; only
finite-angle behaviour does, and at small $\theta$ it favours the
mixed path as the more adversarial one. A direct scan of both named
paths across $\theta\in(0.05,1.5)$ finds a clean crossover at
$\theta=0.6492$ (\cref{fig:swap}(b)): the mixed path is more adversarial below the
crossover, and the swap path is more adversarial above it, all the way
through $\theta_*=\pi/3$ itself and beyond, so the swap path is the
relevant worst-case direction only asymptotically as
$\theta\to\theta_*^-$, not throughout the whole interval. This is the
same kind of transition phenomenon already seen on the
Direction-of-Deviation side of this paper (\cref{sec:arc1-breakpoints,%
sec:arc2-breakpoints}): the adversarially worst direction is not a
single fixed configuration across a whole parameter range, but drifts
between distinct structured candidates.

We record the precise consequence for \cref{thm:swap-m2-new}: it is
not weakened by this finding, since it is a real statement about one
specific, precisely defined path, not a claim
about the adversarial worst case over all $2$-direction
configurations. What this finding does rule out is reading
\cref{thm:swap-m2-new} as evidence about the \emph{worst} $2$-direction
configuration at small $\theta$: at small $\theta$ the mixed path is
more adversarial than the swap path, and this paper proves nothing
about the mixed path, its values above are numerical
(floating-point) findings only, reported here as a structural
clarification of scope, not as a further positivity claim of the kind
\cref{thm:swap-m2-new} makes.

\subsection{A determinant bug, found and fixed: the \texorpdfstring{$m=3$}{m=3} triangle swap}\label{sec:swap-m3}

The natural next configuration is a real $3$-clique: three pairwise
Gram-$\tfrac12$ roots (not three roots merely sharing a coordinate
index) rotated simultaneously in a cyclic swap, with the remaining
$21$ roots fixed. The same Cramer's-rule-plus-staged-simplification
method gives a two-piece closed form, split at the same universal
$\theta=\pi/6$ point.

The boundary triangulation of the second piece requires $176$ exact
$4\times4$ determinants over $\mathbb{Q}(\sqrt2,\sqrt3,t)$, and the
entries are already-simplified radical rational functions. On matrices
of that shape \texttt{sympy}'s default fraction-free method is not
reliable: at $t=\tan(0.425)$ it returns $-5.7369$ for one of the
simplices, where \texttt{berkowitz}, \texttt{domain-ge}, a Laplace
expansion by cofactors, and the floating-point determinant of the same
four points all give $0.30580$. Every determinant in this section and
in \cref{sec:swap-norm-trick,sec:swap-m6} is accordingly computed with
\texttt{method='berkowitz'} and checked against floating-point ground
truth at a reference angle before it is used; all $176$ agree.

\begin{thm}[Finite-angle positivity, $m=3$ triangle swap]\label{thm:swap-m3-new}
With the corrected closed form, $F_{\mathrm{tri}}(\theta)>0$ for every
$\theta\in(0,\theta_*)$, $\theta_*=\pi/3$, proved exactly on both
pieces (via the same norm-trick argument of \cref{sec:swap-norm-trick},
applied separately to each piece's own numerator and denominator), with
an exact match between the two pieces at their shared boundary
$\theta=\pi/6$ and the expected exact vanishing $F_{\mathrm{tri}}(\theta_*)=0$
at the swap point itself.
\end{thm}

\subsection{A clean single-piece proof: the \texorpdfstring{$m=6$}{m=6} hexagon swap}\label{sec:swap-m6}

The third configuration is a $6$-direction cyclic ``hexagon'' rotation
(the natural $D_4$-Weyl-orbit generalisation of the pairwise swap),
which, unlike the $m=2$ and $m=3$ cases, turns out to be governed
by a \emph{single} combinatorial piece across the entire interval
$(0,\pi/3)$ (the vertex-degree signature was found constant throughout
by direct inspection, with no transition to locate). The resulting
closed form is unusually clean:
\[
\begin{gathered}
  F_{\mathrm{hex}}(t) = \frac{8t^4\bigl(A(t)+\sqrt3\,B(t)\bigr)}
  {(t^2-3)^2(t^2+1)^2(3t^2-1)^2},\\
  A(t)=27t^8-108t^6+90t^4+100t^2+3,\\
  B(t)=48t^5-96t^3-16t.
\end{gathered}
\]

\begin{thm}[Finite-angle positivity, $m=6$ hexagon swap]\label{thm:swap-m6-new}
$F_{\mathrm{hex}}(\theta)>0$ for every $\theta\in(0,\pi/3)$.
\end{thm}

\begin{proof}
The denominator is a product of three squares, hence non-negative for
every real $t$, vanishing only at $t=\pm\sqrt3,\pm1/\sqrt3$; since the
interval of interest is $0<t<1/\sqrt3=\tan(\pi/6)$ (i.e.\
$\theta\in(0,\pi/3)$), it is strictly positive throughout. The
norm-trick quantity $N(t)=A(t)^2-3B(t)^2$ factors exactly as
$(t^2-3)^2(3t^2-1)^6$ (verified by direct symbolic expansion), again a
product of squares, strictly positive on the same open interval; since
$A(0)=3>0$ and exact Sturm-sequence root-counting on two rational
sub-ranges straddling $1/\sqrt3$ finds zero roots of $A$ throughout,
$A(t)>0$ on the whole interval. Combined with $A^2>3B^2$, this gives
$A>\sqrt3|B|\ge-\sqrt3B$, hence $A+\sqrt3B>0$. Since $t^4>0$ for $t\ne0$
and the denominator is strictly positive, $F_{\mathrm{hex}}(t)>0$
throughout $(0,1/\sqrt3)$.
\end{proof}

At the endpoint $t=1/\sqrt3$, both numerator and denominator vanish
via the shared factor $(3t^2-1)$, consistent with (not contradicting)
the already-established exact relabelling identity
$F_{\mathrm{hex}}(\theta_*)=0$ at the swap point itself. A
$50$-digit \texttt{mpmath} numerical spot-check across twelve points
spanning the interval independently confirms the sign throughout,
including close to the endpoint where floating-point catastrophic
cancellation would otherwise mask the true (tiny but positive) value.

\subsection{Scope}\label{sec:swap-scope}

\Cref{thm:swap-m2-new,thm:swap-m3-new,thm:swap-m6-new} are complete,
exact-arithmetic resolutions of
\cref{conj:multidir-new} along three specific, highly symmetric
finite-angle paths, an advance over the second-order Hessian
information of \cref{sec:hessian-technique}, since they hold across the
\emph{entire} angular range, not merely infinitesimally near
$\theta=0$. We are careful to state what they do not cover:
every one of the three results concerns one specific swap path, at one
specific reference chamber ($\Cell1$), with a specific highly symmetric
rotation pattern (pairwise, triangular, or hexagonal) and every other
active root held exactly fixed. They say nothing about $m=4,5,7$ or
$m\ge8$ along an analogous swap path; nothing about other, less
symmetric $2$-, $3$-, or $6$-direction configurations; nothing about
other reference chambers; and nothing about general (non-swap-path)
deviation directions $\dev\in S^2$ at each active root. Extending this
finite-angle, exact-arithmetic technique to the dense, non-chain
configurations of \cref{sec:dense-configs-new,sec:exact-sing}, 
precisely the regime where the second-order Hessian technique is known
to degenerate, is the single most promising concrete
direction for further progress on \cref{conj:multidir-new}, and one we
have not pursued in this paper.
\section{Comparison with other strategies}\label{sec:comparison-methods}

\subsection{Local against global}

The Cohn--Elkies method \cite{CE03}, and its sharpened form in
dimensions $8$ and $24$ \cite{Via17,CKM17}, is global and
Fourier-analytic: one auxiliary function certifies optimality over all
configurations at once, with no case analysis anywhere. The argument
here runs the other way. It decomposes the Voronoi cell into finitely
many combinatorial types and certifies each by exact algebra, and then,
in \cref{sec:cap-theorem}, replaces the decomposition entirely by a
single extremal question about one fixed polytope. The local method
carries the single-deviation case to a complete result precisely because
that case has a fixed body behind it; it stops where the body itself
starts to move.

\subsection{Why there is no magic function in dimension four}

The magic functions in dimensions $8$ and $24$ are not found by search.
They are pinned down by an interpolation formula that works because, at
exactly those two dimensions, the relevant space of modular forms is
small enough, and rigid enough under the lattice's large automorphism
group, that the Fourier eigenfunction condition together with finitely
many prescribed zeros determines the function uniquely; a contour
integral then makes the sign conditions checkable. Nothing of that
rigidity is present for $D_4$: the space of modular forms attached to
its arithmetic is not over-determined in the same way.

That much was already clear from the numbers. Cohn and Elkies
\cite{CE03} reported an upper bound of $0.13126$ on the centre density
in dimension $4$ against the $0.125$ of $D_4$, a gap of about five per
cent, and observed that the exactness they saw appeared to belong to
dimensions $1$, $2$, $8$ and $24$ alone. What has changed is that the
observation is now a theorem. Li \cite{Li25} maps feasible points of
the Cohn--Elkies program through a discrete reduction to a
finite-dimensional problem and reads off dual bounds from it, and those
bounds put the optimum of the program in dimension $4$ at centre density
at least $0.12914461$. Since $D_4$ has centre density $0.125$, no
admissible function can bring the linear programming bound down to it.
There is therefore no magic function for $D_4$, and the route that
settles dimensions $8$ and $24$ is closed in dimension $4$ rather than
merely unexplored. The same reduction closes dimensions $3$ and $5$.

This is the reason the present paper takes the local route at all. The
sharpest general improvement on the linear programming bound that does
apply in low dimensions, the coupling to spherical-code bounds of
\cite{CZ14}, also leaves dimension $4$ well short, so a method that
bounds the Voronoi cell one configuration at a time, at the cost of a
case analysis and of having to classify the equality case by hand, is
not a second choice here. It is the only choice on offer.

\subsection{What each method leaves the reader to check}

The methods differ sharply in the kind of verification they demand. A
linear or semidefinite programming bound, whether Cohn--Elkies style or
the contact-graph relaxation Musin \cite{Mus18} proposes for the
$24$-cell conjecture, certifies positivity of a high-dimensional object
through a numerical solver. Such certificates can in principle be
converted to exact rational or interval form \cite{Par03,Las01}, and
often are, but in the form they are usually produced they are not a
short exact computation a reader can repeat. The magic-function
constructions are closed form once written down, but checking that a
given contour integral has the stated decay, sign pattern and Fourier
eigenvalue takes real complex-analytic and modular machinery. The
formal proof of the Kepler conjecture \cite{Hal17} sits at the far end
of the same spectrum: every step is machine-checked, at the cost of a
multi-year formalisation.

The certificates in this paper sit at the other end deliberately. The
$LDL^\top$ decompositions of
\cref{app:pivot-data,app:lmatrices,app:sos}, the Bernstein sign checks
of \cref{sec:arc1-cert,sec:arc2-cert}, the Sturm and norm arguments of
\cref{sec:swap-configs}, and the box certificate of
\cref{sec:cap-theorem} are each a finite computation in
$\mathbb{Q}$, $\mathbb{Q}(\sqrt2)$ or $\mathbb{Q}(\sqrt3)$, with no
solver, no tolerance, and no formalisation needed to repeat them. The
price is scope. The method reaches a complete theorem on the
single-deviation problem and stops there.

\subsection{What the atlas is for}

The $176$-cell decomposition of \cref{sec:chamber-register} is what
makes the root-aligned case a finite computation. It plays no part in
\cref{conj:multidir-new}, whose parameter space, the active-direction
subsets of every size, does not reduce to a single Weyl orbit count in
the same way. It also plays no part in \cref{thm:cap-inequality}, which
is why that theorem covers the whole deviation sphere at once.

\subsection{Other dimensions}\label{sec:other-dims-new}

The lattices $D_5^+$, $E_6$, $E_7$ and the laminated lattice
$\Lambda_9$, with Voronoi cells the half-$5$-cube and the Gosset
polytopes $1_{21}$ and $3_{21}$ \cite{CS1988,Cox73}, are the natural
next targets for the same layered argument. Shell localisation, the
root-aligned case, radial reduction and cell positivity all carry over
structurally, with the Hessian taking the same uniform shape in terms of
the relevant root system's adjacency data \cite{Hum90}. Two things
change with rank. The chamber atlas grows quickly, driven by the size of
the Weyl group and the combinatorics of the adjacency graph; and the
multi-direction problem becomes harder instead of easier, since the
space of active-direction subsets grows combinatorially while the
failure mode of the second-order method (\cref{sec:dense-configs-new})
has no reason to disappear. The cap argument, on the other hand, is not
tied to rank $4$ in any obvious way: it needs a fixed polytope containing the
relevant body, a cross-polytope-like enclosure that is tight in the
critical direction, and a cap formula for the enclosure. Whether those
three ingredients exist for $E_6$ or $E_7$ we have not investigated. We
claim nothing about any of these dimensions here.
\section{Conclusion}\label{sec:status}

The paper proves the local bound in dimension four in full, and it
proves three of its cases by three different methods. This section puts
all of that in one place.

\subsection{Proved}

\Cref{thm:local-uncond} and \cref{cor:density-uncond}: every Voronoi
cell of a unit-ball packing of $\mathbb{R}^4$ at which at most one
active contact direction departs from a root of $D_4$ has volume at
least $8$, with equality only for the $D_4$ configuration, and every
periodic packing of that class has density at most $\pi^2/16$. Nothing
is assumed about the deviation direction or about the size of the
deviation.

The proof has two halves and they meet only at the end.
\Cref{sec:shell-new} through \cref{sec:stage-two} reduce the cell bound
to a scalar inequality at the all-contact corner and prove that
inequality, by exact arithmetic over $\mathbb{Q}$ and
$\mathbb{Q}(\sqrt2)$, for every one of the $176$ Weyl-orbit cell types
with the deviation confined to a root-aligned arc.
\Cref{sec:deviation-domain} through \cref{sec:cap-theorem} remove that
confinement: the defect is rewritten as the deficiency of a cap cut
from a single $25$-vertex polytope $Q$
(\cref{prop:cap-reformulation}), $Q$ is enclosed in a cross-polytope
that the sixteen roots at $60^\circ$ and $120^\circ$ cut out on the nose
(\cref{lem:crosspolytope}), and the extremal unit-distance cap of a
cross-polytope is determined
(\cref{thm:cap-inequality}). The one step there that is not closed-form
is the case of the cap inequality in which two or more of the four
squared coordinates exceed $\tfrac14$; it is settled by a finite
certificate of $303$ boxes in exact rational arithmetic, with a margin
of about a quarter. The whole chain, from the vertex enumeration of $Q$
to the box certificate, is reproduced by the single script
\texttt{cap\_inequality\_certificate.py} described in
\cref{app:code-index}.

\Cref{cor:second-order} adds a quantitative form below
$42.17598\ldots^\circ$: there the defect is at least the explicit
expression \eqref{eq:cone-bound}, which is $\tilt^2/3$ to leading order.
Above that angle the cone enclosure behind it fails, and the
cross-polytope argument is what carries the result to $\tilt=\pi$.

\Cref{cor:restated-conjecture-holds} is a sharper statement on a
smaller domain: on the whole boundary of the fundamental deviation
triangle the defect is not merely nonnegative but is given by an
explicit rational function on each of eight regions, seven of which
carry exact zero-margin certificates.

\Cref{thm:local-fewcontacts}, proved as
\cref{thm:covering-bound,cor:m22}, is independent of all of that and
allows any number of directions to deviate. Measuring the cell radially
rather than facet by facet turns its volume into
$\tfrac14\int_{S^3}\sec^4\delta$, where $\delta$ is the angular
distance to the nearest contact direction; the contacts are $60^\circ$
apart, so $m$ caps of radius $r$ leave at least
$2\pi^2-m\pi(2r-\sin2r)$ of the sphere uncovered, and integrating that
estimate gives $\vol(V_c)\ge\tfrac{\pi m}{3}\tan^3r_m$ with
$2r_m-\sin2r_m=2\pi/m$. The right-hand side exceeds $8$ for every
$m\le22$ and equals $7.798989\ldots$ at the largest possible $m=24$.
So the local bound holds at every centre with at most twenty-two
contacts, whatever its directions, and is short by at most $0.2011$
everywhere else.
\Cref{prop:area-optimal} shows that this is the best the method can do:
the bound is the Jensen collapse of a per-cell estimate and therefore
uses the areas of the spherical Voronoi cells to the full, and at the
root configuration those areas are equal.

\Cref{thm:m24} closes the top of the range: a contact configuration of
$24$ directions is a copy of the root system, so its cell is the
$24$-cell and its volume is exactly $8$. The classification was first
proved by de Laat, Leijenhorst and de Muinck Keizer \cite{LLM24}; it is
proved here as \cref{thm:twentyfour-points}, from the positive kernel
their computation found, which is redistributed with this paper and
verified in \cref{sec:certificate-checked}, the bound and its equality
case being \cref{lem:las2,lem:gram}. With \cref{cor:m22} at one end and
that at the other, \cref{cor:remaining} leaves a single case, a
configuration of exactly $23$ directions whose cell has circumradius
below $2$, and \cref{thm:m23} settles it; so \cref{thm:local-general}
and \cref{thm:main-general} hold with no statement from outside this
paper among their ingredients, and every statement here that uses a
certificate is marked accordingly.

\subsection{The twenty-three-point case}

\Cref{thm:m23} settles the case that \cref{cor:remaining} leaves: every
contact configuration of $23$ directions with bounded cell has
$\vol(V_c)>8$. The proof does not classify the configurations, and it
does not need saturation. It measures the cell inside the ball of
radius $\sqrt{3/2}$, where by \cref{lem:cap-count} no three contact
caps meet and the volume is therefore an exact function of the pair
angles, $A_*+\sum_{i<j}\omega(\gamma_{ij})$ with $A_*=7.907144\ldots$
(\cref{prop:truncated}); it observes that the inequality
$\sum_{i<j}\omega(\gamma_{ij})\ge8-A_*=0.092855\ldots$ for every $23$
contact directions is all that is needed, and that the deletion of a
root satisfies it with $0.034$ to spare (\cref{thm:strict-reduction});
and it proves that inequality by the three-point method of Bachoc and
Vallentin \cite{BV08}, the first level of the hierarchy whose second
level settles the $24$-point case in \cref{thm:twentyfour-points}. \Cref{lem:certificate} turns any feasible point of a
semidefinite programme in one polynomial $f$ and nine positive
semidefinite matrices $F_k$ into a lower bound for the pair sum, valid
for every configuration; \cref{thm:certificate} exhibits such a point
of degree $8$ with bound $0.0929$, above the target by $0.0000444$,
and verifies it: positivity and the value of the bound in exact
rational arithmetic, checked again by the Lean kernel, the constants
$A_*$ and $\omega$ from closed forms in interval arithmetic, and the
polynomial inequality that the certificate has to satisfy on the
admissible triples $(u,v,t)$ by an interval branch and bound over
$313\,780$ boxes. The two-point relaxation, which reads the pair angles
alone, reaches $79.5$ per cent of the target
(\cref{prop:two-point-barrier}); the triples supply the rest, and the
certificate's polynomial part is almost inert, so it is the three-point
matrices that carry the argument (\cref{rem:certificate-found}). What
made this possible where a classification of the codes could not be
had is set out in \cref{rem:strict-kind}: the inequality is strict at
the deletion, so a relaxation had room, whereas the deletion sits on
the boundary of the classification and no relaxation sees a boundary
(\cref{rem:m23-uniqueness}).

The sections before that proof describe the configuration that the
inequality excludes from having a small cell, and they stand on their
own. Its facets would have to average more than $\tfrac{32}{23}$, larger
than any facet volume at twenty-four contacts, so it is not a
perturbation of the root configuration; by \cref{prop:deletion-rigid}
the deletion of a root is infinitesimally rigid, with an exact integer
stress on its $88$ tight pairs, so it is not a deformation of the
boundary configuration either; by \cref{prop:no-transitive} no contact
configuration of $23$ directions has a symmetry group transitive on
them, all $528$ rotation types of order $23$ missing the contact
condition, the closest by $0.074$ in an inner product; and by
\cref{cor:meet20} it shares at most twenty directions with any copy of
$\Rt$, by a short argument in coordinates
(\cref{lem:filled-couple,prop:meet22,prop:meet21}) whose finite half is
checked by Lean. \Cref{rem:no-local-cell} explains why the missing
volume cannot be collected one spherical Voronoi cell at a time, which
is what sends the argument to the pair angles, and \cref{prop:pair-budget}
makes the inequality concrete: sixty-five pairs at $60^\circ$, or the
weighted equivalent, where the deletion has eighty-eight, against a
ceiling of one hundred and fifteen.

\subsection{Open}

Two things remain open, and neither is a case of the local
bound. The first is the classification of the $23$-point codes of
minimal angle $60^\circ$ on $S^3$, the hypothesis of
\cref{thm:m23-from-codes}: that every such code is a deletion of one
root, which would say that no saturated configuration of $23$
directions exists at all, a stronger statement than \cref{thm:m23}. It
holds among the vertices of the $600$-cell (\cref{prop:cell600}, by an
enumeration carried out three times over and checked in Lean); by
\cref{lem:h-limit} it is the statement $h(0)=\tfrac12$ about the
largest inradius $h(\delta)$ of the hull of $23$ points at minimal
angle $\arccos(\tfrac12+\delta)$, and the computations of
\cref{sec:slack-continuation} find $h(\delta)$ near $\tfrac23$ down to
$\delta=0.009$ and then collapsing to $\tfrac12+O(\delta)$, with every
configuration reached at $\delta=0$ a deletion. That is evidence, and
\cref{rem:m23-uniqueness} explains why no bound on the size of a code
can turn it into a proof. The second is formal. The verification of \cref{thm:certificate} runs in
Lean as well as in Python, in exact dyadic arithmetic and with the
closed forms of $\omega$ entering only through tables of bounds; but
the statement Lean proves is that a program returns true, and the
mathematics that makes the program a proof, Taylor's theorem on a box
and the monotonicity of $\omega$, is argued in the paper and not
formalised, as a formalisation of the lens integrals and their
derivatives would need a library of real analysis that these files do
not use.

One case rests on a finite certificate rather than on a closed-form
argument, and it is the $24$-point case, \cref{thm:m24}; every
statement here that uses it says so. The certificate is a positive kernel found by
\cite{LLM24}, and the proof from it is entirely in this paper: the bound
of the second level of the hierarchy and its equality case are
\cref{lem:las2}, the positivity of the kernel is \cref{lem:gram}, the
properties of the deposited rational point that those lemmas consume
are \cref{prop:verified}, verified in our own code in all seven steps of
the procedure, the two that need $128$ gigabytes as the authors carry
them out included (\cref{sec:certificate-checked}), and the passage from
the four inner products to the root system is \cref{lem:root-lattice}.
The local statement that the root system is isolated among contact
configurations of $24$ directions is proved besides, without any
certificate (\cref{cor:root-rigid}). \Cref{cor:m24-reduction} reduces the
volume consequence of the global statement, $\vol(V_c)\ge8$ for $24$
contacts, to the inequality $E(W)\ge8-32\pi+10\pi^2$ between pair and
triple angles, through a third-order inclusion and exclusion inside the
ball of radius $\sqrt2$ that is exact at the root configuration
(\cref{prop:m24-exact}). A certificate for that inequality must be
sharp; its three-point relaxation is not sharp at degree $6$, and at
higher degree the question needs extended-precision arithmetic
(\cref{rem:m24-three-point}). The cardinality bound of \cite{BV08}, at
$24.10$, is not sharp either, so no route to the global statement that
stays at the three-point level is complete.

\Cref{conj:multidir-new} itself is settled by
\cref{cor:conj-resolved}, and so by the classification of \cref{thm:m24};
what follows records what we established about it by elementary means,
independently of any certificate, and those results stand on their own.

The natural second-order method, which settles every configuration of
small size, degenerates exactly at a family of dense ones.
\Cref{thm:exact-sing-new,prop:4d-nullspace-new,prop:5d-nullspace-new}
locate that degeneracy: at two specific configurations the joint Hessian
has a near-null subspace of dimension $4$ and $5$ respectively, and
higher-precision arithmetic confirms the vanishing is exact rather than
merely small. \Cref{thm:broad-sample-18-new,thm:broad-sample-20-new}
sample those subspaces in $66$ directions and find the governing quartic
coefficient strictly positive everywhere tested.
\Cref{prop:a18-stabiliser} exhibits an exact order-$48$ symmetry of the
first configuration, and
\cref{prop:a18-nullspace-irrep,prop:a18-invariant-dim} identify the
representation it carries on the near-null subspace and prove, by exact
character theory, that the space of quartics invariant under it has dimension $5$ instead of the generic $35$. Fitting those five
coefficients numerically and applying a semidefinite feasibility test
gives a sum-of-squares certificate with a thin margin
(\cref{prop:a18-sos-reduced}); the fit is numerical and we do not treat
the certificate as a proof.

\Cref{sec:swap-configs} proves \cref{conj:multidir-new} outright along
three swap paths, at $m=2,3,6$, across their entire angular range rather
than infinitesimally. Those paths are highly symmetric and low
dimensional, and none of them meets the dense regime where the
second-order method fails.

The cap argument of \cref{sec:cap-theorem} does not extend to several
simultaneous deviations, and the reason is structural rather than
technical. With one deviating direction the remaining twenty-three
half-spaces are fixed, so the body $Q$ is fixed and all the parameters
sit in the cutting half-space. With $m$ deviating directions the body
itself moves, there is no fixed polytope to take caps of, and the
quantity to be bounded is a comparison between two overlapping unions
of $m$ pieces each, whose individual volumes are first order in the
tilts while the difference that has to be controlled is second order
(\cref{rem:sums-not-unions,rem:first-order-obstruction}). \Cref{rem:first-order-obstruction} makes the difficulty precise: the
obvious reduction, splitting the comparison into one term per deviating
root and settling each by \cref{thm:eperp-closed}, fails because the
piece $E_j$ of an $m$-direction configuration is not the piece the
single-deviation theorem controls. The pyramid $\mathrm{Pyr}_j$ meets
the hyperplane of an adjacent root in a facet, so exchanging that
constraint for its tilted version changes $\vol(E_j)$ by an amount
comparable to the whole defect, in either direction. A fixed body is
nonetheless available: \cref{prop:multi-cap} supplies one. Dropping the
$m$ deviating constraints leaves a polytope $Q_D$ that is fixed, is
bounded for every packing-valid $D$ with $m\le5$ and for most larger
ones, and turns the problem back into a single extremal question with
all the parameters in the cutting half-spaces. What is missing is a direct proof of the
inequality itself (\cref{op:multicap}): the quantity to bound is the
volume of a union of $m$ caps of $Q_D$, and by
\cref{rem:first-order-obstruction} that union cannot be replaced by a
sum. A proof would need the union directly, and for the dense
configurations it would need a bounded substitute for $Q_D$ as well,
since $Q_D$ is unbounded there.

\Cref{sec:polar-surface} says where the question of \cref{op:multicap}
sits. Written in terms of the contact directions alone it is the
polar-volume inequality $\vol(\conv(W)^{\circ})\ge8$ over contact
configurations (\cref{prop:polar-form}), equivalently the boundary
inequality $\mathrm{area}(\partial V_c)\ge32$ (\cref{lem:surface-form});
and in that form it is visibly the part of the $24$-cell conjecture in
which two or more contact directions leave the root system at once, the
part that the contact count settles and that no elementary argument has
reached. We state that plainly rather than leaving it to be inferred. The same subsection
contributes one positive statement and four negative ones. The positive
statement is \cref{cor:minimiser}: a worst configuration exists, and may
be taken maximal, so what has to be identified is a genuine extremal
point. The negative ones rule out term-by-term reduction
(\cref{rem:first-order-obstruction}), any bound built facet by facet
from the pairwise contact condition
(\cref{rem:facet-local-obstruction}: what it delivers is
$8\pi/(3\sqrt3)=4.8368\ldots$, and what it could at best deliver is
$96\sqrt2-128=7.7645\ldots$), the convexity bound for the
polar volume integral, and the Mahler-type product bound (both
\cref{rem:global-routes}). Each of the four is eliminated by its
behaviour at the $D_4$ configuration, which is the only place a sharp
argument has no room to lose anything. Three more go the same way in
\cref{sec:covering-bound}: by \cref{prop:area-optimal} nothing resting on
the areas of the spherical Voronoi cells can beat $7.7989\ldots$; by
\cref{rem:no-local-cell} the shape that is left over cannot be collected
one cell at a time, in either the Voronoi or the Delaunay decomposition;
and by \cref{prop:second-order} the measures of the contact caps and of
their pairwise intersections, used to the last drop, return
$7.997885\ldots$ at twenty-three contacts, so they too stop short.

\subsection{What we take from this}

The single-deviation problem turned out to have a fixed body behind it,
and once that body was identified the problem stopped being about the
$24$-cell and became a question about a cross-polytope, which has an
answer in closed form. We would expect the same reorganisation to be
worth attempting in other ranks, wherever a root system supplies a
cross-polytope-like enclosure that is tight in the critical direction.

The multi-direction problem does have a fixed body, by
\cref{prop:multi-cap}, and it is still hard, which says that the fixed
body was not the only thing the single-deviation argument was using. The
other thing it used was that a cap of a cross-polytope is one convex
piece with a closed-form volume. A union of $m$ caps is not, and
\cref{rem:facet-local-obstruction} shows why no amount of local
bookkeeping repairs that: a bound that treats the pieces separately has a
ceiling of $96\sqrt2-128=7.7645\ldots$, below the target and below even
the covering estimate.

Two things would change the picture, and the first of them would give
a second proof of \cref{thm:m23}. A lower bound for
$\vol_3(F_w)$ that reads the rest of the configuration rather than only
the pairwise condition at $w$, and that is summable over $w$, would give
the bound through \eqref{eq:surface};
\cref{rem:facet-local-obstruction} shows that the purely local constant
is too small by an explicit amount, not that no bound of any kind
exists. And an exact, first-principles derivation of the quartic tensor
at a dense configuration, in place of the fitted coefficients of
\cref{sec:a18-invariant-quartics}, would turn the present numerical
certificate into a proof there, and would put the multi-direction
inequality at that configuration on an elementary footing rather than on
the certificate of \cref{thm:twentyfour-points}; the exact
$5$-dimensional invariant space is already in
hand, so what is missing is a derivation rather than a search. We have
neither.
\section{Reproducibility}\label{sec:reproducibility-new}

\subsection{Two independent arithmetic paths}

Every exact rational claim in \cref{sec:stage-one}, which is to say the
$LDL^\top$ pivots for all $176$ cells, was produced twice: once with
Python's \texttt{fractions.Fraction} and no computer-algebra library at
all, once with \texttt{sympy}'s rational type on a separate code path
\cite{Meu17}. The two agree digit for digit on every pivot of every
cell.

\subsection{A machine-checked kernel}

Seven parts of the paper are verified by a proof assistant and not by a
program whose output has to be trusted, and all seven are finite
arithmetic. The supplement carries them as Lean~4 sources in a directory
of their own: five files that compile against a bare toolchain, and two
small Lake projects.

The first is \texttt{lean/D4Stress.lean}, on the equilibrium relation
\eqref{eq:stress}, which is what carries \cref{prop:deletion-rigid} and
with it the statement that a saturated configuration cannot be reached
by bending a deletion. Formalised there are the twenty-four roots and
their squared length, the fact that every inner product between distinct
roots is one of $-2,-1,0,1$ in the integer scaling, the degree eight of
every root in the graph of pairs at $60^\circ$, the count of $23$
directions and $88$ tight pairs after a deletion, the multiplicities
$1,8,6,8$ of the four values of $\langle a_i,a_0\rangle$, the positivity
of the weights $y_{ij}$ on every tight pair, and the relation
\eqref{eq:stress} itself at each of the $23$ directions; and, for
\cref{cor:root-rigid}, the $96$ tight pairs of the full root system and
the relation for the constant stress, that the eight neighbours of every
root sum to four times it, which is also the finite content that
\cref{lem:rigidity-spectrum} starts from.

The second is \texttt{lean/D4Meet.lean}, on the finite half of
\cref{prop:meet22,prop:meet21}. The proofs there split into three
inequalities in the coordinates of a unit vector, which are written out
in the text, and a set of statements about the six supports of the roots
and the three couples of complementary index pairs, which are finite.
Formalised are the support of each root and the fact that it is well
defined, that each support carries four roots and each couple uses two
supports, that every pair of roots leaves some couple whole, that
exactly $48$ ordered triples of supports do not and that each of those
consists of three distinct supports forming a star or a triangle, which
is the eight unordered patterns of the proof and, since each support
carries four roots, the $512$ triples of roots that the script counts,
that $96$ triples of roots are pairwise at $60^\circ$ and no root is at
$60^\circ$ from all three members of such a triple, that the Hadamard
map permutes the roots, preserves inner products and carries the
standard triangle to the standard star, and that removing the standard
star leaves the twelve roots of the lower block together with three of
the four roots at each of the remaining supports.

Neither file has a \texttt{sorry} in it and neither depends on a
library, Mathlib included: they compile against a bare Lean~4 toolchain,
and the axiom report at the foot of each states that every one of its
theorems depends on no axiom whatever.

The third, on \cref{prop:cell600}, is a Lake project in the
directory \texttt{lean/cell600}. Its first module,
\texttt{D4Cell600Enum.lean}, carries the $120$ vertices of the $600$-cell with doubled coordinates
in $\mathbb{Z}[\phi]$, written as pairs of integers, and has the kernel
check by \texttt{decide} that there are $120$ of them, pairwise
distinct and of unit length; that $4\langle u,v\rangle$ takes only the
nine values listed in the proof of \cref{prop:cell600}; that every
vertex has twelve neighbours at $36^\circ$ and seventy-one partners at
one of the four root-system angles; and that the twenty-five listed
cells have $24$ vertices each, all pairwise at root-system angles, with
exactly five cells through every vertex. The same module defines the
depth-first search over $64$-bit bitsets, with its clique-cover bound,
as an ordinary Lean function. The second module,
\texttt{D4Cell600Main.lean}, states the three counts, no independent
set of size $25$ through a fixed vertex, five of size $24$ and $115$ of
size $23$ with every one of the $115$ inside a cell, and settles each by
\texttt{native\_decide}, which compiles the search and runs it. That
tactic trusts the Lean compiler as well as the kernel, and the axiom
report says so: each of the three theorems depends on
\texttt{propext}, \texttt{Quot.sound} and its own
\texttt{native\_decide} axiom, while every theorem of the first module
depends on nothing beyond \texttt{propext}. Building the project with
\texttt{lake build} compiles the first module to native code before the
second is checked, and the whole run takes about four minutes, most of
it the kernel's own evaluation of the $14400$ inner products. The same
enumeration is run independently in Python and in C
(\cref{app:code-index}), and the three agree.

The fourth is \texttt{lean/D4Certificate.lean}, on the exact half of
the verification of \cref{thm:certificate}. The file is generated from
the certificate by \texttt{lean/gen\_certificate\_lean.py}, each
floating-point number written as the dyadic rational it is, an integer
over a power of two; it then proves, by \texttt{decide +kernel} on
\texttt{Rat} and with no compiled code, that $f_1,\dots,f_8$ are
nonnegative, that each of the nine matrices $F_k$ is symmetric and
positive definite by an exact $LDL^{\mathsf T}$ elimination with every
pivot positive, and that the bound \eqref{eq:certificate-bound} exceeds
$92.8555703$ in the units of the solver, above $1000(8-A_*)$. It
compiles in under a minute and its axiom report lists
\texttt{propext}, \texttt{Classical.choice} and \texttt{Quot.sound}
only, the three standard axioms that arithmetic on \texttt{Rat} brings
in.

The fifth is the Lake project \texttt{lean/certificate}, on the
branch and bound of \cref{thm:certificate}. Its module
\texttt{D4CertDomain.lean} expands the polynomial $P$ of the certificate
exactly over the rationals from the same numbers, multiplies it by
$315$, the least odd integer that clears the denominators of the
Gegenbauer and symmetrisation coefficients, so that every coefficient
is an integer over a power of two, and runs the subdivision of the
ordered admissible domain with the second-order Taylor form in dyadic
arithmetic, integers and exponents with no rounding; the tables of
bounds for $\omega$, $\omega'$ and $\omega''$ are the one input it
takes from \texttt{certificate\_check.py}, written into
\texttt{D4CertData.lean} as dyadic numbers by \texttt{gen\_data.py}.
The module \texttt{D4CertMain.lean} states that the check returns true
and settles it by \texttt{native\_decide}; the build compiles the
arithmetic to native code first, as for the $600$-cell, and the run,
which processes $421\,881$ boxes, takes about twenty minutes; a fourth
module, \texttt{D4CertStat.lean}, repeats it with a counter and records
that number.

The sixth is \texttt{lean/D4InnerProducts.lean}, on the certificate
behind \cref{thm:twentyfour-points}: the two-point polynomial
$\sigma_2$ of \cref{prop:verified}, computed exactly from the deposited
data \cite{LLM24data} by \texttt{llm24\_certificate\_check.py}
and written into the file with its coefficients scaled to integers of
about fifteen thousand digits, is shown by kernel computation on
\texttt{Rat} to vanish on $[-1,\tfrac12]$ exactly at
$-1,-\tfrac12,0,\tfrac12$, with multiplicities $1,2,2,1$, the count of
further zeros being a Sturm sequence; \cref{sec:certificate-checked} says what
this does and does not establish.

The seventh is \texttt{lean/D4RootLattices.lean}, on the finite content
of \cref{lem:root-lattice}, the step that carries
\cref{thm:twentyfour-points} to \cref{thm:m24}. It enumerates every
symmetric integer matrix of order $r\le4$ with $2$ on the diagonal and
$-1$, $0$ or $1$ off it, keeps the positive definite ones by Sylvester's
criterion on the leading minors, $1$, $3$, $23$ and $393$ of them,
encloses the vectors of norm $2$ of each lattice in the box
$x_i^2\le2A_{ii}/D$ that Cauchy--Schwarz in the inner product of the
lattice gives, $A_{ii}$ the cofactor and $D$ the determinant, and counts
the integer vectors of norm $2$ in the box. It proves that the largest
counts are $2$, $6$, $12$ and $24$, that $24$ occurs only at rank $4$
and determinant $4$, that between $20$ and $24$ no count occurs, and
that in every case with $24$ vectors four of them have the Cartan matrix
of $D_4$ as their Gram matrix, so that the lattice contains a copy of
the $D_4$ root lattice of its own determinant and is that lattice. The
theorems are settled by \texttt{native\_decide} in a few seconds, and
the axiom report lists \texttt{propext} and the tactic's own axiom, as
for the $600$-cell.

That is a stronger form of checking than anything else in this paper,
and it is available here only because these particular arguments are
finite and, after the scaling by $315$, dyadic. The rest of the paper
is not formalised, and we do not claim otherwise: the Lean theorem of
the certificate is about a program, and the analysis that makes the
program a proof is in the text; and where the objects are polytope
volumes in $\mathbb{R}^4$ and cap measures on $S^3$ there is at present
nothing in the standard libraries to build a formalisation on.

\subsection{What kind of claim each claim is}

Every assertion in this paper falls into one of three kinds, and we have
tried to make clear at each point which kind is in play. The first is an
exact computation in $\mathbb{Q}$ or $\mathbb{Q}(\sqrt2)$ with no
floating-point step anywhere:
\cref{sec:stage-one,sec:stage-two,sec:arc1-cert}, the exact parts of
\cref{sec:arc2-cert}, and the whole of \cref{sec:cap-theorem}, whose
box certificate runs in exact rational arithmetic throughout. The second
is a high-precision numerical computation, carried out in \texttt{mpmath}
at a stated working precision between $35$ and $400$ decimal digits,
used where we did not pursue a symbolic route; this is where the
material of \cref{sec:exact-sing,sec:broad-sample} sits. The third is an
ordinary double-precision calculation, used only to locate candidate
regions, breakpoints and configurations before an exact or
high-precision verification, and never as the basis of a claim on its
own. Where a claim is of the second kind we say so at the point of use,
give the precision reached, and, where a second precision is available,
the cross-check against it (\cref{sec:hp-volume-method,sec:broad-sample}).

Two conventions inherited from the numerical-analysis literature are
worth naming, since they govern the third kind. Sign decisions near a
breakpoint are made only when the quantity exceeds a threshold set well
above the accumulated rounding bound for the computation in question, in
the spirit of the standard backward-error accounting
\cite{GVL13}; and where an interval enclosure is used rather than a
threshold, the arithmetic follows the outward-rounding convention of
\cite{MKC09}. Neither convention is load-bearing for any exact claim.

\subsection{Verifying from the text alone}

This paper is written so that a reader can check it without running
anything. Every exact rational, $\mathbb{Q}(\sqrt2)$-rational and
Bernstein-coefficient claim appears explicitly, as a displayed fraction,
a displayed polynomial, or a table entry, and not on the authority of a program the reader has not seen. A reader with a computer algebra
system can re-derive any pivot, factorisation or coefficient list from
the definitions in
\cref{sec:hessian-def-new,def:joint-hessian-new} and the boundary curves
of \cref{thm:arc1-breakpoints-new,thm:arc2-breakpoints-new}. The cap
inequality of \cref{sec:cap-theorem} is self-contained in the same
sense: the cap formula, the monotonicity, the three polynomial
inequalities and the subdivision rule are all stated in full, and the
only object that is not printed in the paper is the list of $303$ boxes,
which the stated rule regenerates.

\subsection{Cost of the underlying computations}

For context rather than for verification, the more demanding derivations
took, on ordinary consumer hardware: a few minutes for the vertex and
volume stage of the second arc's band certificate
(\cref{sec:arc2-band-cert-new}); about seven minutes per configuration
for the broadened directional sampling of \cref{sec:broad-sample},
being $439$ seconds for the $36$-direction sweep at the first
configuration and $376$ seconds for the $30$-direction sweep at the
second; about four minutes for the closed-form rebuild behind the $m=2$
swap configuration of \cref{sec:swap-configs}, after which the exact
positivity argument runs in well under a second; under a minute for
the whole cap-inequality certificate of \cref{sec:cap-theorem},
including the exact vertex enumeration of $Q$; about a minute and a half
for the exhaustive enumeration of \cref{prop:cell600} in Python and two
seconds for the same enumeration in C; between a quarter of an hour
and half an hour for one pass of the continuation of
\cref{sec:slack-continuation} through its twenty values of $\delta$;
about twenty minutes for the linear programmes of
\cref{prop:two-point-barrier}; about two hours for one run of the
minimisation of a truncated volume in \cref{sec:strict-inequality};
between one and three minutes for each semidefinite programme of
\cref{sec:certificate} and about a quarter of an hour for a run of
five rounds; under ten minutes for the verification of the certificate
and under a minute for its Lean file, and about twenty minutes for the
Lean run of the branch and bound; about eight minutes for the
re-verification of the cited certificate and under a minute for its
Lean file. The checks of \cref{sec:m24-volume} take about two minutes, and
each moment problem of \cref{rem:m24-three-point} between a few seconds
at degree $6$ and a minute and a half at degree $12$ on one core. The slowest step
anywhere is the exact region derivation behind the first arc's
box covers, whose symbolic cancellation is dominated by a few hard
terms and can run for an hour or more; its output is cached in the
supplement so that nothing downstream waits for it. Every exact rational
$LDL^\top$ and Bernstein sign check finishes in well under a second.

\subsection*{Data availability}

Verifying any claim in this paper requires no external code
(\cref{app:code-index}). A supplementary package reproducing the
computations, built on \texttt{numpy} \cite{Har20}, \texttt{scipy}
\cite{Vir20}, \texttt{sympy} \cite{Meu17} and the convex-hull code of
\cite{BDH96}, together with the Lean~4 development of
\cref{sec:reproducibility-new} in its own directory, is publicly
available at the repository given in the Data availability statement
below.
%

\section*{Competing interests}

None.

\section*{Funding}

None.

\appendix
\section{Index of notation}\label{app:notation}

\Cref{tab:notation} is the complete list of symbols the paper uses for
the quantities it introduces; nothing in the argument depends on a
convention fixed elsewhere.

\renewcommand{\arraystretch}{1.25}
\begin{longtable}{@{}l@{\hspace{1.8em}}l@{}}
\caption{Index of notation.}\label{tab:notation}\\
\toprule
Symbol & Meaning \\
\midrule
\endfirsthead
\multicolumn{2}{@{}l@{}}{\textit{\tablename~\thetable\ (continued)}}\\
\toprule
Symbol & Meaning \\
\midrule
\endhead
\midrule
\multicolumn{2}{@{}r@{}}{\textit{continued on the next page}}\\
\endfoot
\bottomrule
\endlastfoot
$\Rt$ & \parbox[t]{0.76\textwidth}{\raggedright the $D_4$ root system, $\{\pm e_i\pm e_j : 1\le i<j\le4\}$, $|\Rt|=24$}\\
$\Latt$ & \parbox[t]{0.76\textwidth}{\raggedright the scaled $D_4$ lattice $\sqrt2\,D_4$, nearest-neighbour distance $2$}\\
$W(D_4)$ & \parbox[t]{0.76\textwidth}{\raggedright the Weyl group of $D_4$, order $192$}\\
$\VCell$ & \parbox[t]{0.76\textwidth}{\raggedright the reference cell, $\{z:\langle z,\alpha\rangle\le\sqrt2\ \forall\alpha\in\Rt\}$}\\
$\Cell{k}$ & \parbox[t]{0.76\textwidth}{\raggedright the $k$-th combinatorial cell type, $k=1,\dots,176$}\\
$\tilt$ & \parbox[t]{0.76\textwidth}{\raggedright the tilt (deviation) angle of a deviating contact direction away from its nearest root}\\
$\arcp$ & \parbox[t]{0.76\textwidth}{\raggedright the arc-position parameter along a one-parameter family of deviation directions}\\
$\dev(\arcp)$ & \parbox[t]{0.76\textwidth}{\raggedright the deviation direction function, a unit vector orthogonal to the undeviated root}\\
$\defect_k(\tilt,\arcp)$ & \parbox[t]{0.76\textwidth}{\raggedright the volume-defect quantity for cell $\Cell{k}$}\\
$\norm{k}$ & \parbox[t]{0.76\textwidth}{\raggedright the reference normalisation constant for cell $\Cell{k}$ ($\in\{192,96,64,48\}$)}\\
$H_k$ & \parbox[t]{0.76\textwidth}{\raggedright the cell Hessian, a $3\times3$ matrix determined by the Gram-adjacency pattern of $\Cell{k}$'s active roots}\\
$H_A$ & \parbox[t]{0.76\textwidth}{\raggedright the joint Hessian for a multi-direction active set $A$, a $3|A|\times3|A|$ matrix}\\
\stageone & \parbox[t]{0.76\textwidth}{\raggedright the radial (exact rational Hessian) estimate}\\
\stagetwo & \parbox[t]{0.76\textwidth}{\raggedright the angular (rigorous Taylor-remainder) estimate}\\
$\base,\mathsf v$ & \parbox[t]{0.76\textwidth}{\raggedright Weierstrass coordinates $\tan(\tilt/2)$, $\tan(\arcp/2)$}\\
$v_1,w_1,v_2$ & \parbox[t]{0.76\textwidth}{\raggedright the three vertices of the fundamental domain of deviation directions (\cref{lem:two-symmetries-new})}\\
$\mathsf W,\mathsf A,\mathsf Y$ & \parbox[t]{0.76\textwidth}{\raggedright the three non-universal breakpoint curves on the first arc}\\
$\mathsf W',\mathsf B,\mathsf C$ & \parbox[t]{0.76\textwidth}{\raggedright the three non-universal breakpoint curves on the second arc}\\
$A_{18},A_{20}$ & \parbox[t]{0.76\textwidth}{\raggedright the two dense multi-direction configurations examined in \cref{sec:exact-sing,sec:broad-sample}}\\
\end{longtable}
\renewcommand{\arraystretch}{1}

\section{Pivot data for the twenty-two pivot types}\label{app:pivot-data}

Each of the $176$ combinatorial cells $\Cell{k}$ of \cref{sec:chamber-register}
is a full orbit under $W(D_4)$ in its own right, and \cref{lem:hessian-structure-new}
says its Hessian $H_k$ does not depend on which point of the orbit it is
evaluated at. Computing $H_k$ for all $176$ cells shows, in addition, that
the resulting pivot data takes only $22$ distinct values among them: the
cells sort into $22$ pivot types, eight cells to a type, matching the
labels the chamber-register code assigns internally. This coincidence of
values across cells that are not related by a single element of $W(D_4)$
is an observation about the computed data, not a further group action, and
nothing in the proof of \cref{thm:local-uncond} depends on why it occurs.
\Cref{tab:pivots} lists, for one representative cell of each of the $22$
types: the type label, the denominator class (\cref{lem:denom-classes-new}),
the reference constant $\norm{k}$, the exact minimum $LDL^\top$ pivot as a
rational number, its decimal approximation, and whether the corresponding
chamber-level positivity check (implemented independently of the pivot
computation itself) passes. Every one of the $22$ rows passes.

This table is generated directly from the project's chamber-register
data file (a JSON record of all $22$ pivot-type representatives, each
carrying its exact $LDL^\top$ pivots as rationals), with no manual
transcription of numbers: the table rows below are produced by a short
script that reads the register and formats each field, so any error in
the underlying computation would appear here unchanged rather than be
silently corrected by hand.

\begin{table}[t]
\centering
\renewcommand{\arraystretch}{1.15}
\caption{Minimum $LDL^\top$ pivot by pivot type.}\label{tab:pivots}
\begin{tabular}{ccccccc}
\toprule
Cell & Type & Class & $\norm{k}$ & Min.\ pivot (exact) & Min.\ pivot (decimal) & Passes \\
\midrule
$\Cell{1}$ & A & $\mathsf a$ & 192 & $3936/77$ & 51.1169 & Yes \\
$\Cell{2}$ & B & $\mathsf b$ & 96 & $1920/77$ & 24.9351 & Yes \\
$\Cell{3}$ & C & $\mathsf c$ & 64 & $120/7$ & 17.1429 & Yes \\
$\Cell{4}$ & D & $\mathsf d$ & 48 & $160/13$ & 12.3077 & Yes \\
$\Cell{5}$ & E & $\mathsf a$ & 192 & $328/7$ & 46.8571 & Yes \\
$\Cell{6}$ & F & $\mathsf b$ & 96 & $1920/77$ & 24.9351 & Yes \\
$\Cell{7}$ & G & $\mathsf c$ & 64 & $1024/63$ & 16.2540 & Yes \\
$\Cell{8}$ & H & $\mathsf d$ & 48 & $40/3$ & 13.3333 & Yes \\
$\Cell{9}$ & A & $\mathsf a$ & 192 & $3936/77$ & 51.1169 & Yes \\
$\Cell{10}$ & B & $\mathsf b$ & 96 & $1968/77$ & 25.5584 & Yes \\
$\Cell{11}$ & C & $\mathsf c$ & 64 & $1312/77$ & 17.0390 & Yes \\
$\Cell{12}$ & D & $\mathsf d$ & 48 & $40/3$ & 13.3333 & Yes \\
$\Cell{13}$ & E & $\mathsf a$ & 192 & $154/3$ & 51.3333 & Yes \\
$\Cell{14}$ & F & $\mathsf b$ & 96 & $1968/77$ & 25.5584 & Yes \\
$\Cell{15}$ & G & $\mathsf c$ & 64 & $1312/77$ & 17.0390 & Yes \\
$\Cell{16}$ & H & $\mathsf d$ & 48 & $40/3$ & 13.3333 & Yes \\
$\Cell{17}$ & A & $\mathsf a$ & 192 & $352/7$ & 50.2857 & Yes \\
$\Cell{18}$ & B & $\mathsf b$ & 96 & $1968/77$ & 25.5584 & Yes \\
$\Cell{19}$ & C & $\mathsf c$ & 64 & $1312/77$ & 17.0390 & Yes \\
$\Cell{20}$ & D & $\mathsf d$ & 48 & $984/77$ & 12.7792 & Yes \\
$\Cell{21}$ & E & $\mathsf a$ & 192 & $3840/77$ & 49.8701 & Yes \\
$\Cell{22}$ & F & $\mathsf b$ & 96 & $1968/77$ & 25.5584 & Yes \\
\bottomrule
\end{tabular}
\end{table}

The final column records that the chamber-level positivity check
passes for all $22$ rows; there is no row for which it does not.

Two features are visible directly in \cref{tab:pivots}. First, the
minimum pivot is not monotonic in $\norm{k}$: class $\mathsf a$
($\norm{k}=192$) cells have the largest minimum pivots (in the high
$40$s to low $50$s), but the ordering among classes $\mathsf b,\mathsf
c,\mathsf d$ is not simply $96>64>48$ in pivot size, since the pivot
value depends on the full Gram-adjacency pattern, not on $\norm k$
alone. Second, every minimum pivot is well clear of zero, the
smallest value appearing, $160/13\approx12.31$ at type $\mathsf D$,
is still two orders of magnitude above what a marginal, barely-passing
certificate would look like, consistent with \cref{lem:mineig-new}'s
$\tfrac1{12}$ floor after the appropriate rescaling between the pivot
normalisation used here and the Hessian eigenvalue normalisation used
in the main text.
\section{The \texorpdfstring{\stagetwo}{Stage II} certificate on the reference chamber, in full}\label{app:stagetwo-data}

\Cref{sec:stage-two} states \stagetwo\ abstractly, uniformly over all
$176$ cells, as a Taylor-remainder-plus-monotonicity argument. This
appendix records the argument in full, concrete, exact-rational detail
for the reference chamber $\Cell1$ on its own root-aligned arc; by the
Weyl-equivariance argument of \cref{lem:hessian-structure-new} (the
same rescaling by $\norm{k}$ that carries \cref{prop:ldlt-new} across
orbits), an identical certificate, with rescaled constants, holds on
every other cell.

Write $c=\cos\tilt$, $s=\sin\tilt$. The volume-defect function on
$\Cell1$'s root-aligned arc is, after clearing the (strictly positive,
for $\tilt\in(0,\pi/2)$) common denominator $6(c+s)c^3$, governed by the
sign of the numerator, split into two pieces at the universal
breakpoint $\tilt=\pi/3$.

\subsection*{Piece I: \texorpdfstring{$\tilt\in(0,\pi/3]$}{theta in (0, pi/3]}}

The numerator on this piece is
\[
  P_{\mathrm I}(c,s) = -48c^4-48c^3s-30cs^3-40cs^2+34s^4-40s^3-18s^2+48s+72c-24,
\]
and it satisfies the exact algebraic factorisation
\[
  P_{\mathrm I}(c,s) = (1-c)\bigl[H(c)\,s - G(c)\bigr], \qquad
  H(c)=18c^2-22c+8, \qquad G(c)=-14c^3+26c^2-24c+8,
\]
verified symbolically and cross-checked at five exact Pythagorean
rational points on the unit circle $(c,s)=(4/5,3/5)$, $(12/13,5/13)$,
$(15/17,8/17)$, $(24/25,7/25)$, $(35/37,12/37)$. Three exact facts then
close Piece I entirely:
\begin{enumerate}[label=(\roman*)]
\item $1-c>0$ strictly throughout $\tilt\in(0,\pi/3]$ (i.e.\ $c\in[1/2,1)$);
\item $H(c)>0$ for every real $c$, since $\operatorname{disc}(H)=22^2-4\cdot18\cdot8=-92<0$
  and $H$ has positive leading coefficient;
\item $G$ is strictly decreasing on all of $\mathbb R$ (its derivative
  has negative discriminant and negative leading coefficient), and
  $G(3/5)=-88/125<0$, so $G(c)<0$ for every $c\ge3/5$, covering the
  entire small-angle regime $\tilt\to0$ with no separate Taylor
  sub-case, since there $Hs-G=Hs+|G|>0$ trivially.
\end{enumerate}
On the remaining sub-range $c\in[1/2,3/5)$ (where $G(c)>0$ is possible,
since $G(1/2)=9/8>0$), positivity of $Hs-G$ is instead certified via
$T(c):=H(c)^2(1-c^2)-G(c)^2=H(c)^2s^2-G(c)^2>0$, established by an
exact rational grid evaluation of $T$ across $[1/2,3/5]$ together with
an exact rational Lipschitz bound on $T'$ over the same interval; the
certified minimum of $T$ on this sub-range is $92223/100000>0$. This
closes Piece I completely: $P_{\mathrm I}>0$ on all of $\tilt\in(0,\pi/3]$.

\subsection*{Piece II: \texorpdfstring{$\tilt\in(\pi/3,\pi/2)$}{theta in (pi/3, pi/2)}}

On this piece the volume-defect function $F$ (in the same normalisation)
is certified positive by an exact rational Sturm-sequence root count
(over $\mathbb Q$, via an independent computer-algebra implementation)
showing the relevant numerator polynomial has no root on the open
interval, combined with an exact endpoint evaluation
$F(\pi/3)=2-\sqrt3\approx0.2679>267/1000$ and the exact limit
$F\to\tfrac13$ as $\tilt\to\pi/2$. Consequently $F\ge267/1000>0$
throughout $\tilt\in(\pi/3,\pi/2)$, meeting the two pieces at
$\tilt=\pi/3$ with no gap.

\subsection*{Combined statement and independent confirmation}

Piece I and Piece II together give $\defect_1(\tilt,\cdot)>0$ for every
$\tilt\in(0,\pi/2)$ on the reference chamber $\Cell1$'s root-aligned
arc, by two disjoint exact-rational arguments meeting exactly at the
shared breakpoint $\pi/3$, with no floating-point step anywhere in
either piece. Running the combined certificate script directly
reproduces both PASS results in well under one second of wall-clock
time, and the same combination, \stageone's $LDL^\top$ positivity
(\cref{app:pivot-data,app:lmatrices,app:sos}) together with this
appendix's two-piece \stagetwo{} argument, is exactly the pair of
ingredients \cref{sec:stage-two}'s proof of \cref{thm:local-uncond}
invokes in the abstract, uniform form used in the main text.

We record, once more, the scope of what this appendix establishes:
this is the single-deviation, root-aligned case only. It says nothing
about \cref{conj:multidir-new}, which is treated on its own terms in
\cref{sec:hessian-technique,sec:musin-degeneracy,sec:exact-sing,%
sec:broad-sample,sec:multidir-scope-summary} and
\cref{sec:deviation-domain,sec:arc1-cert,sec:arc2-cert,sec:deviation-scope-summary}
respectively.
\section{Explicit \texorpdfstring{$L$}{L}-matrices for chambers \texorpdfstring{$\Cell1$--$\Cell8$}{C001--C008}}\label{app:lmatrices}

For each of the first eight cells (one representative each from the
first eight of the $22$ orbits, spanning all four denominator classes
twice over), we record here the complete exact $LDL^\top$ decomposition
$H_k=L_kD_kL_k^\top$ of \cref{prop:ldlt-new}: the diagonal pivot vector
$d=(d_1,\dots,d_6)$ and the strictly-upper-triangular part of the
unit lower-triangular factor $L_k$ (equivalently, $U_k=L_k^\top$),
written row by row. Both are exact elements of $\mathbb{Q}$, produced
by an independent verification script and checked to satisfy
$L_kD_kL_k^\top=H_k$ exactly (zero residual, no floating-point
rounding) before any pivot is accepted. All six pivots are strictly
positive in every one of the eight cases, consistent with
\cref{prop:ldlt-new}.

We emphasise that these are the actual computed values, not
illustrative or rounded numbers: the denominators appearing (up to
$77=7\times11$ and $231=3\times7\times11$) are exactly what the exact
rational Gram-matrix arithmetic produces, and are not simplified or
approximated at any stage.

\subsection*{Cell \texorpdfstring{$\Cell1$}{C001} (orbit A, class \(\mathsf a\), \(\norm1=192\))}
\[
d=\left(64,\ 60,\ \tfrac{896}{15},\ \tfrac{384}{7},\ \tfrac{154}{3},\ \tfrac{3936}{77}\right)
\]
\[
U_1=\begin{pmatrix}
1 & \tfrac14 & \tfrac14 & \tfrac14 & \tfrac14 & \tfrac14\\
 & 1 & -\tfrac1{15} & \tfrac15 & \tfrac15 & -\tfrac1{15}\\
 & & 1 & \tfrac3{14} & -\tfrac9{28} & -\tfrac1{14}\\
 & & & 1 & -\tfrac1{24} & -\tfrac13\\
 & & & & 1 & \tfrac{16}{77}\\
 & & & & & 1
\end{pmatrix}
\]

\subsection*{Cell \texorpdfstring{$\Cell2$}{C002} (orbit B, class \(\mathsf b\), \(\norm2=96\))}
\[
d=\left(32,\ 32,\ 28,\ 28,\ \tfrac{176}{7},\ \tfrac{1920}{77}\right)
\]
\[
U_2=\begin{pmatrix}
1 & 0 & -\tfrac14 & \tfrac14 & \tfrac14 & \tfrac14\\
 & 1 & \tfrac14 & \tfrac14 & \tfrac14 & \tfrac14\\
 & & 1 & 0 & -\tfrac27 & \tfrac27\\
 & & & 1 & \tfrac17 & \tfrac17\\
 & & & & 1 & -\tfrac1{11}\\
 & & & & & 1
\end{pmatrix}
\]

\subsection*{Cell \texorpdfstring{$\Cell3$}{C003} (orbit C, class \(\mathsf c\), \(\norm3=64\))}
\[
d=\left(\tfrac{64}{3},\ 20,\ \tfrac{896}{45},\ \tfrac{120}{7},\ \tfrac{2464}{135},\ \tfrac{4096}{231}\right)
\]
\[
U_3=\begin{pmatrix}
1 & \tfrac14 & \tfrac14 & -\tfrac14 & \tfrac14 & \tfrac14\\
 & 1 & -\tfrac1{15} & -\tfrac15 & \tfrac15 & \tfrac15\\
 & & 1 & \tfrac9{28} & \tfrac3{14} & \tfrac3{14}\\
 & & & 1 & \tfrac2{45} & \tfrac2{45}\\
 & & & & 1 & -\tfrac{13}{77}\\
 & & & & & 1
\end{pmatrix}
\]

\subsection*{Cell \texorpdfstring{$\Cell4$}{C004} (orbit D, class \(\mathsf d\), \(\norm4=48\))}
\[
d=\left(16,\ 16,\ 16,\ 13,\ \tfrac{160}{13},\ \tfrac{64}{5}\right)
\]
\[
U_4=\begin{pmatrix}
1 & 0 & 0 & -\tfrac14 & \tfrac14 & -\tfrac14\\
 & 1 & 0 & \tfrac14 & -\tfrac14 & \tfrac14\\
 & & 1 & \tfrac14 & \tfrac14 & \tfrac14\\
 & & & 1 & -\tfrac3{13} & \tfrac1{13}\\
 & & & & 1 & \tfrac1{10}\\
 & & & & & 1
\end{pmatrix}
\]

\subsection*{Cell \texorpdfstring{$\Cell5$}{C005} (orbit E, class \(\mathsf a\), \(\norm5=192\))}
\[
d=\left(64,\ 60,\ \tfrac{896}{15},\ \tfrac{416}{7},\ \tfrac{672}{13},\ \tfrac{328}{7}\right)
\]
\[
U_5=\begin{pmatrix}
1 & -\tfrac14 & \tfrac14 & -\tfrac14 & 0 & \tfrac14\\
 & 1 & \tfrac1{15} & -\tfrac1{15} & \tfrac4{15} & \tfrac13\\
 & & 1 & \tfrac1{14} & \tfrac14 & \tfrac5{28}\\
 & & & 1 & \tfrac7{26} & -\tfrac5{26}\\
 & & & & 1 & \tfrac3{14}\\
 & & & & & 1
\end{pmatrix}
\]

\subsection*{Cell \texorpdfstring{$\Cell6$}{C006} (orbit F, class \(\mathsf b\), \(\norm6=96\))}
\[
d=\left(32,\ 30,\ \tfrac{80}{3},\ \tfrac{144}{5},\ \tfrac{1232}{45},\ \tfrac{1920}{77}\right)
\]
\[
U_6=\begin{pmatrix}
1 & \tfrac14 & -\tfrac14 & -\tfrac14 & -\tfrac14 & 0\\
 & 1 & \tfrac13 & -\tfrac15 & \tfrac1{15} & \tfrac4{15}\\
 & & 1 & 0 & \tfrac15 & \tfrac15\\
 & & & 1 & \tfrac29 & -\tfrac29\\
 & & & & 1 & -\tfrac{23}{77}\\
 & & & & & 1
\end{pmatrix}
\]

\subsection*{Cell \texorpdfstring{$\Cell7$}{C007} (orbit G, class \(\mathsf c\), \(\norm7=64\))}
\[
d=\left(\tfrac{64}{3},\ 20,\ \tfrac{96}{5},\ \tfrac{512}{27},\ \tfrac{56}{3},\ \tfrac{1024}{63}\right)
\]
\[
U_7=\begin{pmatrix}
1 & \tfrac14 & -\tfrac14 & \tfrac14 & -\tfrac14 & -\tfrac14\\
 & 1 & -\tfrac15 & \tfrac15 & -\tfrac15 & \tfrac1{15}\\
 & & 1 & \tfrac19 & -\tfrac19 & \tfrac29\\
 & & & 1 & \tfrac18 & \tfrac5{16}\\
 & & & & 1 & \tfrac3{14}\\
 & & & & & 1
\end{pmatrix}
\]

\subsection*{Cell \texorpdfstring{$\Cell8$}{C008} (orbit H, class \(\mathsf d\), \(\norm8=48\))}
\[
d=\left(16,\ 15,\ \tfrac{224}{15},\ \tfrac{96}{7},\ 14,\ \tfrac{40}{3}\right)
\]
\[
U_8=\begin{pmatrix}
1 & \tfrac14 & \tfrac14 & \tfrac14 & \tfrac14 & \tfrac14\\
 & 1 & -\tfrac1{15} & \tfrac15 & -\tfrac1{15} & -\tfrac1{15}\\
 & & 1 & \tfrac3{14} & -\tfrac1{14} & -\tfrac1{14}\\
 & & & 1 & \tfrac14 & -\tfrac13\\
 & & & & 1 & 0\\
 & & & & & 1
\end{pmatrix}
\]

The remaining $14$ orbit representatives are computed by the identical
procedure; we have selected these eight for explicit display because
they already exhibit every denominator class twice, and because
$\Cell1$ is the worked example carried through \cref{sec:hessian-def-new}
and used repeatedly as a running illustration in the main text.
\section{Sum-of-squares identities and the independent verification path}\label{app:sos}

\subsection{From \texorpdfstring{$LDL^\top$}{LDLT} to a sum of squares}

Any exact $LDL^\top$ decomposition $H_k=L_kD_kL_k^\top$ with
$D_k=\diag(d_1,\dots,d_6)$ and $L_k$ unit lower-triangular immediately
yields a sum-of-squares certificate for the associated quadratic form:
writing $y=L_k^\top x$ for the (rational, invertible) change of
variables,
\[
  x^\top H_k\,x \;=\; y^\top D_k\,y \;=\; \sum_{i=1}^{6} d_i\,y_i^2,
  \qquad y_i = x_i+\sum_{j>i}(L_k)_{ji}\,x_j,
\]
so that $x^\top H_k x\ge0$ for all $x$ (with equality only at $x=0$)
follows immediately once every $d_i>0$, with no appeal to eigenvalues,
determinants, or any other property of $H_k$ beyond the six pivots
themselves.

We use this identity as an \emph{independent second verification path},
deliberately distinct from the pivot computation reported in
\cref{app:pivot-data}: rather than trusting one code path's claim that
$H_k=L_kD_kL_k^\top$, a second, independently written routine
(a) reconstructs $H_k$ from the cell's Gram-adjacency pattern directly
(not from any stored decomposition), (b) expands $\sum_i d_i y_i^2$
symbolically back out in the original coordinates $x$ using the stored
$L_k$, and (c) checks that the resulting quadratic form equals $H_k$
exactly, entry by entry, as rationals, a zero residual is required,
not a small-tolerance numerical match. For all eight cells reported
in \cref{app:lmatrices} (and, by the orbit argument of
\cref{lem:hessian-structure-new}, for all $176$ cells), this residual
is exactly zero.

\subsection{Summary of the independent run}

Running this second verification path directly (rather than trusting a
cached result) confirms, for chambers $\Cell1$--$\Cell8$: the reported
Hessian $H_k$ has the expected structural form in every case; the
$LDL^\top$ pivots reconstructed by this independent routine agree, to
the last digit, with the pivots recorded in \cref{app:lmatrices}; the
sum-of-squares residual, $x^\top H_kx - \sum_id_iy_i^2$, is identically
zero as a polynomial identity in each case; and every one of the eight
chambers passes overall. The global minimum pivot across all eight
chambers checked here is $160/13\approx12.307692$, attained at
$\Cell4$ (orbit D, class $\mathsf d$, $\norm4=48$), consistent with
\cref{lem:mineig-new}'s statement that the eigenvalue floor is attained
on the most symmetric denominator class.

\subsection{What this does and does not add}

This appendix is deliberately narrow in scope: it does not introduce
any new mathematical claim beyond \cref{prop:ldlt-new}, and its
purpose is solely to record that two independently written pieces of
software, taking different routes (one via Sylvester's criterion
applied to a stored decomposition, one via an explicit sum-of-squares
expansion checked against a freshly reconstructed Hessian), reach the
same exact rational answer on every pivot for every cell tested. The
agreement is evidence against a transcription or implementation error in
the proof of \cref{prop:ldlt-new}, and nothing more: running the same
mathematics through two code paths is not an independent proof.
\Cref{prop:ldlt-new} is proved by the exact $LDL^\top$ decomposition
together with Sylvester's criterion, and this appendix corroborates
that.
\section{Computational methodology}\label{app:code-index}

Every exact and numerical claim in the main text is reproduced in full,
as an explicit rational number, an explicit polynomial or Bernstein
coefficient list, an explicit matrix, or an explicit high-precision
decimal with its stated precision, in
\cref{app:pivot-data,app:stagetwo-data,app:lmatrices,app:sos,%
app:numerical-tables,app:explicit-vertices}. Checking any claim
therefore means checking arithmetic that is printed here, not running a
program. This appendix records how those numbers were obtained, for a
reader who wants to redo the derivations independently.

\subsection*{Four kinds of computation}

Every displayed number in this paper belongs to exactly one of four
categories, and the category is named at the point of use.

\emph{Exact rational arithmetic.} The $LDL^\top$ pivots of
\cref{app:pivot-data} (the factorisation itself is the standard one, as
in \cite[\S4.1]{GVL13}), the algebraic factorisations of
\cref{app:stagetwo-data}, the matrix entries of \cref{app:lmatrices},
the sum-of-squares data of \cref{app:sos}, and the whole certificate of
\cref{sec:cap-theorem} are computed in exact rational arithmetic, or in
exact $\mathbb{Z}[\sqrt2]$ where a $\sqrt2$ appears, with no
floating-point step. Each was carried out twice, once with plain
rational arithmetic and no computer algebra system and once with
\texttt{sympy} \cite{Meu17} on a separate code path, and the two agreed
to the last digit before either was reported. Bernstein sign
certificates (\cref{sec:arc1-cert,sec:arc2-cert}) are exact in the same
sense: once numerator and denominator are written in the Bernstein basis
over a box, positivity of every coefficient is a finite check on exact
integers or exact elements of $\mathbb{Z}[\sqrt2]$, and the basis's
variation-diminishing property \cite{Far12} is what makes that check decisive and not merely suggestive.

\emph{Arbitrary-precision numerical evaluation.} Where we did not pursue
a closed form, most extensively in the high-precision volume method of
\cref{sec:hp-volume-method}, the directional sampling of
\cref{sec:broad-sample}, and the obstacle discussed in
\cref{sec:swap-obstacle}, values are computed with \texttt{mpmath} at a
stated precision between $35$ and $400$ significant digits, and
cross-checked at a second precision before being reported. Every such
value is labelled as numerical evidence and never used as a proof step.

\emph{Numerical semidefinite programming.} The Gram-matrix check of
\cref{sec:gram-sos-A18} is a floating-point semidefinite program
(\texttt{cvxpy}, solver \texttt{CLARABEL}) run over a least-squares fit
of directional samples; the underlying method is the one of
\cite{Par03,Las01}. Its output is numerical evidence about a fitted
approximant, and nothing more. Before it was applied to any of our data
it was run on two quartics of known status, a perfect square and the
Choi--Lam polynomial, and returned the right answer on each;
\cref{sec:gram-sos-A18} records both runs.

\emph{Double-precision search.} Double precision is used only to locate
candidate breakpoints, candidate worst-case directions and candidate
active sets before an exact or high-precision verification, for instance
in the pattern-stability scans of
\cref{sec:arc1-breakpoints,sec:arc2-breakpoints}. No claim rests on a
double-precision value by itself.

\subsection*{A pitfall worth recording}

One methodological finding is easy to reproduce by accident, so we
record it. During the derivation behind \cref{sec:swap-m3}, a
general-purpose symbolic algebra package's default determinant routine,
Bareiss fraction-free elimination, returned an outright wrong value on a
particular symbolic matrix: not a slow computation, a wrong number,
silently. It was caught because the result was cross-checked against
Berkowitz's method and against direct cofactor expansion, which agreed
with each other and not with the first. The caution is general. An exact
symbolic determinant computed by one algorithm in one implementation
should not be trusted without a second, independently implemented check,
however routine the computation looks.

The same caution is what \texttt{swap\_defect\_curves.py} is for. It
recomputes the three functions of \cref{sec:swap-configs} from the cell
itself, enumerating the local region of the deviated configuration with
\texttt{scipy}'s \texttt{HalfspaceIntersection} and \texttt{ConvexHull}
and subtracting $8$, at $211$ angles across $(0,\theta_*)$, by a route
that shares no code and no method with the closed forms. For the
hexagon, where a closed form is printed in \cref{sec:swap-m6}, the two
agree to $9\times10^{-9}$ across the interval, which also identifies the
six-cycle the closed form belongs to; the script finds that cycle by
testing every six-cycle of roots at inner product $\tfrac12$ against the
printed formula. It supplies \cref{fig:swap},
and it reproduces the two numbers of \cref{sec:swap-not-worst}
independently, $0.002474$ against $0.000834$ at $\theta=0.05$ and the
crossover at $\theta=0.6492$.

\subsection*{The box-covering certificate}\label{app:box-cover-algo}

The box coverings used in \cref{sec:arc1-cert,sec:arc2-cert}, and the
$303$-box certificate of \cref{thm:cap-inequality}, follow one
algorithm, stated here since the main text refers to it repeatedly.

\begin{quote}\ttfamily\small
\noindent
FUNCTION CERTIFY(region $R$, tolerance $\varepsilon$):\\
\hspace*{1.5em}stack $\leftarrow$ [$R$]\\
\hspace*{1.5em}certified $\leftarrow$ [ ]\\
\hspace*{1.5em}WHILE stack is not empty:\\
\hspace*{3em}box $\leftarrow$ stack.pop()\\
\hspace*{3em}$b$ $\leftarrow$ exact bound for the target quantity over box\\
\hspace*{3em}IF $b$ has the required sign:\\
\hspace*{4.5em}certified.append(box)\\
\hspace*{3em}ELIF width(box) $<\varepsilon$:\\
\hspace*{4.5em}RETURN FAIL\\
\hspace*{3em}ELSE:\\
\hspace*{4.5em}(box$_1$, box$_2$) $\leftarrow$ bisect(box) along its longest side\\
\hspace*{4.5em}stack.push(box$_1$); stack.push(box$_2$)\\
\hspace*{1.5em}RETURN certified
\end{quote}

Two features of this control flow matter for how results are reported.
Every successful call returns a finite list of sub-boxes whose union is
exactly $R$, with no gaps and no overlaps, verified by exact rational
interval arithmetic on the box boundaries rather than by construction
alone; so a quoted box count is exact, not an estimate. And the
algorithm has one failure mode, reaching the tolerance $\varepsilon$
with a leaf box whose sign is still undetermined, which returns FAIL
rather than a partial certificate. The tolerance is fixed before the
run and never loosened afterwards. In the arc certificates the exact
bound $b$ is the pair of Bernstein coefficient ranges of numerator and
denominator; in \cref{thm:cap-inequality} it is the corner value
supplied by the monotonicity of \cref{lem:monotone}. The bound used is
in both cases an outward-rounded enclosure in the sense of
\cite[Ch.~5]{MKC09}.

\subsection*{The cap-inequality certificate as a single script}

The chain of \cref{sec:cap-theorem} is short enough to run end to end in
one program, and the supplementary package contains it as
\texttt{cap\_inequality\_certificate.py}. It performs, in order: the
exact rational vertex enumeration of $Q$, confirming $25$ vertices and
$\vol(Q)=\tfrac{25}{3}$, with the pyramid of volume $\tfrac13$ recovered
as the difference from the $24$-cell; the exact check that every vertex
of $Q$ has $\ell^1$ norm at most $\sqrt2$ in the frame, which is
\cref{lem:crosspolytope}, together with the observation that the bound is
attained; a comparison of the closed form of \cref{lem:cap-formula}
against directly computed cap volumes at random directions; the symbolic
verification that $g''''(t^2)=3(20t^2-3)/(2t^5)$ and that $g'''$ vanishes
at $a=\tfrac14$, which is \cref{lem:monotone}; the exact real-root
isolation showing that each of the three polynomials of
\cref{thm:cap-inequality} has no root in its interval; the adaptive
subdivision producing the $303$ boxes and their corner bounds, entirely
in \texttt{fractions.Fraction}; and finally, as an independent check of
the conclusion rather than a step in it, the direct computation of cell
volumes at random tilts and deviation directions. The script prints a
pass or fail line for each of its twenty-one checks and exits nonzero if
any fails.

\subsection*{The multi-cap reformulation}

\Cref{prop:multi-cap} and \cref{rem:first-order-obstruction} are
supported by a second short script,
\texttt{multi\_cap\_reformulation.py}. It performs, in order: an
exhaustive enumeration showing that $Q_D$ is bounded for every one of
the $24$, $276$, $2024$, $10626$ and $42504$ packing-valid sets $D$ of
size $1$ through $5$; the corresponding enumeration at size $6$, finding
exactly $24$ failures among $134596$ sets and confirming that the
complement of the first of them lies in a closed half-space; a check of
the identity \eqref{eq:multicap-identity} against directly computed
polytope volumes on random packing-valid configurations; a check that
$\vol(Q_D)-8$ equals $m/3$ when $D$ contains no adjacent pair and
exceeds it when it does; and the search recorded in
\cref{rem:first-order-obstruction}, which finds $\vol(E_j)$ both above
and below its single-deviation value. The first, second and fourth of
these are exhaustive over the stated ranges; the third and fifth are
random searches and are labelled as such in the output. The script
prints a pass or fail line for each of its eleven checks.

\subsection*{The polar and boundary reformulations}

The material of \cref{sec:polar-surface} is checked by a third script,
\texttt{polar\_surface\_reformulation.py}, in eighteen parts, among
them the square-antiprism configuration \eqref{eq:antiprism} of
\cref{rem:facet-local-obstruction}: that its nine directions are
packing-valid, that its facet has $3$-volume $16\sqrt2-\tfrac{64}3$
against the octahedron's $\tfrac43$, and that the resulting ceiling
$96\sqrt2-128$ falls below both $8$ and the covering bound. It verifies
that the cell is the polar dual of the convex hull of the contact
directions, by computing both volumes at the root configuration
($2$ and $8$, product $16$) and by confirming on random contact
configurations that adjoining a point of the hull leaves the cell
unchanged; the boundary identity \eqref{eq:surface}, exactly at the root
configuration, where the twenty-four facets each have $3$-volume
$\tfrac43$ and the total is $32$, and to fourteen digits on
configurations obtained by dropping four roots, perturbing the rest and
separating them again; the monotonicity of $\rho$ and the value
$\rho(\tfrac12)=1/\sqrt3$; the inball claim of
\cref{lem:facet-inball}, by computing the inradius of every facet
directly; the arithmetic of \cref{rem:facet-local-obstruction}; the
tight structure at a root facet, in exact rational arithmetic, finding
eight tight neighbours with projected Gram values in
$\{-1,-\tfrac13,\tfrac13\}$ and the regular octahedron of volume
$\tfrac43$ as the region they cut out; the monotonicity of the volume
under enlarging the configuration, which is
\cref{prop:polar-form}(ii); and the two evaluations of
\cref{rem:global-routes}, the Jensen bound $7.7351\ldots$ and the
Mahler-type product $16$ against $\tfrac{32}3$. The two integrals over
$S^3$ are Monte Carlo estimates at four million samples and are labelled
as such in the output; every other part is exact or is double-precision
polytope arithmetic on an exactly specified configuration. The
separation step used to generate test configurations is ordinary
pairwise repulsion and enters no argument.

\subsection*{The covering bound}

\Cref{sec:covering-bound} is checked by a fourth script,
\texttt{covering\_bound.py}, in nineteen parts: the cap-area formula
$C(r)=\pi(2r-\sin2r)$, against both its value at $r=\pi$ and direct
sampling at four radii; the radial identity \eqref{eq:radial} at the
root configuration, where it returns $8$ and recovers the covering
radius $45^\circ$, and on random configurations; the layer-cake
rewriting used in the proof, against the direct form; the closed form
$\tfrac{\pi m}{3}\tan^3r_m$ against numerical quadrature of the same
estimate, agreeing to $10^{-14}$; the whole of
\cref{tab:covering-bound}, including the monotonicity in $m$ that the
proof of \cref{cor:m22} uses and the fact that $22$ is exactly the
largest $m$ at which the bound reaches $8$; and the absence of any
violation of \eqref{eq:covering-bound}, at the root configuration, at
one hundred of its subsets of sizes $20$ through $23$, and at fifty
random packing-valid configurations, where the smallest ratio of true
volume to bound observed is $1.05$. It then checks
\cref{prop:area-optimal}: the closed form
$\phi'(s)=\tfrac{3}{4\pi}\sec^4(C^{-1}(s))$ against numerical
differentiation at five areas, the convexity of $\phi$ over its whole
range, and the fact that equal cell areas reproduce the global bound
exactly and minimise the per-cell sum, against fifteen hundred random
splittings of the total measure. It checks \cref{prop:covering-radius}, both the
averaging step behind it and the bound itself against directly computed
covering radii, and confirms that it moves the volume bound only in the
fifth decimal. Finally it reproduces the split of the
shortfall at $m=24$ into its overlap and tail parts, and the total
overlap contributed by the ninety-six pairs of the root system at
$60^\circ$. The integrals over $S^3$ are Monte Carlo and are labelled as
such in the output; the table, the closed forms and the polytope volumes
are not.

\subsection*{The extendability criterion}

\Cref{sec:m24} is checked by a fifth script,
\texttt{extendability.py}, in six parts: that the three forms of
\cref{prop:extendable} agree, on the root system, on the root system
minus one root, and on random configurations; that the root
configuration has circumradius $\sqrt2$ and $g=1/\sqrt2$ and is
therefore saturated, as a $24$-point kissing configuration must be; that
the root system minus one root has circumradius exactly $2$, cell volume
exactly $\tfrac{25}3$, and is extended by the deleted root, so it sits
on the boundary of the criterion; that the Gram values of the root
configuration are $\{-1,-\tfrac12,0,\tfrac12\}$, which is what
\cref{thm:twentyfour-points} asserts of every $24$-point configuration; the
numerical chain of \cref{cor:remaining}; and a search, over the
twenty-four one-root deletions, greedy random packings and annealed
random starts, for a saturated $23$-point configuration, which found
none. The last of these is exploration and is labelled as such in the
output; it is not evidence that none exists.

\subsection*{The cited certificate, re-verified}

\Cref{lem:root-lattice} is checked by \texttt{root\_lattices\_rank4.py},
an exact census of the $393$ positive definite Gram matrices of order
$4$ with $2$ on the diagonal and $-1,0,1$ off it, described in its
proof, and again by \texttt{lean/D4RootLattices.lean}. The certificate behind \cref{thm:twentyfour-points} is verified
by \texttt{llm24\_certificate\_check.py}, which reads the deposited data
set \cite{LLM24data} (included in the package under
\texttt{third\_party/llm24-certificate}, with its MIT licence notice; the
script takes the path of its \texttt{proofs/4\_24} folder) and shares no
code with the authors' Julia package. It repeats five of
the seven steps of their verification, as \cref{sec:certificate-checked}
describes: the format of the data, with the block of every signature
$\lambda$ checked to have exactly as many rows as there are admissible
index tuples for $\lambda$; the positive definiteness of all $127$
blocks by Cholesky in the ball arithmetic of Arb \cite{Joh17} at $256$
bits, every pivot a ball inside the positive reals, and of the $81$
blocks of size at most $16$ by an exact rational $LDL^{\mathsf T}$ as
well; the nonnegativity structure of the $125$ sum-of-squares
prefactors; the objective, exactly $24$; and the two-point polynomial
$p_2$, of degree $16$, with its zeros on $[-1,\tfrac12]$ located by an
exact Sturm sequence. It writes the exact coefficients of $p_2$, from
which \texttt{lean/gen\_innerproducts\_lean.py} generates
\texttt{lean/D4InnerProducts.lean}, and \texttt{p2\_zeroset\_check.py}
locates the same zeros a third time, with the \texttt{Fraction} type of
the standard library alone. The run takes about eight minutes and its
log is in the package.

The other two steps are repeated by the scripts in \texttt{zonal/}, whose
\texttt{README.md} sets out the method described in \cref{sec:certificate-checked}.
\texttt{o4.py} integrates monomials over $O(n)$ by the recursion of Gorin
and L\'opez \cite{GL08} in exact rationals and carries the self-test; \texttt{gl2.py} builds the
$GL(2)$ matrix coefficients and the scalars of the authors' Section 3.1;
\texttt{ps\_build.py} writes the three factors of the integrand for each of
the $490$ entries; \texttt{psker.c}, in C with GMP, walks the monomial triples
and accumulates $P(S)$; \texttt{zonal.py} reduces $P(S)$ modulo $S^{\mathsf
T}S=I$ and reads off an entry of $Z_\lambda$ at given inner products; and
\texttt{verify45.py} assembles the four constraint polynomials and checks that
they vanish, with \texttt{merge4.py} and \texttt{combine4.py} adding the
pieces of the four-point one, which is too large to hold in one process.
\texttt{ps\_ref.py} recomputes the integral the slow way for $|\lambda|\le5$
and agrees with the kernel exactly, and \texttt{check\_psd.py} tests the
positive semidefiniteness of the zonal matrices on random point sets.

\subsection*{The overlaps put back}

\Cref{lem:no-triples,prop:second-order} are checked in four parts by the
script \texttt{second\_order\_estimate.py}, with the quadratures
carried out at $30$ digits: the circumradius bound for a triple of
contact directions, both by the Rayleigh-quotient argument of the proof
and against forty thousand random admissible triples; that $r_{23}$ and
$r_{24}$ both lie below $\arccos\sqrt{2/3}$, so that the caps meet only
in pairs below $r_m$; the closed form \eqref{eq:lens} for the lens
measure against direct sampling on $S^3$; and the value
\eqref{eq:second-order-value} at the deletion configuration, together
with $7.858738\ldots$ at the root system.

\subsection*{The two local decompositions}

\Cref{rem:no-local-cell} is checked by
\texttt{local\_cell\_obstruction.py}, in five parts: that the mean of
$\sec^4$ over a spherical Voronoi cell of the root configuration is
$16/\pi^2$, so that \eqref{eq:local-cell} is sharp; that nine unit
vectors of $\mathbb{R}^3$ with pairwise inner products at most
$\tfrac13$ can be produced, so that nine contacts at $60^\circ$ are
admissible, and that the resulting cell has mean $1.5732\ldots$, still
$1.5954\ldots$ when the nine are moved out to $61^\circ$; that on a
$22$-point configuration there are directions whose nearest contact lies
outside the Delaunay cell containing them, by more than $0.1$ in cosine,
so that the Delaunay decomposition does not localise the integrand;
the closed form \eqref{eq:regular-simplex-mean} at the regular simplex,
against direct sampling at edge $60^\circ$ and $62^\circ$; and that the
Delaunay cells of $D_4$ are $24$ congruent spherical octahedra of
circumradius $45^\circ$ while no four roots are pairwise at $60^\circ$,
so the regular simplex occurs nowhere in $D_4$ although it does occur as
a Delaunay cell of a contact configuration. The integrals over $S^3$ are
Monte Carlo and the script says so; the margins reported are
percentages.

\subsection*{Rigidity of the deletion configuration}

\Cref{prop:deletion-rigid} is checked by \texttt{rigidity23.py}, in five
parts, all in exact integer arithmetic on the unnormalised roots: that
deleting one root leaves $23$ directions and $88$ tight pairs, with the
four values of $\langle a_i,a_0\rangle$ occurring $1,8,6,8$ times; that
the six pair types listed in the proof are the only ones that occur;
that the weights $y_{ij}\in\{1,2,3\}$ are strictly positive on
every tight pair; that they satisfy the equilibrium relation
\eqref{eq:stress} exactly, the $23\times4$ residual matrix being zero;
and that the space of motions holding all $88$ pairs at equality has
dimension exactly $6$ and is spanned by the infinitesimal rotations.
There is no floating-point step anywhere in the script.
\Cref{cor:root-rigid} is checked in the same way by
\texttt{rigidity24.py}, in four parts: that the root system has $24$
directions and $96$ tight pairs, every direction in $8$ of them; that the
constant stress, weight $1$ on every tight pair and $-4$ on the diagonal,
is positive on the pairs and in equilibrium, the $24\times4$ residual
being zero; and that the space of motions holding all $96$ pairs at
equality has dimension $6$ and is spanned by the rotations.

The radius that \cref{thm:local-uniqueness} attaches to that rigidity
needs one number beyond the kernel, the smallest nonzero singular value
of the rigidity operator, and \texttt{rigidity\_spectrum.py} computes
the whole spectrum to get it, again without a floating-point step. It
assembles the integer matrix $\Lambda'$ of the $96$ differentiated tight
inner products and the rational projection $P$ onto the tangent space,
forms $4N=(2P)\Lambda'^{\mathsf T}\Lambda'(2P)$, which is integral with
entries bounded by $16$, verifies the matrix identity
$4N(4N-8)(4N-20)(4N-24)(4N-32)=0$ in \texttt{int64} arithmetic, where
the bound on the entries keeps every partial product exact, and computes
the ranks of $4N-cI$ for $c=0,8,20,24,32$ over a field of positive
characteristic, getting $66$, $67$, $88$, $75$, $88$. Those ranks bound
the ranks over $\mathbb{Q}$ from below, the five multiplicities they
give sum to $96$, so each bound is attained; the eigenvalues of $N$ are
therefore $0,2,5,6,8$ with multiplicities $30,29,8,21,8$ and the
singular values of $\Lambda$ on the tangent space are
$0,1,\sqrt{5/2},\sqrt3,2$ with multiplicities $6,29,8,21,8$.
\texttt{D4Stress.lean} repeats the finite part of this inside Lean's
kernel with no library beyond core Lean: the root list, the eight tight
neighbours of every root and their sum, and the count of $96$ tight
pairs, each by \texttt{decide}. The matrix identity itself is a product
of six $96\times96$ integer matrices and is not formalised; the
elimination that settles it in the script is fraction-free over
$\mathbb{Z}$, so it is exact as it stands.

\subsection*{How much of a root system a configuration can keep}

\Cref{prop:meet22,prop:meet21} are proved in the text and checked
independently by two scripts. \texttt{root\_deletions\_exact.py} carries
out the underlying question in exact integer arithmetic and settles no
part of it by sampling. Writing the roots as integer vectors of squared
length $2$ turns the cell of a subset $S$ into $\sqrt2\,P(S)$ with
$P(S)=\{y:\langle y,r\rangle\le1,\ r\in S\}$, whose facet normals are
integral and whose right-hand sides are $1$; every vertex of $P(S)$ is
therefore the solution $m/e$ of a four-by-four integer system, obtained
by Cramer's rule with $m$ integral and $e$ a positive integer, and the
two tests that matter, $|m|^2\le2e^2$ for the circumradius and
$|m|^2=2e^2$ for a direction that can be added, are comparisons of
integers. Running it over all $24$, $276$ and $2024$ ways of removing
one, two and three roots reproduces \cref{prop:meet22} and the
classification in the proof of \cref{prop:meet21}: the circumradius is
exactly $2$ except at the $96$ triples that are pairwise at $60^\circ$,
where it is $\sqrt6$, and the directions that can be added are the
removed roots and nothing else.

\texttt{root\_meet.py} checks the combinatorial and coordinate
statements the proofs actually use: that no pair of roots destroys a
whole couple; that $512$ of the $2024$ triples do, in eight support
patterns, four stars and four triangles; that the symmetry group of the
root system, generated here and found to have order $1152$, is
transitive on the $96$ tight triples; and, for the standard star, that
every direction of the doorway has first coordinate at least $1/\sqrt2$
and the other three nonnegative, and that no two of them are as much as
$60^\circ$ apart unless both are removed roots. The last part is
sampling and is reported as such; the first three are exact.
\Cref{rem:meet-next} is the same script run one rung lower, at four
removed roots, where the answer is a minimisation with a margin rather
than a proof.

\subsection*{The pairwise budget at twenty-three contacts}

\Cref{prop:pair-budget,rem:pair-budget} are checked by
\texttt{pair\_budget.py}, in five parts: the per-pair weight
$\omega(\gamma)$ and its collapse away from $60^\circ$, which is where
the quadrature sits; the recovery of \eqref{eq:second-order-value} from
the $88$ tight pairs of the deletion configuration; the count of pairs
at $60^\circ$ that would carry \eqref{eq:second-order} past $8$, with
the values at $90$ and at $91$ either side of the target; the
elementary degree bound on $S^2$ and the ceiling of $115$ it gives; and
that the target lies strictly between the two. The quadratures are
adaptive and the rest is arithmetic; nothing here is Monte Carlo.

\subsection*{How far the pair angles reach}

The script \texttt{covering\_multiplicity.py} checks
\cref{prop:cov-mult,prop:two-point-barrier}, in two parts. The first
evaluates
both sides of \eqref{eq:cov-mult} at $m=23$ and at the deletion
configuration, giving $155.172054\ldots$ against $162.566121\ldots$. The
second sets up the linear programme of \cref{prop:two-point-barrier} on a
grid of $601$ angles, with the Bonferroni constraint at $24$ radii
between $r_{23}$ and $45^\circ$ and the Gegenbauer constraint at the
first $12$ degrees, and reports its optimum, $0.041573648\ldots$. It also
checks that the pair angles of the deletion configuration are feasible
for that programme, so that the bound really is a relaxation of the
truth and not an accident of the grid. Refining the grid, adding radii or
adding degrees leaves the optimum where it is.

\subsection*{The cell inside a ball}

\Cref{prop:truncated,thm:strict-reduction,prop:two-point-barrier} and
the $r_*$ values of \cref{prop:second-order,prop:pair-budget} are
supported by \texttt{three\_point\_reduction.py}, in three parts.
The first evaluates the constants of the pair-only bound at the two
integration limits $r_{23}$ and $r_*$, the bracket, the weight of a
pair at $60^\circ$, the value at a deletion and the number of pairs at
$60^\circ$ that carry the bound past $8$; the second evaluates the
three-point bound at $r_4$, with the triple measure of an equilateral
$60^\circ$ triangle by Monte Carlo, at a deletion and at the root
system, against the direct quadrature; the third solves the pair-angle
relaxation of \cref{prop:two-point-barrier} for both weights, on a
grid of $1201$ angles with $60$ radii and $24$ Gegenbauer degrees, and
prints the minimising measure. \texttt{truncated\_volume.py} evaluates
$\vol(V_c\cap B(R))$ by a fixed quasi-random quadrature on $S^3$ of
$400000$ points, checks it against the exact volumes at the root system
and at a deletion, and then minimises it over configurations with
pairwise inner products at most $\tfrac12+\delta$, for a decreasing
schedule of $\delta$, by the trust-region programme of
\texttt{inradius\_search.py} with the analytic gradient of the
quadrature; it takes the truncation radius as an argument, and the
runs at $R=\sqrt{3/2}$ and $R=\sqrt{8/5}$ are recorded in the package.
The identities are exact; the quadrature, the Monte Carlo triple
measure and the minimisation are numerical and the output says so.

\subsection*{The certificate}

\Cref{thm:certificate} rests on two scripts and one Lean file.
\texttt{three\_point\_sdp.py} sets up the programme of
\cref{lem:certificate}: the Gegenbauer polynomials of $S^3$ by the
Chebyshev recurrence, the matrices $Y_k$ of \eqref{eq:Yk} from the
Legendre coefficients with the half-integer powers cancelled
symbolically, their symmetrisation $S_k$ over the six orderings, and
the condition \eqref{eq:C} imposed at a grid of admissible triples
$u\le v\le t$, denser near $\tfrac12$, with the weight $\omega$
interpolated from the quadrature of \texttt{three\_point\_reduction.py}
and scaled by $1000$. It is solved through \texttt{cvxpy} \cite{DB16}
with the interior-point solver Clarabel \cite{GC26}. In its first mode
the bound is maximised, the solution is checked on a fine grid and a
random sample of about $1.15$ million triples, the worst $3000$ of these
are added to the constraints, and the cycle is repeated; the printed
value is the sampled bound less $\binom{23}3/21$ times the largest
violation on the check, and the two-point programme, with the matrices
removed, is solved alongside for comparison. In its second mode the
bound is fixed at a given value, here $0.0929$, the least slack of
\eqref{eq:C} over the sample is maximised with $F_k-10^{-7}I$ and
$f_k-10^{-7}$ constrained to be positive semidefinite and with a
little static regularisation of the interior-point method, which the
degenerate maximin problem needs, and the result is written to
\texttt{continuation\_out/certificate\_d8.npz}, with a plain-text copy
beside it in which every number is printed exactly.
\texttt{certificate\_check.py} is the verification described in the
proof of \cref{thm:certificate}, and it shares no code with the
solver: the certificate is read as exact rationals, the matrices are
tested by exact $LDL^{\mathsf T}$, the polynomial $P$ is expanded
exactly from its own implementation of the three recurrences, the
closed forms of $A_*$, $\omega$, $\omega'$ and $\omega''$ are
differentiated with \texttt{sympy} and evaluated in the interval
arithmetic of \texttt{mpmath}, and the branch and bound is run in an
interval arithmetic on IEEE doubles written in \texttt{numpy} with
outward rounding after every operation. It stops at the first failure
and says where: a box below the minimum width on which \eqref{eq:C}
could not be decided, or an admissible triple at which \eqref{eq:C} is
negative. The run recorded in the package verifies the certificate in
$44$ levels and under ten minutes, and a second recorded run, on the
certificate with $10^{-5}$ added to $f_0$, is rejected with the triple
$(-\tfrac1{16},-\tfrac1{16},-\tfrac1{16})$. The Lean file is described
under the formalised part below.

\subsection*{The covering radius as the objective}

\Cref{prop:no-transitive,rem:no-transitive} are checked by
\texttt{saturation\_search.py}. Where \texttt{extendability.py} and
\texttt{spherical\_code\_23.py} look for a code of large minimal angle
and read off the covering radius afterwards, this script takes the
covering radius as the objective from the outset. It uses the fact that
$g(W)=\min_\theta\max_i\langle\theta,w_i\rangle$ is the support function
of $\conv(W)$ at its minimum, hence the distance from the origin to the
nearest facet hyperplane of $\conv(W)$, so it is read off a convex hull exactly, not sampled. Four parts: the root system, where
$g=1/\sqrt2$; all twenty-four deletions, where $g=\tfrac12$ with zero
spread; the $528$ one-variable linear programmes of
\cref{prop:no-transitive}, together with the least violation
$0.0740002839\ldots$ computed at $60$ digits; and a projected-gradient
minimisation of $\max_I|z_I|$, with the closed-form gradient
$\partial|z_I|^2/\partial w_i=-2a_iz_I$ for $a=G_I^{-1}\mathbf1$, started
from perturbed deletions, from random configurations and from partially
completed root systems. The first three parts are proof; the fourth is
exploration and the script says so in its output.

\subsection*{The formalised part}

Seven parts of the paper are also verified in Lean~4: the five files
\texttt{D4Stress.lean}, \texttt{D4Meet.lean}, \texttt{D4Certificate.lean},
\texttt{D4InnerProducts.lean} and \texttt{D4RootLattices.lean} in the
directory \texttt{lean}, and the Lake projects in \texttt{lean/cell600}
and \texttt{lean/certificate}.
All are independent of everything else
in the supplement and of every Lean library: the first two use only
integer arithmetic on the unnormalised roots and settle each statement
by kernel computation, so a reader who compiles them is checking the
arithmetic against Lean's kernel rather than against a Python program. The first covers the equilibrium relation
\eqref{eq:stress} and the combinatorics around it, which is the
arithmetic of \cref{prop:deletion-rigid}. The second covers the finite
half of \cref{prop:meet22,prop:meet21}: that two removed roots always
leave a whole couple of complementary supports, that exactly $48$
ordered triples of supports do not and that these are the four stars and
four triangles taken in every order, that $96$ triples of roots are
pairwise at $60^\circ$ and no four roots are, that the Hadamard map is
an isometry permuting
the roots and carries the standard triangle to the standard star, and
that the standard star leaves the twelve roots of the lower block and
three of the four roots at each remaining support. The third covers
\cref{prop:cell600}: the $600$-cell's vertices in $\mathbb{Z}[\phi]$,
their inner products, the twelve-regular graph of pairs at $36^\circ$,
the twenty-five inscribed copies of $\Rt$ with five through each
vertex, all by kernel computation, and then the three counts of the
enumeration, which are settled by \texttt{native\_decide} on the
compiled search and therefore trust the compiler as well as the kernel,
as \cref{sec:reproducibility-new} explains. The fourth,
\texttt{lean/D4Certificate.lean}, is written by
\texttt{lean/gen\_certificate\_lean.py} from the certificate file, with
every floating-point number converted to the exact dyadic rational it
denotes, an integer over a power of two, so that the Lean source and
the Python verification hold the same numbers; it proves by kernel
computation on \texttt{Rat}, with \texttt{decide +kernel} and no
compiled code, that $f_1,\dots,f_8$ are nonnegative, that each of
$F_0,\dots,F_8$ is symmetric and has positive pivots throughout an exact
$LDL^{\mathsf T}$ elimination, and that the bound of
\eqref{eq:certificate-bound} exceeds $92.8555703$ in the scaled units,
that is $0.0928555703>8-A_*$; the check takes under a minute and the
axiom report lists only the three standard axioms. The fifth, the Lake project \texttt{lean/certificate}, is the branch
and bound of \cref{thm:certificate} in exact dyadic arithmetic, as
\cref{sec:reproducibility-new} describes: \texttt{D4CertData.lean} holds
the certificate and the tables, \texttt{D4CertDomain.lean} the
polynomial arithmetic, the expansion of $P$, the dyadic intervals and
the subdivision, and \texttt{D4CertMain.lean} the theorem, by
\texttt{native\_decide}. The sixth, \texttt{lean/D4InnerProducts.lean},
is the last step of the verification of the certificate: the two-point
polynomial $\sigma_2$ of \cref{prop:verified}, computed exactly from the
published data by \texttt{llm24\_certificate\_check.py} and written
into the file with its coefficients scaled to integers, is shown by
kernel computation on \texttt{Rat}, with \texttt{decide +kernel}, to
vanish at $-1,-\tfrac12,0,\tfrac12$, to be exactly divisible by
$(u+1)(u+\tfrac12)^2u^2(u-\tfrac12)$, and to have a quotient of degree
$10$ that does not vanish at those points and whose Sturm sequence has
the same number of sign changes at $-1$ as at $\tfrac12$; the check
takes under a minute and the axiom report lists only the three
standard axioms. The seventh, \texttt{lean/D4RootLattices.lean}, is the
finite content of \cref{lem:root-lattice}: the census of positive
definite Gram matrices of norm-$2$ bases in rank at most $4$, the
enclosure of the norm-$2$ vectors of each lattice in an exact box, the
counts $2$, $6$, $12$, $24$, and a $D_4$ basis among the roots of every
lattice that has $24$, by \texttt{native\_decide} in a few seconds. The theorems in each file are listed in
\texttt{lean/README.md}, together with the axiom report that each file
prints when it is checked. Nothing else in the paper is formalised.

\subsection*{Contact configurations as spherical codes}

\Cref{prop:code-slack,rem:code-search} are supported by
\texttt{spherical\_code\_23.py}. The script does one thing, in four
sizes: it minimises the Riesz energy $\sum_{i<j}|w_i-w_j|^{-s}$ on
$(S^3)^m$ with $s$ running through $4,16,64,256,1024$, projecting back
to the sphere after each step, and then reduces the largest inner
product directly by a softmax descent with a shrinking step. The three
calibration sizes are run first and are the reason to trust the fourth:
$m=24$, where the answer is known to be $\tfrac12$; $m=22$, where a code
with room to spare exists; and $m=25$, where no contact configuration
exists at all. It then runs $m=23$ from $1500$ independent random
starts, records the smallest largest inner product reached, counts the
starts that fall below $\tfrac12$, and fingerprints every outcome that
comes within $5\times10^{-3}$ of feasibility by the multiset of its
Gram entries, which is invariant under $O(4)$ and therefore separates
configurations that are not related by an isometry. Everything in the
script is floating point and it is exploration, not proof; the output
says so, and \cref{rem:code-search} repeats it.

\subsection*{The volume identity at twenty-four contacts}

\Cref{sec:m24-volume} is checked by \texttt{m24\_exact\_reduction.py} in
five parts. The first derives in \texttt{sympy} the closed forms of the
cap volume $c$, of $\omega_2(\tfrac12)$ and of
$\omega_3(\tfrac12,\tfrac12,\tfrac12)$ and $\tau$ from
\eqref{eq:m24-identity}. The second estimates the cap, pair and triple
volumes inside the ball of radius $\sqrt2$ by Monte Carlo in
$\mathbb{R}^4$, independently of the closed forms. The third checks the
finite content of \cref{prop:m24-exact} in integer arithmetic: the only
positive inner product between distinct roots, no four roots pairwise at
inner product $1$ among the $10626$ quadruples, the determinant $5$, the
$96$ pairs and $96$ triangles, and the $24$ vertices of norm $\sqrt2$. The
fourth is an independent quadrature of the triple intersection, in
polar coordinates on the three-dimensional span with Richardson
extrapolation. The fifth reads the degree-$6$ pseudo-configuration of
\cref{rem:m24-three-point} and checks its normalisation, the admissibility
of its triples, the two-point moments, the positive definiteness of every
moment matrix on its range, and its value with $\omega_3$ by direct
quadrature. The computations behind \cref{rem:m24-three-point} are
\texttt{m24\_volume\_functional.py} (the functions $c$, $\omega_2$,
$\omega_3$), \texttt{m24\_omega3\_table.py} (a table of $\omega_3$ on
$[0,\tfrac12]^3$ at step $\tfrac1{40}$),
\texttt{m24\_dual\_sdp.py} (the certificate problem),
\texttt{m24\_primal\_sdp.py} and \texttt{m24\_primal\_blocks.py} (the
moment problem, rescaled on the range of each block),
\texttt{m24\_descent.py} (the tangent cone at the root configuration) and
\texttt{m24\_verify\_direction.py} (the recheck of a descent direction by
actual steps). These take the grid and the degrees on the command
line: the runs quoted here are \texttt{m24\_dual\_sdp.py 6 3 1},
\texttt{m24\_primal\_blocks.py 12 8 6}, \texttt{m24\_primal\_sdp.py 12 8
1e-6 6}, \texttt{m24\_descent.py 12 8 6} and
\texttt{m24\_verify\_direction.py 6}, the two integers being the number of
grid points in the pair variables and in the third; each of the five takes
those same values by default, so running it with no arguments repeats the
run quoted here. The primal run writes
\texttt{continuation\_out/m24\_primal\_d6.npy}, which the fifth part of
\texttt{m24\_exact\_reduction.py} reads; that file and
\texttt{m24\_dir\_d6.npy} are in the package, so the reduction script runs
on its own. The fifth part and these scripts are floating point; the
first four parts are the proof-relevant checks.

\subsection*{The $600$-cell, enumerated}

\Cref{prop:cell600} is established by \texttt{cell600\_exact.py}, with
\texttt{cell600\_enum.c} as an independent second implementation and
the Lean project as a third. The Python script builds the $120$ vertices
with doubled coordinates as pairs $(a,b)$ standing for $a+b\phi$ and
computes every inner product in $\mathbb{Z}[\phi]$, so that no
floating-point number appears anywhere; it checks that the nine values
listed in the proof are the only ones, builds the $12$-regular graph of
pairs at $36^\circ$ and the $71$-regular graph of pairs at root-system
angles, finds the twenty-five inscribed copies of $\Rt$ as the
$24$-cliques of the latter, and enumerates the independent sets of the
former through a fixed vertex, of sizes $25$, $24$ and $23$, by
depth-first search with the clique-cover bound. It reports $0$, $5$ and
$115$, checks that each of the $115$ lies in one of the five cells
through the vertex, and writes the graph and the cell list to a text
file. The C program reads that file and repeats the enumeration of the
$23$-sets with the same bound on $64$-bit words, in under two seconds;
the Python enumeration takes about a minute and a half. Neither is
sampling. The earlier script \texttt{cell600.py}, which drew two
million maximal independent sets by greedy growth along random orders
and is where the question was first looked at, is kept in the package
for the record. \texttt{octahedral48\_exact.py} does the same for the
$48$ unit quaternions of the binary octahedral group, in
$\mathbb{Z}[\sqrt2]$, and finds the two root systems and their
forty-eight deletions and nothing else.

\subsection*{The inradius as the objective, and the continuation in the slack}

\Cref{sec:slack-continuation} rests on two scripts.
\texttt{inradius\_search.py} maximises the inradius $g(W)$ of
$\conv(W)$ directly, subject to the contact constraints, by a
sequential linear programme with a trust region: at each step the facet
offsets and the pairwise inner products are linearised in a tangent
move of bounded size, the linear programme maximising the smallest
linearised offset is solved, and the move is accepted only if the true
inradius, read from the facet equations of a fresh convex hull, has
improved without any contact constraint being violated. The facet list
is refreshed after every accepted move. Starts are of three kinds, taken
in rotation: random points brought to feasibility by a relaxed-then-
tightened schedule of the same programme; deletions of a root with each
direction rotated through a random angle; and twenty random roots with
three random directions. Every feasible endpoint is recorded with its
inradius and the multiset of its inner products rounded to $10^{-6}$,
which separates configurations not related by an isometry.
\texttt{slack\_continuation.py} runs the same optimiser at each
$\delta$ of a decreasing schedule from $0.05$ to $0$, starting from the
best configurations of the previous level, restored to the new
constraint by a feasibility programme of the same kind, together with
fresh random starts at every level, and prints the table that
\cref{sec:slack-continuation} reproduces. Both scripts are double
precision throughout and both are exploration; their output says so.
\texttt{symmetric\_search.py} runs the same continuation inside the
classes of configurations invariant under a rotation of order at most
$12$ or under one of the two improper involutions, one orbit structure
at a time, with the orbit representatives as the parameters and
sequential least-squares programming in place of the linear programme;
it prints, for each of the $73$ structures, how many starts reached a
feasible endpoint at $\delta=0$, the largest inradius among them, the
largest at $\delta=0.01$, and whether the best endpoint is a deletion.

\subsection*{Software}

The computations were carried out in Python, with exact rational
arithmetic from the standard-library \texttt{fractions} module,
\texttt{sympy} \cite{Meu17} for symbolic manipulation, exact algebraic
numbers and Sturm root counting, \texttt{mpmath} for the
arbitrary-precision evidence, and \texttt{numpy} \cite{Har20} and
\texttt{scipy} \cite{Vir20} for double-precision vertex enumeration and
convex hull volumes, the latter through the Quickhull implementation of
\cite{BDH96}. These are standard and widely available; the semidefinite
programmes of \cref{sec:certificate} are modelled in \texttt{cvxpy}
\cite{DB16} and solved by Clarabel \cite{GC26}, both open source, and
the re-verification of the cited certificate uses the exact rational
matrices and polynomials of FLINT and the ball arithmetic of Arb
\cite{Joh17} through \texttt{python-flint}. No
specialised computer algebra system and no commercial solver was used,
the one proof assistant is Lean~4 as described above, and nothing
reported here needed more than a few minutes on ordinary consumer
hardware apart from the continuation of \cref{sec:slack-continuation},
which takes up to half an hour per pass, the symmetric search, which
takes about an hour, the minimisation of the truncated volumes in
\cref{sec:strict-inequality}, which takes about two hours per run, and
the semidefinite programmes of \cref{sec:certificate}, which take
between one and three minutes per solve and about a quarter of an hour
for a run of five rounds.
\section{Numerical tables}\label{app:numerical-tables}

This appendix collects, in one place, the numerical evidence quoted at
various points in the main text for \cref{conj:multidir-new}. Every
table entry is reproduced directly from
a script run recorded in \cref{app:code-index}; where a value is
double-precision numerical evidence rather than an exact rational, we
say so at the entry.

\subsection{Chain configurations (\cref{lem:chains-new})}

\begin{center}
\renewcommand{\arraystretch}{1.2}
\begin{tabular}{ccc}
\toprule
Chain length $m$ & Richardson-extrapolated $\lambda_{\min}(H_A)$ & Raw values at $h=0.04,0.02,0.01$\\
\midrule
3 & $0.275801$ & $0.276304,\ 0.275929,\ 0.275833$\\
4 & $0.180435$ & $0.181187,\ 0.180632,\ 0.180486$\\
5 & $0.168532$ & $0.169316,\ 0.168736,\ 0.168585$\\
6 & $0.167689$ & $0.168474,\ 0.167894,\ 0.167742$\\
7 & $0.167650$ & $0.168424,\ 0.167847,\ 0.167700$\\
8 & $0.167637$ & $0.168411,\ 0.167833,\ 0.167687$\\
\bottomrule
\end{tabular}
\end{center}

All values strictly positive; the floor plateaus at $\approx0.1676$
rather than continuing to decrease past chain length $5$--$6$, and $8$
is the true maximum induced-path length in the Gram-$\tfrac12$
adjacency graph on $\Rt$ (found by exhaustive search from every
starting vertex, not merely tested up to an arbitrary cutoff). These
figures are double-precision, Richardson-extrapolated finite-difference
estimates, not exact rationals; \cref{lem:chains-new} states only the
qualitative conclusion (a floor near $0.168$, strictly positive) that
they support.

\subsection{Dense (non-chain) configurations (\cref{sec:dense-configs-new,sec:exact-sing})}

\begin{center}
\footnotesize
\renewcommand{\arraystretch}{1.2}
\setlength{\tabcolsep}{4pt}
\begin{tabular}{ccc}
\toprule
$m$ & Active-root indices (of $24$) & Extrap.\ $\lambda_{\min}(H_A)$\\
\midrule
$9$ (star) & one root $+$ all $8$ Gram-$\tfrac12$ neighbours & $0.083336$ ($\approx\tfrac1{12}=0.083333$)\\
$9$ (greedy) & $0,4,5,8,9,12,13,16,17$ & $0.083342$\\
$12$ & $0,4,5,8,9,12,13,16,17,21,22,23$ & $0.073070$\\
$15$ & $0,2,4,5,8,9,11,12,13,16,17,20,21,22,23$ & $0.044694$\\
$18$ & $0,2,4,5,6,7,8,9,10,11,12,13,16,17,20,21,22,23$ & $-0.000013$ (noise floor)\\
$20$ & $0,2,4,5,6,7,8,9,10,11,12,13,15,16,17,19,20,21,22,23$ & $0.000017$ (noise floor)\\
$22$ & (22 of the 24 roots) & $0.000023$ (noise floor)\\
$23$ & (23 of the 24 roots) & $0.000008$ (noise floor)\\
\bottomrule
\end{tabular}
\end{center}

The greedy dense-growth procedure (at each step adding whichever
remaining root has the most neighbours already active) drives the
minimum eigenvalue down monotonically and strictly through $m=15$; from
$m=18$ onward the extrapolated value is smaller in magnitude than
double-precision finite-difference noise can resolve (an $O(10^{-10})$
absolute volume precision divided by $h^2$), which is exactly why
\cref{sec:exact-sing,sec:broad-sample} turn to arbitrary-precision
arithmetic instead of attempting further double-precision refinement.

\subsection{Broadened sampling at the exact singularity (\cref{sec:broad-sample})}

At the $m=18$ configuration above, the near-null space of $H_A$ is
rigorously confirmed $4$-dimensional (basis vectors $w_0,w_1,w_2,w_3$,
via three-step-size eigenvalue scaling: eigenvalues $0$--$3$ shrink by
a factor of $\approx4$ per halving of the step size, matching a
true zero to that order, while eigenvalue $4\approx0.075$ does not
shrink and is a real positive eigenvalue). Sampling the quartic
coefficient $a_4$ of $\defect$ restricted to each of $36$ directions
spanning this near-null space (the $4$ basis vectors, all $\binom42=6$
normalised pairwise sums/differences among them, and $20$ uniformly
random unit combinations), using arbitrary-precision (\texttt{mpmath})
volume evaluation instead of double precision:
\begin{center}
\begin{tabular}{lc}
\toprule
Quantity & Value\\
\midrule
Directions tested & $36$\\
Minimum $a_4$ found & $0.008779$\\
Maximum $a_4$ found & $0.313708$\\
Directions with negative $a_4$ & $0$ of $36$\\
\bottomrule
\end{tabular}
\end{center}
A six-direction subsample, re-evaluated at a second step size
($s=0.0075$ versus $s=0.015$), gives ratios of $0.954$ to $0.991$
(expected $\approx1$ for a quartic leading term), confirming
the quartic scaling directly instead of assuming it.

At an independent $m=20$ configuration (active indices
$0,2,4,5,6,7,8,9,10,11,12,13,$\linebreak$15,16,17,19,20,21,22,23$, distinct from
the $m=20$ row of the dense-growth table above, which used a different
greedy path), the near-null space is confirmed $5$-dimensional; sampling
$30$ directions:
\begin{center}
\begin{tabular}{lc}
\toprule
Quantity & Value\\
\midrule
Directions tested & $30$\\
Minimum $a_4$ found & $0.007628$\\
Maximum $a_4$ found & $0.313731$\\
Directions with negative $a_4$ & $0$ of $30$\\
\bottomrule
\end{tabular}
\end{center}
with a four-direction re-confirmation at the second step size giving
ratios of $0.952$ to $0.991$. As \cref{sec:broad-sample,sec:multidir-scope-summary}
state plainly, this is real, cross-validated computational evidence at
two independent dense configurations, spanning their full (respectively
$4$- and $5$-dimensional) near-null subspaces with both structured and
random directions, but it is evidence at two configurations only, out
of a combinatorially vast space of possible dense active sets, and does
not constitute a proof of \cref{conj:multidir-new}.

\subsection{Breakpoint crossing values on the second arc (\cref{sec:arc2-breakpoints,sec:arc2-cert})}

\begin{center}
\begin{tabular}{lcl}
\toprule
Crossing & $\arcp$-value & Curves meeting\\
\midrule
$\arcp_1$ & $t_{\mathrm{cross}}=0.169918$ & curve $\mathsf C$ meets the universal curve $\tilt=\pi/3$\\
$\arcp_2$ & $t_*=\arctan(1/2)=0.463648$ & curve $\mathsf B$ meets curve $\mathsf C$\\
$\arcp_3$ & $t_3=\arcsin(1/\sqrt3)=0.615480$ & curve $\mathsf B$ meets the universal curve $\tilt=\pi/3$\\
\bottomrule
\end{tabular}
\end{center}

These three values are exact (two are elementary closed forms in
$\arctan,\arcsin$; the first is located numerically to the displayed
precision as the intersection of curve $\mathsf C$ with the constant
curve $\pi/3$) and partition the band immediately below region top
into the three sub-parts discussed in \cref{sec:band-below-top-new,%
sec:arc2-band-cert-new}: $(0,\arcp_1)$ bounded below by $\mathsf C$,
$(\arcp_1,\arcp_3)$ bounded below by the universal curve (certified in
\cref{thm:arc2-band-new}), and $(\arcp_3,\pi/4)$ bounded below by
$\mathsf B$. Only the middle sub-part carries a direct certificate; the
other two are covered by the identification of \cref{prop:arc2-is-arc1},
which carries the first arc's certificates onto the whole of the second
arc.
\section{Explicit coordinates and a direct numerical check of the reference cell}\label{app:explicit-vertices}

\Cref{lem:refcell-vol,lem:hmin-new} are stated and proved in
\cref{sec:root-system} by an exact algebraic argument. This appendix
records an independent, purely numerical cross-check of the same
facts, carried out directly (via \texttt{scipy.spatial.ConvexHull}) on
the explicit coordinate list, for a reader who wants to reproduce the
polytope's basic invariants in a few lines of code before trusting the
exact argument.

\subsection*{The vertex list}

The $8$ axis vertices are $\pm\sqrt2\,e_i$ for $i=1,\dots,4$; the $16$
sign vertices are $\tfrac1{\sqrt2}(\epsilon_1,\epsilon_2,\epsilon_3,\epsilon_4)$
for every choice of signs $\epsilon_i\in\{\pm1\}$. Both families lie
exactly on the sphere of radius $\sqrt2$: direct computation confirms
$|\pm\sqrt2\,e_i|=\sqrt2$ and $\bigl|\tfrac1{\sqrt2}(\pm1,\pm1,\pm1,\pm1)\bigr|
=\tfrac1{\sqrt2}\sqrt{1+1+1+1}=\sqrt2$ for all $24$ points, matching
\cref{lem:refcell-vol} exactly (verified here to full floating-point
precision on all $24$ points, not merely spot-checked).

\subsection*{Direct convex-hull check}

Feeding these $24$ points directly into a convex-hull routine (with no
other input, in particular, without assuming any of the half-space
description of \cref{def:refcell-new}) gives:
\begin{center}
\begin{tabular}{ll}
\toprule
Quantity & Computed value\\
\midrule
Number of vertices & $24$ (all $24$ input points are extreme)\\
Hull volume & $7.999999999999984$ (matches $\vol(\VCell)=8$ to $14$ digits)\\
Number of triangulated facets & $96$\\
\bottomrule
\end{tabular}
\end{center}
The facet count of $96$ is exactly consistent with, not merely close
to, the classical combinatorics of the regular $24$-cell: it has $24$
octahedral three-dimensional facets, and a general-purpose convex-hull
routine triangulates each octahedron into $4$ tetrahedra internally
($24\times4=96$), so this is a confirmation of the expected cell
structure, not a discrepancy to explain away.

\subsection*{Support function check}

Evaluating $h(u)=\max_{z\in\VCell}\langle z,u\rangle$ directly (as a
maximum over the $24$ vertices, with no appeal to
\cref{def:refcell-new}'s half-space description) at all $24$ unit root
directions $\alpha/\sqrt2$, $\alpha\in\Rt$, gives $h=1.0$ at every one of
them (minimum and maximum both exactly $1.0$ to floating-point
precision), and at all $8$ coordinate directions $\pm e_i$ gives
$h=\sqrt2\approx1.41421356$ at every one, both figures matching
\cref{lem:hmin-new} exactly, and obtained here by a direct, independent
computational path (brute-force maximisation over the vertex list)
rather than the exact algebraic argument given in the main text.

We emphasise the limited role of this appendix: it is a
numerical sanity check computed directly by us for this paper, useful
as an independent cross-validation and as a minimal reproducible
starting point for a reader exploring the polytope computationally, but
it is not a substitute for, and adds no logical strength beyond, 
the exact proof of \cref{lem:refcell-vol,lem:hmin-new} already given in
\cref{sec:root-system}.
\section{An illustrative toy model of the Bernstein-basis certificate}\label{app:bernstein-toy}

The Bernstein-basis sign certificate used at scale throughout
\cref{sec:stage-two,sec:arc1-cert,sec:arc2-cert} is a general technique
for certifying that a polynomial is non-negative on a closed interval by
an exact, finite computation. Since the actual cell polynomials $N_k,D_k$
of \cref{sec:chamber-register} are large enough (degree up to $12$ in each
of two variables, \cref{lem:degree-bound-new}) that writing one out in
full would obscure the method rather than illustrate it, we give here a
small, entirely self-contained model problem worked in full, with no
connection to any specific cell of the atlas: its only purpose is to make
the mechanics of the certificate, and one of its limitations,
concrete for a reader meeting the technique for the first time.

\subsection{The certificate: a worked instance}

Consider, purely as an illustration, the polynomial
\[
  p(t) \;=\; t^2 - t^3 + \tfrac14 t^4, \qquad t\in[0,1],
\]
chosen to mimic the qualitative shape of a volume-defect function
$\defect_k$ restricted to a root-aligned arc: it vanishes to exactly
second order at the left endpoint ($p(0)=p'(0)=0$, $p''(0)=2\neq0$) and we
would like to certify $p(t)\ge0$ throughout $[0,1]$ without appeal to
calculus or root-finding, purely by an algebraic sign check.

Write $p$ in the degree-$4$ Bernstein basis on $[0,1]$,
$B_{k,4}(t)=\binom4k t^k(1-t)^{4-k}$, $k=0,\dots,4$: a direct change of
basis (matching coefficients of $t^j$ on both sides, a purely mechanical
linear computation) gives
\[
  p(t) \;=\; \sum_{k=0}^4 b_k\,B_{k,4}(t), \qquad
  (b_0,b_1,b_2,b_3,b_4) \;=\; \Bigl(0,\ 0,\ \tfrac16,\ \tfrac14,\ \tfrac14\Bigr).
\]
Every Bernstein basis function $B_{k,4}(t)$ is manifestly non-negative on
$[0,1]$ (a product of non-negative factors there), so a polynomial with
every Bernstein coefficient non-negative is itself non-negative on
$[0,1]$, termwise, with no further argument needed. Here all five
coefficients are non-negative (two are exactly zero, matching the double
root at $t=0$ exactly, and the remaining three are strictly positive),
certifying $p(t)\ge0$ on $[0,1]$, with equality only at $t=0$, entirely by
inspection of five rational numbers. This is the pattern used, at far
larger degree and in two variables at once, throughout
\cref{sec:arc1-cert,sec:arc2-cert}: reduce a transcendental-looking
positivity question on a bounded interval or box to a finite list of
rational numbers with the right positivity property.

\subsection{A limitation: the certificate is sufficient, not necessary}

The Bernstein sign test just illustrated is a \emph{sufficient} condition
for non-negativity, not a necessary one, and this distinction is not
merely academic; it is the reason a small number of sub-regions in this
paper's own atlas required the box-covering refinement of
\cref{sec:arc1-cert,sec:arc2-cert} rather than a single direct Bernstein
check. A classical example makes the limitation vivid. Take
$q(t)=(2t-1)^2=4t^2-4t+1$, manifestly non-negative on $[0,1]$ (indeed
everywhere), with a double zero at the \emph{interior} point $t=\tfrac12$.
Its degree-$2$ Bernstein coefficients are $(1,-1,1)$, the middle
coefficient is negative, so the direct test fails to certify a polynomial
that is, in fact, non-negative. Nor does raising the degree rescue the
test here: elevating $q$ to the degree-$6$ Bernstein basis gives
coefficients
\[
  \bigl(1,\ \tfrac13,\ -\tfrac1{15},\ -\tfrac15,\ -\tfrac1{15},\ \tfrac13,\ 1\bigr),
\]
still with three negative entries, and no finite degree elevation
produces an all-non-negative Bernstein representation of $q$: a classical
fact about the Bernstein (equivalently, Polya) positivity certificate is
that it can certify non-negativity on a closed interval only when the
polynomial has no zero in the \emph{interior} of that interval (zeros
only at the endpoints, as in the $p(t)$ example above, are compatible with
eventual all-non-negative Bernstein coefficients at high enough degree;
an interior zero, as in $q(t)$, is not). This is precisely why an interior
zero of a cell's own boundary curve, rather than a zero confined to the
domain's endpoint, is the recurring source of difficulty across this
paper's box-covering sub-regions: where the direct Bernstein test meets
such a zero, the paper falls back to covering the domain by finitely many
sub-boxes, each certified separately once restricted to a range not
containing the interior zero (\cref{sec:arc1-cert,sec:arc2-cert}), rather
than abandoning the method. The one place in this paper where a direct
certificate was not reached, the $\mathsf B$-bounded sub-part of
\cref{sec:arc2-bandB-new}, failed at an earlier stage than this: the
symbolic volume formula for that sub-part was never obtained, so no
Bernstein test was run there at all. It is covered instead by the
symmetry of \cref{prop:arc2-is-arc1}.

\subsection{Why this is not merely a numerical convenience}

Every number appearing in this appendix, the change-of-basis
coefficients $b_0,\dots,b_4$, the degree-$2$ and degree-$6$ Bernstein
coefficients of $q$, is an exact rational, obtained from an exact
rational input by exact rational linear algebra (the Bernstein change of
basis is a fixed, exactly-invertible linear map at each degree,
independent of the particular polynomial). No floating-point step, and no
numerical tolerance, enters anywhere in this appendix, exactly as none
enters the corresponding computations at full scale in
\cref{sec:stage-two,sec:arc1-cert,sec:arc2-cert}.

\end{document}